\documentclass[10pt,twoside,a4paper,abstract=on]{scrartcl}

\usepackage{authblk}

\usepackage{geometry}
\usepackage{multicol}
\usepackage{multirow}

\usepackage[utf8]{inputenc}
\usepackage{hyperref}
\usepackage{caption}

\usepackage{todonotes}

\usepackage[toc,title,page]{appendix}

\usepackage{amsmath,amsfonts,amssymb,amsthm}
\usepackage{bbm}

\usepackage{enumitem}

\usepackage{tikz}  

\usepackage{algorithm}
\usepackage{algpseudocode}

\newtheorem{Theorem}{Theorem}[section]

\newtheorem{Proposition}[Theorem]{Proposition}
\newtheorem{Lemma}[Theorem]{Lemma}
\newtheorem{Assumption}[]{Assumption}

\newtheorem{Remark}[Theorem]{Remark}

\usepackage[a-1b]{pdfx}

\usepackage{graphicx}
\usepackage{subcaption}

\usepackage[maxbibnames=6]{biblatex}
\definecolor{electricultramarine}{rgb}{0.25, 0.0, 1.0}
\definecolor{ikb}{rgb}{0.0, 0.18, 0.65}
\definecolor{green(colorwheel)(x11green)}{rgb}{0.16, 0.5, 0.0}

\newcommand*{\fzcst}[1]{\relax\ifmmode\text{\textcolor{green(colorwheel)(x11green)}{\sout{\ensuremath{#1}}}}\else\textcolor{green(colorwheel)(x11green)}{\sout{#1}}\fi}

\usepackage{todonotes,varwidth}

\title{
	Dynamical Low-Rank Filters \\
	for Data Assimilation}
\author[1]{Yoshihito Kazashi}
\author[2]{Youssef Marzouk}
\author[3]{Fabio Nobile}
\author[3]{Fabio Zoccolan}
\affil[1]{Department of Mathematics,
	University of Manchester, Oxford Road, Manchester, M13 9PL, UK. email: y.kazashi@manchester.ac.uk}
\affil[2]{Department of Aeronautics and Astronautics, Massachusetts Institute of Technology, Cambridge, MA 02139, USA. email: ymarz@mit.edu}
\affil[3]{Institut de Mathématiques, École Polytechnique Fédérale de Lausanne, 1015 Lausanne, Switzerland. email: fabio.nobile@epfl.ch, fabio.zoccolan@epfl.ch}
\date{}

\begin{document}
   
   \maketitle
	
	\begin{abstract}
		We propose dynamical low-rank (DLRA) filters for data-assimilation problems based on stochastic differential equations (SDEs).
		We first consider a Kalman filter applied to a Euler-Maruyama discretization of the SDE and observation processes, and derive a DLRA Kalman-type filter, named DLRA-JMCO, by minimizing jointly the error on the mean and covariance. We then obtain a time continuous DLRA-JMCO filter as the continuous limit of the time discrete one, and introduce a corresponding particle approximation. We also propose a strategy to enrich the DLRA-JMCO filter with a simple Gaussian process to capture the neglected dynamics in the orthogonal to the DLRA subspace. This correcting term allows the DLRA subspace evolution to be directly influenced by the observation operator. 
		
		We finally propose a particle filter version of the DLRA-JMCO filter. Numerical experiments assess the performance of the proposed filters.		
    \end{abstract}
    %\tableofcontents
    
    %\listoftodos

	\section{Introduction}
	\addcontentsline{toc}{section}{Introduction}
	
    Numerical simulations of real-world systems rely on mathematical models derived from physical principles and engineering studies. Their accuracy, however, is limited by several sources of uncertainty, including model error and the intrinsic randomness affecting natural and engineering systems. Quantifying these uncertainties often requires probabilistic approaches, which can be significantly more computationally demanding than deterministic analyses.
    
    Data assimilation addresses these limitations by combining model predictions with observational data to improve state estimates and reduce uncertainty. A common approach is to formulate the problem as a \emph{Bayesian inverse problem}, treating both the system state and the observations as random variables. The resulting filtering procedure typically consists of two steps: a \emph{prediction} step, in which the state evolves according to the chosen model, and an \emph{analysis} step, in which observational data are used to update the state distribution through conditioning. Since both the dynamics and the observations may be affected by noise, stochastic differential equations (SDEs) naturally arise in this framework.
    
    Data assimilation plays a central role in many applications, including weather forecasting and oceanography \cite{carrassi2018data,bennett2005inverse,kalnay2003atmospheric}, as well as automotive engineering \cite{huang2025improved}, robotics \cite{thrun2005probabilistic}, and finance \cite{javaheri2003filtering,kellerhals2013financial,kulikova2022estimation,wells2013kalman}.   
    Its applicability, however, is often limited by the computational cost of \emph{high-dimensional systems}, where both the state and the observations may require large representations, for instance through particle-based approximations. 
    
    A common remedy is \emph{reduced-order modelling} (ROM), in which the original system is replaced by a lower-dimensional surrogate that is cheaper to evaluate while retaining high accuracy when the effective dynamics are intrinsically low-dimensional \cite{cao2007reduced,cheng2023machine,hoteit2002simplified,popov2022multifidelity}. ROM can reduce both the physical and stochastic complexity of the problem, potentially decreasing the number of particles required when the solution is concentrated near a low-dimensional subspace.
    Nevertheless, standard ROM strategies may be inadequate for data assimilation, particularly when based on a fixed reduced subspace. Such a subspace cannot adapt to evolving, possibly chaotic dynamics or to the information introduced by incoming observations.
    
   To overcome the limitations of fixed subspace ROMs, an appealing approach is offered by dynamical low-rank approximation (DLRA) techniques. DLRA is a surrogate method characterized by time-dependent components, and, hence, the procedure can be computed completely on-the-fly. In this regard, the main low-dimensional subspace individuated by the DLRA is allowed to evolve over time, without any dependence on a full-order solution. This property is very appealing for filtering procedures, as DLRA can manage to adapt to unexpected changes due to the chaotic dynamics and to enrich its bases incorporating the observed high-dimensional data, while always remaining low-dimensional, and, hence, cheap to compute. \emph{The aim of this work is to provide suitable DLRA filtering techniques for data-assimilation problems modeled by stochastic differential equations} (SDEs), which can also be deployed in nonlinear or non-Gaussian problems.
    
    DLRA was first introduced by \cite{koch2007dynamical} in the context of matrix ordinary differential equations (ODEs), which continues to remain a very active theme of research \cite{carrel2024randomised,ceruti2022rank,ceruti2024robust,lubich2014projector,}. Thanks to its useful benefits, this framework was rapidly studied for deterministic \cite{bachmayr2021existence,dektor2025interpolatory} and random partial differential equations \cite{musharbash2015error,musharbash2018dual,kazashi2021stability}, and only recently a well-posed framework for SDEs has been studied \cite{bao2026exponential,kazashi2025dynamical,kazashi2026locally} through dynamically orthogonal approximation \cite{sapsis2009dynamically} (see \cite{cao2018stochastic} for another approach). Some preliminary results of DLRA applied to data assimilation have been proposed in \cite{lu2021bayesian,sapsis2013blending, sapsis2013interaction, sondergaard2013dataI, sondergaard2013dataII}, where the dynamically orthogonal method - a possible instantiation of DLRA - was tested to oceanography or chaotic systems and merged with other techniques like Gaussian Mixtures, too. In that framework, the assimilation was performed at discrete time.    
    Recently, in \cite{nobile2025dynamicallowrankapproximationskalman,nobile2026dynamicallowrankensemblekalman}, the DLRA-for-SDEs framework from \cite{kazashi2025dynamical} was applied to derive Kalman-type filters. In detail, in \cite{nobile2025dynamicallowrankapproximationskalman} a DLRA Kalman-Bucy filter is derived for linear drift, a procedure that naturally generalizes to an ensemble algorithm and can be effectively applied to real-world problems, such as hemodynamics \cite{nobile2026dynamicallowrankensemblekalman}.
    
    In the context of SDEs, DLRA relies on a DO formulation of the form $X^{\mathrm{DLRA}}_t=U^{\top}_tY_t$, where $U_t$ is a deterministic $d \times k$ matrix with orthonormal columns satisfying $\dot{U}_tU^{\top}_t=0$ and $Y_t$ is a suitable stochastic process, both having $k \ll d$ independent components \cite{kazashi2025dynamical}. 
    If, on the one hand, the DLRA equations for SDEs present numerical and theoretical difficulties due to the heavy coupling of the deterministic and stochastic bases and the strong dependence on the law of the process, on the other, due to this particular structure, the fully-discretized reduced solution has the advantage to be a noisy interacting particle system. This interpretation turns out to be very suitable for building ensemble filtering methods to deal with nonlinear and/or non-Gaussian problems. 
	
	First, we propose a DLRA-type filter that minimizes jointly the mean and the covariance at each time step. In the case of linear drift and diffusion, this procedure can be interpreted as a type of reduced Kalman-Bucy filter (KBF), different to the one derived in \cite{nobile2025dynamicallowrankapproximationskalman}. We propose its construction in Section \ref{sec: DLRA for JMCO}, where the Algorithm will be denoted as \emph{DLRA JMCO} filter. Later, via a similar procedure, we advance an akin filtering method aimed to track the locally irrelevant, but maybe not negligible over time, complementary directions to the subspace. It turns out that in this method the main subspace can also evolve depending on the observation operator. This property ensures more adaptiveness to the surrogate evolution, especially in the case of discrete observations obtained via a time-dependent observation operator or when the subspace individuated by the observation is not aligned with the one of the state dynamics. We highlight the derivation of this filtering procedure in Section \ref{sec: DLRA compl JMCO}. Another relevant feature of all of the proposed techniques is the versatility of being applicable to both discrete or continuous-time data framework.
	
    Finally, we advance another filtering procedure that can be exploited when dealing with strongly nonlinear problems: we propose a reduced-type particle filter, where the prediction step is obtained through a DLRA update, and the analysis one acts through a branching procedure whose weights depend on the likelihood estimator, as proposed in \cite{bain2009fundamentals}. Also in this case, the subspace of DLRA does not evolve between the prediction and the analysis step. To overcome this problem, we propose a computationally cheap procedure that exploits the usual DLRA projection approximability assumption \cite{lubich2014projector,kazashi2021stability,kazashi2025dynamicalpartI}. We summarize the whole procedure in Section \ref{sec: compl DLRA PF}. For the sake of readability, we will postpone most of the proofs presented in these Sections to the Appendix at the end of the article.

    Numerical experiments showing the potentiality of all these methods will be presented in Section \ref{sec: numerics filt}.
	
	\section{Setting}\label{sec: setting}
	Let $\left( \Omega, \mathcal{F}, \mathbb{P}, (\mathcal{F}_t)_{t \geq 0} \right)$ be a filtered complete probability space with the usual conditions \cite[Remark 6.24]{schilling2021brownian}.
	For $T>0$, we consider the following continuous-time data assimilation problem based on SDE models for the state and observation processes, respectively:
	\begin{equation}\label{eq: DA SDE}
		\begin{aligned}
			\mathrm{d}X_t &= A(t, X_t) \, \mathrm{d}t + Q^{\frac{1}{2}} \, \mathrm{d}W_t, \quad t>0, \  X_0 \text{ given} \\
			\mathrm{d}Z_t &= H(t) X_t \, \mathrm{d}t + R^{\frac{1}{2}} \, \mathrm{d}B_t, \quad t>0, \ Z_0 =0,
		\end{aligned}
	\end{equation}
	where at time $t$ the random vector $X_t=\left(X_1(t), \ldots, X_d(t)\right)^{\top}$ belongs to $L^2(\Omega, \mathbb{R}^d)$, whereas $Z_t=\left(Z_1(t), \ldots, Z_h(t)\right)^{\top}$ belongs to $L^2(\Omega, \mathbb{R}^h)$, $A : [0, \infty) \times \mathbb{R}^d \to \mathbb{R}^d$ is the (possibly nonlinear) drift, $H : [0, \infty) \to \mathbb{R}^{d \times h}$ is a linear observation operator, $W$ and $B$ are real $d$-dimensional and $h$-dimensional $(\mathcal{F}_t)$-Brownian motions, respectively,
	i.e., $W(t) = \left(W_1(t), \ldots , W_d(t) \right)^{\top}$ and $B(t) = \left(B_1(t), \ldots , B_h(t) \right)^{\top}$, and $Q \in \mathbb{R}^{d \times d}$, $R \in \mathbb{R}^{h \times h}$ are symmetric positive semidefinite matrices which we will refer to as the noise
	covariances of the state and the observation, respectively. We assume that $\{W(t)\}_{t\in [0,T]}$, $\{B(t)\}_{t\in [0,T]}$, and $X_0$ are independent.
	
	In the context of data assimilation, $X_t \in \mathbb{R}^d$ is the “noisy’’ \emph{state} that we are trying to infer,
	whereas $Z_t \in \mathbb{R}^h$ is the “noisy’’ \emph{observation} or \emph{data}. We assume $h \leq d$, i.e.\ the dimension of the available information is always smaller or equal of the state dimension. The goal of data assimilation and filtering is to estimate the law of the state $X_t$ conditional to the observation process $\{Z_s\}_{0\leq s \leq t}$ up to time $t$. In general, the filtering problem \eqref{eq: DA SDE} does not admit a closed-form solution without additional assumptions. A notable exception is the linear-Gaussian setting. In particular, when the drift $A$ is linear and the initial condition is Gaussian, the filtered process \eqref{eq: DA SDE} remains Gaussian. The \emph{Kalman--Bucy filter} then provides evolution equations for its conditional mean $m_t$ and covariance $C_t$, which completely characterize the filtered distribution \cite[Theorem 8.1]{law2015data}:
	\begin{equation}\label{eq: kalman-bucy}
		\begin{aligned}
			\mathrm{d} \mathbb{E}[X_t]= & A(t)  \mathbb{E}[X_t] \mathrm{d}t + C_t H(t)^{\top}  R^{\frac{1}{2}} \left( \mathrm{d}Z_t - H(t)\mathbb{E}[X_t] \mathrm{d}t\right) \\
			\mathrm{d} C_t = & \left( A(t) C_t + C_t A(t)^{\top} + Q -  C_t H(t)^{\top}  R^{-1} H(t)C_t \right) \mathrm{d}t.
		\end{aligned}
	\end{equation}
	Equivalently, one can construct a filtered process $X_t^f$ whose Gaussian distribution has mean and covariance as in \eqref{eq: kalman-bucy}. Such process can be generalized to the case of a non-linear drift and reads 
	\begin{equation}\label{eq: kalman-bucy nonlinear}
		\mathrm{d}X_t^f = A(t, X_t^{f}) \mathrm{d}t + Q^{\frac{1}{2}} d\widehat{W}_t + C_t^{f} H(t)^{\top} R^{-1} (\mathrm{d}Z_t - H(t)X_t^{f} \mathrm{d}t - R^{-\frac{1}{2}}\mathrm{d}\widehat{B}_t),
	\end{equation}
	where $\widehat{W}_t$ and $\widehat{B}_t$ are Brownian motions with the same law as $W_t$ and $B_t$ and $C_t^{f}$ is the covariance matrix of $X_t^f$. We should point out, however, that if the drift is non-linear, the law of $X_t^f$ does not coincide with the exact filtered distribution of \eqref{eq: DA SDE}. In this respect, \eqref{eq: kalman-bucy nonlinear} is an inconsistent approximation of the filtered distribution. However, it is widely used in practice because a particle approximation of \eqref{eq: kalman-bucy nonlinear} can be easily implemented and often provides acceptable results in applications.
	
	For easiness of presentation, in this work we restrict to SDEs in \eqref{eq: DA SDE} with additive noise. However, the treatment proposed in this work can be extended with no difficulty to the case of multiplicative-noise type dynamics of the state and/or observation.
	
	The running assumption in this work is that the dimension of the state is high, i.e.\ $d \gg 1$, however its (filtered) distribution remains reasonably concentrated around a relatively low dimensional subspace at each time $t$, thus justifying the use of model order reduction techniques and, particularly, DLRA.
	\subsection{DLRA for SDEs}
	We briefly summarize the main elements of DLRA for SDEs in order to highlight the framework in which our filtering procedure will be built. In \cite{kazashi2025dynamical}, DLRA was formulated via the so-called dynamically orthogonal fields equations, a time-depedent reduced order model made by two components, a deterministic orthogonal basis and a stochastic (linearly independent) one, both allowed to evolve in time \cite{sapsis2009dynamically}. In this work, we follow the slightly more complex three-term formulation in \cite{sapsis2009dynamically}, which isolates the evolution of the mean. In more details, for any realization $\omega \in \Omega$ the DO solution is built as the following sum of the mean and a linear combination of $k$ members:
	\begin{equation}\label{eq: DO}
		X^{\mathrm{DLRA}}_t(\omega) = m^{\mathrm{DLRA}}_t + \sum_{i = 1}^{k} (U^{i}_t)^{\top}Y^{i}_t(\omega)=  m^{\mathrm{DLRA}}_t + U_t^{\top}Y_t(\omega), \quad t \geq 0,
	\end{equation}
	where $m^{\mathrm{DLRA}}_t$ is the mean of the DLRA at time $t$, $U^{\top}_t = \left((U^{1}_t)^{\top}, \dots, (U^{k}_t)^{\top} \right)\in  \mathbb{R}^{d \times k}$ is the deterministic basis, and $Y_t = \left(Y^{1}_t, \dots, Y^{k}_t\right) \in \left(L^2_0(\Omega)\right)^k$ is the centered stochastic basis, respectively, and $k \ll d$, is the rank of the approximation. Notice that the covariance of $X^{\mathrm{DLRA}}_t$ is simply given by $C^{\mathrm{DLRA}}_{t}= U_t^{\top} \mathbb{E}[Y_t Y_t^{\top}]U_t$. In this work we consider this three-term DLRA version, at the place of a two-component ones of \cite{kazashi2025dynamical}, because it allows us to have an immediate connection with the first and second (centered) moments of the full-order solution associated to \eqref{eq: DA SDE}.
	
	For a general Itô SDE of the following type
	\begin{equation}\label{eq:SDE-diff}
		\mathrm{d}X_t = a(t,X_t)\mathrm{d}t+b(t,X_t)\mathrm{d}W_t, \quad \text{ for all } t \in (0, T], \quad X_0 \ \text{given},
	\end{equation} 
	where $a\colon[0,\infty)\times\mathbb{R}^{d}\to\mathbb{R}^{d}$ 
	and $b\colon[0,\infty)\times\mathbb{R}^{d}\to\mathbb{R}^{d\times d}$ are suitable drift and diffusion operators, the dynamic equations governing the evolution of the triple $(m^{\mathrm{DLRA}}_t, U_t, Y_t)$, i.e.\ for the mean, the deterministic and stochastic bases, are given by~\cite{kazashi2025dynamical,sapsis2009dynamically}:
	\begin{equation}\label{3-terms DLR-conditions}
		\begin{aligned}
			\dot{m}^{\mathrm{DLRA}}_t & = \mathbb{E}[a(t,U_{t}^{\top}Y_t)]\\
			C_{Y_t}\dot{U}_t & = \mathbb{E}[Y_t \left(a(t,U_{t}^{\top}Y_t)-\mathbb{E}[a(t,U_{t}^{\top}Y_t)]\right)^{\top}]\left(I_{d \times d} - P_{U_t} \right), \\
			\mathrm{d}Y_t  & = U_t \left(a(t,U_{t}^{\top}Y_t)-\mathbb{E}[a(t,U_{t}^{\top}Y_t)]\right) \mathrm{d}t + U_t b(t,U_{t}^{\top}Y_t)\mathrm{d}W_t,
		\end{aligned}
	\end{equation}
	where $C_{Y_t}:= \mathbb{E}[Y_t Y_t^{\top}]$ is the Gram matrix for $Y_t$ and $P_{U_t}=U_t^{\top}U_t$ is the Euclidean orthogonal projector onto the rows of $U_t$.
	From \eqref{3-terms DLR-conditions} one can notice that the evolution of $m_t$, $U_t$, and $Y_t$ is highly coupled, depending also on the law of the process itself via the expectation present in the equation of the mean, the deterministic basis, as well as in the Gram matrix $C_{Y_t}$. These issues make the well-posedness and numerical analysis of DLRA not trivial. A study on well-posedness of DO equations under different conditions of $a$, $b$ can be found in \cite{bao2026exponential,kazashi2025dynamical,kazashi2026locally}, whereas the analysis of time and stochastic discretization schemes can be found in \cite{kazashi2025dynamicalpartI} and \cite{kazashi2026dynamicalpartII}, respectively·
	
	In the context of this work, we say that a process $Z \in L^2(\Omega,\mathbb{R}^d)$ is of rank $k$ if $\mathrm{rank}(\mathbb{E}[ZZ^{\top}])=k$ (for more details see e.g.\ \cite[Section 2.2]{kazashi2025dynamical}). Furthermore, let us denote with the superscript $\mathring{}$ a centered random variable, i.e.\ $\mathring{Z} = Z - \mathbb{E}[Z]$.  Then the DLR approximation $X^{\mathrm{DLRA}}$ of \eqref{eq: DA SDE} is a rank-$k$ process and satisfies the following McKean-Vlasov type equation
	\begin{equation}\label{eq:eq-manifold DO 3 terms} 
		\begin{aligned}
			\mathrm{d}X^{\mathrm{DLRA}}_{t} =& \Bigg(\mathbb{E}\left[a(t, X^{\mathrm{DLRA}}_t)\right]+ \bigl(\bigl(I_{d\times d}-\mathcal{P}_{\mathcal{U}( X^{\mathrm{DLRA}}_t)}\bigr)\left[\mathcal{P}_{\mathcal{Y}(X^{\mathrm{DLRA}}_t)}\left[\mathring{a}(t, X^{\mathrm{DLRA}}_t)\right]\right]\\
			&+\mathcal{P}_{\mathcal{U}(X^{\mathrm{DLRA}}_t)}\left[\mathring{a}(t, X^{\mathrm{DLRA}}_t)\right]\Bigg) \mathrm{d}t+\mathcal{P}_{\mathcal{U}(X^{\mathrm{DLRA}}_t)}b(t, X^{\mathrm{DLRA}}_t)\,\mathrm{d}W_{t},
		\end{aligned}
	\end{equation}
	where $\mathcal{P}_{\mathcal{U}(X^{\mathrm{DLRA}}_{t})} v= U_{t}^{\top}U_{t} v =P_{U_{t}}v$, $v \in \mathbb{R}^d$ and $\mathcal{P}_{\mathcal{Y}(X^{\mathrm{DLRA}}_{t})}w = \mathbb{E}[w Y^{\top}_{t}]C^{-1}_{Y_{t}}Y_{t}=P_{Y_{t}} w$, $w \in L^2(\Omega)$, are the projectors onto the range of the covariance of $X^{\mathrm{DLRA}}_{t}$, i.e.\ $\mathrm{Im}(\mathbb{E}[\mathring{X}_t^{\mathrm{DLRA}}\,\cdot\,]) = \mathrm{span}\{(U^{1}_t)^{\top}, \dots, (U^{k}_t)^{\top}\}$, and the corange $\mathrm{Im}\bigl(((\mathring{X}_t^{\mathrm{DLRA}})^{\top}\cdot\,)\bigr) = \mathrm{span}\{Y^1_{t}, \dots, Y^k_{t}\}$, respectively \cite{kazashi2025dynamical}. From \eqref{eq:eq-manifold DO 3 terms}, we define the orthogonal projector onto the tangent space of rank $k$ processes at the point $X_t^{\mathrm{DLRA}}$ as 
	\begin{equation}\label{eq: proj tang X}
		P_{X_t^{\mathrm{DLRA}}}[\cdot] :=\bigl(I_{d\times d}-\mathcal{P}_{\mathcal{U}(X_t^{\mathrm{DLRA}})}\bigr)\mathcal{P}_{\mathcal{Y}(X_t^{\mathrm{DLRA}})}[ \, \cdot \, ]+\mathcal{P}_{\mathcal{U}(X_t^{\mathrm{DLRA}})} [\, \cdot \,],
	\end{equation}
	so that \eqref{eq:eq-manifold DO 3 terms} can be written more compactly as  
	\begin{equation*}
		\mathrm{d}X^{\mathrm{DLRA}}_{t} = \left(\mathbb{E}\left[a(t, X^{\mathrm{DLRA}}_t)\right]+ P_{X^{\mathrm{DLRA}}_t}\left[\mathring{a}(t, X^{\mathrm{DLRA}}_t)\right]\right) \mathrm{d}t+\mathcal{P}_{\mathcal{U}(X^{\mathrm{DLRA}}_t)}b(t, X^{\mathrm{DLRA}}_t)\,\mathrm{d}W_{t}.
	\end{equation*}
	
	\subsection{Notation}
	For a vector $x \in \mathbb{R}^{d}$, its Euclidean norm is denoted by $|x|$. For a matrix $A\in \mathbb{R}^{m \times n}$, $| A |$ and $\| A \|_{\mathrm{F}}$ denote its spectral and Frobenius norm, respectively. 
	The symbol $\operatorname{Tr}$ indicates the trace of a squared matrix, i.e.\ $\operatorname{Tr}(A) = \sum_{i=1}^{m} A_{ii}$ for $A \in \mathbb{R}^{m \times m}$. We cite hereafter some useful properties of the trace operator that we will use in this work. For all matrices $A,X,B$ with suitable dimensions one has (see e.g.\ \cite{petersen2008matrix}):
	\begin{equation}\label{eq: matrix der}
		\begin{aligned}
			\frac{\partial}{\partial X} \operatorname{Tr}(XA) = A^{\top},\quad
			\frac{\partial}{\partial X} \operatorname{Tr}(AXB) &= A^{\top} B^{\top}, \quad
			\frac{\partial}{\partial X} \operatorname{Tr}(AX^{\top}B) = BA,  \\
			\frac{\partial}{\partial X} \operatorname{Tr}(X^{\top}A) &= A, \quad 
			\frac{\partial}{\partial X} \operatorname{Tr}(AX^{\top}) = A,
		\end{aligned}
	\end{equation}
	where the derivative of a scalar function $X \mapsto f(X) \in \mathbb{R}$ with respect to a matrix $X$ of general dimension, $\frac{\partial f}{\partial x_{ij}}$, can be identified with a matrix of the same dimension as $X$, whereas the derivative of a matrix valued function $X \mapsto A(X)$ of general dimension with respect to the input matrix $X$, $\frac{\partial A_{kl}}{\partial X_{ij}}$, can be identified with a $4$-th order tensor.
	Given two $(m \times n)$-matrices $A,B$, the Frobenius scalar product is defined as
	$$\left\langle A, B
	\right\rangle_{\mathrm{F}} := \operatorname{Tr}(A^{\top}B) =  \operatorname{Tr}(B^{\top}A).$$
	
	We will consider discretizing all the stochastic quantities via a Monte Carlo method. For instance, the expectation $\mathbb{E}$ will be approximated by the Monte Carlo average $\overline{\mathbb{E}}$ with the following properties: we denote $\overline{\mathbb{E}}[ab^{\top}]= \frac{ \sum_{i=1}^{M} a^{i}(b^{i})^{\top} }{M}$ for any real matrices $a =[a^i]_{i=1,\dots,M}\in \mathbb{R}^{j \times M}$, $b =[b^i]_{i=1,\dots,M}\in \mathbb{R}^{\ell \times M}$, where $a^i \in \mathbb{R}^{j}$, $b^i \in \mathbb{R}^{\ell}$ are columns of $a$ and $b$, respectively, for positive integers $j,\ell$. 
	
	We define the Stiefel matrix manifold $\mathrm{St}(k, d)$ as the set of all $(k \times d)$-dimensional matrices whose rows are orthonormal, i.e.\
	$$\mathrm{St}(k, d):=\left\{F \in \mathbb{R}^{k \times
		d} \ : \ F F^{\top} = I_{k \times k}\right\}.$$
	
	An element $F \in \mathrm{St}(k, d)$ identifies a $k$-dimensional subspace of $\mathbb{R}^d$; actually infinitely many element of $\mathrm{St}(k, d)$ identify the same subspace of $\mathbb{R}^d$.
	
	\section{DLRA-JMCO filter}\label{sec: DLRA for JMCO}
	In this section, we derive a DLRA-filter that minimizes locally in time the joint error on mean and covariance between the DLRA surrogate and the filtered full-order process \eqref{eq: kalman-bucy nonlinear} when started from the DLRA surrogate and advanced for a small time interval. Namely, this error involves the error between the mean of the surrogate and the one of the full filtered equation measured in Euclidean norm plus the error between their covariances measured in Frobenius norm. In detail, given two processes $\mathcal{X}, \mathcal{Z}$ with respective mean $m_x$, $m_z$, and covariances $C_x$, $C_z$, we consider the following discrepancy between $\mathcal{X}$ and $\mathcal{Z}$:
	\begin{equation}\label{eq: JMCO}
		\sqrt{ | m_x - m_z|^2 + \| C_x - C_z \|_{\mathrm{F}}^2 },
	\end{equation}
	which quantifies the mismatch on the mean and covariance between $\mathcal{X}$ and $\mathcal{Z}$. 
	
	First, we describe the general idea behind our derivation of the DLRA-filtered equations, which was inspired by the derivation of DLRA for matrix ODEs in \cite{koch2007dynamical} (and exploited for SDEs in \cite{kazashi2026approaches}), postponing its detailed construction to the remaining of this section. 
	
	In \cite{koch2007dynamical}, the DLRA for matrix ODEs was derived in the following fashion. A time-dependent surrogate $X^{\mathrm{DLRA}}$ was sought via minimizing the error between the derivative of the surrogate and the right-hand side of the equation computed in the point $X^{\mathrm{DLRA}}_t$, keeping $X^{\mathrm{DLRA}}$ constrained to the manifold of rank-$k$ matrices. In that context, differentiability of the surrogate was crucial in order to express its time derivative approximation as a linear combination of components lying in the tangent space of the manifold of rank-$k$ matrices at the point $X^{\mathrm{DLRA}}_t$. 
	
	In the context of Itô SDEs, the differentiability property is not satisfied and this fact brings an additional difficulty in deriving suitable DLRA equations.
    To overcome the difficulty of non-differentiability of Itô processes, we minimize \eqref{eq: JMCO} for a time-discretized version of the filtered equations and then recover time continuous filtered DLRA equations by formally taking the limit $\Delta t \to 0$. We mention that the work \cite{nobile2025dynamicallowrankapproximationskalman} follows a different construction, closer to \cite{kazashi2025dynamical}, based on a consistency argument and representation of Itô processes to obtain DLRA equations directly projecting \eqref{eq: kalman-bucy nonlinear}. However, their obtained DLRA filter does not have the minimization property \eqref{eq: JMCO} and differs from the one we obtain here via time discretization.
	
	In our derivation, we consider a reasonable time discretization of the DLRA and of the full-order filtered model \eqref{eq: kalman-bucy nonlinear}.  With an abuse of notation, we denote the filtered DLRA at time mesh $t_n$ by $X_n^{\mathrm{DLRA}}$ and suppose this process has already assimilated all the observations up to $t_n$, i.e.\ $\{Z_j\}_{0\leq j \leq n}$. Then, suppose that at time $t_n$ the full-order filtered solution $X_n^{f}$ equals the DLRA surrogate, i.e.\ $X_n^{f}=X_{n}^{\mathrm{DLRA}}$. Then, the reasonable updating equations that produce $X_{n+1}^{\mathrm{DLRA}}$, i.e.\ the DLRA at time $t_{n+1}$, are derived minimizing the error \eqref{eq: JMCO} between $X_{n+1}^{\mathrm{DLRA}}$ and the full-order filtered solution $X_{n+1}^{f}$ after one time step. The complete procedure is explained in detail hereafter.

	\subsection{Setting for time-discrete DLRA}
    Consider a time-grid
	$\Delta = \{ t_i,\; i \in \{1,\dots,N\} : t_0 = 0,\; t_N = T,\; \text{and}\;
	t_i - t_{i-1} = \Delta t_n,\; i \in \{1,\dots,N\} \}$.
	From now on, we consider a time discrete DLR process $\{X_n^{\mathrm{DLRA}}\}_{n=0}^N$ of rank $k$, whose distribution approximates the filtering distribution of the full model \eqref{eq: DA SDE}. Such DLRA is composed by 3 terms and reads $X_n^{\mathrm{DLRA}} = m_{n}^{\mathrm{DLRA}} + U_{n}^{\top} Y_{n}$, whose components are
	\begin{itemize}
		\item the mean $m^{\mathrm{DLRA}}_n= \mathbb{E}[X_{n}^{\mathrm{DLRA}}] \in \mathbb{R}^{d}$, which is a deterministic vector;
		\item the deterministic basis $U_n \in \mathbb{R}^{k \times d}$, which has orthonormal rows, i.e.\ $U_n U_n^\top = I_{k \times k}$ for all $n$;
		\item the stochastic basis $Y_n \in L^2(\Omega,  \mathbb{R}^k)$, adapted to the filtration $\mathcal{F}_{t_n}$, having zero mean, and with linearly independent components $Y_n^1, \dots, Y_n^k$ in $L^2(\Omega,  \mathbb{R})$ for all $n$.
	\end{itemize}
	Given the triple $(m^{\mathrm{DLRA}}_n, U_n, Y_n)$ at time $t_n$, we seek reasonable updates $(\Delta m^{\mathrm{DLRA}}_n, \Delta U_n, \Delta Y_n)$ such that the discrete DLRA at time $t_{n+1}$,
	\begin{equation}\label{eq: DLRA 3 terms}
	X_{n+1}^{\mathrm{DLRA}} = m_{n+1}^{\mathrm{DLRA}} + U_{n+1}^{\top} Y_{n+1},
	\end{equation}
	approximating the filtered distribution at time $t_{n+1}$, is characterized by the updates
	\begin{align}\label{eq: DLRA updates JMCO}
		m_{n+1}^{\mathrm{DLRA}} &= m_n^{\mathrm{DLRA}} + \Delta m_n^{\mathrm{DLRA}}, \\
		U_{n+1} &= U_n + \Delta U_n,
		\qquad \Delta U_n U_n^\top = 0, \label{eq: discrete gauge} \\
		Y_{n+1} &= Y_n + \Delta Y_n, \label{eq: DLRA updates JMCO fin}
		\qquad \Delta Y_n \text{ adapted},
	\end{align}
	where the increments satisfy the following order conditions 
	\begin{equation}\label{eq: updates}
		\begin{aligned}
	\Delta m_n &= c_n \Delta t + d_n \Delta Z_n, \quad	\Delta U_n = e_n \Delta t + f_n \Delta Z_n, \\
	\Delta Y_n &= g_n \Delta t + h_n \Delta \widehat{W}_n + i_n \Delta Z_n +  l_n \Delta \widehat{B}_n
	\end{aligned}
	\end{equation}
	for suitable $c_n \in \mathbb{R}^d$, $d_n  \in \mathbb{R}^{ d \times h}$, $e_n \in \mathbb{R}^{ d \times k}$, $f_n  \in \mathbb{R}^{ d \times k \times  h}$, $g_n \in L^2(\Omega, \mathbb{R}^k)$, $h_n \in L^2(\Omega, \mathbb{R}^{k \times m})$, $i_n  \in L^2(\Omega, \mathbb{R}^{d \times h})$, $l_n \in L^2(\Omega, \mathbb{R}^{k \times h})$, where $\Delta \widehat{W}_n$ and $\Delta \widehat{B}_n$ are independent Brownian increments. The ansatz \eqref{eq: updates} is consistent with the derivation of the DLRA in \cite{kazashi2025dynamical} and filtered DLRA in \cite{nobile2025dynamicallowrankapproximationskalman}. 
	
	Notice that condition \eqref{eq: discrete gauge} resembles the usual gauge condition for the continuous DLRA \cite{kazashi2025dynamical,koch2007dynamical}, and it is enforced in order to ensure a unique representation of the deterministic basis (and, hence, of the DLRA) in the limit of $\Delta t_n \to 0$, choosing one element in the equivalence class of elements of $\mathrm{St}(k,d)$ identifying the same $k$-dimensional subspace of $\mathbb{R}^d$. In this regard, we refer to \eqref{eq: discrete gauge} as the \textit{discrete gauge condition}. Once we succeed in obtaining reasonable discrete equations for the triple $(\Delta m^{\mathrm{DLRA}}_n, \Delta U_n, \Delta Y_n)$, we formally pass to the limit $\Delta t_n \to 0$ to gather reasonable DLRA differential equations describing a continuous-time reduced order model for \eqref{eq: kalman-bucy nonlinear}. Then, \eqref{eq: discrete gauge} transforms into the usual (continuous-time) gauge condition $\dot{U}_t U_t^{\top}=0$. %For the sake of notation, we denote the centered DLRA as 
	%$$\mathring{X}_n^{\mathrm{DLRA}}= X_n^{\mathrm{DLRA}} - \mathbb{E}[X_n^{\mathrm{DLRA}}] = X_n^{\mathrm{DLRA}} - m_n^{\mathrm{DLRA}} = U_n^{\top} Y_n.$$
	
	\subsection{Discretization of the filtering system associated to \eqref{eq: DA SDE}}

	We consider an Euler-Maruyama discretization of the full-order filtered process \eqref{eq: kalman-bucy nonlinear} in the interval $[t_n,t_{n+1}]$, starting at $X_n^{f} = X_n^{\mathrm{DLRA}}$. Then $X_{n+1}^{f}$ reads
	\begin{equation}
		\label{eq: EM Kalman Nonlinear}
		\begin{aligned}
		X_{n+1}^{f}= &X_n^{\mathrm{DLRA}}+ A(t_n, X_n^{\mathrm{DLRA}})\Delta t_n + Q^{\frac{1}{2}} \Delta \widehat{W}_n \\
		&+ C_n^{\mathrm{DLRA}} H(t_n)R^{-1} \left( \Delta Z_n - H(t_n)X_n^{\mathrm{DLRA}}-R^{\frac{1}{2}} \Delta \widehat{B}_n \right)
		\end{aligned}
	\end{equation}
where $\Delta Z_{n} = Z(t_{n+1})-Z(t_{n})$, $\Delta \widehat{W}_n = \widehat{W}(t_{n+1})-\widehat{W}(t_{n}) \sim \mathcal{N}(0, \Delta t_n I_{d \times d})$ and $\Delta \widehat{B}_n = \widehat{B}(t_{n+1})-\widehat{B}(t_{n}) \sim \mathcal{N}(0, \Delta t_n I_{h \times h})$ with $\Delta \widehat{W}_n$, $\Delta \widehat{B}_n$ independent of each other and also independent of $X_n^{\mathrm{DLRA}}$. Then, the filtered mean $m_{n+1} = \mathbb{E}[X_{n+1}^{f} \mid \{Z_\ell\}_{\ell=0,\dots,n}, \Delta Z_n]$ reads 
\begin{equation}\label{eq: EM cond mean - JMCO}
	\begin{aligned}
		m_{n+1}=& m_n^{\mathrm{DLRA}}
		+ \mathbb{E}\!\left[A(t_n,X_n^{\mathrm{DLRA}})\right]\Delta t_n \\
		&+ C_{n}^{\mathrm{DLRA}}\,H(t_n)^{\top}\,R^{-1}
		\left( \Delta Z_{n}  - H(t_n)\, m_n^{\mathrm{DLRA}} \Delta t_n \right) 
	\end{aligned}
\end{equation}
whereas, discarding terms of order higher than $O(\Delta t_n)$, the filtered covariance $C_{n+1} = \mathbb{E}\left[\left(X_{n+1}^{f}-m_{n+1}\right) \left(X_{n+1}^{f}-m_{n+1}\right)^{\top}\right]$ reads
 \begin{equation}\label{eq: EM cond cov - JMCO}
 	\begin{aligned}
	C_{n+1}
	= & 
	C_n^{\mathrm{DLRA}}
	+ \mathbb{E}\!\left[\mathring{X}_n^{\mathrm{DLRA}}\mathring{A}(t_n,X_n^{\mathrm{DLRA}})^{\!\top}\right]\Delta t_n + \mathbb{E}\!\left[ \mathring{A}(t_n,X_n^{\mathrm{DLRA}})(\mathring{X}_n^{\mathrm{DLRA}})^{\!\top}\right]\Delta t_n
	+  Q\,\Delta t_n \\
	&-  
	C_n^{\mathrm{DLRA}}  H(t_n)^{\!\top} R^{-1} H(t_n)
	C_n^{\mathrm{DLRA}}\Delta t_n.
\end{aligned}
\end{equation}

	\subsection{Finding DLRA updates via a minimization problem}
	As anticipated earlier, to derive reasonable discrete DLRA updates, we follow a strategy akin to the continuous-time derivation of DLRA equations in \cite{koch2007dynamical}: we minimize the RMSE  
		\begin{equation}\label{eq: RMSE DLRA}
		\left| m_{n+1} - m_{n+1}^{\mathrm{DLRA}} \right|^2
		+ \left\| C_{n+1} - C_{n+1}^{\mathrm{DLRA}} \right\|_{\mathrm{F}}^2
	\end{equation}	
	computed between $(m_{n+1}, C_{n+1})$ from \eqref{eq: EM cond mean - JMCO}-\eqref{eq: EM cond cov - JMCO} and the respective DLR approximation $(m_{n+1}^{\mathrm{DLRA}}, C_{n+1}^{\mathrm{DLRA}})$ from \eqref{eq: DLRA updates JMCO}-\eqref{eq: DLRA updates JMCO fin}, where $C_{n+1}^{\mathrm{DLRA}}$ at time $t_{n+1}$ is given by
	\begin{equation}\label{eq: C_{n+1}}
		\begin{aligned}
	C_{n+1}^{\mathrm{DLRA}}
	= & \mathbb{E}\!\left[ (U_{n+1}^\top Y_{n+1})(U_{n+1}^\top Y_{n+1})^\top \mid  \Delta Z_{n}, \{Z_\ell\}_{\ell=0,\dots,n} \right] \\
	= & U_{n+1}^\top\, \mathbb{E}\!\left[ Y_{n+1} Y_{n+1}^\top \mid  \Delta Z_{n}, \{Z_\ell\}_{\ell=0,\dots,n}\right] U_{n+1}
	=  (U_n + \Delta U_n)^\top\, C_{Y_{n+1}}\, (U_n + \Delta U_n),
	\end{aligned}
	\end{equation}
	where  
	\begin{equation}\label{eq: C_y}
	C_{Y_{n+1}} = \mathbb{E}\!\left[ (Y_n + \Delta Y_n)(Y_n + \Delta Y_n)^\top \mid  \Delta Z_{n}, \{Z_\ell\}_{\ell=0,\dots,n} \right]
	\end{equation}
	is the Gramian of the stochastic basis $Y_{n+1}$ at time $t_{n+1}$. For the sake of computation, we consider that relation \eqref{eq: C_y} is obtained by the covariance at the previous step $t_n$ plus an increment, namely
	\begin{equation}\label{eq: C_y 2}
	C_{Y_{n+1}} = C_{Y_n} + \Delta C_{Y_n},
	\end{equation}
	with $\Delta C_{Y_n} \in \mathbb{R}^{k \times k}$.

	We then look for the triple $(\Delta m_n^{\mathrm{DLRA}}, \Delta U_n, \Delta C_{Y_n})$
	that minimizes
	\begin{equation}\label{eq: prob DLRA JMCO}
		\begin{aligned}
	\min_{ \substack{\Delta m_n^{\mathrm{DLRA}},
			\Delta U_n , \\
			\Delta C_{Y_n},
	\text{ satisfying \eqref{eq: discrete gauge}-\eqref{eq: updates}} }}
	\left| m_{n+1} - (m_n^{\mathrm{DLRA}} + \Delta m_n^{\mathrm{DLRA}}) \right|^2 
 + \left\| C_{n+1} - (U_n + \Delta U_n)^\top (C_{Y_{n}} + \Delta C_{Y_n}) (U_n + \Delta U_n) \right\|_{\mathrm{F}}^2.
			\end{aligned}
	\end{equation}
	and then construct $\Delta Y_n$ so that $X_{n+1}^{\mathrm{DLRA}}$ approximates pathwise $X_{n+1}^{f}$ from \eqref{eq: EM Kalman Nonlinear}.
	
	To simplify the calculation and since we are ultimately interested in the continuous limit as $\Delta t_n \to 0$, we simplify the expression of $C_{n+1}^{\mathrm{DLRA}}$ by retaining only first order terms in $\Delta t_n$:
	\begin{equation}\label{eq: DLRA C mute}
		\begin{aligned}
	C^{\mathrm{DLRA}}_{n+1}  = & (U_n + \Delta U_n)^\top\, \left(C_{Y_{n}} + \Delta C_{Y_{n}}\right)\, (U_n + \Delta U_n) \\
	=& C_n^{\mathrm{DLRA}}
	+  U_n^{\!\top} C_{Y_n} \Delta U_n
	+ U_n^{\!\top} \Delta C_{Y_n} U_n
	+ U_n^{\!\top} \Delta C_{Y_n} \Delta U_n \\
	&+ \Delta U_n^{\!\top} C_{Y_n} U_n
	+ \Delta U_n^{\!\top} C_{Y_n}  \Delta U_n
	+ \Delta U_n^{\!\top} \Delta C_{Y_n} U_n +  \Delta U_n^{\!\top} \Delta C_{Y_n} \Delta U_n \\
	\approx &C_n^{\mathrm{DLRA}}
	+ U_n^{\!\top} C_{Y_n} \Delta U_n
	+ U_n^{\!\top} \Delta C_{Y_n} U_n
	+  \Delta U_n^{\!\top} C_{Y_n} U_n +  \Delta U_n^{\!\top} C_{Y_n}  \Delta U_n.
	\end{aligned}
	\end{equation}
	
 	We are now ready to derive the equations for $m^{\mathrm{DLRA}}_{n}$ and $C^{\mathrm{DLRA}}_{n}$ that satisfy relation \eqref{eq: prob DLRA JMCO} as a first order approximation in $O(\Delta t_n)$.
	
	\begin{Proposition}[DLRA-JMCO filter]\label{prop: triple JMCO}
	The triple $(\Delta m_{n}^{\mathrm{DLRA}} ,\, \Delta U_n,\, \Delta Y_n)$ satisfying
	\begin{equation}\label{eq: JMCO triple}
		\begin{aligned}
	\Delta m_{n}^{\mathrm{DLRA}}  :=& \mathbb{E}\big[ A(t_n, X_n^{\mathrm{DLRA}}) \big]\Delta t_n
	+ U_n^{\top}C_{Y_n}U_n H(t_n)^\top R^{-1} (\Delta Z_n - H(t_n)\mathbb{E}[X_n^{\mathrm{DLRA}}]\Delta t_n) \\
	\Delta U_n :=& C_{Y_n}^{-1}\, \mathbb{E}\!\left[ Y_n\, \mathring{A}(t_n, X_n^{\mathrm{DLRA}})^\top \right] P_{U_n}^\perp \Delta t_n
	+ C_{Y_n}^{-1} U_n Q\, P_{U_n}^\perp \Delta t_n \\
	\Delta Y_n :=& U_n\, \mathring{A}(t_n, X_n^{\mathrm{DLRA}})\Delta t_n
	+ U_n Q^{\frac{1}{2}}\Delta \widehat{W}_n
	-\, C_{Y_n} U_n H(t_n)^{\top} R^{-1} H(t_n) U_n^\top Y_n\, \Delta t_n \\
	&- C_{Y_n} U_n H(t_n)^{\top} R^{-\frac{1}{2}}\Delta \widehat{B}_n,
	\end{aligned}
	\end{equation}
    with $C_{Y_n} = \mathbb{E}[Y_nY_n^{\top}]$, $\Delta \widehat{W}_n$, $\Delta \widehat{B}_n$ Brownian increments independent of $X_n^{\mathrm{DLRA}}$, and $ P_{U_n}= U_n^{\top}U_n$ the projector onto the rows of $U_n$,
	solves \eqref{eq: prob DLRA JMCO} at first–order in $\Delta t_n$.
	\begin{proof}
		The proof can be found in Appendix \ref{app: proof DLRA JMCO}.
	\end{proof}
	\end{Proposition}
	
	\begin{Remark}[On the dependence on $\Delta Z_n$ for \eqref{eq: JMCO triple}]
		One can notice that in \eqref{eq: JMCO triple} only the mean depends on the observation $\Delta Z_n$. This fact is not surprising: indeed, by taking a Kalman approximation of the mean and covariance of the filtered process associated to \eqref{eq: DA SDE}, we know a priori that only the mean $m_{n+1}$ depends on $\Delta Z_n$, (see equations \eqref{eq: kalman-bucy})  which does not enter in the minimization of the covariance when we seek equations for $\Delta U_n$ and $\Delta C_{Y_n}$.
	\end{Remark}

\subsection{Formal continuous limit of DLRA-JMCO filter}
Proposition \ref{prop: triple JMCO} suggests the following discrete DLRA equations to compute the surrogate at time $t_{n+1}$ given the previous solution $(m_{n}^{\mathrm{DLRA}}, U_n, Y_n)$ at time $t_n$
	\begin{equation}\label{eq: DLRA JMCO}
		\begin{aligned}
	m_{n+1}^{\mathrm{DLRA}}
	= &m_n^{\mathrm{DLRA}}
	+ \mathbb{E}[A(t_n, X_n^{\mathrm{DLRA}})]\, \Delta t_n
	+ C_n^{\mathrm{DLRA}} H(t_n)^\top R^{-1} (Z_{n+1} - Z_{n} - H(t_n) m_n^{\mathrm{DLRA}}\, \Delta t_n),\\
	U_{n+1}
	= &U_n
	+ C_{Y_n}^{-1} \mathbb{E}[Y_n\, \mathring{A}(t_n, X_n^{\mathrm{DLRA}})^\top]
	\, P_{U_n}^\perp\, \Delta t_n
	+ C_{Y_n}^{-1} U_n Q\, P_{U_n}^\perp\, \Delta t_n. \\
	Y_{n+1}
	= &Y_n
	+ U_n \mathring{A}(t_n, X_n^{\mathrm{DLRA}})\, \Delta t_n
	+ U_n Q^{\frac{1}{2}}\, \Delta \widehat{W}_n \\
	&- C_{Y_n}\, U_n H(t_n)^\top R^{-1} H(t_n) U_n^\top Y_n\, \Delta t_n
	- C_{Y_n}\, U_n H(t_n)^\top\, R^{-\frac{1}{2}}\, \Delta \widehat{B}_n.
     \end{aligned}
	\end{equation}
	
	Taking the limit for the step-size $\Delta t_n$ going to $0$ in \eqref{eq: DLRA JMCO}, we retrieve the following differential equations for the continuous-time modes of the sought DLRA filtering procedure:
	\begin{align}\label{eq: m_t JMCO}
		\mathrm{d} m_t^{\mathrm{DLRA}}
		=& \mathbb{E}\!\left[ A\!\left(t, X_t^{\mathrm{DLRA}}\right) \right] \mathrm{d}t
		+U_t^{\top} C_{Y_t}U_t\, H(t)^\top R^{-1}
		\bigl( \mathrm{d}Z_t - H(t)\, m_t^{\mathrm{DLRA}}\, \mathrm{d}t \bigr),\\
		\mathrm{d} U_t
		=& C_{Y_t}^{-1}
		\Bigl(
		\mathbb{E}\!\left[ Y_t\, \mathring{A}(t, X_t^{\mathrm{DLRA}})^{\top} 
		  \right]
		+ U_t Q
		\Bigr)
		P_{U_t}^{\perp}\, \mathrm{d}t, \label{eq: U_t JMCO}
		\\
		\mathrm{d} Y_t
		=&
		\Bigl(
		U_t \mathring{A}(t, X_t^{\mathrm{DLRA}})
		- U_t U_t^{\top} C_{Y_t}\, U_t H(t)^\top R^{-1} H(t)\, U_t^\top Y_t
		\Bigr)\mathrm{d}t \nonumber
		\\
		&
		+ U_t Q^{\frac{1}{2}}\, \mathrm{d}\widehat{W}_t
		- U_t U_t ^{\top}C_{Y_t}\, U_t H(t)^\top R^{-\frac{1}{2}}\, \mathrm{d}\widehat{B}_t \nonumber
		\\
		=&
		U_t\left(\mathring{A}(t, X_t^{\mathrm{DLRA}})
		- U_t ^{\top} C_{Y_t}U_t H(t)^\top R^{-1} H(t) U_t^\top Y_t\right) \mathrm{d}t \label{eq: Y_t JMCO}
		\\
		&\qquad
		+ U_t Q^{\frac{1}{2}}\, \mathrm{d}\widehat{W}_t
		- U_t (U_t ^{\top} C_{Y_t}U_t ) H(t)^\top R^{-\frac{1}{2}}\, \mathrm{d}\widehat{B}_t \nonumber,
	\end{align}	
	where $\widehat{W}_t$ and $\widehat{B}_t$ are independent Brownian motions.
	Then, using It\^o's formula, $X^{\mathrm{DLRA}}$ satisfies the following equation:
	\begin{equation*}
	\begin{aligned}
		\mathrm{d} X_t^{\mathrm{DLRA}}
		&= \mathrm{d} m_t^{\mathrm{DLRA}} + \mathrm{d}(U_t^\top Y_t)
		\\
		&= \mathrm{d} m_t^{\mathrm{DLRA}}
		+ (\mathrm{d} U_t^\top)\, Y_t
		+ U_t^\top\, \mathrm{d} Y_t, \\
					\end{aligned}
			\end{equation*}
which translates into 
			\begin{equation}\label{eq: X_t JMCO}
			\begin{aligned}
		\mathrm{d}X_t^{\mathrm{DLRA}}
		&=
		\mathbb{E}\!\left[A\!\left(t, X_t^{\mathrm{DLRA}}\right)\right]\mathrm{d}t
		+ C_t^{\mathrm{DLRA}} H(t)^\top R^{-1}\bigl(\mathrm{d}Z_t - H(t)\mathbb{E}\!\left[ X_t^{\mathrm{DLRA}}\right] \mathrm{d}t\bigr)
		\\
		&\quad
		+ P_{U_t}^\perp P_{Y_t}
		\Bigl(
		\mathring{A}(t,X_t^{\mathrm{DLRA}})
		+ Q\, U_t^\top C_{Y_t}^{-1} Y_t
		\Bigr)\mathrm{d}t
		\\
		&\quad
		+ P_{U_t}
		\Bigl(
		\mathring{A}(t, X_t^{\mathrm{DLRA}})
		\Bigr)\mathrm{d}t
		\\
		&\quad
		- P_{U_t} C_t^{\mathrm{DLRA}} H(t)^\top R^{-1} H(t) U_t^\top Y_t \, \mathrm{d}t
		+ P_{U_t} Q^{\frac{1}{2}}\, \mathrm{d}\widehat{W}_t
		- P_{U_t} C_t^{\mathrm{DLRA}} H(t)^\top R^{-\frac{1}{2}} \, \mathrm{d}\widehat{B}_t
		\\
		&=
		\mathbb{E}\!\left[A\!\left(t, X_t^{\mathrm{DLRA}}\right)\right]\mathrm{d}t
		+ P_{X_t^{\mathrm{DLRA}}}
		\Bigl(
		\mathring{A}(t, X_t^{\mathrm{DLRA}})
		\Bigr)\mathrm{d}t
		\\
		&\quad
		+ P_{U_t}^\perp Q\, U_t^\top C_{Y_t}^{-1} Y_t \, \mathrm{d}t
		+ P_{U_t} Q^{\frac{1}{2}}\, \mathrm{d}\widehat{W}_t
		\\
		&\quad
		+ C_t^{\mathrm{DLRA}} H(t)^\top R^{-1}
		\Bigl( \mathrm{d}Z_t - H(t) X_t^{\mathrm{DLRA}} \mathrm{d} t - R^{\frac{1}{2}} \mathrm{d}\widehat{B}_t \Bigr),
   \end{aligned}
	\end{equation}
   	where in the last line we use the fact that
    $U_t U_t^\top = I_{k \times k}$, and, hence, $P_{U_t} C_t^{\mathrm{DLRA}} = C_t^{\mathrm{DLRA}}$, for all $t \geq 0$.

\subsection{The numerical algorithm}\label{sec: DLRA JMCO alg}
Given the way we derived the system \eqref{eq: m_t JMCO}-\eqref{eq: U_t JMCO}-\eqref{eq: Y_t JMCO}, it is reasonable to think that equations \eqref{eq: DLRA JMCO} represent a possible time-discrete approximation. However, \eqref{eq: DLRA JMCO} is a fully explicit procedure in a Kalman-type update, and it is well-known that this structure can present numerical instabilities \cite{blomker2018strongly}. In this regard, we propose a semi-implicit staggered discretization, instead.
      
	Notice that \eqref{eq: X_t JMCO} can be interpreted as a composition of two distinct mechanisms. The 
	first one (the first two lines in the last equality in \eqref{eq: X_t JMCO}) corresponds to a DLR approximation of the original SDE \eqref{eq:eq-manifold DO 3 terms}, with the extra term $P_{U_t}^\perp Q\, U_t^\top C_{Y_t}^{-1} Y_t \, \mathrm{d}t$ which still belongs to the tangent space to the manifold of rank $k$ process at the point $X_t^{\mathrm{DLRA}}$. This extra term has been already considered in \cite{cao2018stochastic} and \cite{kazashi2026approaches} and was obtained precisely by deriving the continuous DLRA equations as a formal limit of time-discretized low-rank processes, as opposed to the argument followed in \cite{kazashi2025dynamical}.  The second mechanism (the last line in \eqref{eq: X_t JMCO}) is a
	filtering correction. Therefore, the first can be thought as a prediction update, whereas the second can be considered as an analysis step, exactly as in a Kalman type filter.
	
    We propose to discretize the prediction step via the DLR Projector Splitting for SDEs, as it shows remarkable stability properties \cite{kazashi2025dynamicalpartI,kazashi2026dynamicalpartII}. This leads to the following scheme for the predicted quantities $(\widehat{m}_{n+1}^{\mathrm{DLRA}}, \widehat{U}_{n+1},\widehat{Y}_{n+1})$ with intermediate step $(\widetilde{U}_{n+1}, \widetilde{Y}_{n+1})$:
	\begin{equation}
		\begin{aligned}\label{eq: DLRA JMCO prediction}
			\widehat{m}_{n+1}^{\mathrm{DLRA}}
			=&  m_n^{\mathrm{DLRA}}
			+  \mathbb{E}\!\left[ A\!\left(t_n, X_n^{\mathrm{DLRA}}\right) \right]\Delta t_n \\
			\widetilde{Y}_{n+1}
			=&  Y_n
			+ U_n  \left(\mathring{A}(t_n, X_n^{\mathrm{DLRA}})\right) \Delta t_n + U_n Q^{\frac{1}{2}} \Delta \widehat{W}_n, \\
		\widetilde{U}_{n+1} 
			=& U_n+ C_{\widetilde{Y}_{n+1}}^{-1}
			\Bigl(
			\mathbb{E}\!\left[ \widetilde{Y}_{n+1}\, \left(\mathring{A}(t_n, X_n^{\mathrm{DLRA}}) \right)^{\top}\right]
			+ U_n Q
			\Bigr)
			P_{U_n}^{\perp}\, \Delta t_n
		\end{aligned}
	\end{equation}
after which we apply an orthonormalization procedure $\widehat{U}_{n+1}^{\top}\widehat{Y}_{n+1}=\widetilde{U}_{n+1}^{\top}\widetilde{Y}_{n+1}$ with $\widehat{U}_{n+1}\widehat{U}_{n+1}^{\top} = I_{k\times k}$ to recover an orthonormal deterministic basis.
	Next, we propose the following analysis step 
	\begin{equation}\label{eq: DLRA JMCO analysis}
		\begin{aligned}
			m_{n+1}^{\mathrm{DLRA}}
			= &\widehat{m}_{n+1}^{\mathrm{DLRA}}
			+ \widehat{U}_{n+1}^{\top}C_{\widehat{Y}_{n+1}}\widehat{U}_{n+1} H(t_n)^\top R^{-1} (Z_{n+1} - Z_n - H(t_n)m_{n+1}^{\mathrm{DLRA}} \Delta t_n)\\
			U_{n+1}
			= & \widehat{U}_{n+1}\\
			Y_{n+1}
			= &\widehat{Y}_{n+1}
	        - C_{\widehat{Y}_{n+1}}\, \widehat{U}_{n+1} H(t_n)^\top R^{-1} H(t_n) \widehat{U}_{n+1}^{\top} Y_{n+1}\, \Delta t_n
			- C_{\widehat{Y}_{n+1}}\, \widehat{U}_{n+1} H(t_n)^{\top}\, R^{-\frac{1}{2}}\, \Delta \widehat{B}_n,.
		\end{aligned}
	\end{equation}
	In \eqref{eq: DLRA JMCO analysis}, we use the just-computed prediction triple $(\widehat{m}_{n+1}, \widehat{U}_{n+1}, \widehat{Y}_{n+1})$ as well an implicit approximation of the Kalman-gain type terms to determine $m_{n+1}^{\mathrm{DLRA}}$ and $Y_{n+1}$. Notice that we do not need an additional orthonormalization step after the analysis update as the deterministic basis does not change. 

	Equations \eqref{eq: DLRA JMCO prediction}-\eqref{eq: DLRA JMCO analysis} could be used also in the case of discrete observations, in which the time interval between observations is larger than the time step needed to follow the dynamics. In this case, one will perform multiple time steps of \eqref{eq: DLRA JMCO prediction} before incurring in the assimilation step.
	\begin{Remark}[Subspace $U$ independent of the observations]
		Relation \eqref{eq: DLRA JMCO analysis} highlights that the evolution from time $t_n$ to time $t_{n+1}$ of the deterministic basis $U$ is independent of the observation process. This can be seen as a limitation of the currently developed filtering procedure.% and the complemented formulation proposed in Section \ref{sec: DLRA compl JMCO} addresses precisely this point.
	\end{Remark}
	
	To complete this section, we present a fully discrete version of \eqref{eq: DLRA JMCO prediction}-\eqref{eq: DLRA JMCO analysis}, in which all expectations are approximated by a Monte-Carlo estimator based on a particle ensemble which is propagated by the dynamics. We denote by $\overline{\mathbb{E}}$ the Monte Carlo estimator. The fully discrete algorithm is therefore a noisy interacting particle system. We call this Algorithm \emph{DLRA Joint-Mean-Covariance-Optimized Filter} (see Algorithm \ref{alg: DLRA JMCO Filter algorithm}).
	
	\begin{algorithm}[!h]
		\caption{DLRA-JMCO Filter (Continuous time observations)}\label{alg: DLRA JMCO Filter algorithm}
		\begin{flushleft}
			\textbf{Input}: initial data  ($m_0$, $U_0$, $Y_0$), ensemble size $M$, sequence of time steps $\{\Delta t_n\}_{n \in \{1, \dots, N\}}$, observations $\{Z_n\}_{n\in \{1, \dots, N\}}$
			
			\textbf{Output:} filtered approximation $\{X^{\mathrm{DLRA},i}_n= m_n^{\mathrm{DLRA}} + U_n^{\top}Y^{i}_n\}_{n=0,\ldots, N}$ for $i=1,\dots,M$. 
		\end{flushleft}
		\begin{algorithmic}[1]
			
			\ForAll {$n \in  \{0, \ldots, N-1\}$} 
			
			\State Generate $M$ Brownian increments $\Delta \widehat{W}_n^{i} \overset{\mathrm{i.i.d.}}{\sim} \mathcal{N}(0,\Delta t_n I_{m \times m})$, $\Delta \widehat{B}_n^{i} \overset{\mathrm{i.i.d.}}{\sim} \mathcal{N}(0,\Delta t_n I_{h \times h})$  with $i \in \{1,\dots,M\}.$
			
			\State Compute $\widehat{m}_{n+1}^{\mathrm{DLRA}}
			= m_n^{\mathrm{DLRA}}
			+ \left(\frac{1}{M} \sum_{i=1}^{M}A(t_n,  m_n^{\mathrm{DLRA}} + U_n^{\top}Y_n^{i})\right)\, \Delta t_n.$
			
			\State Compute $\widetilde{Y}^{i}_{n+1} = Y_n^{i} + U_n \mathring{A}\!\left(t_n, m_n^{\mathrm{DLRA}} + U_n^{\top}Y_n^{i}\right) \Delta t_n +  U_n Q^{\frac{1}{2}} \Delta \widehat{W}_n^{i}$, for all $i \in \{1,\dots,M\}.$
			
				\State Assemble $\overline{C}_{\widetilde{Y}_{n+1}} = \frac{1}{M-1} \sum_{i=1}^{M}\widetilde{Y}_{n+1}^{i}(\widetilde{Y}^{i}_{n+1})^{\top}$
			\State Compute $\widetilde{U}_{n+1}$ by solving
			$$\overline{C}_{\widetilde{Y}_{n+1}} \widetilde{U}_{n+1}
			= \overline{C}_{\widetilde{Y}_{n+1}}  U_n
			+\left( \frac{1}{M-1} \sum_{i=1}^{M}\widetilde{Y}_{n+1}^{i}\, \mathring{A}(t_n,  m_n^{\mathrm{DLRA}} + U_n^{\top}Y_n^{i})^\top \right)
			\, P_{U_n}^\perp\, \Delta t_n
			+ U_n Q\, P_{U_n}^\perp\, \Delta t_n. $$
			\State Reorthonormalize the deterministic basis: find   $(\widehat{U}_{n+1}, \widehat{Y}_{n+1})$ such that:
			\begin{equation*}
				 \widehat{U}^{\top}_{n+1} \widehat{Y}_{n+1}^{i}=\widetilde{U}_{n+1}^{\top} \widetilde{Y}_{n+1}^{i}, \ \forall i \in 1, \dots,M  \text{ and } \widehat{U}_{n+1}\widehat{U}^{\top}_{n+1} = I_{k\times k}.
			\end{equation*}
			via $(\widehat{U}^{\top}_{n+1}, R_{n+1}) = \texttt{QR}(\widetilde{U}^{\top}_{n+1})$ and $\widehat{Y}_{n+1}^{i} = R_{n+1} \widetilde{Y}^{i}_{n+1}$, for all $i \in \{1,\dots,M\}.$
			
			\State Recentering: $\widehat{Y}_{n+1}^{i} = \widehat{Y}_{n+1}^{i} - \frac{1}{M} \sum_{j= 1}^{M}  \widehat{Y}_{n+1}^{j}$  
			for all $i$ and $\widehat{m}_{n+1} = \widehat{m}_{n+1} + \widehat{U}_{n+1}^{\top} \left(\frac{1}{M} \sum_{j= 1}^{M}   \widehat{Y}_{n+1}^{j}\right)$
			
			\State With $\Delta Z_n = Z_{n+1}-Z_n$, set $U_{n+1} =\widehat{U}_{n+1}$ and solve for $m_{n+1}^{\mathrm{DLRA}}$ and $Y_{n+1}^{i}$
			\begin{equation*}
				\begin{aligned}
					&(I_{d \times d}+ \widehat{U}_{n+1}^{\top} K^H_n \Delta t_n) m_{n+1}^{\mathrm{DLRA}}
					= \widehat{m}_{n+1}^{\mathrm{DLRA}}
					+ \widehat{U}_{n+1}^{\top} \overline{C}_{\widehat{Y}_{n+1}} \widehat{U}_{n+1} H(t_n)^\top R^{-1}  \Delta Z_n,\\
					&(I_{k \times k} + K^H_n \widehat{U}^{\top}_{n+1} \Delta t_n )Y_{n+1}^{i}
					= \widehat{Y}_{n+1}^{i}
					- \overline{C}_{\widehat{Y}_{n+1}} \, \widehat{U}_{n+1}H(t_n)^{\top} R^{-\frac{1}{2}}\, \Delta \widehat{B}_n^{i}, i =1,\dots,M
				\end{aligned}
			\end{equation*}
		for $K^H_n := \overline{C}_{\widehat{Y}_{n+1}} \, \widehat{U}_{n+1} H(t_n)^\top R^{-1} H(t_n)$
			\State Recentering: $Y_{n+1}^{i} = Y_{n+1}^{i} - \frac{1}{M} \sum_{j= 1}^{M}  Y_{n+1}^{j}$  
			for all $i$ and $m_{n+1}^{\mathrm{DLRA}} = m_{n+1}^{\mathrm{DLRA}} + U_{n+1}^{\top} \left(\frac{1}{M} \sum_{j= 1}^{M}   Y_{n+1}^{j}\right)$
			\EndFor
		\end{algorithmic}
	\end{algorithm}
	
	Algorithm \ref{alg: DLRA JMCO Filter algorithm} has been stated for continuous-time observations. However, it can be easily modified to account for discrete observations arriving at time intervals of length $\Delta t_a > \Delta t_n$, as described below. In the case of observations acquired in discrete time, we suppose $\Delta t_n$ is chosen as an integer fraction of $\Delta t_a$ and the data are obtained as
	\begin{equation}\label{eq: discr obs cont}
		\Xi_{n+1} = H(t_{n+1})X_{n+1} + \eta_n,
	\end{equation}
	 whenever $t_{n+1}$ is an observation time, with $\eta_n \sim \mathcal{N}(0,\widetilde{R} )$. Notice that $\Xi_{n+1}$ corresponds to the increment $\Delta Z_{n} = Z(t_{n+1})- Z(t_{n})$ is the continuous time observation case $(\Delta t_a = \Delta t_n)$. Likewise, the noise covariance $\widetilde{R}$ should correspond to $R \Delta t_a$. 
	
	The prediction step performed at every time step $\Delta t_n$ by a step of the DLR Splitting for SDEs. If $\Delta t_n \mod \Delta t_a = 0$, the analysis step implies solving the following equations for $m_{n+1}^{\mathrm{DLRA}}$ and $Y_{n+1}^{i}$:
	\begin{equation*}
		\begin{aligned}
			(I_{d \times d}+ \widehat{U}_{n+1}^{\top}K^H_{n+1}) m_{n+1}^{\mathrm{DLRA}}
			= &\widehat{m}_{n+1}^{\mathrm{DLRA}}
			+ \widehat{U}_{n+1}^{\top} \overline{C}_{\widehat{Y}_{n+1}} \widehat{U}_{n+1} H(t_{n+1})^\top \widetilde{R}^{-1}  \Xi_{n+1},\\
			(I_{k \times k} + K^H_{n+1} \widehat{U}^{\top}_{n+1} )Y_{n+1}^{i}
			= &\widehat{Y}_{n+1}^{i}
			- \overline{C}_{\widehat{Y}_{n+1}} \, \widehat{U}_{n+1}H(t_{n+1})^{\top} \, \widetilde{R}^{-\frac{1}{2}}\,  \xi_i, 
		\end{aligned}
	\end{equation*}
	for all $i \in \{1,\dots,M\}$, with $K^H_{n+1} := \overline{C}_{\widehat{Y}_{n+1}} \, \widehat{U}_{n+1} H(t_{n+1})^\top \widetilde{R}^{-1} H(t_{n+1})$ and where $(\xi_i)_{i=1,\dots,M}$ are M i.i.d.\ standard normal random variables. The deterministic basis does not change in the analysis step, i.e.\: $U_{n+1} = \widehat{U}_{n+1}.$
	If, on the other, $\Delta t_n \mod \Delta t_a \neq 0$, then no assimilation occurs and we have
	\begin{equation*}
		\begin{aligned}
			m_{n+1}^{\mathrm{DLRA}}
			= \widehat{m}_{n+1}^{\mathrm{DLRA}}, \quad U_{n+1} = \widehat{U}_{n+1}, \quad 
			Y_{n+1}^{i}
			= \widehat{Y}_{n+1}^{i}, \quad \text{for all } i \in \{1,\dots,M\}.
		\end{aligned}
	\end{equation*}

	\begin{Remark}[Invertibility of $\overline{C}_{Y_n}$]
		In Algorithm \ref{alg: DLRA JMCO Filter algorithm}, the Gramian of $\overline{C}_{\widetilde{Y}_n}$ is implicitely assumed invertible at all $n$ and for all $M$. In order to guarantee invertibility of $\overline{C}_{\widetilde{Y}_n}$, one may consider regularization techniques depending on the sample size $M$. For more details, the reader may refer to \cite{kazashi2026dynamicalpartII}.
	\end{Remark}

	\subsection{The affine-drift case}
	
	We now specialize the continuous-time DLRA equation \eqref{eq: X_t JMCO}, derived from the DLRA-JMCO filter, to the affine-drift case, i.e.\
	\begin{equation*}
	A(t,x)=A(t)x+f(t).
	\end{equation*}
	Since the observation model is linear, if the initial (DLRA) state is Gaussian, the law of \(X_t^{\mathrm{DLRA}}\) remains Gaussian and is then characterized by the first two moments, i.e.\ its mean and covariance. In this case, the evolution of \(X_t^{\mathrm{DLRA}}\) can be replaced by evolution equations for these two moments, avoiding the need to evolve a noisy-particle system. This leads to a reduced-order Kalman--Bucy filter formulation for the continuous-time DLRA dynamics.
	
	In case of affine drift, the equation \eqref{eq: m_t JMCO} of the first moment of the DLRA-JMCO becomes
	\begin{equation}\label{eq: m KBF DLRA}
		\mathrm{d} m_t^{\mathrm{DLRA}} = \bigl(A(t)\, m_t^{\mathrm{DLRA}} + f(t)\bigr)\, \mathrm{d}t +C_t^{\mathrm{DLRA}} H(t)^\top R^{-1} \bigl( \mathrm{d}Z_t - H(t) m_t^{\mathrm{DLRA}} \, \mathrm{d}t \bigr).
	\end{equation}
	
	To derive the equation of the covariance, we exploit the usual DLRA-covariance relation, i.e.\
	\begin{equation}\label{eq: C_t DLRA}
	C_t^{\mathrm{DLRA}} = U_t^\top \, \mathbb{E}\!\left[ Y_t Y_t^\top \right] U_t,
	\end{equation}
	and Itô's formula to find
	\begin{equation*}
		\begin{aligned}
		\mathrm{d} C_t^{\mathrm{DLRA}}
		& =
		\mathrm{d}
		\bigl(
		U_t^\top \, \mathbb{E}[ Y_t Y_t^\top ] \, U_t
		\bigr) \\
		&=
		(\mathrm{d} U_t^\top)\, \mathbb{E}[ Y_t Y_t^\top ]\, U_t
		+
		U_t^\top \, \mathrm{d}\mathbb{E}[ Y_t Y_t^\top ] \, U_t
		+
		U_t^\top \, \mathbb{E}[ Y_t Y_t^\top ]\, \mathrm{d} U_t.
	\end{aligned}
\end{equation*}
Notice that in case of affine drift, the equation \eqref{eq: U_t JMCO} for the deterministic modes simplifies to
\begin{equation*}
	\mathrm{d} U_t = \left(U_t A(t)^{\top} + C_{Y_t}^{-1} U_t Q \right) P_{U_t}^\perp  \mathrm{d} t
\end{equation*}
and, hence, one obtains
		\begin{equation}\label{eq: C_t}
		\begin{aligned}
		\mathrm{d} C_t^{\mathrm{DLRA}} &=
		\Big(P_{U_t}^\perp \left(A(t)C_t^{\mathrm{DLRA}} + Q P_{U_t} \right)
		+
		U_t^\top \, \mathrm{d}\mathbb{E}[ Y_t Y_t^\top ]\, U_t
		+
		\left(C_t^{\mathrm{DLRA}} A(t)^\top + P_{U_t} Q\right)\, P_{U_t}^\perp \,\Big) \mathrm{d}t.
	\end{aligned}
		\end{equation}
	where, as it can be seen from \eqref{eq: C_t DLRA}, $P_{U_t}$ is the orthogonal projection onto the range of $C_t^{\mathrm{DLRA}}$.
	To determine $\mathrm{d} C_{Y_t} = \mathrm{d} \mathbb{E}[ Y_t Y_t^\top ]$, first let us recall that $Y_t$ follows the relation \eqref{eq: Y_t JMCO}. Then, using again Itô's formula
	\begin{equation*}
		\mathrm{d} \left(Y_t Y_t^{\top}\right)= (\mathrm{d}Y_t) Y_t^{\top} +  Y_t (\mathrm{d}Y_t)^{\top} + \mathrm{d} \langle Y_t, Y_t \rangle_t,
\end{equation*}
where $ \langle Y_t, Y_t \rangle_t$ denotes the quadratic covariation of $Y$, we have
	\begin{equation}\label{eq: C_Y_t}
	\begin{aligned}
		\mathrm{d}\mathbb{E}[ Y_t Y_t^\top ]
		&=
		\bigl(
		U_t A(t) U_t^\top C_{Y_t}
		+ C_{Y_t} U_t A(t)^\top U_t^\top
		- 2\, U_t C_t^{\mathrm{DLRA}} H(t)^\top R^{-1} H(t) C_t^{\mathrm{DLRA}} U_t^\top
		+ U_t Q U_t^\top
		\\
		&\qquad
		+ U_t C_t^{\mathrm{DLRA}} H(t)^\top R^{-1} H(t) C_t^{\mathrm{DLRA}} U_t^\top
		\bigr)\, \mathrm{d}t
		\\
		&=
		\bigl(
		U_t A(t) U_t^\top C_{Y_t}
		+ C_{Y_t} U_t A(t)^\top U_t^\top
		- U_t C_t^{\mathrm{DLRA}} H(t)^\top R^{-1} H(t) C_t^{\mathrm{DLRA}} U_t^\top
		+ U_t Q U_t^\top
		\bigr)\, \mathrm{d}t .
	\end{aligned}
	\end{equation}
	Inserting expression \eqref{eq: C_Y_t} into \eqref{eq: C_t} and using the orthogonality of $U_t$, one gets
	\begin{equation}\label{eq: C KBF DLRA}
	\begin{aligned}
		\mathrm{d}C_t
		=& \big(P_{U_t}^\perp (A(t)C_t^{\mathrm{DLRA}} + QP_{U_t})\, \;+\; P_{U_t} A(t) C_t^{\mathrm{DLRA}} \;+\; C_t^{\mathrm{DLRA}} A^{\top}(t) P_{U_t} \\
		&\quad -\; P_{U_t} C_t^{\mathrm{DLRA}} H(t)^\top R^{-1} H(t) C_t^{\mathrm{DLRA}} P_{U_t}
		\;+\; P_{U_t} Q P_{U_t}
		\;+\; (C_t^{\mathrm{DLRA}}A(t)^{\top}+ P_{U_t}Q ) P_{U_t}^{\perp}\Big)
	\, \mathrm{d}t \\
	= &\Big(
	A(t)\, C_t^{\mathrm{DLRA}} \;+\; A(t)^{\top} C_t^{\mathrm{DLRA}}
	\;+\; Q -  P_{U_t}^{\perp} Q\, P_{U_t}^{\perp}-\; C_t^{\mathrm{DLRA}} H(t)^\top \, R^{-1} H(t)\, C_t^{\mathrm{DLRA}}
	\Big)\, \mathrm{d}t.
	\end{aligned}
	\end{equation}
	where in the last line we used the fact that
	\begin{equation*}
	P_{U_t}\, C_t^{\mathrm{DLRA}} = C_t^{\mathrm{DLRA}},
	\end{equation*}
	as $P_{U_t}$ is the projection onto the range of $C_t^{\mathrm{DLRA}}$. 
	
	For implementation purposes, the whole procedure translates into a system of three equations defining the mean, the reduced covariance, and the subspace that defines the transformation of the reduced covariance in the ambient space:
		\begin{equation}\label{eq: DLRA KBF}
		\begin{aligned}
			\mathrm{d} m_t^{\mathrm{DLRA}} =& \bigl(A(t)\, m_t^{\mathrm{DLRA}} + f(t)\bigr)\, \mathrm{d}t +U_t^{\top}C_{Y_t} U_t H(t)^\top R^{-1} \bigl( \mathrm{d}Z_t - H(t) m_t^{\mathrm{DLRA}} \, \mathrm{d}t \bigr)\\
			\mathrm{d}C_{Y_t}
			=&
			\bigl(
			U_t A(t) U_t^\top C_{Y_t}
			+ C_{Y_t} U_t A(t)^\top U_t^\top
			- C_{Y_t} U_t H(t)^\top R^{-1} H(t) U_t^\top C_{Y_t}
			+ U_t Q U_t^\top
			\bigr)\, \mathrm{d}t \\
			\mathrm{d} U_t
			=& \bigl( U_t A(t)^{\top} + C_{Y_t}^{-1} U_t Q
			\Bigr)
			P_{U_t}^{\perp}\, \mathrm{d}t,
		\end{aligned}
	\end{equation}
    which is equivalent of the following 2-equations system, exploiting the fact that $P_{U_t} C_t^{\mathrm{DLRA}} = C_t^{\mathrm{DLRA}}$,
    \begin{equation}\label{eq: DLRA KBF 2eq}
    	\begin{aligned}
    		\mathrm{d} m_t^{\mathrm{DLRA}} =& \bigl(A(t)\, m_t^{\mathrm{DLRA}} + f(t)\bigr)\, \mathrm{d}t +C_{t}  H(t)^\top R^{-1} \bigl( \mathrm{d}Z_t - H(t) m_t^{\mathrm{DLRA}} \, \mathrm{d}t \bigr)\\
    		\mathrm{d}C_t^{\mathrm{DLRA}}
    		= & \Big(
    		A(t)\, C_t^{\mathrm{DLRA}} \;+\; A(t)^{\top} C_t^{\mathrm{DLRA}}
    		\;+\;  Q -  P_{C_t^{\mathrm{DLRA}}}^{\perp} Q\, P_{C_t^{\mathrm{DLRA}}}^{\perp}\\
    		&-\; C_t^{\mathrm{DLRA}} H(t)\, R^{-1} H(t)\, C_t^{\mathrm{DLRA}}
    		\Big)\, \mathrm{d}t,
    	\end{aligned}
    \end{equation}
    where $P_{C_t^{\mathrm{DLRA}}}$ is the orthogonal projector onto the range of $C_t^{\mathrm{DLRA}}$.
	\begin{Remark}[Differences with the DLRA-KBF proposed in \cite{nobile2025dynamicallowrankapproximationskalman}]
	  In \cite{nobile2025dynamicallowrankapproximationskalman} a DLRA-type KBF was derived according to the DLRA formulation for SDEs proposed in \cite{kazashi2025dynamical}. The main difference between the DLRA-KBF proposed in \cite{nobile2025dynamicallowrankapproximationskalman} and the one defined by \eqref{eq: DLRA KBF} is the presence of a diffusion component in the computation of the deterministic basis. This new term, which originates from deriving the filtered DLRA equations on the discrete process first and then taking the limit as $\Delta t \to 0$, can be beneficial to better track possible noise-dominated components in the evolution dynamics, as well as to tackle nearly ``low-rank" problems where the diffusion is responsible of determining some components of the dominant $k$-dimensional subspace. In this regard, \eqref{eq: DLRA KBF} is a good starting point to build ``enriched"-noise-type DLRA filters to better approximate nearly low-rank problems, where we can use an artificial ad-hoc noise to capture the DLRA projection error and influence the evolution of the subspace $U$, possibly making it depend on the observation operator. Notice that the DLRA-KBF filter from \cite{nobile2025dynamicallowrankapproximationskalman} and our filter \eqref{eq: DLRA KBF} coincide when $U_t \in \mathrm{Im}(Q)^{\perp}$. Proofs of well-posedness and useful properties of \eqref{eq: DLRA KBF} can be derive similarly to the treatment in \cite{nobile2025dynamicallowrankapproximationskalman}, and we refer to this reference for more theoretical details.
	\end{Remark}

\section{An enriched DLRA-JMCO filter}\label{sec: DLRA compl JMCO}

In this section we want to propose an ``enriched'' DLRA in order to better capture the possible non-negligible directions that lie outside the subspace constructed by the standard DLRA prediction strategy and that might not be negligible over the time integration. 

A further motivation for an enrichment strategy is given by the terms appearing in the equation of the deterministic basis \eqref{eq: DLRA JMCO analysis}, which is not influenced by the observation operator. This is rather unfortunate, because we expect this behavior when considering the best $k$-approximation of the standard Ensemble Kalman Filter in the full-dimensional space. In fact, the observation operator $H$ can change the eigenvalues of the filtered covariances, and in turn can change the order of its eigenvectors, and, hence, of its best $k$-dimensional subspace. This phenomenon is even more emphasized when $H$ is time dependent and produces observations that lies in a subspace different from the one where the DLRA actually lives.
For DLRA, the issues to address are essentially the following two: how to build a complemented strategy that it is still inexpensive to perform (which we will address at the end of this section) and how to estimate this complemented component.

Indeed, relating to the latter issue, we know that one of the benefits of DLRA is its on-the-fly computation, which allows us to compute a surrogate cheaply without the need of an expensive offline phase and to have access to the true solution of the system. In this regard, we make use of the following approximation assumption, which does not involve the true solution in \eqref{eq: DA SDE} and
has been considered in similar forms in the DLRA literature \cite{lubich2014projector,kazashi2021stability,kazashi2025dynamicalpartI,nobile2026high}, where error bounds comparing DLRA schemes with time-convergent numerical methods typically depend explicitly on this quantity.

\begin{Assumption}[$\varepsilon$-error]\label{ass: eps error}
	Let $\varepsilon:=\varepsilon(k) \geq0$  be defined as
	\begin{equation}\label{eps bound projection}
		\varepsilon^2 := \sup_{\substack{t \in [0,T] \\ Z \in L^2(\Omega,\mathbb{R}^d) \text{ of rank } k \\
				\text{ with } ·\mathbb{E}[|Z|^2] \leq \tilde{K}}} \max\{\mathbb{E}[|\mathring{A}(t,Z) -  P_{Z}\mathring{A}(t, Z) |^2], \ \mathbb{E}[\|  \mathcal{P}_{\mathcal{U}(Z)}^{\perp}Q^{\frac{1}{2}}\|^2_{\mathrm{F}} ] \},
	\end{equation}
	where $\tilde{K}:=\min \{ K > 0 \text{ such that } \max\{\mathbb{E}[\sup\limits_{0 \leq t \leq T} |X_t^{\mathrm{DLRA}}|^2], \mathbb{E}[\max\limits_{0 \leq n \leq N} |X_n^{\mathrm{DLRA}}|^2] \} \leq K \}$, $P_{Z}[ \ \cdot\ ]:= (I_{d\times d}-\mathcal{P}_{\mathcal{U}(Z)})[\ \cdot \ ]\mathcal{P}_{\mathcal{Y}(Z)}+\mathcal{P}_{\mathcal{U}(Z)}[\ \cdot \ ]$ is an orthogonal projector, whereas $P_\mathcal{U(Z)}$ and $\mathcal{P}_{\mathcal{Y}(Z)}$ denote the projection onto range and corange of $Z$, respectively, and $\mathring{A}= A - \mathbb{E}[A]$ denotes the centered drift.
\end{Assumption}

Assumption \ref{ass: eps error} can be locally estimated when employing the DLR Projector Splitting scheme for SDEs \cite{kazashi2025dynamicalpartI} with the respective estimator
\begin{equation}\label{eq: eps_n}
	\varepsilon^2_n = \max\{\mathbb{E}[|\mathring{A}(t,X_n^{\mathrm{DLRA}}) -  P_{U_n^{\top}\widehat{Y}_{n+1}}\mathring{A}(t, X_n^{\mathrm{DLRA}}) |^2], \ \mathbb{E}[\|  \mathcal{P}_{U_n}^{\perp}Q^{\frac{1}{2}}\|^2_{\mathrm{F}} \}.
\end{equation}
Evaluating \eqref{eq: eps_n} is still cheap from the computational point of view if the problem is effectively low-rank, and, hence, characterized by the drift $A$ and the diffusion $Q$ which are low-rank themselves. Furthermore, notice that $\varepsilon_n$ can be computed with quantities
already stored in the elaboration of the algorithm.  In the derivation of the DLRA-JMCO in Proposition \ref{prop: triple JMCO}, we make evolve the DLRA $X_{n+1}^{\mathrm{DLRA}}$ and the full-order model with the same initial point $X_n=X_{n}^{\mathrm{DLRA}}$ in a prediction-then-analysis procedure. It is then natural to query whether an enriched predicted DLRA can produce a similar JMCO filter that better approximate the true filtering distribution. We pursue this analysis by first exploiting the local error between the predicted DLRA $\widehat{X}_{n+1}^{\mathrm{DLRA}}$ and the full-order model $\widehat{X}_{n+1}$ as function of $\varepsilon_n$ and then deriving complemented updates $(\Delta m_n^{\mathrm{DLRA},c}, \Delta U_n^{c}, \Delta Y_n^{c})$ for this enriched DLRA filter, which are denoted by the superscript ``$c$'' to indicate their complemented nature, again as minimization of Kalman updates.

Once this local error has been quantified, we enrich the DLRA dynamics prior to the analysis step by introducing an independent centered diffusion term, whose variance is chosen to commensurate the estimated accumulation of local errors over time, as we will explain in details below.

In practice, we estimate the local error dependent of \eqref{eq: eps_n} comparing the DLR Projector Splitting scheme for SDEs \cite{kazashi2025dynamicalpartI} with a Euler-Maruyama method for the true predicted update. Under standard assumption of the drift and the diffusion, both methods are strongly convergent to the DLRA and the true solution, respectively, yielding the consistency of this estimate. This result is summarized in Lemma \ref{lem: approx error}.

\begin{Lemma}[Local approximation error]\label{lem: approx error}
	Fix $n$ such that $1 \leq n \leq N$. Assume that at time $t_{n+1}$ $\widehat{X}_{n+1}$ is the Euler Maruyama solution where $X_{n} = X_{n}^{\mathrm{DLRA}}$, whereas $\widehat{X}^{\mathrm{DLRA}}_{n+1}$ denotes the DLR Projector Splitting scheme for SDEs. Assume that the drift $A$ is Lipschitz continuous with respect to the state variable uniform in $t$ with constant $C_{\mathrm{Lip}}$. Then, one has
	\begin{equation}\label{eq: Lambda_n}
		\begin{aligned}
			\mathbb{E}\!\left[\left| \widehat{X}_{n+1} - \widehat{X}^{\mathrm{DLRA}}_{n+1} \right|^2 \right]
			\leq 3\varepsilon^2 (2+k)\Delta t_n
			=: \Lambda(\varepsilon,\Delta t_n)= \Lambda_{n+1},
		\end{aligned}
	\end{equation}
	and
	\begin{equation}\label{eq: cov Lambda_n}
		\begin{aligned}
			\left\| \widehat{C}_{n+1} - \widehat{C}_{n+1}^{\mathrm{DLRA}} \right\|_{\mathrm{F}}
			\leq \Lambda_{n+1},
		\end{aligned}
	\end{equation}
	where $\widehat{C}_{n+1} = \mathbb{E}\left[\widehat{X}_{n+1} \widehat{X}_{n+1}^{\top}\right]$ and $\widehat{C}_{n+1}^{\mathrm{DLRA}} = \mathbb{E}\left[\widehat{X}^{\mathrm{DLRA}}_{n+1} (\widehat{X}^{\mathrm{DLRA}}_{n+1})^{\top}\right]$.
	\begin{proof}
		The proof can be found in Appendix \ref{app: proof DLRA-JMCO enriched}.
	\end{proof}
\end{Lemma}

Lemma \ref{lem: approx error} also gives an estimate on the error discrepancy between the covariances of the predicted DLRA and the predicted full-order solution, which is of order $O(\Delta t)$. With this result, we consider the following ansatz for our enriched DLRA: the continuous true dynamics of the state in \eqref{eq: DA SDE} can be actually well approximated by the DLRA dynamics and a diffusion noise living in the orthogonal component of the subspace of the DLRA, namely
\begin{equation}\label{eq: enrichment ansatz}
	\begin{aligned}
	X_t =& X_{t}^{\mathrm{DLRA}} + P_{U_{t}}^{\perp} D_t \xi_t \\
	=&  X_{t}^{\mathrm{DLRA}} + \mathcal{E}_t
	\end{aligned}
\end{equation}
where $D_t : [0, \infty) \times  \mathbb{R}^{d} \to \mathbb{R}^{d \times s}$ is jointly measurable and  $\xi_t$ is a $s$-dimensional $(\mathcal{F}_t)$-Brownian motion, independent of $W_t,$ $B_t$, $\widehat{B}_t$, and $X_0$, where the dimension $s\in \mathbb{N}$ is at most the dimension of the orthogonal complement subspace to $U_t$, i.e.\ $s \leq d -k$, and $\left\| \mathbb{E}[D_t D_t^{\top}] \right\|_{\mathrm{F}} \leq 3\varepsilon^2 (2+k) $. Notice that if $s$ is comparable to the magnitude of the rank $k$, then estimating the term $\mathcal{E}_t = P_{U_t}^{\perp} D_t \xi_t$ should still be affordable from the computational point of view. 

The ansatz \eqref{eq: enrichment ansatz} is motivated by the following two observations. As the mean of the DLRA follows the same equation of the mean of the full-order Kalman update, we can further assume that the DLRA tracks perfectly this latter mean, and therefore the added orthogonal component $\mathcal{E}_t$ is centered. In addition, the covariance of the complemented diffusion term $D_t$ has value to commensurate the integration over time of the local error of the continuous version of \eqref{eq: cov Lambda_n}. Moreover, notice that in this
perspective if the full-order solution has rank $k$ and the DLRA matches its evolution
exactly at each time $t$, then $\varepsilon = 0$ and, hence, $D_t$ is the
null-matrix, i.e.\ there is no need to complement the dynamics.

After discretization, \eqref{eq: enrichment ansatz} becomes
\begin{equation}\label{eq: enrichment ansatz discr}
	\begin{aligned}
	\widehat{X}_{n+1} =& \widehat{X}_{n+1}^{\mathrm{DLRA}} + \mathcal{E}_{n+1} \\
	=& \widehat{X}_{n+1}^{\mathrm{DLRA}}  + P_{\widehat{U}_{n+1}}^{\perp} D_{n}  \xi_n,
		\end{aligned}
\end{equation}
where we use the superscript $ \ \widehat{} \ $ to denote the evolution as a prediction step.
Then, 
\begin{equation}\label{eq: pred E_n+1}
	\begin{aligned}
	\widehat{E}_{n+1}:= &\mathbb{E}[\mathcal{E}_{n+1} \mathcal{E}_{n+1}^{\top}] \\
	=& P_{\widehat{U}_{n+1}}^{\perp} \mathbb{E}\left[D_{n} \xi_n \xi_n^{\top}D_{n}^{\top}\right] P_{\widehat{U}_{n+1}}^{\perp}=  P_{\widehat{U}_{n+1}}^{\perp} \mathbb{E}\left[D_{n} D_{n}^{\top}\right] P_{\widehat{U}_{n+1}}^{\perp} t_n 
	\end{aligned}
\end{equation}
denotes the predicted covariance of the orthogonal component, and the norm of its covariance can be bound by
\begin{equation}\label{eq: overline{Lambda}}
	\begin{aligned}
	\| \widehat{E}_{n+1} \|_{\mathrm{F}} \leq \| \mathbb{E}\left[D_{n} D_{n}^{\top}\right] t_n \|_{\mathrm{F}} \leq  & 3\varepsilon^2 (2+k) \Delta t_{n} \\
	= & \sum_{\ell=0}^{n} 3\varepsilon^2 (2+k) \Delta t_{\ell} = \sum_{\ell=0}^{n} \Lambda_{\ell} =: \overline{\Lambda}_{n} \\
	\end{aligned}
\end{equation}
 We want to exploit \eqref{eq: enrichment ansatz} in order to derive new DLRA equations that can evolve also thanks to the observation operator, but that preserve their computational advantages.

\iffalse
the Gramian $\widehat{C}_{n+1}$ of the full order solution $\widehat{X}_{n+1}$ admits the following decomposition
\begin{equation}\label{eq: predic enrich cov 1}
	\begin{aligned}
		\widehat{C}_{n+1} &= \widehat{C}_{n+1}^{\mathrm{DLRA}} + \widehat{E}_{n+1}\\
		&= \widehat{C}_{n+1}^{\mathrm{DLRA}}
		+ \widehat{E}_{n}
		+ P_{\widehat{U}_{n+1}}^\perp
		\mathbb{E}\left[D_nD_n^{\top}\right] 
		P_{\widehat{U}_{n+1}}^\perp \Delta t \\
		&= \widehat{C}_{n+1}^{\mathrm{DLRA}}
		+ 
		\left(\sum_{\ell=1}^{n+1}
		P_{\widehat{U}_{\ell}}^\perp \mathbb{E}\left[D_\ell D_\ell^{\top}\right]  P_{\widehat{U}_{\ell}}^\perp\right) \Delta t.
	\end{aligned}
\end{equation}
We can give an upper bound of the estimate of 
Assuming $\{D_\ell\}_{\ell=1,\ldots,n}$ as conditionally independent, then the norm of $\mathcal{E}_n$ is bounded from above by
\begin{equation}\label{eq: var E_n+1}
	\begin{aligned}
		\left| \widehat{E}_{n+1} \right|  = & \left|
		\left(\sum_{\ell=1}^{n+1}
		P_{\widehat{U}_{\ell}}^\perp \mathbb{E}\left[D_\ell D_\ell^{\top}\right]  P_{\widehat{U}_{\ell}}^\perp\right) \Delta t
		\right| \\
		\leq & 
		\sum_{\ell=1}^{n+1} \left|
		P_{\widehat{U}_{\ell}}^\perp \mathbb{E}\left[D_\ell D_\ell^{\top}\right]  P_{\widehat{U}_{\ell}}^\perp\right|^2 \Delta t \\
		\leq &
		\sum_{\ell=1}^{n+1}
		\mathbb{E}\left[D_\ell D_\ell^{\top}\right]  \Delta t
		= \sum_{\ell=1}^{n+1}\Lambda_\ell
		= : \overline{\Lambda}_{n}.
	\end{aligned}
\end{equation}
\fi
Before going into computations, it is worth highlight two technical details. The former concerns the magnitude of $\overline{\Lambda}_n$ defined in \eqref{eq: overline{Lambda}}. Indeed, $\overline{\Lambda}_\ell$ is not said a priori to be negligible in the order of $O(\Delta t)$ and therefore its presence may modify the evolution of $\Delta U_n^c$ and $\Delta Y_n^c$.
Our assumption $P_{\widehat{U}_{n+1}}^{\perp} D_{n} \xi_n$ is a centered isotropic Gaussian living in the orthogonal component of $\widehat{U}_{n+1}$, i.e.\ $P_{\widehat{U}_{n+1}}^{\perp} D_{n} \xi_n \sim \mathcal{N}(0,\frac{\overline{\Lambda}_{n}}{d-k}P_{\widehat{U}_{n+1}}^\perp)$, which we recall that it is independent of the initial condition $X_0$ and of the Brownian motions $\{W_t\}_{t\geq 0}$ and $\{B_t\}_{t\geq 0}$ (and hence of the stochastic basis $Y$). If more information concerning the orthogonal complement is available, one could confine this enrichment term in a smaller subspace belonging to the orthogonal complement of $\widehat{U}_{n+1}$. Then, from \eqref{eq: enrichment ansatz} the full-order filtered covariance would be expressed by the following Kalman update
\begin{equation}\label{eq: cov enriched}
	\begin{aligned}
		C_{n+1}
		&= \widehat{C}_{n+1}
		- \widehat{C}_{n+1} H^\top R^{-1} H \widehat{C}_{n+1}\Delta t \\
		&= \left(\widehat{C}_{n+1}^{\mathrm{DLRA}} + \widehat{E}_{n+1} \right)
		- \left(\widehat{C}_{n+1}^{\mathrm{DLRA}} + \widehat{E}_{n+1} \right)
		H^\top(t_n) R^{-1} H(t_n)
		\left(\widehat{C}_{n+1}^{\mathrm{DLRA}} + \widehat{E}_{n+1}\right)\Delta t
	\end{aligned}
\end{equation}
in a first order of $\Delta t$, and the enriched mean reads as
\begin{equation}\label{eq: mean enriched}
	\begin{aligned}
		m_{n+1}
		&= \widehat{m}_{n+1}^{\mathrm{DLRA}} 
		+ \widehat{C}_{n+1} H(t_n)^\top R^{-1} \left( \Delta Z_n - H(t_n) \widehat{m}_{n+1}^{\mathrm{DLRA}} \Delta t \right)\\
		&= \widehat{m}_{n+1}^{\mathrm{DLRA}} 
		+ 	\left(\widehat{C}_{n+1}^{\mathrm{DLRA}} + \widehat{E}_{n+1}\right) H^\top R^{-1} \left( \Delta Z_n - H(t_n) \widehat{m}_{n+1}^{\mathrm{DLRA}} \Delta t \right).\\
	\end{aligned}
\end{equation}
We seek new updates $(\Delta m_n^{\mathrm{DLRA},c},\Delta U_n^c,\Delta Y_n^c)$ minimizing \eqref{eq: mean enriched} and \eqref{eq: cov enriched} to obtain an enrich DLRA strategy in the following proposition.
\begin{Proposition}\label{prop: compl DLRA-JMCO}
	The updates $(\Delta m_n^{\mathrm{DLRA},c},\Delta U_n^c,\Delta Y_n^c)$ which solves
	\begin{equation}\label{eq: DLRA JMCO enriched min}
		\begin{aligned}
			\min_{ \substack{
				\Delta m_n^{\mathrm{DLRA},c},\Delta U_n^c, \\
				\Delta Y_n^c,
				\text{ satisfying \eqref{eq: DLRA updates JMCO}-\eqref{eq: DLRA updates JMCO fin}} }}
			&\left|m_{n+1}-(m_n+\Delta m_n^{\mathrm{DLRA},c})\right|^2 \\
			&\quad
			+ \left\|
			C_{n+1}
			- (U_n+\Delta U_n^c)^\top
			(C_{Y_n}+\Delta C_{Y_n}^c)
			(U_n+\Delta U_n^c)
			\right\|_{\mathrm{F}}^2
		\end{aligned}
	\end{equation}
	with $C_{n+1}$ and $m_{n+1}$ defined in \eqref{eq: cov enriched} and \eqref{eq: mean enriched}, respectively, and $\Delta C_{Y_n}^c=\mathbb{E}\!\left[\Delta Y_n^c(\Delta Y_n^c)^\top\right]$, are given by
	\begin{equation*}
		\begin{aligned}
			\Delta m_n^{\mathrm{DLRA},c}
			:=&\ \mathbb{E}\!\left[A(t_n,X_n^{\mathrm{DLRA},c})\right]\Delta t
			+ \left( C_n^{\mathrm{DLRA}} + P_{U_n}^{\perp} \mathbb{E}\left[D_n D_n^{\top}\right] P_{U_n}^{\perp} \right) H^\top R^{-1}
			\left(\Delta Z_n-H(t_n)m_n^{\mathrm{DLRA}}\Delta t\right) \\
			\Delta U_n^c
			:=&\ C_{Y_n}^{-1}
			\left(
			\mathbb{E}\!\left[Y_n \mathring{A}(t_n,X_n^{\mathrm{DLRA},c})^\top\right]
			+ U_nQ
			\right)
			P_{U_n}^\perp\Delta t \\
			&\quad
			- U_n H^\top(t_n) R^{-1}H(t_n) P_{U_n}^{\perp} \mathbb{E}\left[D_n D_n^{\top}\right] P_{U_n}^\perp \Delta t \\
			\Delta Y_n^c
			:=&\ U_n \mathring{A}(t_n,X_n^{\mathrm{DLRA},c})\Delta t
			+ U_nQ^{\frac{1}{2}}\Delta W_n
			-  C_{Y_n}U_nH(t_n)^\top R^{-1}H(t_n) U_n^\top Y_n\Delta t \\
			&\quad
			- C_{Y_n}U_n H(t_n)^\top R^{-\frac{1}{2}}\Delta\widehat{B}_n.
		\end{aligned}
	\end{equation*}
	\begin{proof}
		The proof can be found in Appendix \ref{app: proof DLRA-JMCO enriched}.
	\end{proof}
\end{Proposition}
Denoting the Kalman update of the covariance of $\widehat{E}_{n+1}$ as 
$$E_{n+1} = \left( I_{d \times d} - \widehat{E}_{n+1} H^\top(t_n) R^{-1} H(t_n) \right) \widehat{E}_{n+1},$$ 
then we update the accumulate approximation error as
\begin{equation}\label{eq: lambda}
	\overline{\Lambda}_{n} = \lambda_{\mathrm{max}} (E_{n+1}) = \lambda_{\mathrm{max}} \left( \widehat{E}_{n+1} -  \widehat{E}_{n+1} H^\top(t_n) R^{-1} H(t_n)\widehat{E}_{n+1}\right),
\end{equation}
where $\lambda_{\mathrm{max}}(D)$ indicates the largest eigenvalue of a matrix $D$, and the next update of the discrepancy will be computed as $\overline{\Lambda}_{n+1} = \overline{\Lambda}_{n} + \Lambda_{n+1}$, i.e.\ adding the new local error. Notice that as $\widehat{E}_{n+1} \succeq E_{n+1}$, then the quantity $\overline{\Lambda}_{n}$ after the analysis can be smaller than its prediction counterpart, in compliance to the fact that the information available diminish the model error also in the direction of propagation of the complement.

For the sake of computability, we approximate the drift $A(t_n,X_n^{\mathrm{DLRA},c}$ in \eqref{eq: DLRA JMCO enriched min} by $A(t_n,X_n^{\mathrm{DLRA}}$, assuming that this discrepancy is negligible and its small variation can be captured by the analysis step on $\widehat{E}_{n+2}$ at the next step (see Remark \ref{rmk: DLRA enriched drift}).

Notice that the update of the deterministic modes in \eqref{eq: DLRA JMCO enriched min} can now evolve with respect to the observation operator $H$. If $h$, the number of observed coordinates, is small or we know a priori that $\mathcal{E}_{n+1}$ remains confined in a small set of directions, then the term $\Delta U_n^{c}$ is still computationally affordable. Furthermore, as proposed in the analysis updates for $\Delta m_n^{\mathrm{DLRA}}$ and $\Delta Y_n$ in \eqref{eq: DLRA JMCO analysis}, we can evolve $U_{n+1}$ in a staggered way as follows
\begin{equation*}
	\begin{aligned}
		C_{\widetilde{Y}_{n+1}} \widetilde{U}_{n+1}
		=& C_{\widetilde{Y}_{n+1}}  U_n +
		\left(
		\mathbb{E}\!\left[\widetilde{Y}_{n+1} \mathring{A}(t_n,X_n^{\mathrm{DLRA}})^\top\right]
		+  U_nQ
		\right)
		P_{U_n}^\perp\Delta t, \\
		U_{n+1}=& \widehat{U}_{n+1}
		- U_{n+1}H^\top(t_n) R^{-1}H(t_n) P_{U_n}^{\perp} \mathbb{E}\left[D_n D_n^{\top}\right] P_{U_n}^{\perp} \Delta t, \\
		& - C_{Y_n}^{-1} P_{U_n}^{\perp} \mathbb{E}\left[D_n D_n^{\top}\right] P_{U_n}^{\perp} ) H^\top(t_n) R^{-1}H(t_n) P_{U_n}^{\perp} \mathbb{E}\left[D_n D_n^{\top}\right] P_{U_n}^\perp \Delta t \\
		 \text{ with } &\widehat{U}_{n+1}^{\top}, R_{n+1} = \texttt{QR}(\widetilde{U}_{n+1}^{\top}),
	\end{aligned}
\end{equation*}
and $m_{n+1}^{\mathrm{DLRA}}$ and $Y_{n+1}$ as well
\begin{equation}
	\begin{aligned}
	\widehat{m}_{n+1} = & m_{n} + \mathbb{E}\!\left[A(t_n,X_n^{\mathrm{DLRA},c})\right]\Delta t, \\
	m_{n+1}^{\mathrm{DLRA}} = & \widehat{m}_{n+1}
	+ \left( C_n^{\mathrm{DLRA}} + P_{U_n}^{\perp} \mathbb{E}\left[D_n D_n^{\top}\right] P_{U_n}^{\perp} \right) H^\top R^{-1}
	\left(\Delta Z_n-H(t_n)m_{n+1}^{\mathrm{DLRA}}\Delta t\right),  \\
	\widetilde{Y}_{n+1} = &Y_n + U_n A(t_n,X_n^{\mathrm{DLRA}}) \Delta t + U_n Q^{\frac{1}{2}} \Delta W_n, \\
	Y_{n+1} = & \widehat{Y}_{n+1} - C_{ \widehat{Y}_{n+1} } \widehat{U}_{n+1} H^\top(t_n) R^{-1}H(t_n) H(t_n) \widehat{U}_{n+1}^\top Y_{n+1} \Delta t \\
	& - C_{ \widehat{Y}_{n+1} } \widehat{U}_{n+1} H(t_n)^\top R^{-\frac{1}{2}}\Delta\widehat{B}_n, \text{ with }  \widehat{Y}_{n+1} = R_{n+1} \widetilde{Y}_{n+1}.
	\end{aligned}
\end{equation}
We recall that when we ask for $P_{\widehat{U}_{n+1}}^{\perp} D_{n} \xi_n \sim \mathcal{N}(0,\overline{\Lambda}_{n}P_{\widehat{U}_{n+1}}^\perp)$, then we have $P_{U_n}^{\perp} \mathbb{E}\left[D_n D_n^{\top}\right] P_{U_n}^{\perp} =  \overline{\Lambda}_{n} P_{U_n}^{\perp}$.

It remains to discuss how practically implement the correction $\mathcal{E}_n$ in order to maintain the whole computational cost low also when we are not considering an analysis step as the one of an ensemble Kalman procedure, for instance when it is needed to sample from this distribution, as in the case of the particle filter that we will see in Section \ref{sec: compl DLRA PF}.
In this regard, we have to make explicit $\mathcal{E}_n$ to add to the predicted $\widehat{X}_{n+1}^{\mathrm{DLRA}}$. We propose to compute $s$ eigenvectors of $P_{U_n}^{\perp}\mathbb{E}[D_n D_n^{\top}]P_{U_n}^{\perp}$, namely $q_1,\ldots,q_s$ and build the matrix
\begin{equation*}
	\begin{aligned}
		Q=[q_1,\ldots,q_s]\in\mathbb{R}^{d\times s}.
	\end{aligned}
\end{equation*}
Then we set
\begin{equation*}
	\begin{aligned}
	\mathcal{E}_n^{i}=\sqrt{\frac{\overline{\Lambda}_n}{s}}\,Q\xi_i,\qquad \text{with}\qquad \xi_i\sim\mathcal{N}(0,I_{s\times s})
	\end{aligned}
\end{equation*}
and, hence
\begin{equation}\label{eq: compl X_hat}
	\begin{aligned}
		\widehat{X}_{n+1}^{\mathrm{DLRA},c,i}=\widehat{X}_{n+1}^{\mathrm{DLRA},i}+\mathcal{E}_n^{i}.
	\end{aligned}
\end{equation}

Of course, in order for this procedure to be computationally affordable, $s$ should be of the same magnitude as the rank $k$. As $\mathcal{E}_n$ enters only in the computation of $\Delta U_n$ multiplied by the term $H^\top(t_n) R^{-1} H(t_n)$, then in case of isotropic noise $\mathcal{E}$ another option is to obtain the first $s$ singular vectors of
\begin{equation*}
	\begin{aligned}
		(I-U_n^\top U_n)H(t_n)^{\top},
	\end{aligned}
\end{equation*}
whose procedure is still affordable if $h$ is small. As already mentioned, the purpose of this additional stochastic component is to account for uncertainty that is not captured by the low-rank dynamics: if further structural information about the underlying system is available, both $\xi_n$ and the covariance structure of the complementary dynamics can be modified accordingly to incorporate this information.
\begin{Remark}[Approximation of the drift in \eqref{eq: DLRA JMCO enriched min}]\label{rmk: DLRA enriched drift}
	Notice that equations \eqref{eq: DLRA JMCO enriched min} have the drift $A$ which is computed using the (non enriched) DLRA $X_{n}^{\mathrm{DLRA}}$ at time $t_n$ instead of the complemented one $X_{n}^{\mathrm{DLRA},c}$. This is pursued by assuming that the discrepancy between the two quantities can be tracked by the approximability assumption. To keep into account of this error, one can consider another diffusion-type complemented term in \eqref{eq: enrichment ansatz}. A possible complement can be described by an Ornstein–Uhlenbeck process whose scalar coefficient of the drift and variance are this error in the drift and the accumulate error $\overline{\Lambda}_{n}$, respectively. Then, one would pursue the minimization problem \eqref{eq: DLRA JMCO enriched min} also with respect to these parameters. This research direction will be pursued in future works.
\end{Remark}

%\section{Further approaches: (Complementary) DLRA Particle Filter}\label{sec: compl DLRA PF}
\section{A DLRA Particle Filter}\label{sec: compl DLRA PF}
In this section, we propose a particle filter for DLRA for SDEs, which we call \emph{DLRA particle filter} (DLRA PF). 
The adopted strategy reflects the usual continuous-time particle filter for SDEs proposed 
in \cite{bain2009fundamentals}. Indeed, we are seeking a Bayesian-type 
reduced-order model algorithm where the prediction step is performed using
the usual DLRA for SDE and the analysis step is performed via
an empirical measure whose weights are computed from the likelihood functional dependent on the observation process $\mathrm{d}Z_t$, and a possible branching procedure at fixed intervals of time,
as proposed in \cite{bain2009fundamentals}.  For the sake of clarity we briefly recall how the approximated particle filter 
proposed in  \cite{bain2009fundamentals} works for \eqref{eq: DA SDE}. Assume that we have at our disposal
$M$ particles $(X^i)_{i = 1,\dots,M}$ with independent 
initial positions $X_0^i$ and initial equal weights $\lambda_0^{i}= \frac{1}{M}$ for all $i = 1,\dots,M$. In the time 
interval $[t_n, t_{n+1}]$, the particles evolve according to the general SDE
\begin{equation*}
X^i_t = X^i_{t_n} + \int_{t_n}^t A(s, X^i_s) \, \mathrm{d}s 
+ \int_{t_n}^t Q^{\frac{1}{2}} \, \mathrm{d}W_s^{i},
\end{equation*}
or a suitable time-integration scheme.
To each particle one associates a time-dependent positive weight $\lambda^i(t)$ obtained through the following relation
\begin{equation}\label{eq: weights}
\widehat{\lambda}^i(t) 
:= \exp \left(
\int_{t_n}^{t} \left( H(s) X^i_s \right)^{\top} \, R^{-1}\mathrm{d} Z_s
- \frac{1}{2}
\int_{t_n}^{t} \left| R^{-\frac{1}{2}}H(s) X^i_s \right|^{2} \, \mathrm{d}s
\right)\lambda_n^{i}
\end{equation}
and the final weights at time $t_{n+1}$ are obtained through a normalization procedure, i.e.\
\begin{equation*}
\lambda^i_{n+1} =
\frac{\widehat{\lambda}^i(t_{n+1})}{\sum_{j=1}^{M} \widehat{\lambda}^j(t_{n+1})},
\qquad
\text{for all } i = 1,\dots,M.
\end{equation*}

At certain time steps, a branching process on 
$(X^i)_{i = 1,\dots,M}$ is applied and all the weights are reinitialized 
to $\frac{1}{M}$. The branching procedure \cite[Section 9.2.1]{bain2009fundamentals} is
fundamental in the construction of the particle filter as it
guarantees the convergence of the empirical measure induced by the particle filter to the one
of the true filtered process at the Monte Carlo rate of $\frac{1}{\sqrt{M}}$ \cite{bain2009fundamentals} and prevent weights collapse. For the sake of completeness, we detail the branching algorithm in Appendix \ref{app: proof DLRA PF}.

We aim to build a similar algorithm that performs the analysis
step in the low-rank subspace.
Namely, starting from the triple $(m_n, U_n, (Y^{i}_n)_{i})$, mean, deterministic, and stochastic
basis of the DLRA $X_n^{\mathrm{DLRA}}$ at time $t_n$, where each $(Y^{i}_n)_{i}$ has associated a normalized weight $\lambda_n^{\mathrm{DLRA,i}}$, we obtain the prediction triple of rank $k$
$(\widehat{m}_{n+1}^{\mathrm{DLRA}}, \widehat{U}_{n+1}, \widehat{Y}_{n+1})$ via the DLRA for SDEs, for instance using the DLR Projector
Splitting for SDEs \cite{kazashi2026approaches}, i.e.\
\begin{equation*}
\begin{aligned}
	\widehat{m}_{n+1}^{\mathrm{DLRA}} &= m_n + \mathbb{E}[A(t_n, X_n^{\mathrm{DLRA}})]\Delta t_n \\
	\widetilde{Y}_{n+1} &= Y_n + U_n \left( \mathring{A}(t_n, X_n^{\mathrm{DLRA}}) \right) \Delta t_n + U_n Q^{\frac{1}{2}} \Delta W_n \\
	\widetilde{U}_{n+1} &= U_n + C_{\widetilde{Y}_{n+1}}^{-1} \left( \mathbb{E}[\widetilde{Y}_{n+1} \left( \mathring{A}(t_n, X_n^{\mathrm{DLRA}}) \right)^{\top}] + U_n Q \right) P_{U_n}^\perp \Delta t_n \\
\end{aligned}
\end{equation*}
with $\widehat{X}_{n+1}^{\mathrm{DLRA}} = \widehat{m}_{n+1}^{\mathrm{DLRA}} + \widetilde{U}_{n+1}^{\top} \widetilde{Y}_{n+1}$ and $\widehat{U}_{n+1}, R_{n+1} = \texttt{QR}(\widetilde{U}_{n+1})$, $\widehat{Y}_{n+1} = R_{n+1}\widetilde{Y}_{n+1}$. Then, we perform the analysis
step similarly to \eqref{eq: weights} where we obtain the unnormalized weights
as
\begin{equation*}
	\widetilde{\lambda}_{n+1}^{\mathrm{DLRA}, i}  = \lambda_n^{\mathrm{DLRA}, i} \exp\left( (H(t_n)\widehat{X}_{n+1}^{\mathrm{DLRA}, i})^{\top} R^{-1}\Delta Z_n - \frac{1}{2} |R^{-\frac{1}{2}}H(t_n) \widehat{X}_{n+1}^{\mathrm{DLRA}, i}|^2 \Delta t_n \right),
\end{equation*}
and, hence, the normalized ones as $\widehat{\lambda}_{n+1}^{\mathrm{DLRA}, i} = \frac{\widetilde{\lambda}_{n+1}^{\mathrm{DLRA}, i} }{\sum_{j=1}^{M} \widetilde{\lambda}_{n+1}^{\mathrm{DLRA}, j}}$ for $i=1,..,M$. If $t_{n+1}$ is not a branching time, we simply set 
\begin{equation*}
	m_{n+1} = \widehat{m}_{n+1}, \quad U_{n+1} = \widehat{U}_{n+1}, \quad Y_{n+1} = \widehat{Y}_{n+1}, \quad \lambda_{n+1} = \widehat{\lambda}_{n+1}.
\end{equation*}	
On the other hand, if $t_{n+1}$ is a branching time, then a branching procedure is performed as proposed in \cite[Section 9.2.1]{bain2009fundamentals}
\begin{equation}\label{eq: branch}
X_{n+1}^{\mathrm{DLRA}}, (\lambda_{n+1}^{\mathrm{DLRA}, i})_{i=1,..,M} = \mathrm{branching}\left( \widehat{X}^{\mathrm{DLRA}}_{n+1}, (\lambda_{n+1}^{\mathrm{DLRA}, i} )_{i=1,..,M} \right),
\end{equation} 
where the weights are reset to $(\lambda_{n+1}^{\mathrm{DLRA}, j})_{j=1,..,M} = \frac{1}{M}$, for all $j=1,...,M$, after procedure \eqref{eq: branch}. Notice that this procedure can be performed on the reduced particles $\widehat{Y}^{i}_{n+1}$ instead of the $d$-dimensional ones $\widehat{X}^{\mathrm{DLRA},i}_{n+1}$, enabling computational time savings.
In either case, this procedure gives us the following empirical filtered distribution at $t_{n+1}$:
\begin{equation*}
	\pi_{n+1}^{\mathrm{DLRA}, M}(\mathrm{d}x)=\frac{1}{M}\sum_{i=1}^M 	\lambda_{n+1}^{\mathrm{DLRA}, i} \delta_{	X_{n+1}^{\mathrm{DLRA}, i}}(\mathrm{d}x).
\end{equation*}

\begin{Remark}[Choice of the branching times]
		The distance $\Delta t_{\mathrm{branch}}$ between two consecutive branching times is a parameter of choice and is regulated by a tradeoff. Choosing $\Delta t_{\mathrm{branch}}$ very large may lead to degeneracy of the particle ensemble and slow convergence. On the other hand, branching at the same mesh points of the dynamics may lead to large variance of the ensemble.
\end{Remark}

Like the DLRA-JMCO filter and similar to \cite[Section 7]{Vidlickovathesis}, the range of the deterministic basis
obtained by the analysis step coincides with the one of the prediction step.

\begin{Lemma}\label{lem: invar ran}
	For all $n$, the bases $\widehat{U}_{n+1}, U_{n+1}$ produced by the DLRA PF
	satisfy the following property
	\begin{equation*}
		 \mathrm{Ran}(U_{n+1}^{\top})  \subseteq \mathrm{Ran}(\widehat{U}_{n+1}^{\top}).
	\end{equation*}
	\begin{proof}
		The proof can be found in Appendix \ref{app: proof DLRA PF}.
	\end{proof}
\end{Lemma}

Lemma \ref{lem: invar ran} highlights a property that can be dangerous for filtering:
the deterministic basis is not affected by the observations and nor by the observation
operator in the analysis step. From this perspective, one can observe that in a regime
of continuous-time observations this might not be an issue. Indeed, for very small time steps, although $U_{n+1}$ does not depend on $\Delta Z_n$, $Y_{n+1}$ does, and so $U_{n+2}$ at the next time step. Hence the subspace can adjust fast to the incoming information. However, this may not be the case in the regime of discrete
observations, where the assimilation time-step can be much longer than the one used to evolve the state. If the subspace informed by the observation operator differs considerably from the one propagated by the DLRA, this may result in a significant error in predicting the filtering distribution.
%For instance, this is the case of a time-dependent operator $H$ such that it tracks observations in a different system of reference. We propose then
%a modification of the above algorithm, aiming at enriching the DLRA in order
%to obtain a range of $U_{n+1}$ different to the one of the deterministic basis $\widehat{U}_{n+1}$  %of
%the prediction step.

\subsection{Complementary DLRA-PF}
Similarly to the discussion in Section \ref{sec: DLRA compl JMCO}, we want to build a modified DLRA-PF algorithm such that the deterministic
basis of the analysis step can differ a priori  from the one of  the prediction
step. To derive such a filter, we exploit the action of the likelihood in directions orthogonal to the predicted deterministic basis. This property will allow us to make the range of the surrogate
depend on the obtained observation data. Again, we point out that this modified
dynamics needs to remain computationally cheap, so that
the advantage of DLRA is not lost. In this regard, we propose a complementation of the prediction DLRA as stated in \eqref{eq: compl X_hat}, which we recall here
\begin{equation*}
	\begin{aligned}
		\widehat{X}_{n+1}^{\mathrm{DLRA}, c, i}
		=
		\widehat{X}_{n+1}^{\mathrm{DLRA},i}
		+ \mathcal{E}_{n+1}^{i},
	\end{aligned}
	\end{equation*}
where $\mathcal{E}_{n+1}^{i} = \sqrt{\frac{\overline{\Lambda_{n+1}}}{s}} Q \xi^{i}$, with $\xi^{i} \sim  \mathcal{N}(0, I_{s \times s})$ and $P_{\widehat{U}_{n+1}}^{\perp} Q = Q$, where $Q \in \mathbb{R}^{d \times s}$.

We approximate the filtering distribution inferred by \eqref{eq: compl X_hat} via a \emph{Gaussian mixture distribution}. Indeed for each particle, we suppose that our accumulated model error is added a posteriori on the predicted DLRA $\widehat{X}_{n+1}^{\mathrm{DLRA},i}$ as Gaussian perturbation of the $i$-th particle, i.e.\ the measure induced by the enriched DLRA is a Gaussian mixture which reads
\begin{equation}\label{eq: GMM DLRA}
	\pi_{	\widehat{X}_{n+1}^{\mathrm{DLRA}, c, i}} = \sum_{i=1}^{M} \lambda_{n}^{\mathrm{DLRA}, i} \mathcal{N}(\widehat{X}_{n+1}^{\mathrm{DLRA},i}, \widehat{E}_{n+1})
\end{equation}
where the weights $(\lambda_{n}^{\mathrm{DLRA}, i})_{i=1,\dots, M}$ satisfy $\sum_{i = 1}^M\lambda_{n}^{\mathrm{DLRA}, i}=1$.

Then, after the analysis step on \eqref{eq: compl X_hat}-\eqref{eq: GMM DLRA}, the updated equations for the mean of each particle $m_{n+1}^{\mathrm{DLRA},i}$ and the weights $\lambda_{n}^{\mathrm{DLRA}, i}$ are the following
\begin{equation}\label{eq: kalman update full GMM}
	\begin{aligned}
	m_{n+1}^{\mathrm{DLRA},i} := & \widehat{X}_{n+1}^{\mathrm{DLRA},i}+ \widehat{E}_{n+1} H(t_n)^{\top} R^{-1} \left( \Delta Z_n - H(t_n)\widehat{X}_{n+1}^{\mathrm{DLRA},i} \Delta t\right) \\
	\lambda_{n}^{\mathrm{DLRA}, i} :=& \frac{\widehat{\lambda}_{n+1}^{\mathrm{DLRA}, i}}{\sum_{j=1}^{M} \widehat{\lambda}_{n+1}^{\mathrm{DLRA}, j}}, \text{ with update}  \\
	\widehat{\lambda}_{n+1}^{\mathrm{DLRA}, i}  = &\lambda_n^{\mathrm{DLRA}, i} \exp\left( (H(t_n)X_{n+1}^{\mathrm{DLRA}, c, i})^{\top} R^{-1}\Delta Z_n - \frac{1}{2} |R^{-\frac{1}{2}}H(t_n)X_{n+1}^{\mathrm{DLRA}, c, i}|^2 \Delta t_n \right) \\
	E_{n+1} = & \widehat{E}_{n+1} - \widehat{E}_{n+1}  H(t_n)^{\top} R^{-1}H(t_n) \widehat{E}_{n+1} \Delta t,
	\end{aligned}
\end{equation}

One can notice that in \eqref{eq: kalman update full GMM} the update of the particles is not belonging to the $k$-rank manifold due to the influence of $\widehat{E}_{n+1}$, whose range belongs to $\mathrm{Im}(\widehat{U}_{n+1}^{\perp})$. Therefore, we want to derive the DLRA triple $(m^{\mathrm{DLRA}, c}_{n+1}, U_{n+1}^{c}, Y_{n+1}^{c})$ in order to restart the filtering procedure.
First, notice that the update of the mean of the particles in \eqref{eq: kalman update full GMM} can be rewritten as
\begin{equation}\label{eq: expanded dlra gmm}
	\begin{aligned}
	m_{n+1}^{\mathrm{DLRA},i} := & \widehat{X}_{n+1}^{\mathrm{DLRA},i}+ \widehat{E}_{n+1} H(t_n)^{\top} R^{-1} \left( \Delta Z_n - H(t_n)\widehat{X}_{n+1}^{\mathrm{DLRA},i} \Delta t\right) \\
	 =& \widehat{m}_{n+1} + \widehat{U}_{n+1}^{\top} \widehat{Y}_{n+1}^{i} + \widehat{E}_{n+1} H(t_n)^{\top} R^{-1} \left( \Delta Z_n - H(t_n)\left(\widehat{m}_{n+1} + \widehat{U}_{n+1}^{\top} \widehat{Y}_{n+1}^{i}\right) \Delta t \right) \\
	\end{aligned}
\end{equation}
Via applying the expectation to relation \eqref{eq: expanded dlra gmm}, one obtains the update for the mean:
\begin{equation}\label{eq: DLRA GMM mean}
	m_{n+1}^{\mathrm{DLRA}} = \widehat{m}_{n+1} + \widehat{E}_{n+1} H(t_n)^{\top} R^{-1} \left( \Delta Z_n - H(t_n)\widehat{m}_{n+1} \Delta t \right).
\end{equation}

In order to find an update for $U_{n+1}^{c}$ and $Y_{n+1}^{c,i}$, by ansatz we suppose that $U_{n+1}^{c} = \widehat{U}_{n+1} + A P_{\widehat{U}_{n+1}}^{\perp}$ for a matrix $A \in \mathbb{R}^{k \times d}$.
Then subtracting \eqref{eq: expanded dlra gmm} from \eqref{eq: DLRA GMM mean}, one has
\begin{equation}\label{eq: interm GMM uy}
	(U_{n+1}^{c})^{\top}Y_{n+1}^{c,i} =\widehat{U}_{n+1}^{\top} \widehat{Y}_{n+1}^{i} + \widehat{E}_{n+1} H(t_n)^{\top} R^{-1} \left( \Delta Z_n - H(t_n)\left(\widehat{U}_{n+1}^{\top} \widehat{Y}_{n+1}^{i}\right) \Delta t \right) \\
\end{equation}
Via multiplying relation \eqref{eq: interm GMM uy} on the left by $ \widehat{U}_{n+1}$, using the orthogonality of $ \widehat{U}_{n+1}$,  $\widehat{U}_{n+1}P_{\widehat{U}_{n+1}}^{\perp}=0$, and $\widehat{U}_{n+1} E_{n+1} =0$, one obtains that $Y^{c,i}_{n+1} = \widehat{Y}_{n+1}^{i}$ for all $i=1,\dots,M$. Finally, multiplying \eqref{eq: interm GMM uy} on the right by $(\widehat{Y}_{n+1}^{i})^{\top}$, taking the expectation, and dividing by $\mathbb{E}\left[\widehat{Y}_{n+1}\widehat{Y}_{n+1}^{\top}\right]$, one finally gets
\begin{equation}\label{eq: DLRA GMM U}
	U_{n+1}^{c} = \widehat{U}_{n+1}^{\top}  - \widehat{U}_{n+1}  H(t_n)^{\top} R^{-1} H(t_n) \widehat{E}_{n+1} \Delta t.
\end{equation}

Notice that even though the equation of the $Y$ does not change between the prediction and the analysis, actually its distribution does due to the update of the weights $\lambda_{n}^{\mathrm{DLRA}, i}$ associated to each particle $Y_{n}^{i}$. Therefore, each component of the DLRA is updated during the analysis step.
Our whole procedure algorithm is described in Algorithm \ref{alg: DLRA PF}, where we highlight again that we compute our drift $A(t,x)$ only in the DLRA point $X_{n}^{\mathrm{DLRA}}$, assuming that the discrepancy between $A(t_n, X_{n}^{\mathrm{DLRA}})$ and $A(t_n, X_{n}^{\mathrm{DLRA},c})$ can be captured by analysis on the model error. The effectiveness of the DLRA relies on the assumption that the discrepancy \(\varepsilon\) defined in \eqref{eq: eps_n} between the drift and diffusion coefficients and their projections onto the tangent space of the rank-\(k\) manifold remains sufficiently small along the numerical solution.
\begin{algorithm}
	\caption{ (Complemented) DLRA Particle Filter}\label{alg: DLRA PF}
	
	\begin{flushleft}
		\textbf{Input}: initial data ($m_0$, $U_0$, $Y_0$), number of samples $M$, $(\lambda_0^{\mathrm{DLRA}, i} )_{i=1,..,M}=\frac{1}{M}$ with $i \in \{1,\dots,M\}.$
		
		\textbf{Output:} approximation  $\{X^{\mathrm{DLRA},c,i}_n\}_{n=0,\ldots, N}$ for $i=1,\dots,M$. 
	\end{flushleft}
	
	\begin{algorithmic}[1]
		\ForAll {$n \in \{0, \ldots, N-1\}$}
		\State Generate  $\Delta W_n^{i} \sim \mathcal{N}(0,\Delta t_n I_{m \times m})$  with $i \in \{1,\dots,M\}.$
		\State Compute $\widehat{m}_{n+1}^{\mathrm{DLRA}}
		= m_n^{\mathrm{DLRA}}
		+ \left(\frac{1}{M} \sum_{i=1}^{M}A(t_n,  m_n^{\mathrm{DLRA}} + U_n^{\top}Y_n^{i})\right)\, \Delta t_n.$		
		\State Compute $\widetilde{Y}^{i}_{n+1} = Y_n^{i} + U_n \mathring{A}\!\left(t_n, m_n^{\mathrm{DLRA}} + U_n^{\top}Y_n^{i}\right) \Delta t_n +  U_n Q^{\frac{1}{2}} \Delta W_n^{i}$ 
		\State Assemble $\overline{C}_{\widetilde{Y}_{n+1}} = \frac{1}{M-1} \sum_{i=1}^{M}\widetilde{Y}_{n+1}^{i}(\widetilde{Y}_{n+1}^{i})^{\top}$
		and compute $\widetilde{U}_{n+1}$: 
		$$\overline{C}_{\widetilde{Y}_{n+1}} \widetilde{U}_{n+1}
		= \overline{C}_{\widetilde{Y}_{n+1}}  U_n
		+  \frac{1}{M-1} \sum_{i=1}^{M-1}\left(\widetilde{Y}_{n+1}^{i}\, \mathring{A}(t_n,  m_n^{\mathrm{DLRA}} + U_n^{\top}Y_n^{i})^\top\right)
		\, P_{U_n}^\perp\, \Delta t_n
		+ U_n Q\, P_{U_n}^\perp\, \Delta t_n. $$
		\State Reorthonormalize the deterministic basis: find   $(\widehat{U}_{n+1}, \widehat{Y}_{n+1})$ such that:
		\begin{equation*}
			\widehat{U}^{\top}_{n+1} \widehat{Y}_{n+1}=\widetilde{U}_{n+1}^{\top} \widetilde{Y}_{n+1}, \quad \widehat{U}_{n+1}\widehat{U}^{\top}_{n+1} = I_{k\times k}.
		\end{equation*}
		via $(\widehat{U}^{\top}_{n+1}, R_{n+1}) = \texttt{QR}(\widetilde{U}^{\top}_{n+1})$ and $\widehat{Y}_{n+1}^{i} = R_{n+1} \widetilde{Y}^{i}_{n+1}$, for all $i \in \{1,\dots,M\}.$
		\State Recentering: $\widehat{Y}_{n+1}^{i} = \widehat{Y}_{n+1}^{i} - \frac{1}{M} \sum_{j= 1}^{M}  \widehat{Y}_{n+1}^{j}$, $\widehat{m}_{n+1} = \widehat{m}_{n+1} + \widehat{U}_{n+1}^{\top} \left(\frac{1}{M} \sum_{j= 1}^{M}   \widehat{Y}_{n+1}^{j}\right)$

		\If{ Complemented}
		\State Compute $\varepsilon_n$ according to \eqref{eq: eps_n},  set $\Lambda_{n+1}$ as \eqref{eq: Lambda_n} and $\overline{\Lambda}_{n}=\sum_{\ell=1}^{n+1}\Lambda_{\ell}$, build $\xi_n^{i} =\sqrt{\frac{\overline{\Lambda}_n}{s}}\,Q\zeta_i,$ with $\zeta_i\sim\mathcal{N}(0,I_{s\times s})$, and set 
		\begin{equation*}
			X_{n+1}^{\mathrm{DLRA}, c, i}
		=
		\widehat{m}_{n+1}^{\mathrm{DLRA}} + \widehat{U}_{n+1}^{\top} \widehat{Y}_{n+1}^{i}
		+ P_{\widehat{U}_{n+1}}^\perp \xi_n^{i}
    \end{equation*}
    	\State Recentering: $\xi_{n}^{i} = \xi_{n}^{i} - \frac{1}{M} \sum_{j= 1}^{M}  \xi_{n}^{j}$, $\widehat{m}_{n+1} = \widehat{m}_{n+1} +   P_{\widetilde{U}_{n+1}}^\perp  \left(\frac{1}{M} \sum_{j= 1}^{M}  \xi_{n}^{j}\right) $
    	\State Compute (and then set $\widetilde{U}_{n+1}= U_{n+1}$)
    \begin{equation*}
    	\begin{aligned}
    	m_{n+1}^{\mathrm{DLRA}} = &\widehat{m}_{n+1} + \widehat{E}_{n+1} H(t_n)^{\top} R^{-1} \left( \Delta Z_n - H(t_n)\widehat{m}_{n+1} \Delta t \right), \quad \\
    	U_{n+1} = & \widehat{U}_{n+1}^{\top}  - \widehat{U}_{n+1}  H(t_n)^{\top} R^{-1} H(t_n) \widehat{E}_{n+1} \Delta t. \\
    	\end{aligned}
    \end{equation*}
    
		\Else
		\begin{equation*}
		X_{n+1}^{\mathrm{DLRA}, c, i}
		=
		\widehat{m}_{n+1}^{\mathrm{DLRA}} + \widehat{U}_{n+1}^{\top} \widehat{Y}_{n+1}^{i}
	\end{equation*}
			\EndIf		
		 \State Compute $\lambda_{n+1}^{\mathrm{DLRA}, i} = \frac{\widehat{\lambda}_{n+1}^{\mathrm{DLRA}, i} }{\sum_{j=1}^{M} \widehat{\lambda}_{n+1}^{\mathrm{DLRA}, j}}$ for $i=1,..,M$ with
		 	\begin{equation*}
		 			\widehat{\lambda}_{n+1}^{\mathrm{DLRA}, i}  = \lambda_n^{\mathrm{DLRA}, i} \exp\left( (H(t_n)X_{n+1}^{\mathrm{DLRA}, c, i})^{\top}R^{-1} \Delta Z_n - \frac{1}{2} |R^{-\frac{1}{2}}H(t_n)X_{n+1}^{\mathrm{DLRA}, c, i}|^2 \Delta t_n \right)
		 	\end{equation*}
		 and set the empirical filtering measure to 
		 
		 $$\pi_{n+1}^{\mathrm{DLRA}, c ,M}(\mathrm{d}x)=\frac{1}{M}\sum_{i=1}^M 	\lambda_{n+1}^{\mathrm{DLRA}, i} \delta_{	X_{n+1}^{\mathrm{DLRA}, c, i}}(\mathrm{d}x).$$
		
		\If {$\Delta t_{\mathrm{branch}} \mod \Delta t_n$ $=$ $0$}
		
		 $m_{n+1}^{\mathrm{DLRA}}+(\widehat{U}_{n+1}^{\top}\widehat{Y}_{n+1}^{i})_{i=1,..,M},
		 (\lambda_{n+1}^{\mathrm{DLRA}, i})_{i=1,..,M} = \mathrm{branching}\left( 	X_{n+1}^{\mathrm{DLRA}, c, i}, (\lambda_{n+1}^{\mathrm{DLRA}, i} )_{i=1,..,M} \right)$
			\EndIf

			\State Reorthonormalization: $(U_{n+1}, Y_{n+1})$ such that $U^{\top}_{n+1} Y_{n+1} =\widehat{U}_{n+1}^{\top} \widehat{Y}_{n+1}$, $U_{n+1}U_{n+1}^{\top} = I_{k\times k}$.
		via $(U_{n+1}^{\top}, R_{n+1}) = \texttt{QR}(\widehat{U}^{\top}_{n+1})$ and $Y_{n+1}^{i} = R_{n+1} \widehat{Y}^{i}_{n+1}$, for all $i \in \{1,\dots,M\}.$
		\EndFor
	\end{algorithmic}
\end{algorithm}

We conclude this section by showing the equations for the assumption of isotropic noise in the orthogonal complement in a first order of $O(\Delta t)$. The prediction triple $(\widehat{m}_{n+1}, \widehat{U}_{n+1}, \widehat{Y}_{n+1})$ is updated by the usual DLR Projector Splitting for SDEs \eqref{eq: DLRA JMCO prediction}, whereas the analysis step reads as 
\begin{equation}
	\begin{aligned}
	m_{n+1}^{\mathrm{DLRA}} = &\widehat{m}_{n+1} + \overline{\Lambda}_{n}P_{U_n}^{\perp} H(t_n)^{\top} R^{-1} \left( \Delta Z_n - H(t_n)\widehat{m}_{n+1} \Delta t \right) \\
	U_{n+1}= & \widehat{U}_{n+1}  - \overline{\Lambda}_{n} U_{n}  H(t_n)^{\top} R^{-1} H(t_n)  P_{U_n}^{\perp} \Delta t \\
	Y_{n+1}^{i} =& \widehat{Y}_{n+1}^{i}, \quad \lambda_{n+1}^{\mathrm{DLRA}, i} = \frac{\widehat{\lambda}_{n+1}^{\mathrm{DLRA}, i} }{\sum_{j=1}^{M} \widehat{\lambda}_{n+1}^{\mathrm{DLRA}, j}} \\
\widehat{\lambda}_{n+1}^{\mathrm{DLRA}, i}  = &\lambda_n^{\mathrm{DLRA}, i} \exp\left( (H(t_n)X_{n+1}^{\mathrm{DLRA}, c, i})^{\top}R^{-1} \Delta Z_n - \frac{1}{2} |R^{-\frac{1}{2}}H(t_n)X_{n+1}^{\mathrm{DLRA}, c, i}|^2 \Delta t_n \right).
	\end{aligned}
\end{equation} 
Notice that if $H \in \mathbb{R}^{d \times h}$ with $h \ll d$, then $H(t_n)P_{U_n}^{\perp}$ is still computationally affordable. Therefore, in case we do not have any further information on the orthogonal complement, a basis isotropic noise assumption may be again advantageous.

\begin{Remark}[The case of discrete-time observations]
	The strategy presented for the DLR PF and its possible complementation are presented for the case of continuous-time observations. We propose also in this case a variant for the case of discrete-time observations. In case of observations that follows \eqref{eq: discr obs cont}, when $\Delta t_n \mod \Delta t_a = 0$, where $\Delta t_a$ is time-step of assimilation of the observations,  our DLRA weights are given by 
	\begin{equation*}
		\widehat{\lambda}_{n+1}^{\mathrm{DLRA}, i}  = \lambda_n^{\mathrm{DLRA}, i} \exp\left( -\frac{1}{2}  (Z_{n+1} - H(t_{n+1})X_{n+1}^{\mathrm{DLRA}, c, i})^{\top} R^{-1} (Z_{n+1} - H(t_{n+1})X_{n+1}^{\mathrm{DLRA}, c, i})\right).
	\end{equation*}
	Otherwise, the analysis step does not take place, and the DLRA updates remain the same of the prediction step.
\end{Remark}

\section{Numerical experiments}\label{sec: numerics filt}

In this section, we assess the performance of the proposed filtering algorithms on some numerical test cases. As reference algorithms to compare our DLRA Filters with, we consider a \emph{particle filter} \cite[Section 9]{bain2009fundamentals} (PF) and a \emph{continuous-time ensemble Kalman filter} for stochastic differential equations \cite[Result 8.4]{law2015data} (EF SDE). These algorithms will be computed with a much larger number of particles than the DLRA ones. The observation $Z$ will be generated from a single trajectory $X^{\ast}$, simulated with the same time step also used by the filters. 

We will consider several types of errors, all being averaged over the time interval $[T_0,T]$ which disregards an initial transient dynamics. This way of computing errors is standard practice in data assimilation, as filtering procedures need a warm-up time to eliminate the possible additional bias due to their initialization. 
To analyze the joint-mean and covariance error considered in Section \ref{sec: DLRA for JMCO}, we will plot the error between the conditional mean of the DLRA filters and the one of the reference PF or EF SDE, as well as the error between their conditional covariances, as relative error quantities normalized with respect to the value of the reference algorithms PF or EF SDE, respectively.
In details, we consider the \emph{averaged-over-time relative error between the mean} $(m_s)_{s \in [0, T]}$ and the \emph{covariance} $(C_s)_{s \in [0, T]}$ of the DLRA filters versus the reference ones $(m_s^{\mathrm{reference}})_{s \in [0, T]}$ and $(C_s^{\mathrm{reference}})_{s \in [0, T]}$, namely
\begin{equation}\label{eq: avg rel mean and cov}\tag{E.1}
	\frac{1}{T-T_0}	\int_{T_0}^{T}\frac{|m^{\mathrm{reference}}_s-m_s|}{|m^{\mathrm{reference}}_s|} \mathrm{d}s, \quad \frac{1}{T-T_0}\int_{T_0}^{T}\frac{\|C^{\mathrm{reference}}_s-C_s\|_{\mathrm{F}}}{\|C^{\mathrm{reference}}_s\|_{\mathrm{F}}} \mathrm{d}s,
\end{equation}
respectively.

Furthermore, we will also analyze the relative Euclidean error between the true state $X^{\ast}_t$ generating the observations $Z$, and the conditional mean of a filtering algorithm (for instance $m_{t}^{\mathrm{DLRA}}$ for any DLRA-type filter), because the latter is considered as the best estimate of the true solution given the observations.  In the caption of the error plots, the conditional mean $m_s$ will be denoted by $\overline{\widetilde{X}}(t)$ and we will plot the \emph{averaged-over-time relative root-mean square error (RMSE)}, which reads as
\begin{equation}\label{eq: avg rel RMSE}\tag{E.2}
	\frac{1}{T-T_0}\int_{T_0}^{T}\frac{|X^{\ast}_s- \overline{\widetilde{X}}_s |}{|X^{\ast}_s|} \mathrm{d}s,
\end{equation}
where the (component-wise) root-mean square error is usually defined as $\frac{|X^{\ast}_s- \overline{\widetilde{X}}_s |}{\sqrt{d}}$. Furthemore, in the legend of the plots inherent to \eqref{eq: avg rel RMSE}, ``True'' will indicate the true state $X^{\ast}$. 

The DLRA-filters are all discretized in time using the DLR PS SDE integrator, for instance as described in Algorithm \ref{alg: DLRA JMCO Filter algorithm}, which features good stability and convergence properties. On the other hand, the reference filters (PF and EF SDE) are discretized in time using an Euler-Maruyama, and an implicit treatment of the innovation term for the EF SDE.

In the legends of our plots, the DLRA-JMCO filter will be indicated simply by \emph{DLR}, the algorithm proposed in \cite{nobile2025dynamicallowrankapproximationskalman} by \emph{DLR Trad}, standing for ``traditional DLR filter'', the DLRA particle filter will be denoted as \emph{DLR PF}. 

%\iffalse
\subsection{Linear problem: Stochastic Advection-Diffusion-Reaction Model}

In this section, we analyze the performance of the DLRA-JMCO filter in the context of linear drift. We consider a noisy one-dimensional advection-diffusion-reaction equation with low-rank additive noise (white in time)
\begin{equation}\label{ex: SADR DA}
	\begin{cases}
		\begin{aligned}
			\partial_t u(x,t,\omega)&= L u(x,t,\omega) + \sum_{i=1}^{m} \phi_i(x) \dot{W}_i(\omega), \quad (x,t,\omega) \in [0,L] \times [0,T] \times \Omega, \\
			u_{0}(x,\omega)&=  \sum_{i = 1}^m \cos\left( \frac{\pi i x}{ L} \right) \mathcal{N}_{i}(0,1), \quad (x,\omega) \in [0,L] \times \Omega,\\
			\partial_x u(0,t, \omega)&= \partial_x u(1,t,\omega) = 0, \quad (t,\omega) \in [0,T] \times \Omega,
		\end{aligned}
	\end{cases}
\end{equation}
where in \eqref{ex: SADR DA} the linear operator $L$ is given by $L u(x,t,\omega) = a \partial_x^2 u(x,t,\omega) - v \partial_x u(x,t,\omega) + r u(x,t,\omega)$,  $\{\mathcal{N}_{i}(0,1)\}_{i=1,\dots,m}$ are i.i.d.\ normal random variables independent of the one-dimensional Brownian motions $\{W_i
\}_{i=1,\dots,m}$, $a=0.05$, $v= 0.05$, and $r = 0.02$. The spatial domain is $[0,L]$ with $L=5$, with Neumann boundary conditions at the endpoints, while the temporal domain is $[0,T]$, with $T=20$.  Concerning the noise term, we set $\phi_{i}(x) = \cos(\frac{i \pi f x}{L})$ for $i=1,\cdots, m$ and $f=5$. We discretize in space \eqref{ex: SADR DA} by second order centered finite differences, with a first order upwinding treatment of the advection term, using a uniform grid with mesh size $\Delta x = 0.1$. This results in a system of SDEs of dimension $d=51$ and we denote by $U^{\Delta x}_t$ the state vector containing the spatially discretized components of $u$ at time $t$. The time discretization to generate the true state $u^{\ast}$ and the observation is an Euler-Maruyama with uniform mesh size $\Delta t =5\cdot 10^{-3}$. We choose $m=12$ and we employ a DLR approximation of rank $k=12$. 

We consider the following process for the observation
\begin{equation}\label{eq: SADR obs process 1}
	\mathrm{d} Z_{t} = H U^{\Delta x}_t \mathrm{d}t + R^{\frac{1}{2}} \mathrm{d}B_t
\end{equation}
We observe the state in $10$ equidistant grid points (discarding the last endpoint), hence the dimension $h$ of the observation vector $Z_{t}$ is $h=10$. The covariance of the observation noise is $R= 0.01 \cdot I_{h \times h}$.
The number of simulated paths for the reference algorithms is $M=10^5$.  Furthermore, to compute all the average-in-time errors we consider a warm-up time of $T_0=2$.

As we are considering a linear S(P)DE problem with a Gaussian initial condition, we know that the measure induced by the problem is Gaussian at all times and an ensemble Kalman filter for SDEs \cite{law2015data} converges to true filtering distribution as the number of particles goes to infinity. As the sum of the dimensions of the noise and the initial condition is low, we are expecting a good performance of the DLRA-JMCO filter even with relatively small ranks. Concerning the DLR PF, we employ a branching procedure with a time-step $\Delta t_{\mathrm{branch}} = 4 \Delta t$.

\begin{figure}[!h]
	\centering
	\includegraphics[scale=0.046]{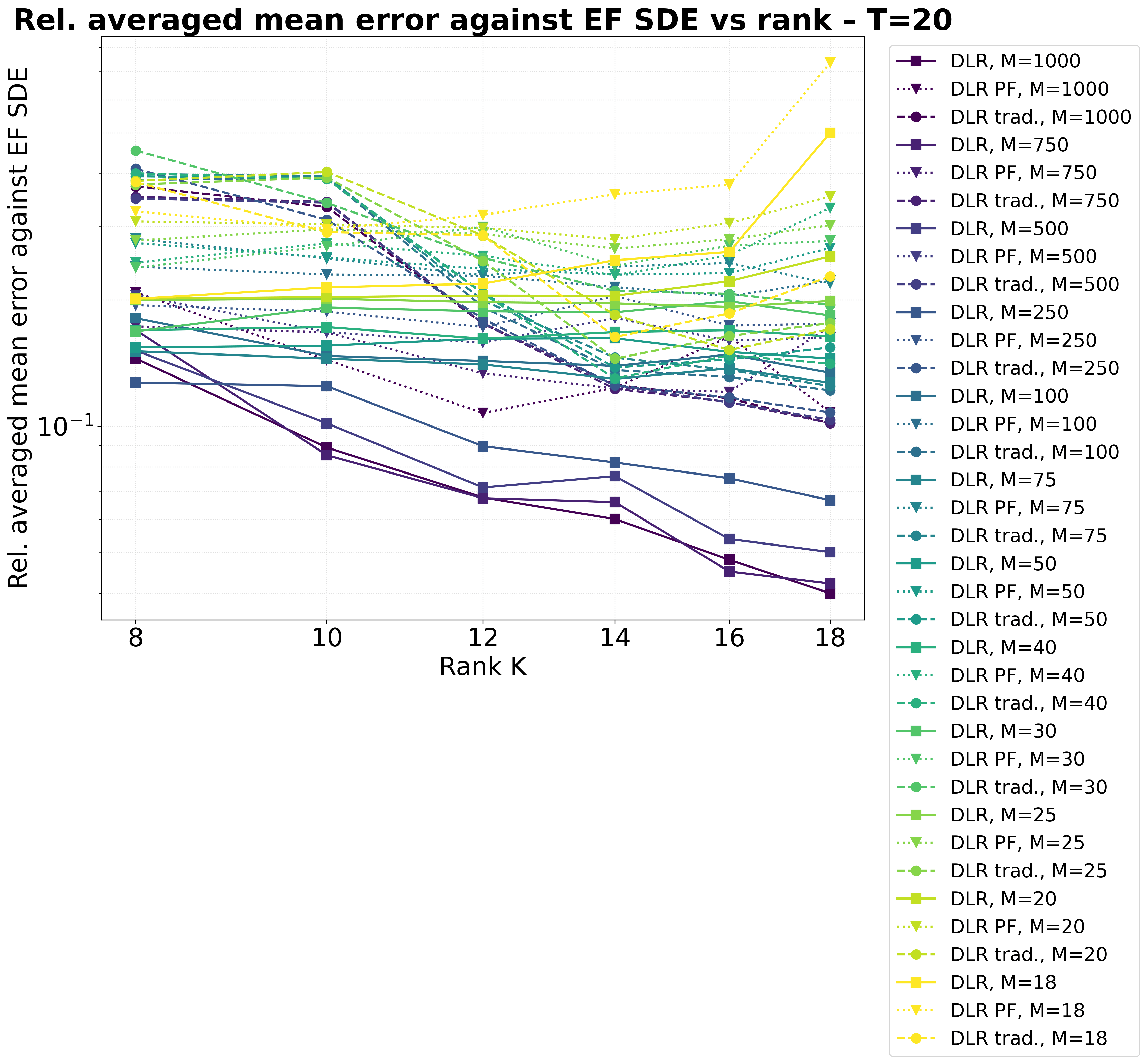}
	\includegraphics[scale=0.046]{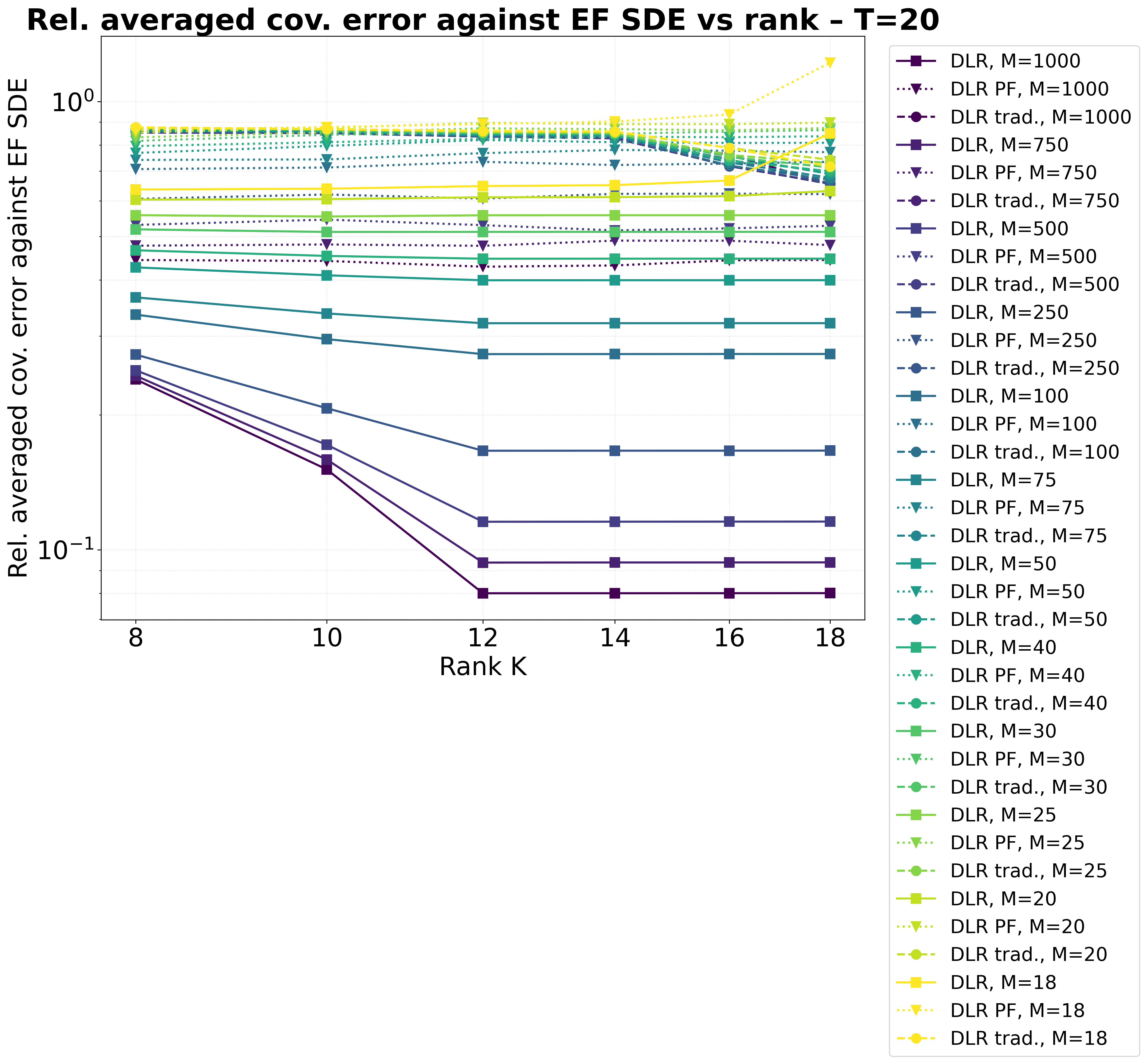} 
	\caption{Rel.\  averaged-over-time errors against the ensemble Kalman filter (EF SDE) as a function of the rank $K$ and for various sample sizes $M$, (left) mean, (right) covariance. Problem \eqref{ex: SADR DA} with continuous observations \eqref{eq: SADR obs process 1}, $\Delta t_{branch}= 4\Delta t$, $d=51$.}
	\label{fig: SADR CDCO EF SDE rank}
\end{figure}

\begin{figure}[!h]
	\centering
	\includegraphics[scale=0.046]{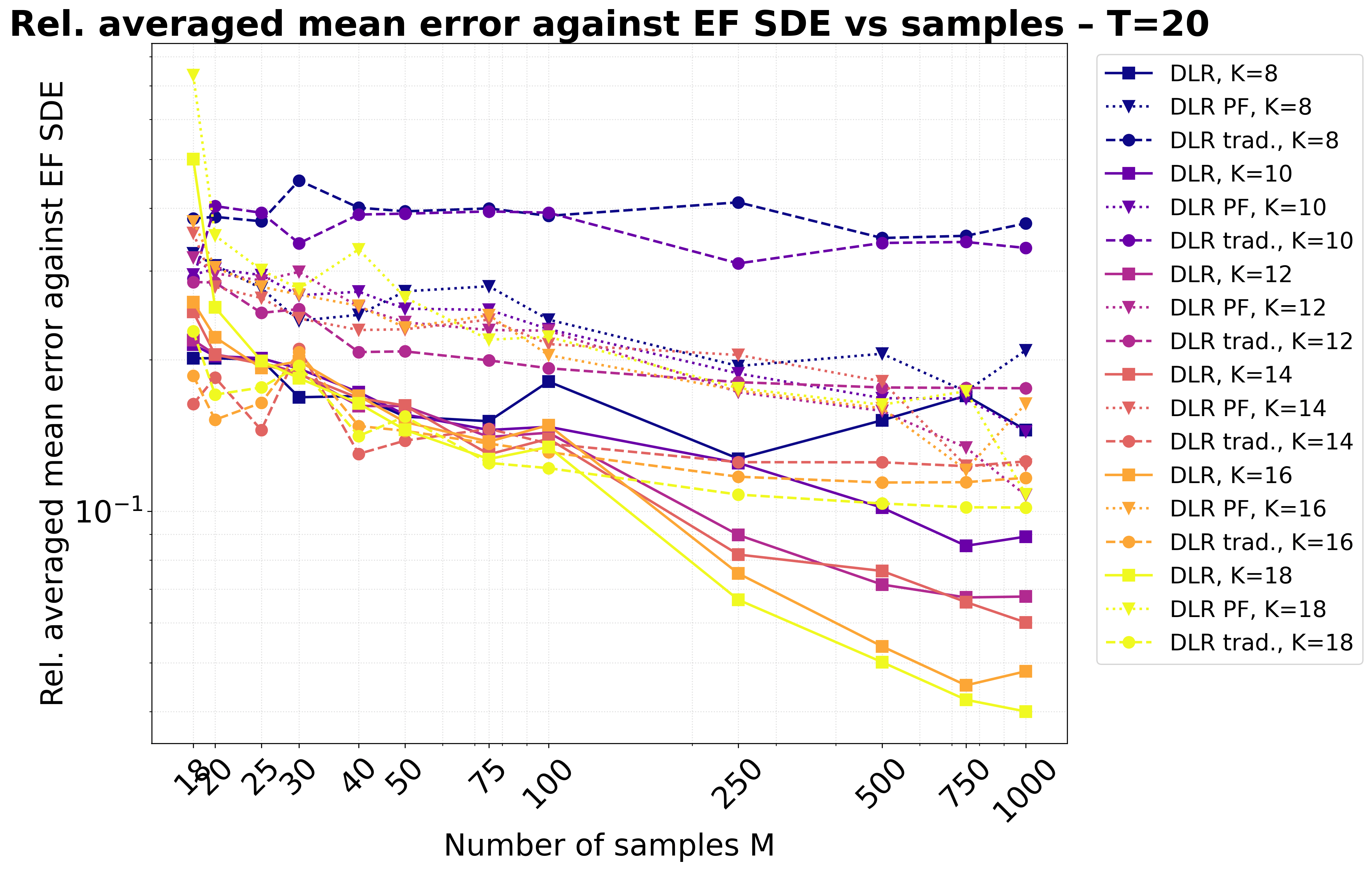}
	\includegraphics[scale=0.046]{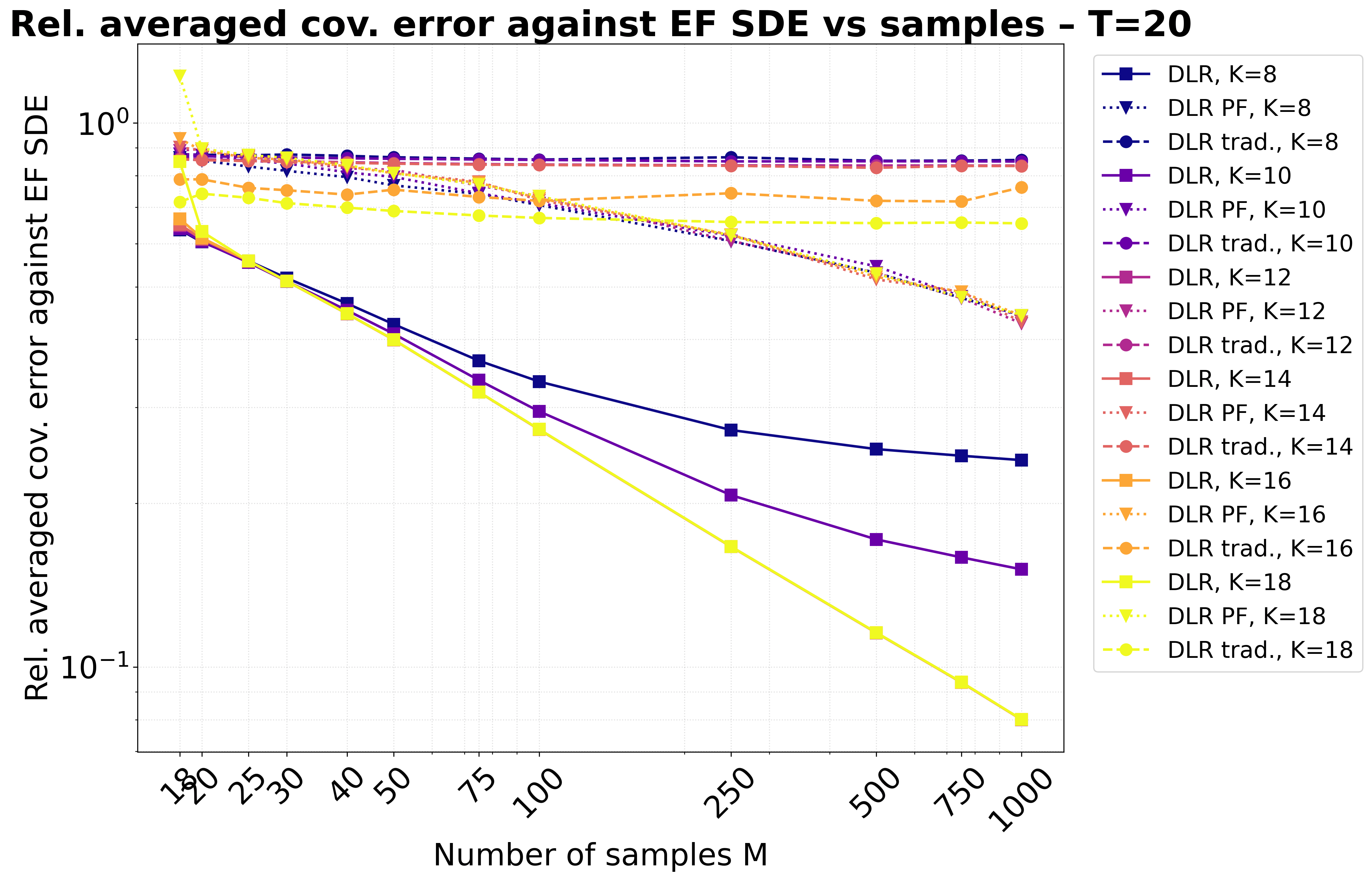} 
	\caption{Rel.\  averaged-over-time errors against the ensemble Kalman filter (EF SDE) as a function of the sample size $M$ and for various ranks $K$, (left) mean, (right) covariance. Problem \eqref{ex: SADR DA} with continuous observations \eqref{eq: SADR obs process 1}, $\Delta t_{branch}= 4\Delta t$, $d=51$.}
	\label{fig: SADR CDCO EF SDE samples}
\end{figure}

\begin{figure}[!h]
	\centering
	\includegraphics[scale=0.044]{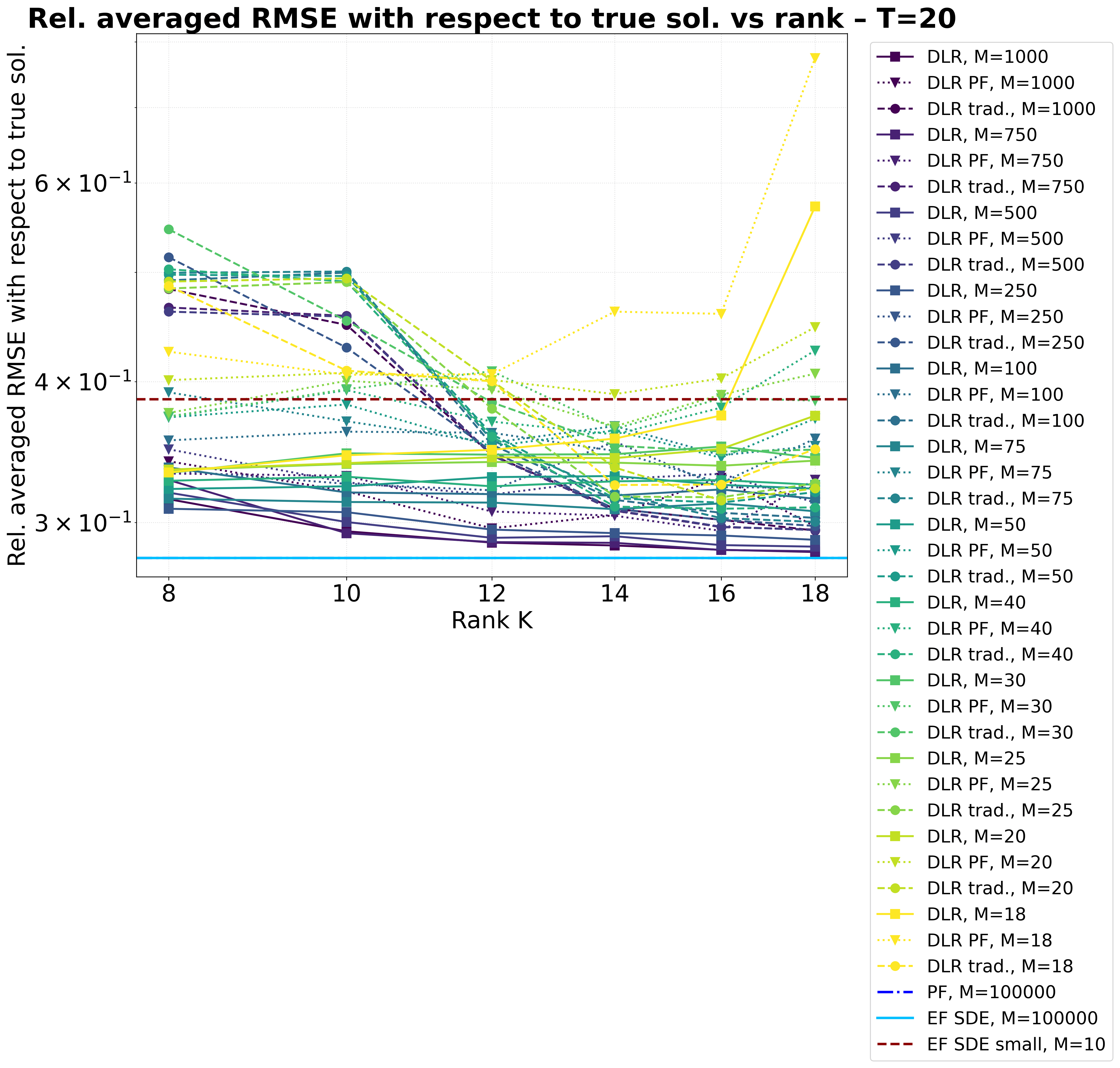}
	\includegraphics[scale=0.044]{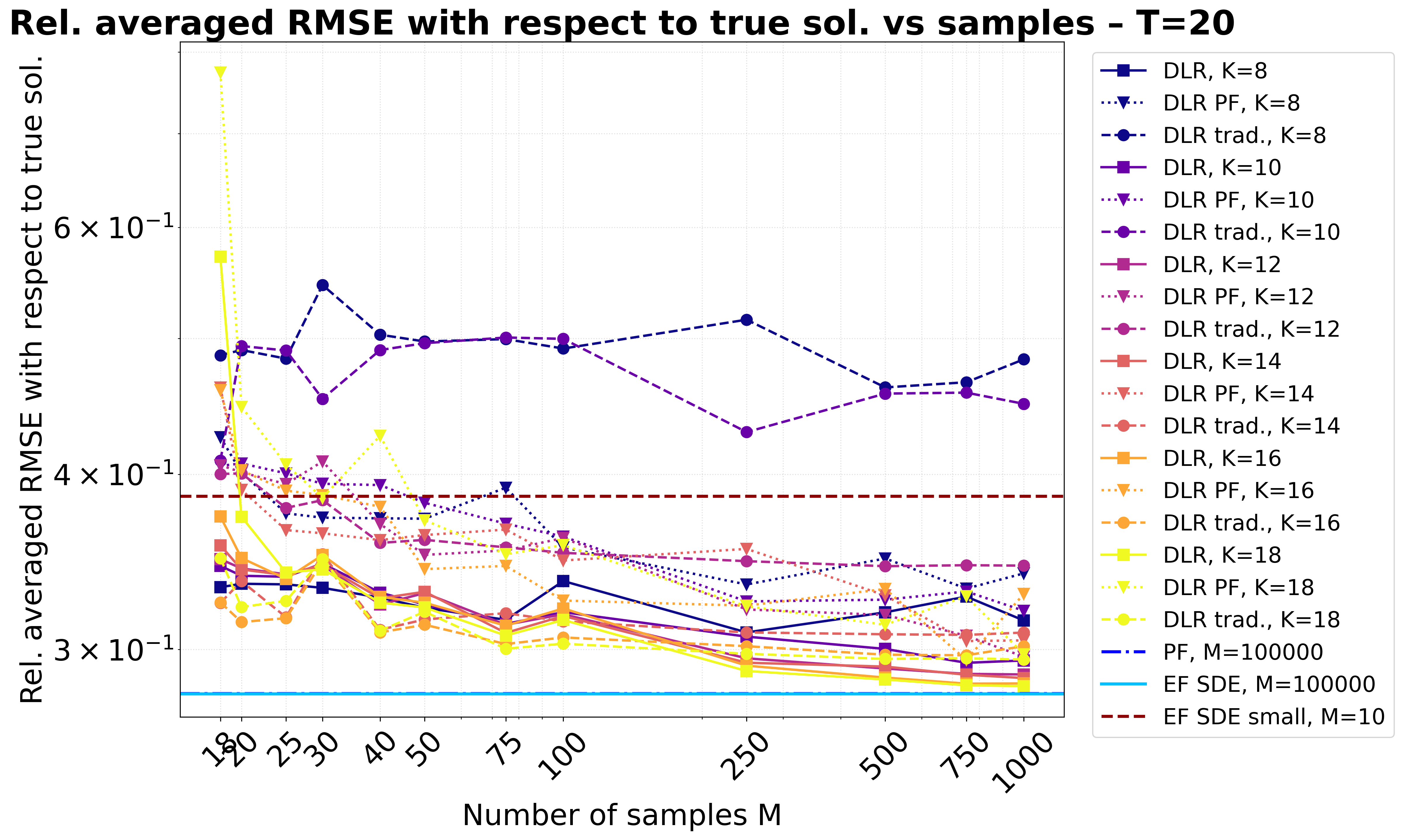} 
	\caption{Rel.\  averaged-over-time \emph{RMSE}, as a function of the rank $k$ (left) and sample size $M$ (right). Problem \eqref{ex: SADR DA} with continuous observations \eqref{eq: SADR obs process 1}, $d=51$, $\Delta t_{branch}= 4\Delta t$.}
	\label{fig: SADR CDCO RMSE}
\end{figure}

For the case of continuous-time observations, in Figures \ref{fig: SADR CDCO EF SDE rank} and \eqref{fig: SADR CDCO EF SDE samples}, we plot the error on the mean and the covariance with respect to the EF-SDE filter, for different ensemble sizes and DLRA ranks. We see that overall DLR JMCO filter has the best performance, outperforming DLR Trad., whereas the DLRA PF filter needs more and more particles to reach the same accuracy of the DLR JMCO. In Figure \ref{fig: SADR CDCO EF SDE rank}, the error on the covariance stagnates at $K=12$, suggesting that it is dominated by the Monte-Carlo error. Indeed, in Figure \eqref{fig: SADR CDCO EF SDE samples} we see that the error is still decreasing in $M$ for high ranks, indicating that we need more samples to fully resolve the DLRA main subspace. On the other hand, the error on the mean as well as the RMSE with respect to the true solution shown in Figure \ref{fig: SADR CDCO RMSE} present a decreasing behavior with respect to both the rank and the sample size, where in the latter we are approaching the performance of the reference EF SDE. Concerning the RMSE in Figure \ref{fig: SADR CDCO RMSE}, one can also notice that an EF SDE with a small number of particles ($M=10$) is outperformed by all the DLRA JMCO with a similar number of particles ($M=18$) but with smaller rank $K=8,10$ and hence maximum dimension that the surrogate can capture. Furthermore, we also notice that for ranks higher than $K=12$, the errors of the DLRA tend to increase when a small number of samples is used: this is a stability issue due to the fact that we need a sufficient number of samples to well represent the subspace and avoid degeneracy of the covariance matrix that has to be inverted in the equation of the deterministic modes \eqref{eq: DLRA JMCO}.

\begin{figure}[!h]
	\centering
	\includegraphics[scale=0.046]{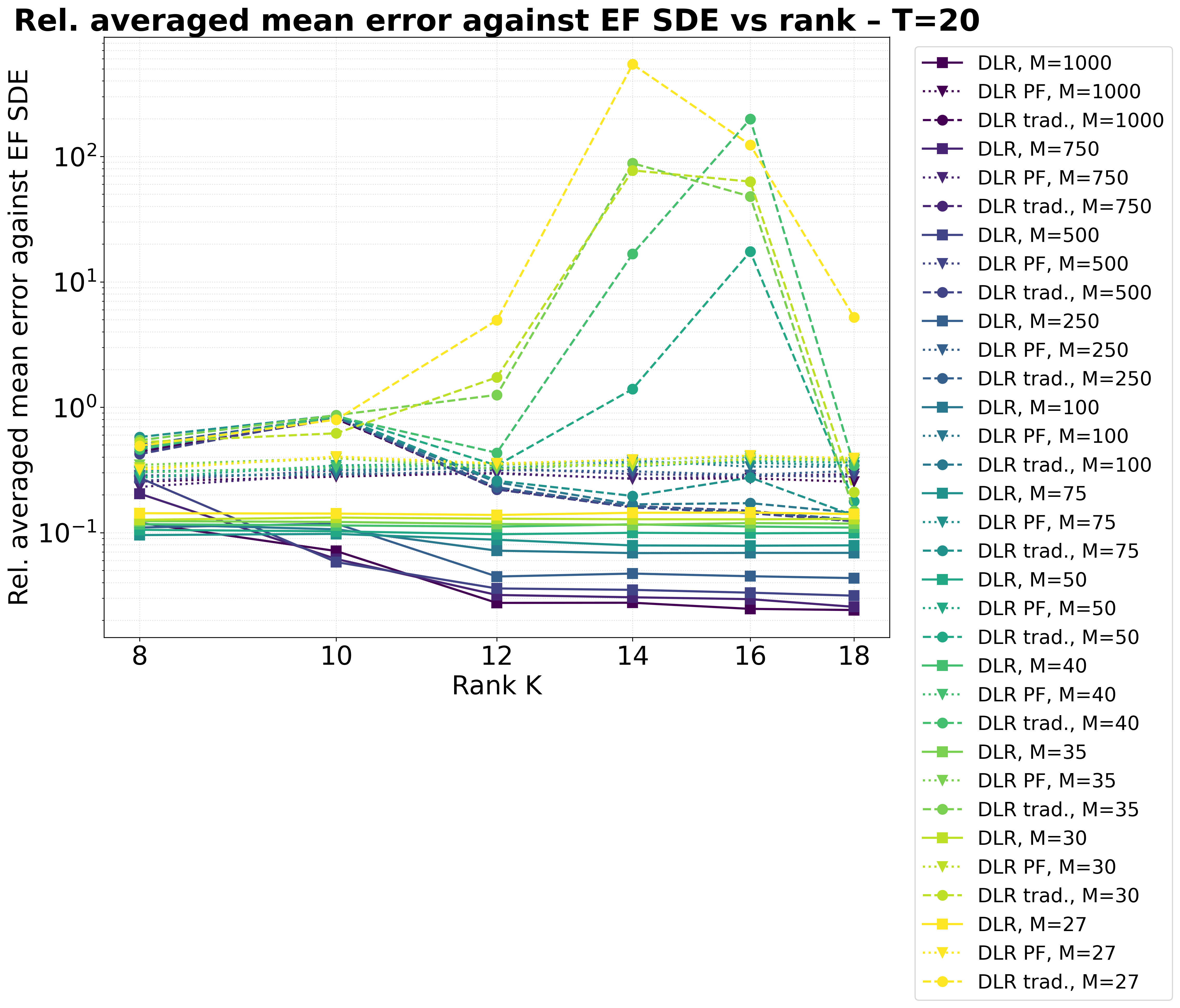}
	\includegraphics[scale=0.046]{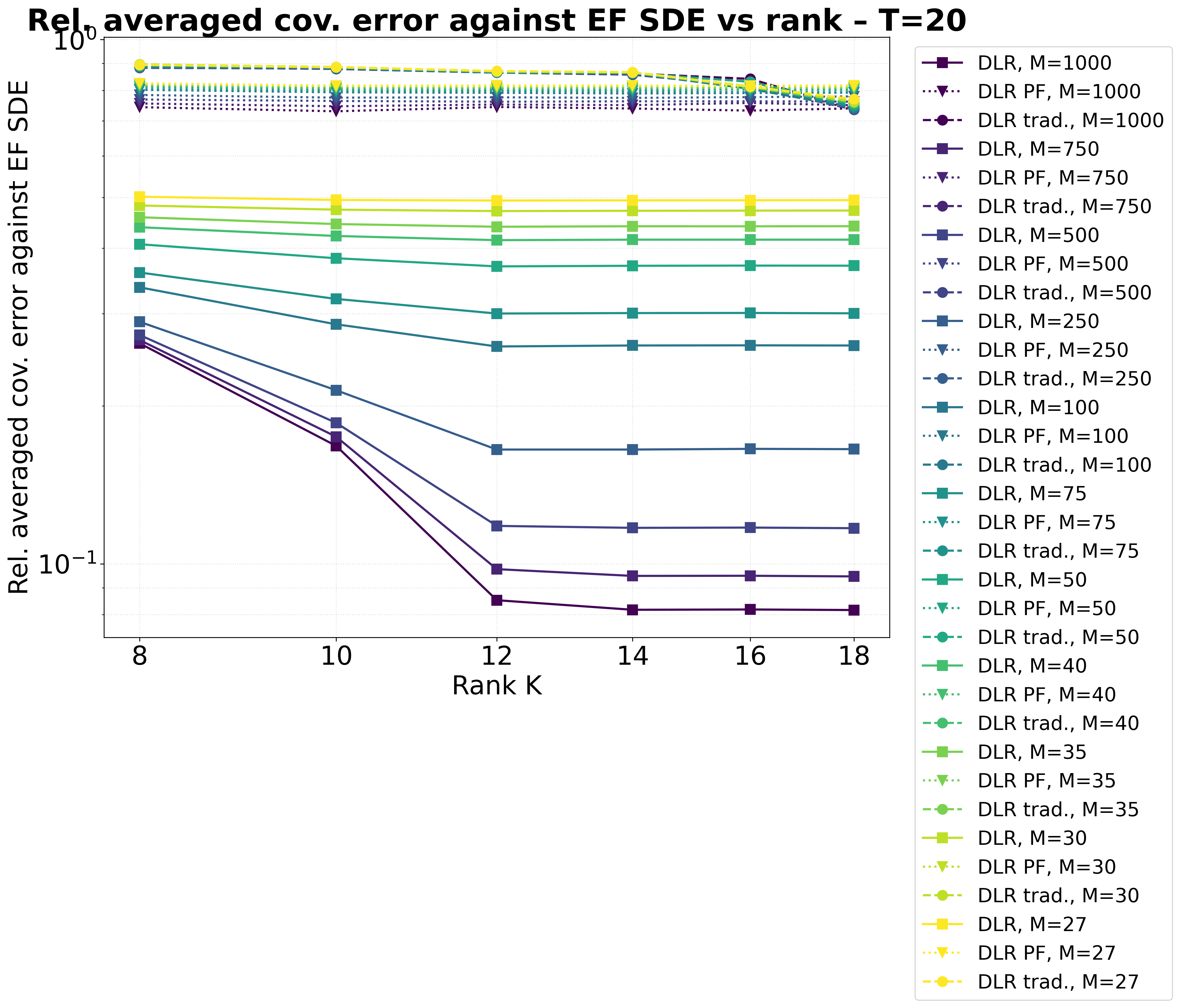} 
	\caption{Rel.\  averaged-over-time errors against the ensemble Kalman filter (EF SDE) as a function of the rank $K$ and for various sample sizes $M$, (left) mean, (right) covariance. Problem \eqref{ex: SADR DA} with discrete observations from the observation process \eqref{eq: SADR obs process 1},  $\Delta t_a = 10 \Delta t$, $d=51$.}
	\label{fig: SADR CDDO EF SDE rank}
\end{figure}

\begin{figure}[!h]
	\centering
	\includegraphics[scale=0.046]{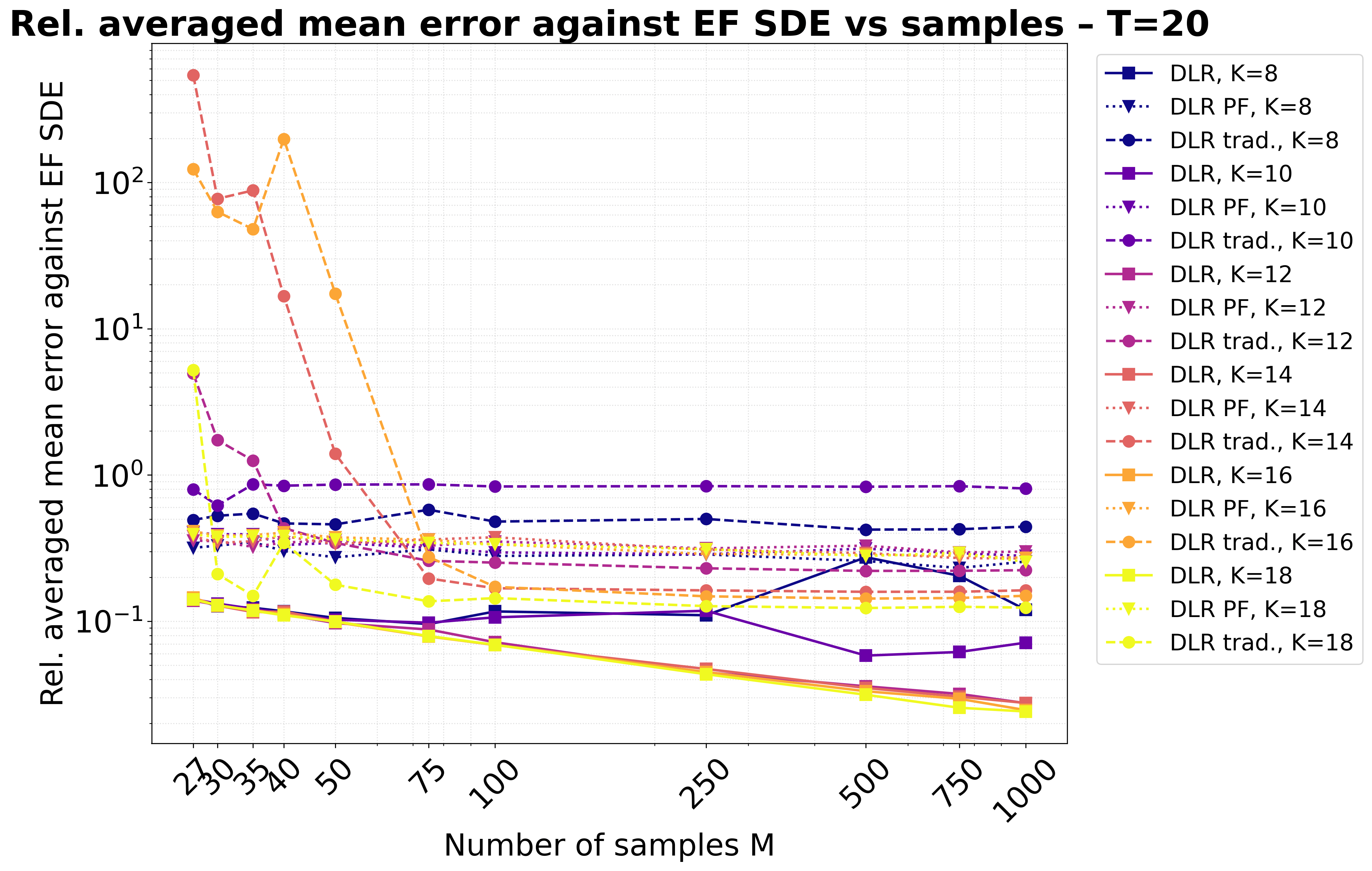}
	\includegraphics[scale=0.046]{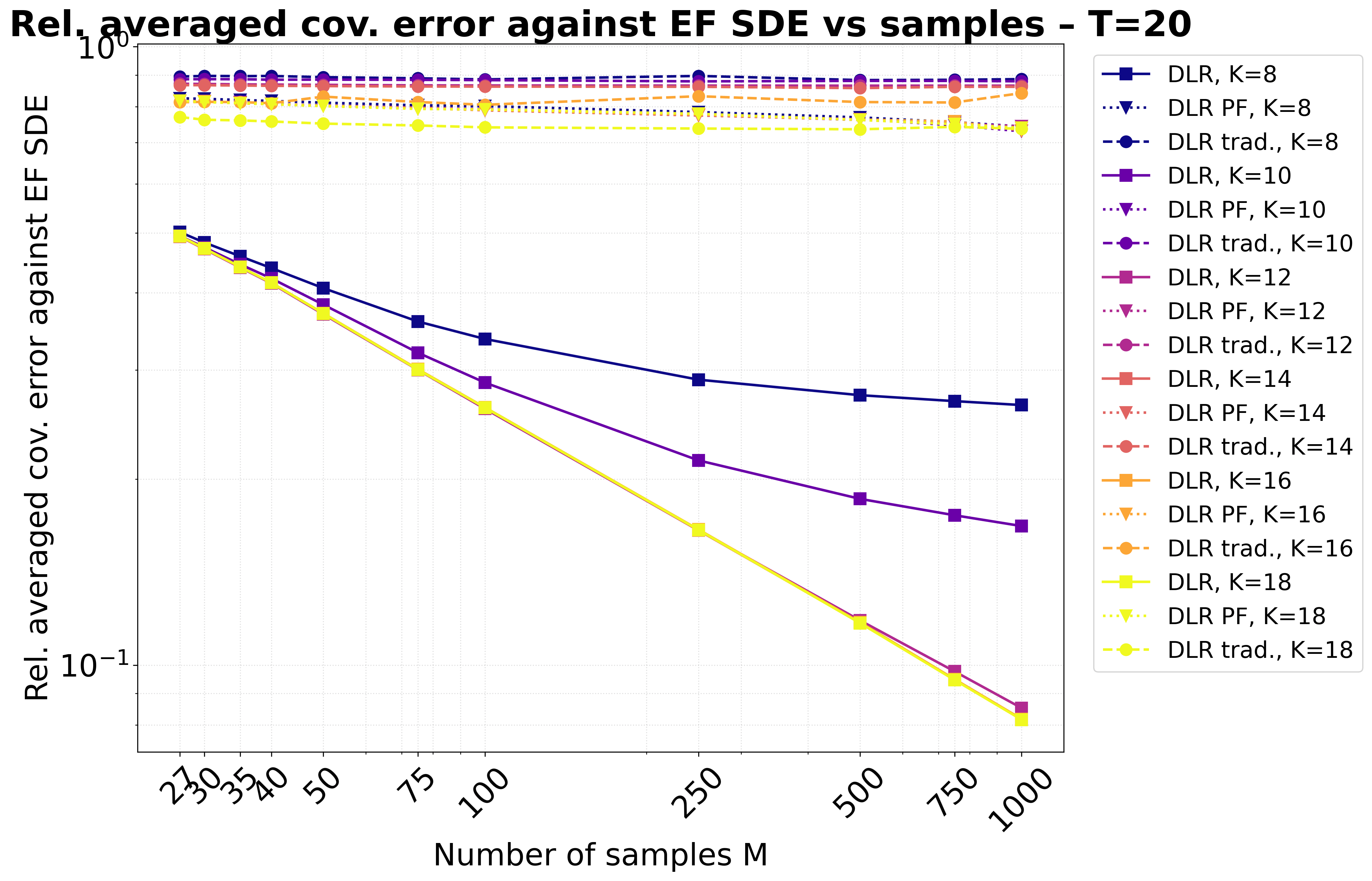} 
	\caption{Rel.\  averaged-over-time errors against the ensemble Kalman filter (EF SDE) as a function of the sample size $M$ and for various ranks $K$, (left) mean, (right) covariance. Problem \eqref{ex: SADR DA} with discrete observations obtained from the observation process \eqref{eq: SADR obs process 1},  $\Delta t_a = 10 \Delta t$, $d=51$.}
	\label{fig: SADR CDDO EF SDE samples}
\end{figure}

\begin{figure}[!h]
	\centering
	\includegraphics[scale=0.044]{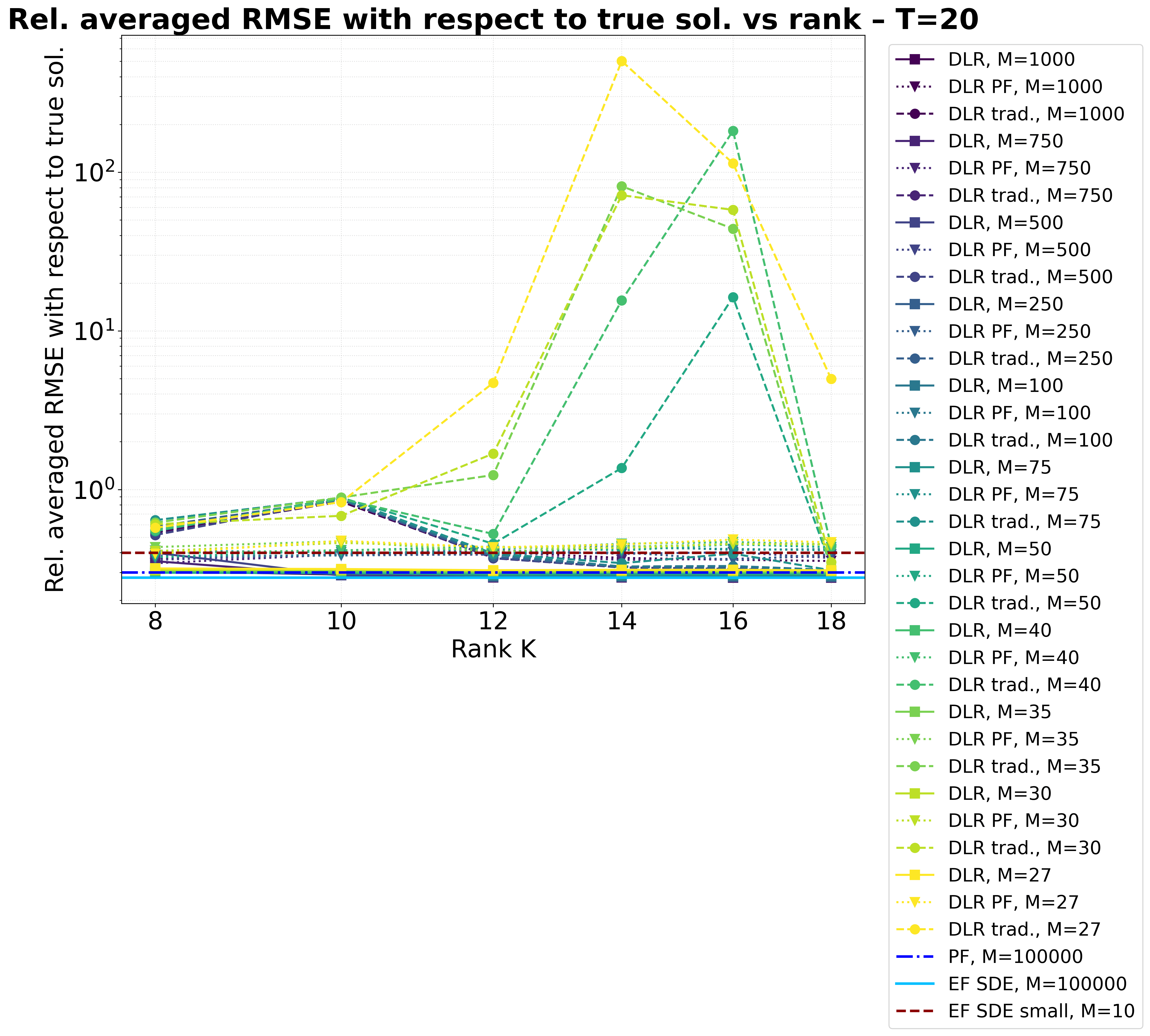}
	\includegraphics[scale=0.044]{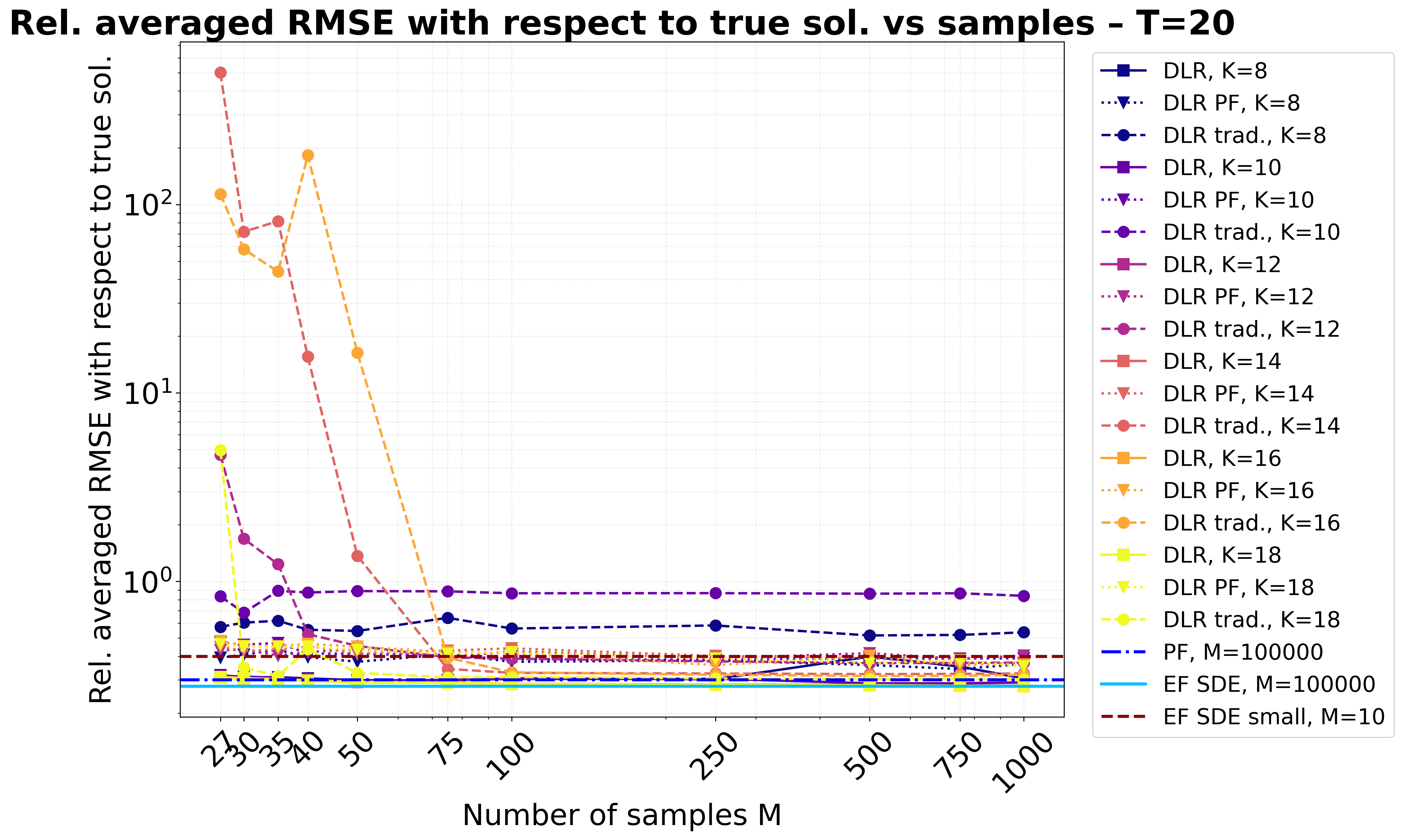} 
	\caption{Rel.\  averaged-over-time \emph{RMSE}, as a function of the rank $k$ (left) and sample size $M$ (right). Problem \eqref{ex: SADR DA} with discrete observations from the observation process \eqref{eq: SADR obs process 1},  $\Delta t_a = 10 \Delta t$, $d=51$.}
	\label{fig: SADR CDDO RMSE}
\end{figure}

In Figures \ref{fig: SADR CDDO EF SDE rank} and \eqref{fig: SADR CDDO EF SDE samples}, we compare the error on the mean and covariance with respect to the reference EF SDE in the case of discrete in time observations, where we assimilate them at a time-step $\Delta t_a = 10 \Delta t$. One can notice once again that the DLRA JMCO filter outperforms the other DLRA filters. 
DLRA-JMCO shows better performance than DLR Trad (and DLR PF) in all types of errors thanks to its additional noise-based term in the evolution of the covariance. The same behavior can be observed in the RMSE plots in Figure \ref{fig: SADR CDDO RMSE}.

\subsection{Nonlinear problem: Lorenz '63}\label{sec: Lorenz 63}
We test our DLR filter with a first example of SDEs with nonlinear drift. In this section, we simulate a dynamics described by a Lorenz '63 problem, firstly presented in \cite{lorenz2017deterministic}, where we add a (multiplicative) stochastic forcing term:
\begin{equation}\label{ex:lorenz63}
	\mathrm{d}
	\begin{pmatrix}
		X_1\\
		X_2\\
		X_3
	\end{pmatrix}
	=
	\begin{pmatrix}
		\sigma (X_2 - X_1)\\
		X_1(\rho - X_3) - X_2\\
		X_1X_2 - \beta X_3
	\end{pmatrix}
	\mathrm{d}t
	+
	b(X_t)\,\mathrm{d}W_t,
	\qquad t\in[0,T],
\end{equation}
where $W_t=(W_{1}(t),W_{2}(t))^\top$ is a two-dimensional standard Wiener process and the state-dependent diffusion matrix
$b(X_t)\in\mathbb{R}^{3\times 2}$ appearing in \eqref{ex:lorenz63} is defined as
\begin{equation}\label{noise l63}
	b(X_t)
	=
	\begin{pmatrix}
		0 & 0\\
		\varepsilon_\rho X_1 & 0\\
		0 & -\varepsilon_\beta X_3
	\end{pmatrix}.
\end{equation}
In particular, we choose $\varepsilon_\rho = 1$ and $\varepsilon_\beta = \frac{1}{2}$.
Hence, a rank-$2$ multiplicative noise is applied to the deterministic Lorenz-63 dynamics, which can be interpreted as accounting for unresolved temporal fluctuations in the physical parameters of the Lorenz '63 model. In particular, the stochastic term $\varepsilon_\rho X_1\mathrm{d}W_{1}(t)$ corresponds formally to a random perturbation of the parameter $\rho$, which controls the thermal forcing of the system, whereas the term $-\varepsilon_\beta X_3 \mathrm{d}W_{2}(t)$ represents fluctuations in the parameter $\beta$, associated with the damping of the third state variable. 

Moreover, in \eqref{ex:lorenz63} the parameters are chosen as follows: $\sigma=10$, $\rho=28$, and $\beta = \frac{8}{3}$. This choice allows to have a chaotic dynamics already in the deterministic regime, and the addition of noise implies extra difficulty in tracking the state, making the problem interesting for data assimilation purposes \cite{lorenz2017deterministic}.  For the sake of numerical stability, we have used a tamed-Euler method \cite{hutzenthaler2012strong} to generate the state $X^{\ast}$, the observations, the reference, and the DLR filters. Although the Lorenz ’63 system is not high-dimensional, it serves as a useful benchmark for assessing whether DLRA can accurately track the chaotic behavior of a system whose stochastic dynamics evolve mainly within a low-dimensional subspace.

The observation process follows the relation
\begin{equation}\label{eq: lorenz 63 obs process 1}
	\mathrm{d} Z_{t} = X_t \mathrm{d}t + R^{\frac{1}{2}} \mathrm{d}B_t,
\end{equation}
i.e.\ we observe the full state polluted by an additive noise with covariance $R=0.01 \cdot I_{k \times k}$.  Notice that, as we are in the case of fully observed state, we can take the rescaled observation differences $X_{n+1}^{\mathrm{obs}}=\frac{\Delta Z_n}{\Delta t_n}$ themselves, which are known quantities, as a reference trivial filter to estimate the true state from the observations. Indeed, a backward discretization of \eqref{eq: lorenz 63 obs process 1} reads $\Delta Z_n = X_{n+1} \Delta t + \Delta B_n$, where the noise of the observation $\Delta B_n$ is unknown, and we can simply ignore it and take $X_{n+1}^{\mathrm{obs}}=\frac{\Delta Z_n}{\Delta t_n}$ as our best estimate of the state. The expected RMSE of such filter is 
\begin{equation}\label{eq: obs. noise ref}
	\begin{aligned}
		\mathbb{E}\left[\int_{T_0}^{T} |X_s^{\mathrm{obs}}-X_s | \mathrm{d}s\right] \approx & 	\mathbb{E}\left[\sum_{n=N_0}^{N} |X_{n+1}^{\mathrm{obs}}-X_{n+1} | \Delta t\right] \\
		=& \mathbb{E}\left[\sum_{n=N_0}^{N} |R^{\frac{1}{2}} \Delta B_n | \Delta t\right] = C | R^{\frac{1}{2}}| N \sqrt{\Delta t} = C \frac{T-T_0}{\sqrt{\Delta t}}|R^{\frac{1}{2}}|.
	\end{aligned}
\end{equation}
with $C = \mathbb{E}\left[\frac{| \Delta B_n|}{\sqrt{\Delta t}} \right]$.
In our plots, we will refer to this quantity as the \emph{observation noise} (\emph{Obs. noise} in the legends of the plots). In case of continuous-time observation, if a filtering strategy manages to have a lower RMSE than \eqref{eq: obs. noise ref}, then we will say that our filter is \emph{tracking the trajectory of the true state of reference} as it is effectively filtering out the observation noise using the model.

The final time is chosen as $T=20$, whereas the number of paths simulated for the reference algorithms is $M=10^{5}$. As another reference, we run a particle filter with only $M=10$ particles. We consider a uniform time step $\Delta t = 5 \cdot 10^{-3}$. The chosen rank for all the DLRA algorithms is $k=2$, which is the same rank of the multiplicative noise $b$. The initial condition $X_0$ has been built as follows: a random variable $\gamma=\left(\gamma_1,\gamma_2,\gamma_3\right)$, whose components have uniform values between $1$ and $2$, i.e.\ $\gamma_i \sim \mathcal{U}(1, 2)$, independent of $W_t$ and of $B_t$, has been sampled. Then, it is projected in a subspace $U$ of dimension equal to $2$, randomly generated and independent of all the other quantities, i.e.\ $X_0 \sim P_U \gamma$.  With this setup, we expect that a DLRA with rank $k=2$ can well track the dynamics of the problems.
In all the average-in-time errors, the warm-up time is $T_0=2$.

\begin{figure}[!h]
	\centering
	\includegraphics[scale=0.3]{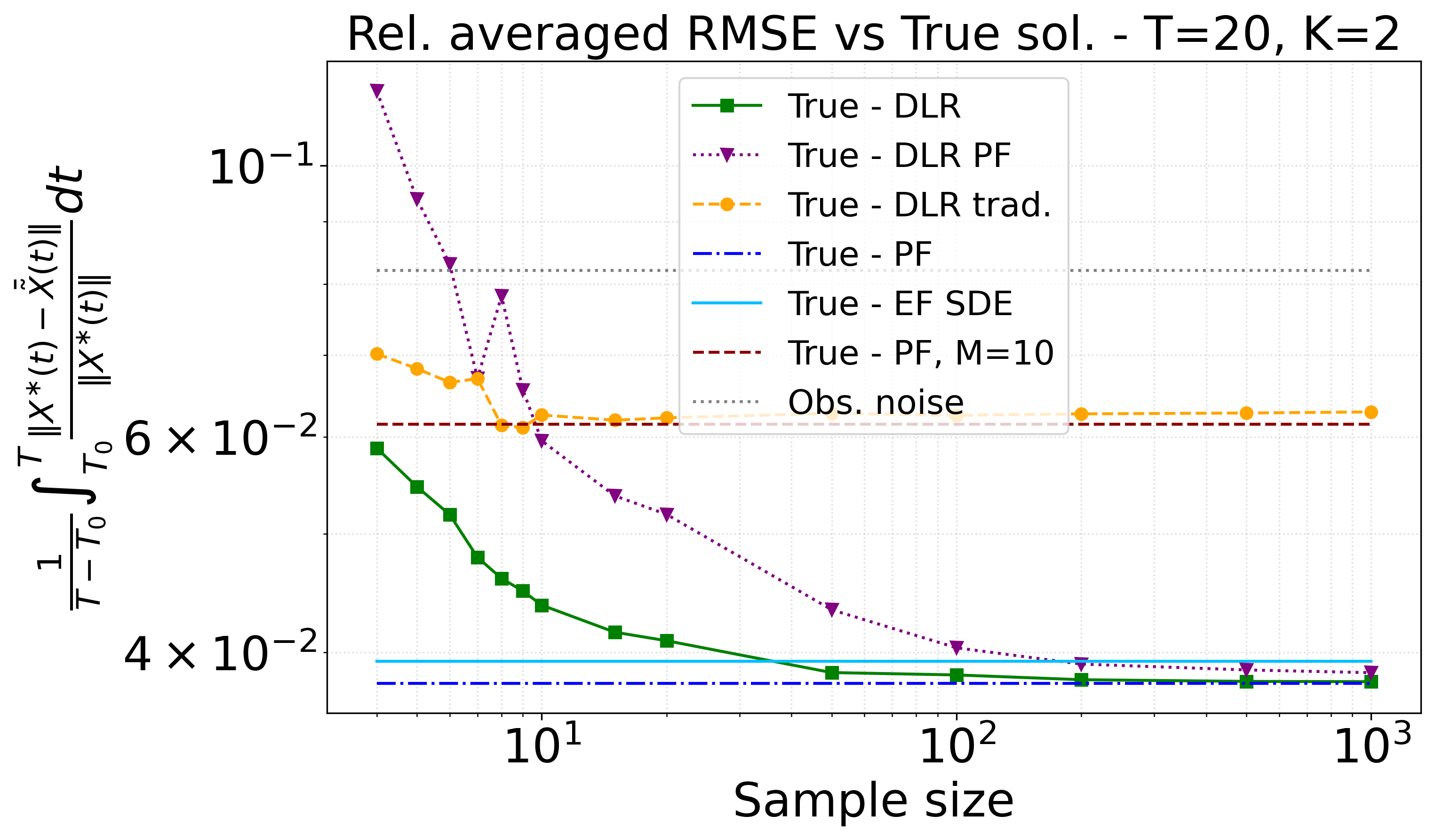}
	\caption{Rel.\  averaged-over-time \emph{RMSE}, as a function of the sample size $M$. Problem \eqref{ex:lorenz63} with continuous-time observations, transient time $T_0=2$, $k=2$, reference computed with $M=10^5$.}
	\label{fig: Lorenz63 RMSE vs True sol - co}
\end{figure}

For the continuous-time-observation case, the employment of the branching procedure is occuring at intervals of time $\Delta t_{\mathrm{branch}} = 4 \Delta t$.
In Figure \ref{fig: Lorenz63 RMSE vs True sol - co} we see how the various filters perform concerning the relative average-in-time \emph{RMSE against the True solution \eqref{eq: avg rel RMSE}} in the case of continuous-time observations. We see that both the DLR JMCO and the DLR PF follow almost perfectly the PF of reference for large number of samples, whereas the DLR Trad.\ stagnates. On the other hand, the PF with $M=10$ samples has a huger error with respect to the DLR JMCO filter employed with the same number of particles.

\begin{figure}[!h]
	\centering
	\includegraphics[scale=0.3]{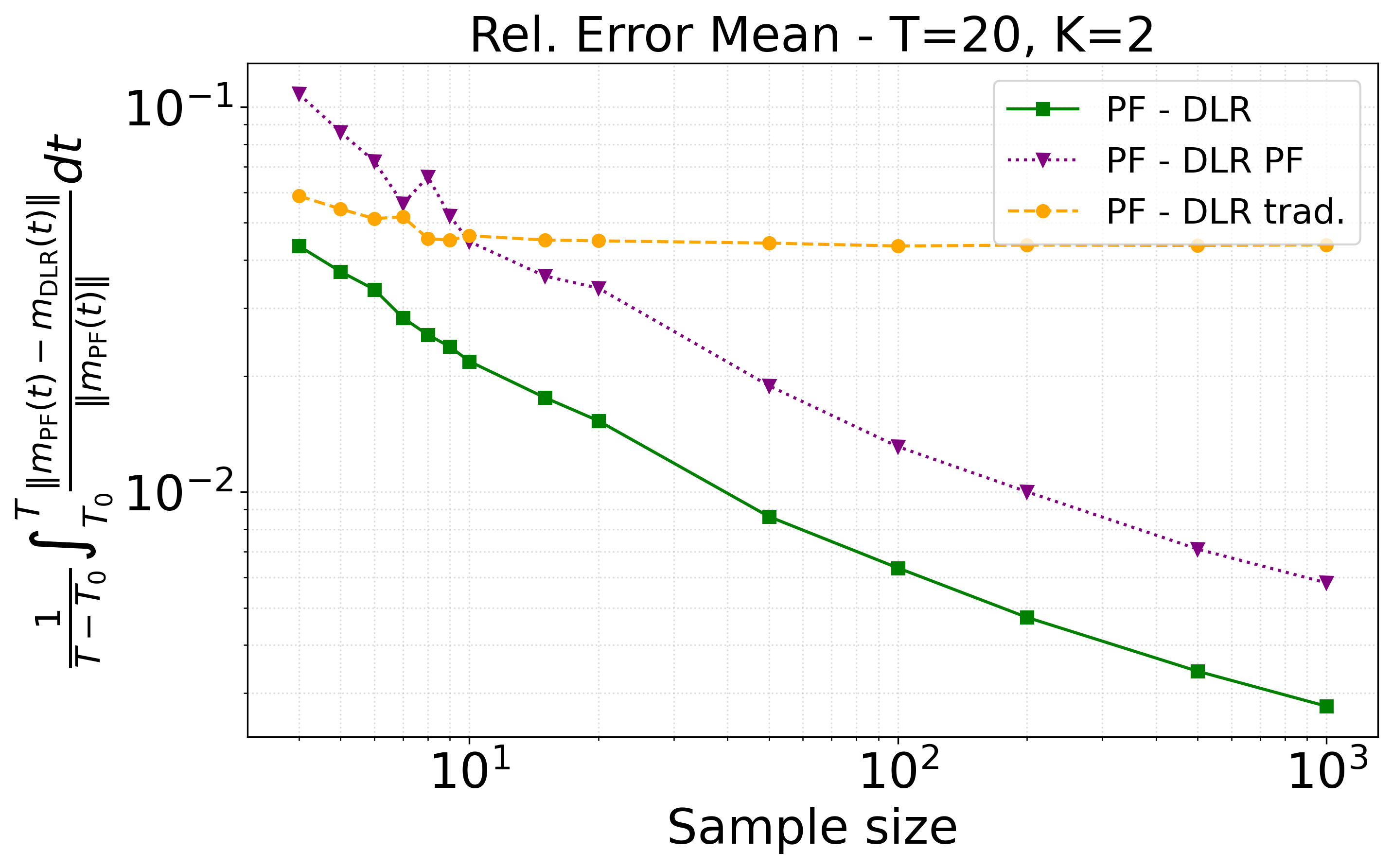}
	\includegraphics[scale=0.3]{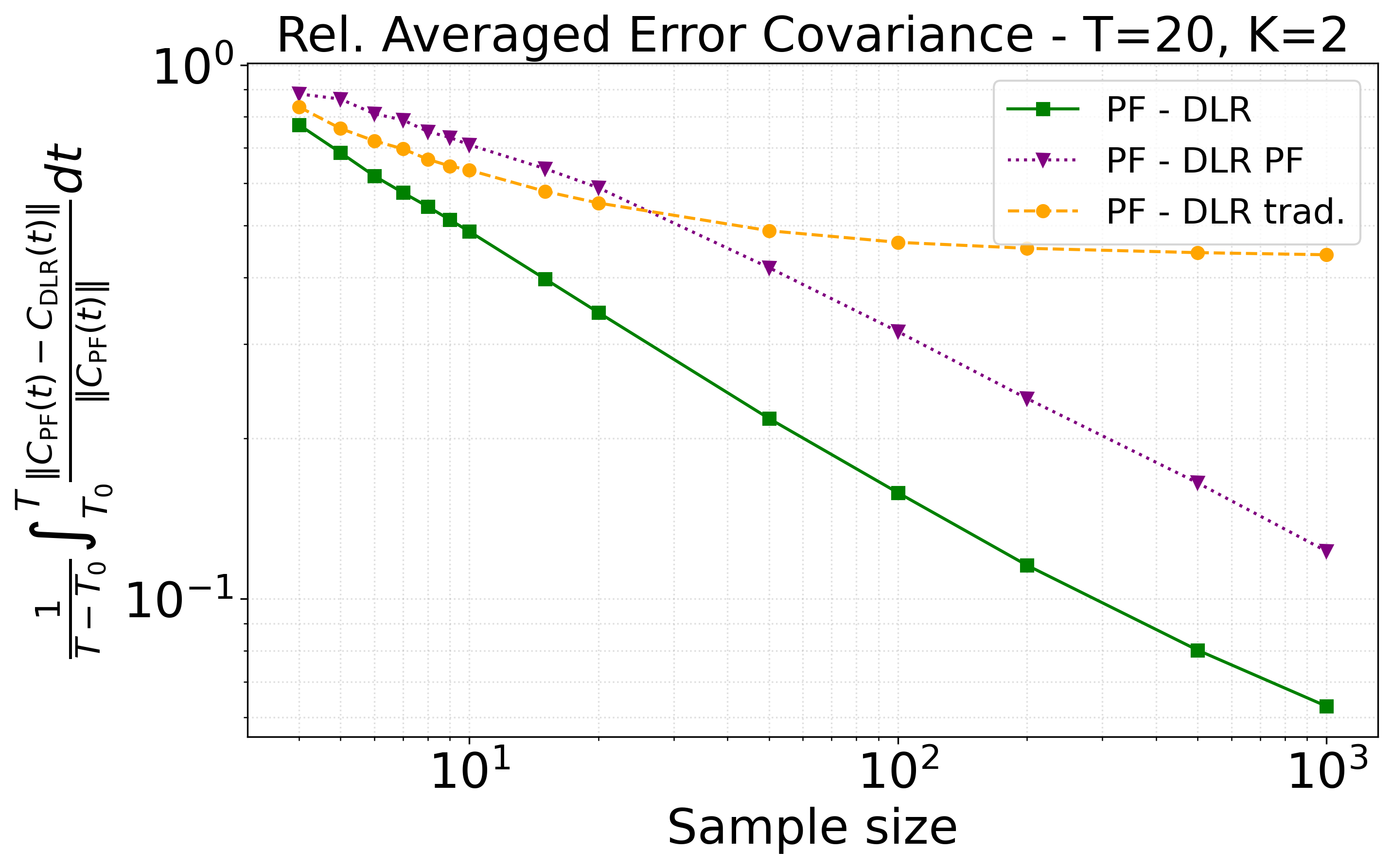}
	\caption{Rel.\  averaged-over-time \emph{errors against PF reference} as a function of the sample size $M$, (left) mean, (right) covariance. Problem \eqref{ex:lorenz63} with continuous-time observations, transient time $T_0=2$, $k=2$, reference computed with $M=10^5$, with $\Delta t_{branch}= 4\Delta t$.}
	\label{fig: lorenz63 mean - cov PF co}
\end{figure}

As we are dealing with a chaotic nonlinear system, our reference algorithm is the PF. In this case, the DLR PF emulates its full-order counterpart, showing convergence for the first two moments with respect to the number of particles (see Figure \ref{fig: lorenz63 mean - cov PF co}). However, the DLRA-JMCO filter still has the better accuracy and, quite surprisingly, seems to converge to the PF as $M \to \infty$. 

We repeat now the same simulations, this time changing the branching time of DLR PF and the PF. Figure \eqref{fig: Lorenz63 - co, branch 2} and Figure \eqref{fig: Lorenz63 - co, branch 1} are the same as Figure \ref{fig: lorenz63 mean - cov PF co} but with $\Delta t_{\mathrm{branch}} = 2 \Delta t$ and  $\Delta t_{\mathrm{branch}} = \Delta t$, respectively.
For a large number of employed samples $M$, DLR PF definitely shows better performance when taking the smallest $\Delta t_{\mathrm{branch}} = \Delta t$, improving its accuracy in all types of errors. However, the DLR JMCO still wins against DLR-PF and seems to be a better filtering choice for this problem.
\begin{figure}[!h]
	\centering
	\includegraphics[scale=0.19]{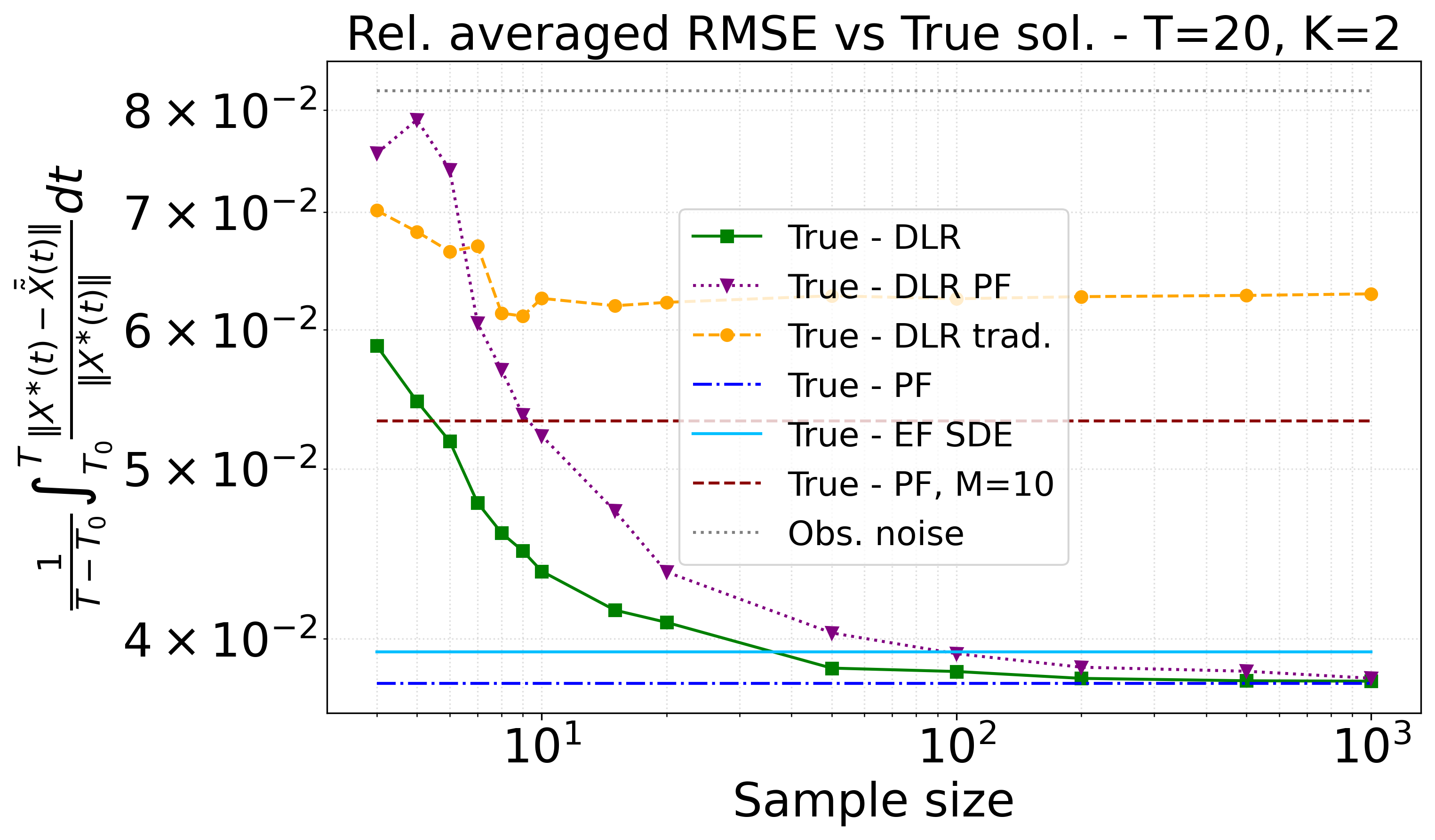}
	\includegraphics[scale=0.19]{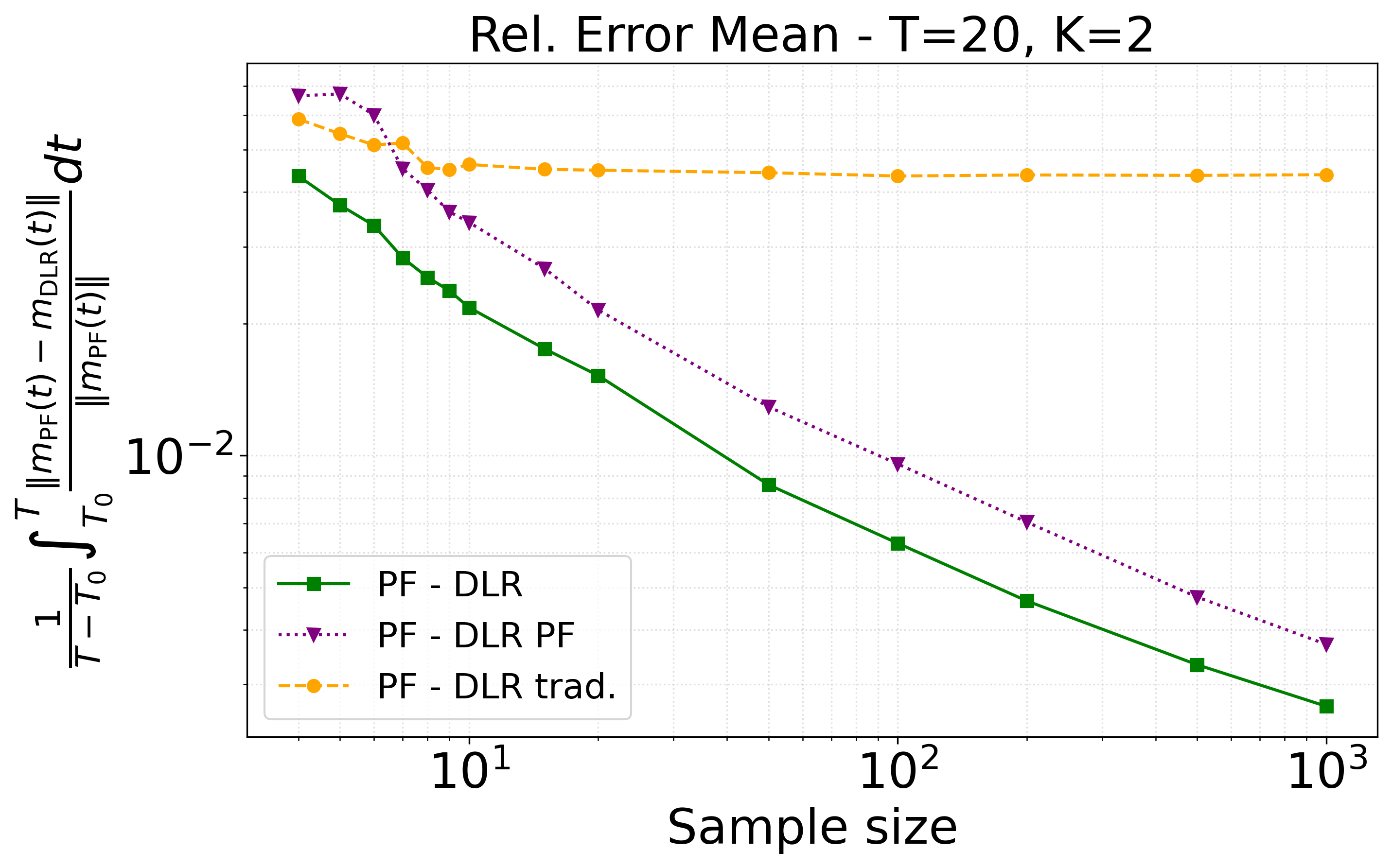}
	\includegraphics[scale=0.19]{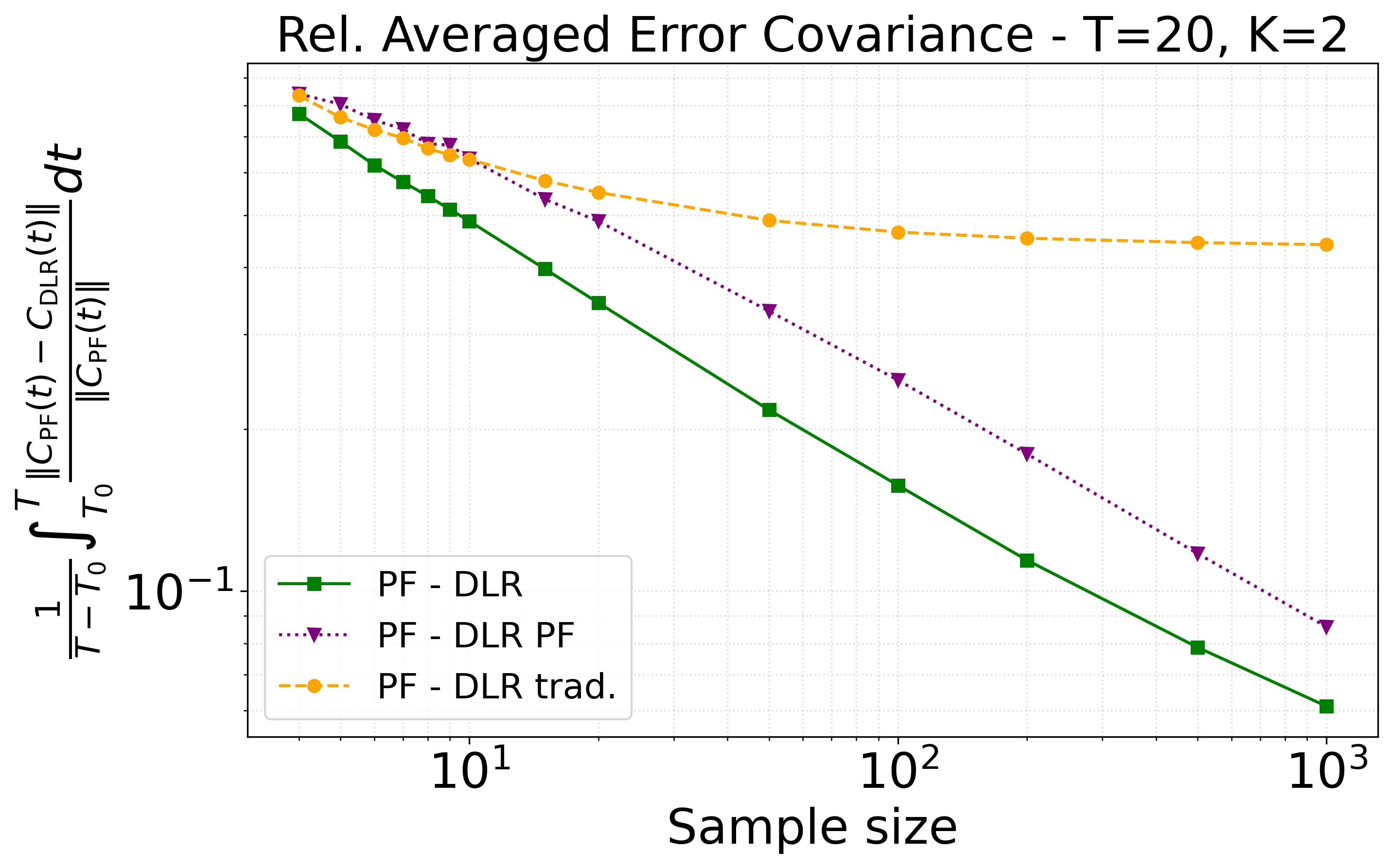}
	\caption{Rel.\  averaged-over-time \emph{RMSE} (left), \emph{errors vs mean of PF} (center) and \emph{errors vs covariance of PF} (right) as a function of the sample size $M$.
		Problem \eqref{ex:lorenz63} with continuous-time observations, transient time $T_0=2$, $k=2$, reference computed with $M=10^5$, with $\Delta t_{branch}= 2\Delta t$.}
	\label{fig: Lorenz63 - co, branch 2} 
\end{figure}

\begin{figure}[!h]
	\centering
	\includegraphics[scale=0.19]{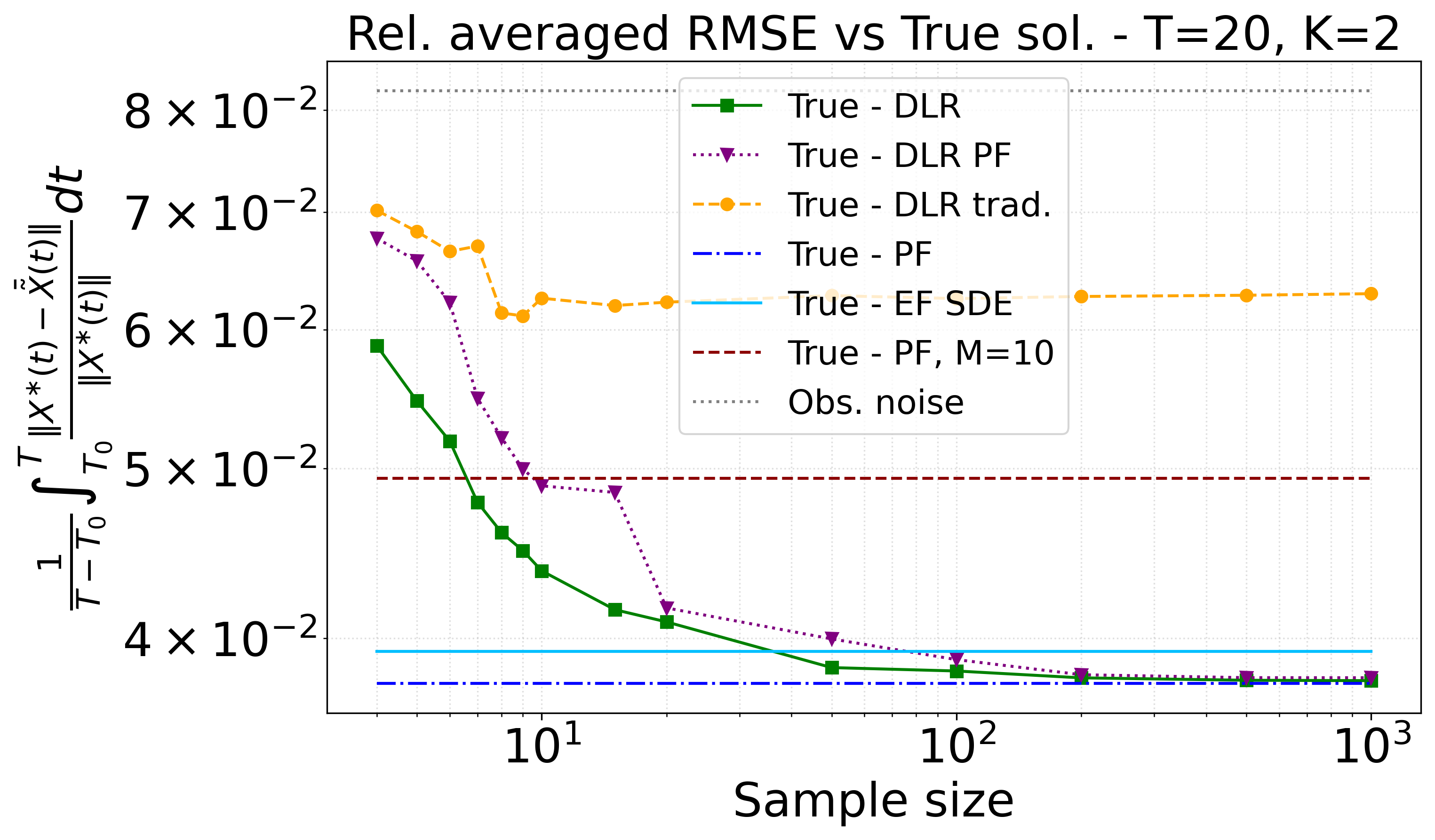}
	\includegraphics[scale=0.19]{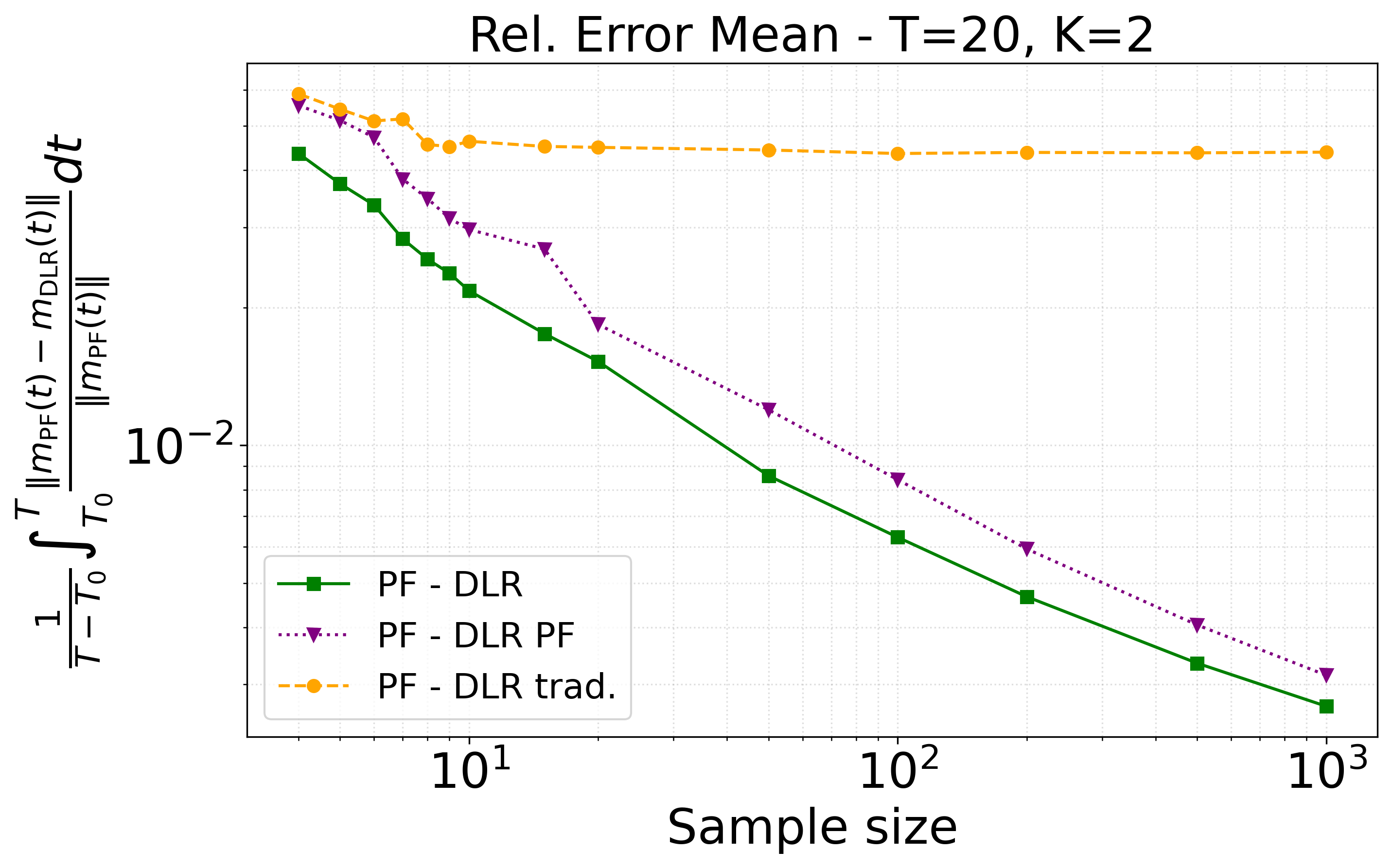}
	\includegraphics[scale=0.19]{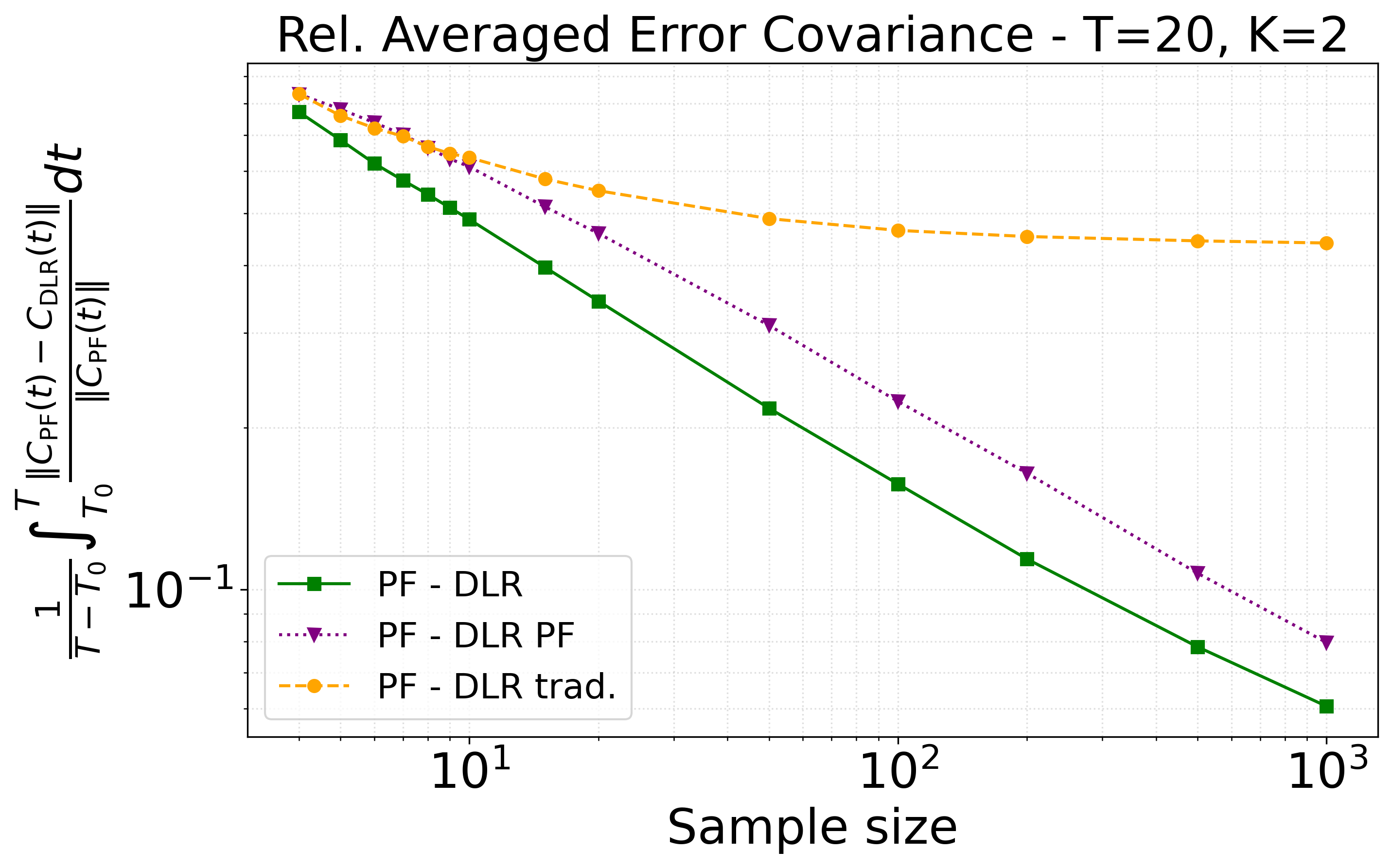}
	\caption{Rel.\  averaged-over-time \emph{RMSE} (left), \emph{errors vs mean of PF} (center) and \emph{errors vs covariance of PF} (right) as a function of the sample size $M$.
		Problem \eqref{ex:lorenz63} with continuous-time observations, transient time $T_0=2$, $k=2$, reference $M=10^5$, with $\Delta t_{branch}= \Delta t$.}
	\label{fig: Lorenz63 - co, branch 1} 
\end{figure}

\begin{figure}[!h]
	\centering
	\includegraphics[scale=0.28]{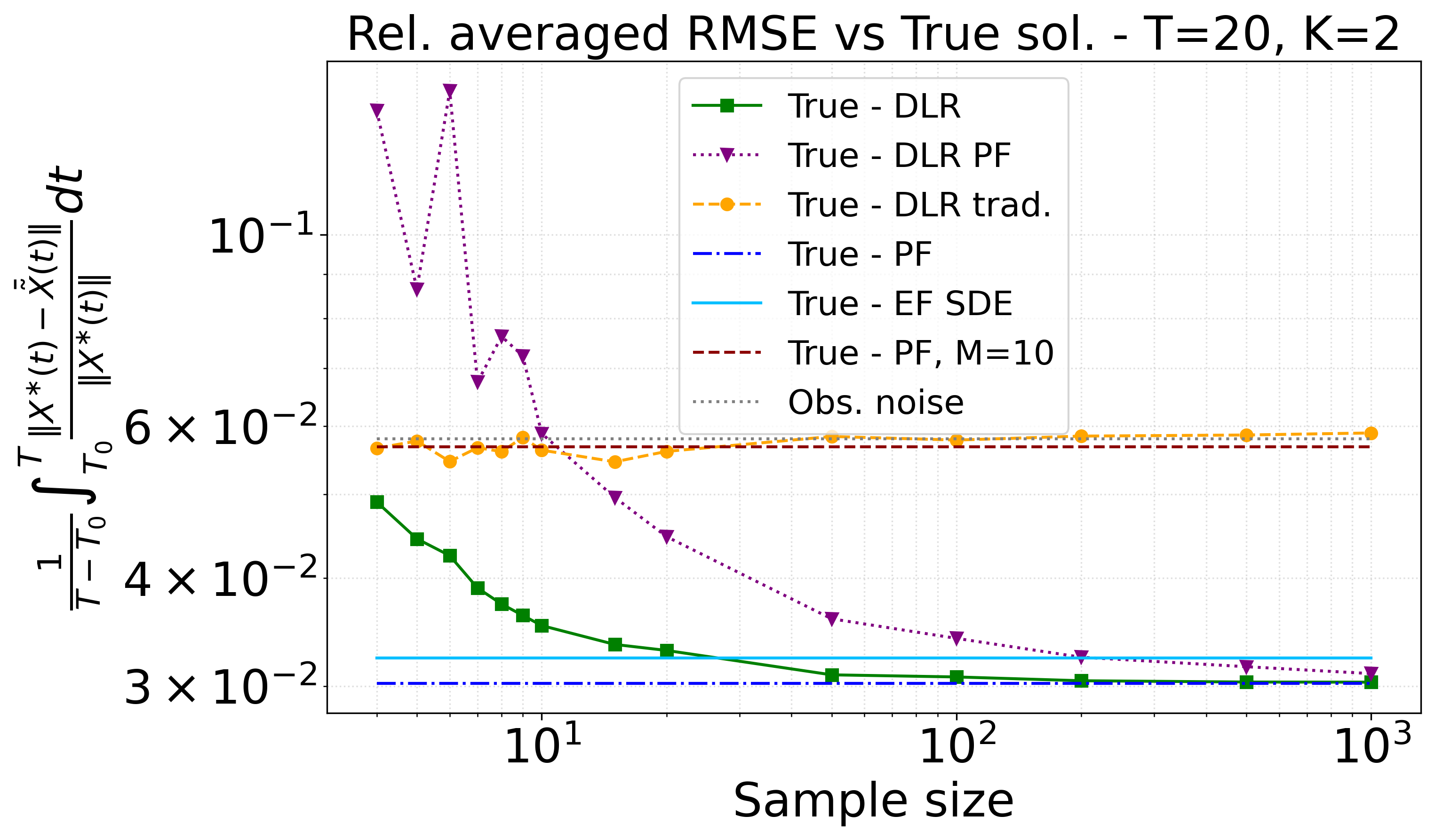}
	\includegraphics[scale=0.28]{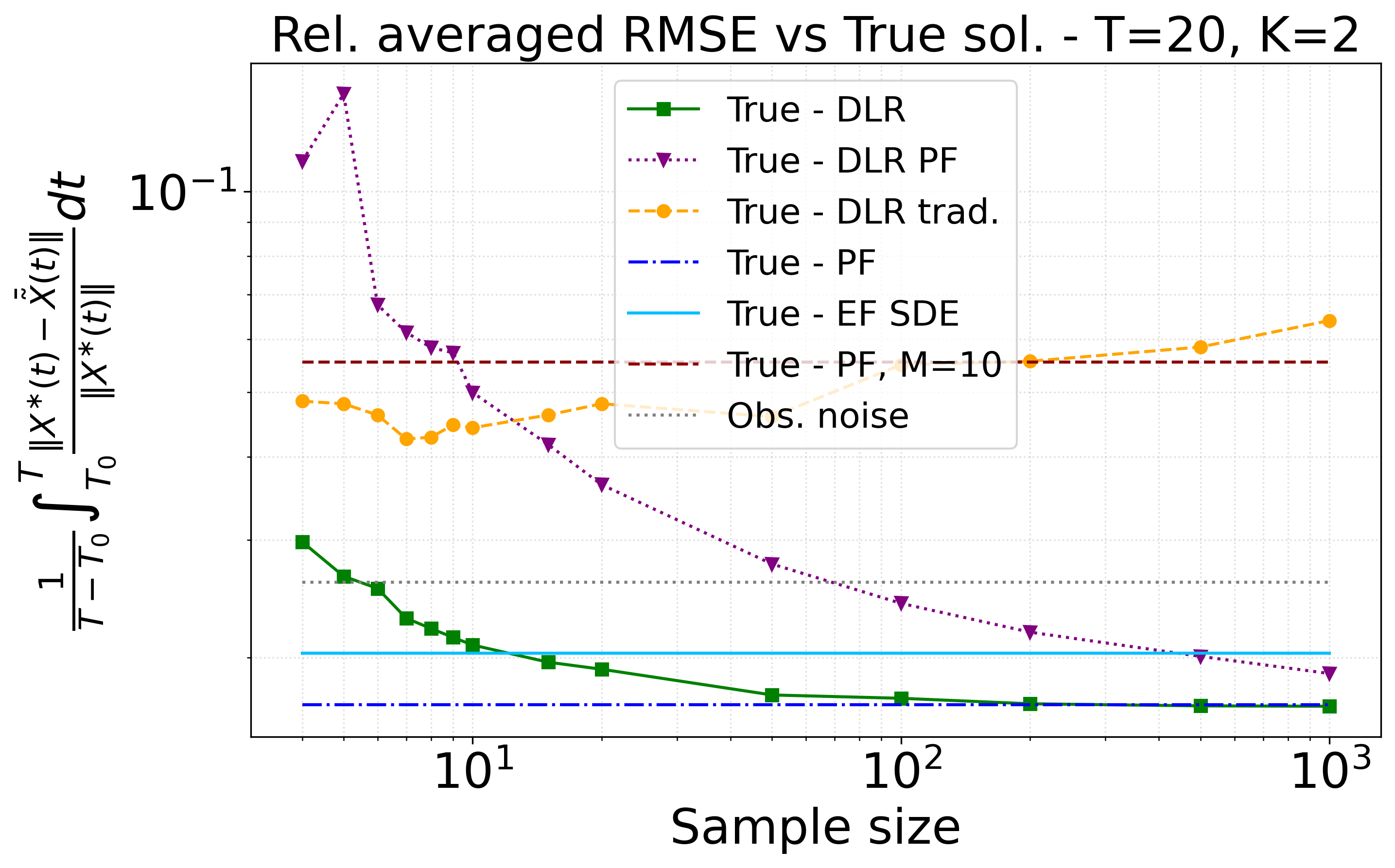}
	\caption{Rel.\  averaged-over-time \emph{RMSE} as a function of the sample size $M$. Problem \eqref{ex:lorenz63} with continuous-time observations for $R=0.005$ (left) and $R=0.001$ (right), transient time $T_0=2$, $k=2$, reference $M=10^5$, $\Delta t_{branch}= 4\Delta t$.}
	\label{fig: Lorenz63 RMSE vs True sol - co R div}
\end{figure}

Still in the case of continuous-time observations, we simulate the same problem this time with smaller observation noise covariance $R=0.005$ and $R=0.001$ for $\Delta t_{branch}= 4\Delta t$. This change increases the accuracy of the trivial filter as a reference filter. However, the DLRA JMCO still manages to well track the true state with a relatively small number of particles, as we see in Figure \ref{fig: Lorenz63 RMSE vs True sol - co R div}.

\begin{figure}[!h]
	\centering
	\includegraphics[scale=0.3]{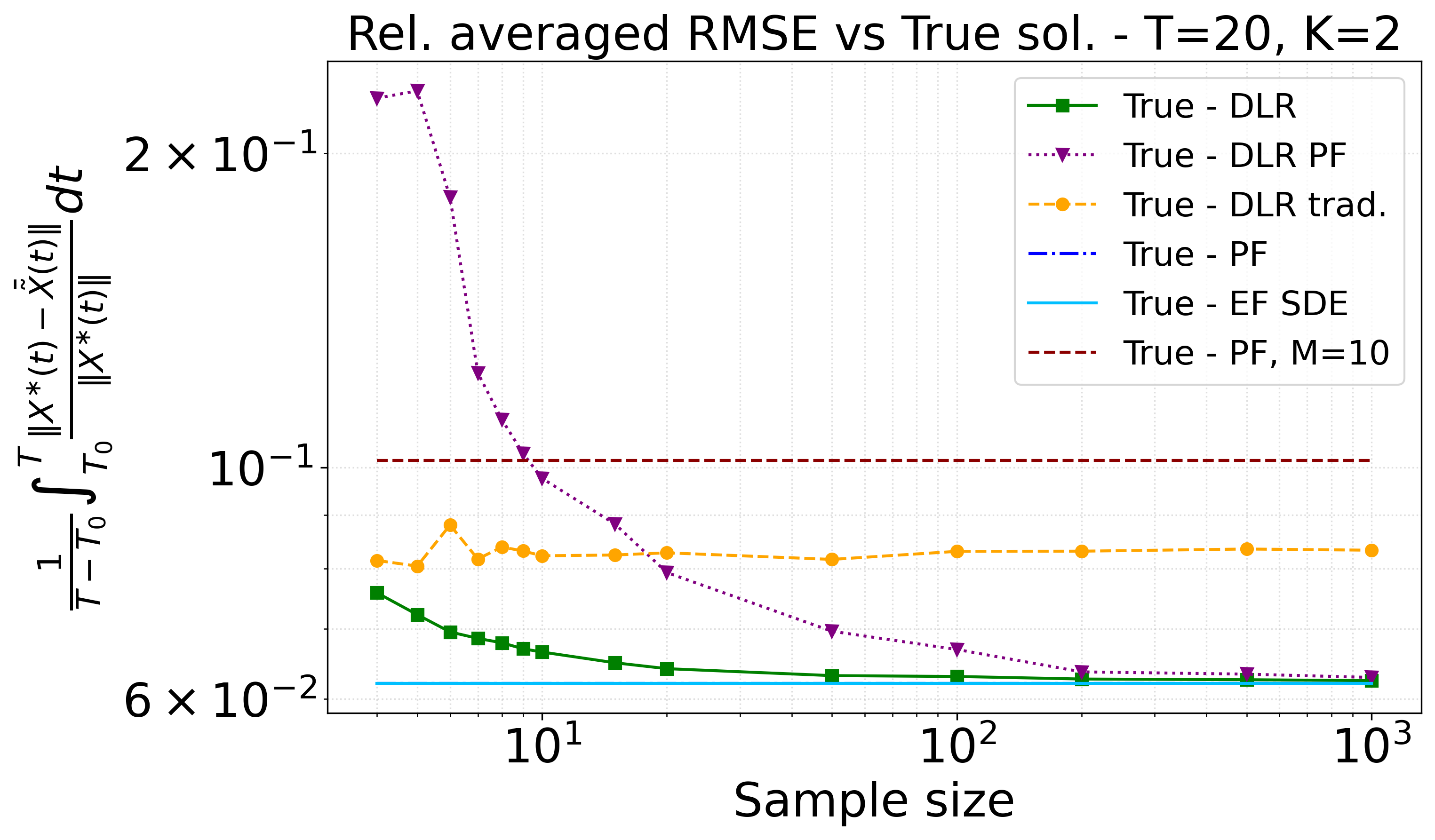}
	\caption{Rel.\  averaged-over-time \emph{RMSE} as a function of the sample size $M$. Problem \eqref{ex:lorenz63} with discrete-time observations, transient time $T_0=2$, $k=2$, reference computed with $M=10^5$, $\Delta t_a = 10 \Delta t$.}
	\label{fig: Lorenz63 RMSE vs True sol - do}
\end{figure}

\begin{figure}[!h]
	\centering
	\includegraphics[scale=0.3]{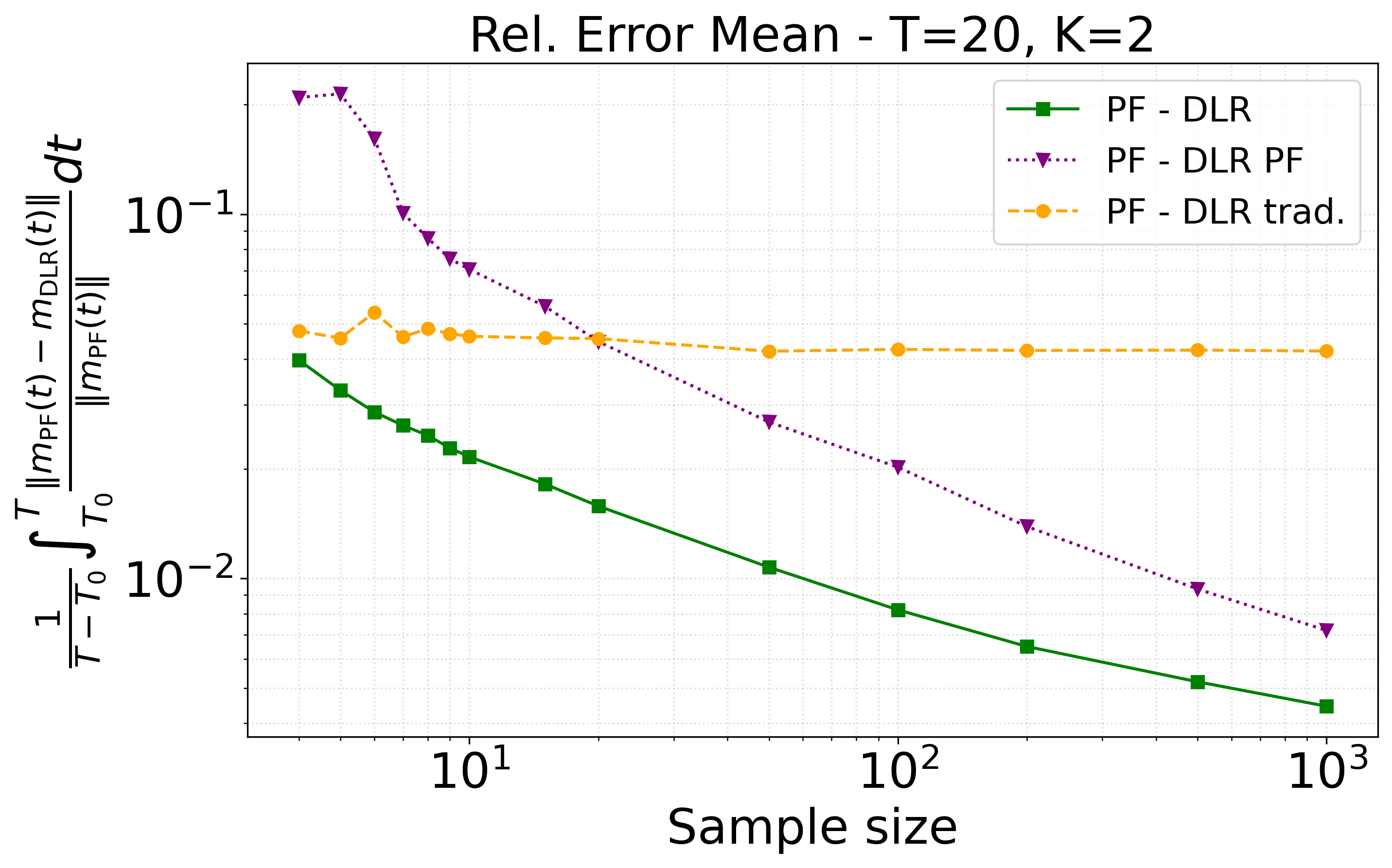}
	\includegraphics[scale=0.3]{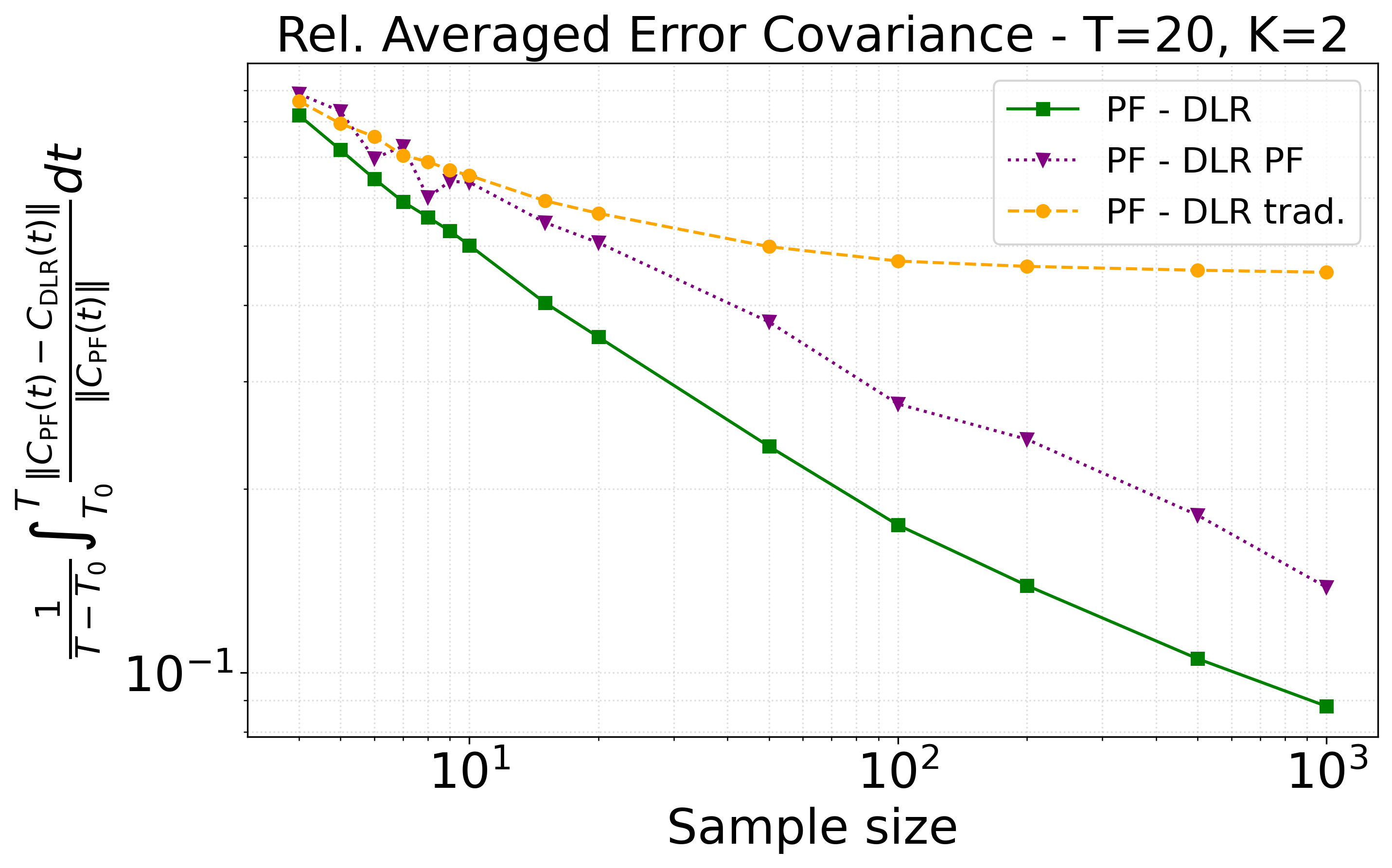}
	\caption{Rel.\  averaged-over-time \emph{errors against PF reference} as a function of the sample size $M$: (left) mean, (right) covariance. Problem \eqref{ex:lorenz63} with discrete-time observations, transient time $T_0=2$, $k=2$, reference computed with $M=10^5$, $\Delta t_a = 10 \Delta t$.}
	\label{fig: lorenz63 mean - cov PF do}
\end{figure}

We run a simulation also for the case of discrete-time observations for a reference number of particles $M=10^5$, where we assimilate at a time-step $\Delta t_a = 10 \Delta t$. In Figures \ref{fig: Lorenz63 RMSE vs True sol - do} and \ref{fig: lorenz63 mean - cov PF do}, we see that the trend concerning the filters is similar to the case of continuous-time observations: the DLR JMCO filter shows good performance, whereas the DLR PF needs a larger number of particles to obtain a good approximation in the RMSE against the True solution, as well as in the error on the first two moments against the PF. Again, the DLR Trad.\ stagnates, probably because it cannot well track the noise of the true dynamics. 

\subsection{Nonlinear problem: Lorenz '96}

\begin{figure}[!h]
	\centering
	\includegraphics[scale=0.35]{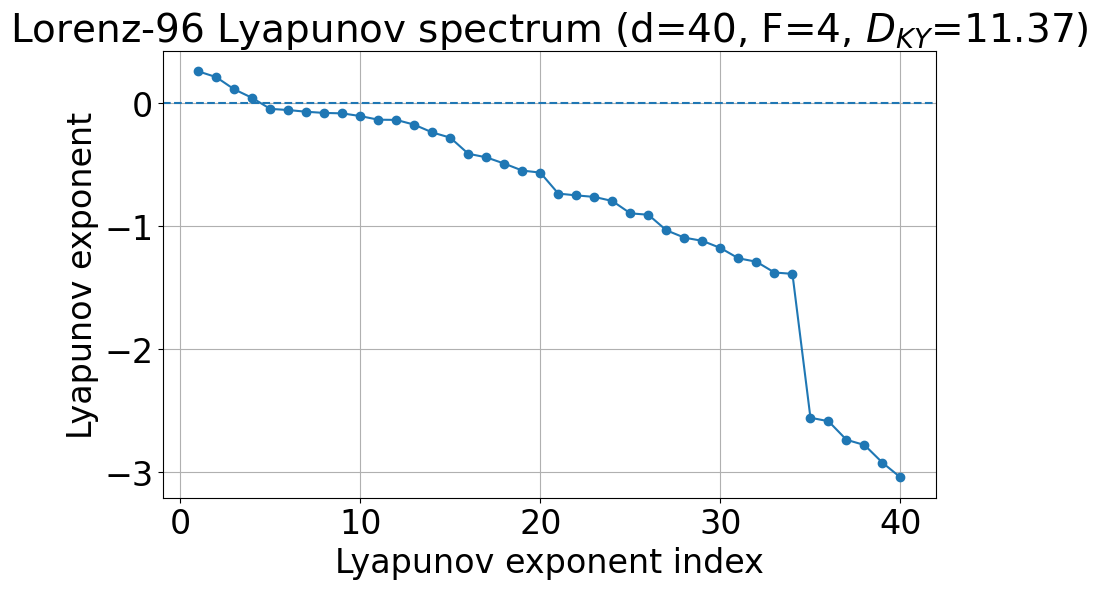}
	\caption{Lyapunov exponents for the deterministic version of \eqref{ex:lorenz96} for $F=4$, $T=20$.}
	\label{fig: Lorenz96 D_KY}
\end{figure}

In this subsection, we consider the Lorenz '96 system, first proposed in
\cite{lorenz1996predictability}, which provides another well-known example
of chaotic dynamics. This problem is described by the following stochastic system
\begin{equation}\label{ex:lorenz96}
	\begin{aligned}
		\mathrm{d}X_i
		&=
		\left[
		X_{i-1}(X_{i+1}-X_{i-2})
		-X_i+F
		\right]\mathrm{d}t
		+
		\bigl(b(X_t)\,\mathrm{d}W_t\bigr)_i,
		\qquad i=1,\ldots,d,\\
		X_0&=X_d,\qquad
		X_{d+1}=X_1,\qquad
		X_{-1}=X_{d-1},
		\qquad t\in[0,T].
	\end{aligned}
\end{equation}
Depending on the value of the forcing parameter $F$ and on the dimension
$d$, the deterministic version of \eqref{ex:lorenz96}, obtained by
setting $b\equiv0$, can exhibit chaotic behavior
\cite{karimi2010extensive,lorenz1996predictability}. In this case, the deterministic version is characterized by the presence of a chaotic attractor of a certain dimension $D_{KY}$ dependent on the number of nonzero Lyapunov exponents $(\lambda_i)_{i=1,\dots,d}$. The dimension $D_{KY}$ is said the Kaplan-Yorke dimension and its formula is given by
\begin{equation*}
	D_{KY} = j^{*} + \frac{\sum_{i = 1}^{j^{*}} \lambda_i}{|\lambda_{j^{*}+1}|},
\end{equation*}
where $j^{*}$ is the largest exponent such that $\sum_{i = 1}^{j^{*}} \lambda_i \geq 0$ \cite{karimi2010extensive}. We consider $d=40$ and a value $F=4$ in this section. We illustrate the trend of the Lyapunov exponents of \eqref{ex:lorenz96} without noise in Figure \ref{fig: Lorenz96 D_KY}, where we have $D_{KY} \approx 11$.

To further account for unresolved spatially structured fluctuations, we
introduce a state-dependent multiplicative noise $b$ acting along selected
Fourier modes of the Lorenz '96 ring, which we describe next. We consider the set
of wave numbers $\mathcal{K}=\{1,\ldots,8\}$ which results in a stochastic forcing of rank $m=2|\mathcal{K}|=16$
since each wave number is represented by a cosine and a sine mode.

For each $k\in\mathcal{K}$, we define the discrete Fourier modes
\begin{equation}
	\phi_k^{c}(j)
	=
	\sqrt{\frac{2}{d}}
	\cos\left(\frac{2\pi k j}{d}\right),
	\qquad
	\phi_k^{s}(j)
	=
	\sqrt{\frac{2}{d}}
	\sin\left(\frac{2\pi k j}{d}\right),
	\qquad j=0,\ldots,d-1.
\end{equation}
Collecting these modes gives the matrix $B= \left( \phi_1^c \ \phi_1^s  \
\cdots \phi_8^c \  \phi_8^s \ \right) \in\mathbb{R}^{d\times m}$.
The intensity of the stochastic perturbation decreases with increasing
wave number. More precisely, for $k\in\mathcal{K}$ we define $\varepsilon_k= 0.5\exp\left[-0.1(k-1)\right]$.
The amplitude of the noise depends on the instantaneous projection of
the state onto these Fourier modes. Defining
\begin{equation}
	\widehat{X}
	=
	B^\top(X-F\mathbf{1}),
\end{equation}
we introduce the saturated modal amplitudes
\begin{equation}
	c(X)
	=
	C_{\max}
	\tanh\left(
	\frac{\widehat{X}}{A_{\max}}
	\right),
	\qquad
	C_{\max}=1,
\end{equation}
where the hyperbolic tangent (tanh) is applied componentwise. Consequently,
each modal amplitude satisfies $|c_\ell(X)|\leq C_{\max}$,
$\ell=1,\ldots,m.$
Notice that, since the selected Fourier modes have nonzero wave number,
they are orthogonal to the constant vector $\mathbf{1} \in \mathbb{R}^{d}$, and, hence, 
$B^\top(F\mathbf{1})=0.$ In the end,
the state-dependent diffusion matrix in \eqref{ex:lorenz96} is
defined as
\begin{equation}\label{eq:l96_mult_noise}
	b(X)
	=
	B 	\operatorname{diag}
	\left(
	\varepsilon_1,\varepsilon_1,
	\ldots,
	\varepsilon_8,\varepsilon_8
	\right)
	\operatorname{diag}
	\left[
	C_{\max}
	\tanh\left(
	\frac{B^\top(X-F\mathbf{1})}{C_{\max}}
	\right)
	\right]
	\in\mathbb{R}^{d\times m}.
\end{equation}

The exponential decay of the coefficients
$\varepsilon_k$ assigns larger stochastic perturbations to the
large-scale, low-wave-number components and progressively weaker
perturbations to smaller spatial scales. The hyperbolic tangent further
prevents individual modal amplitudes from becoming arbitrarily large.

The final time of the simulation is chosen as $T=20$. We employ a uniform time step $\Delta t = 5 \cdot 10^{-3}$ and number of paths $M=1.5 \cdot 10^5$ for the reference algorithms. The initial condition $X_0$ is constructed similarly to Section \ref{sec: Lorenz 63}: we sample a vector with independent component uniformly as $\mathcal{U}(1, 2)$, also independent of $W_t$ and of $B_t$ and we projected them to a subspace of dimension $k$.  Again, we implement the drift via a tamed Euler for the sake of numerical stability \cite{hutzenthaler2012strong}, for the true state, the observation, the DLR and the reference filters. 

The observation operator is defined as
\begin{equation}\label{eq: Lorenz '96 obs process 1}
	\mathrm{d} Z_{t} = X_t  \mathrm{d}t + R^{\frac{1}{2}} \mathrm{d}B_t,
\end{equation}
where we observe the full state perturbed by an additive noise of covariance $R=0.01 \cdot I_{40\times40}$, and, hence, the number of observations is $h=40$.

\begin{figure}[!h]
	\centering
	\includegraphics[scale=0.044]{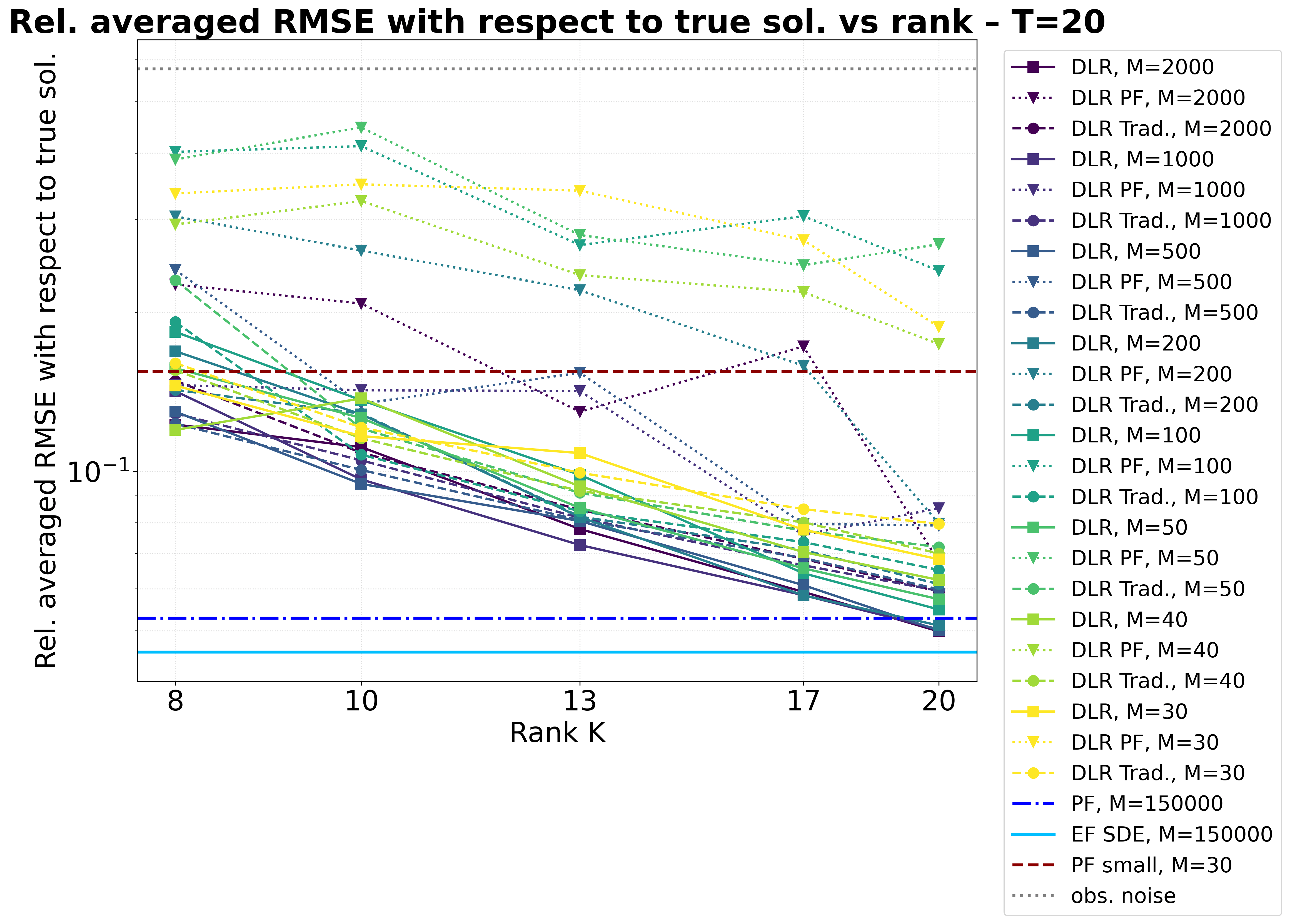}
	\includegraphics[scale=0.044]{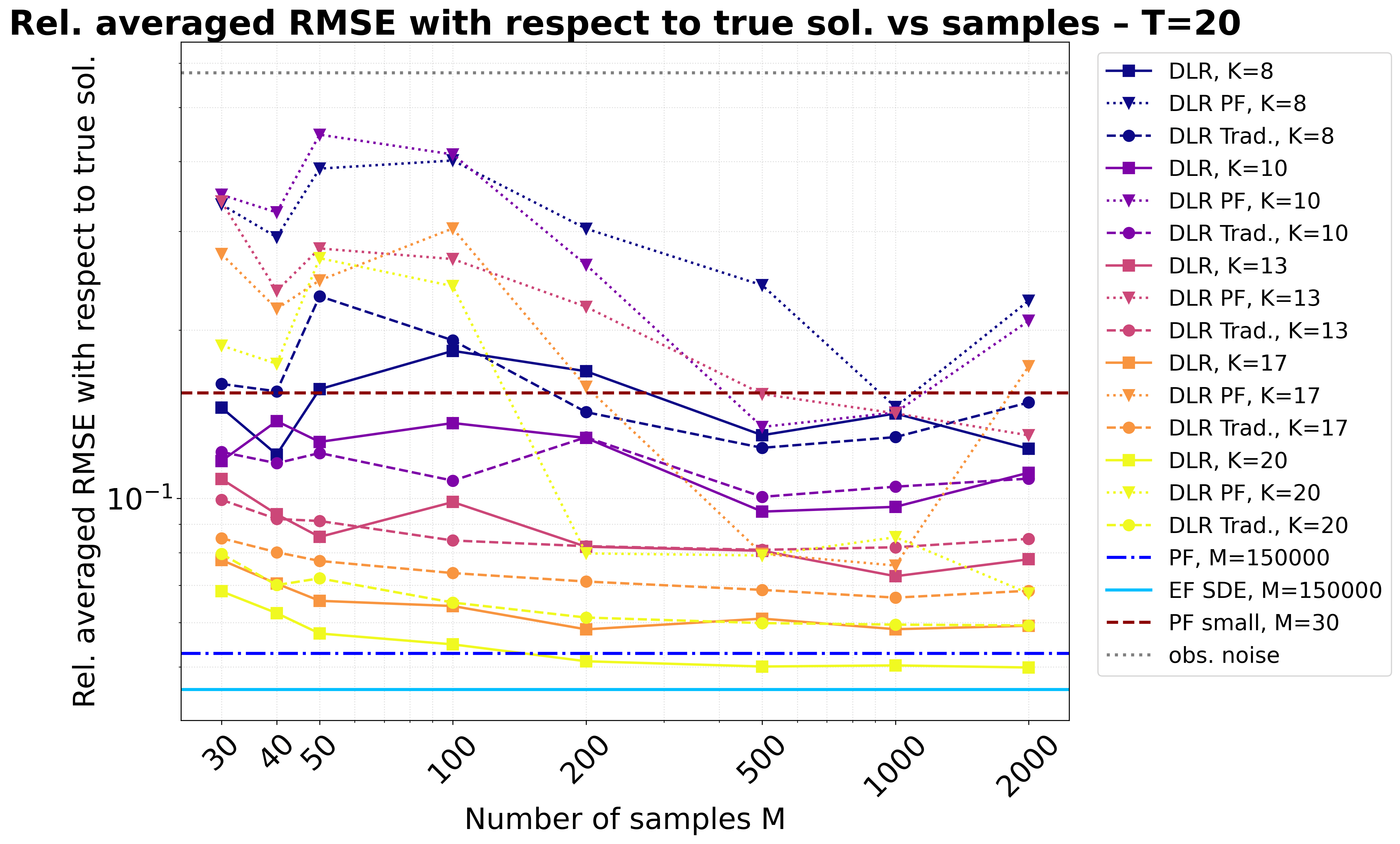}
	\caption{Rel.\  averaged-over-time \emph{RMSE} as a function of the rank $K$ (left) and as a function of the samples sizes $M$ (right). Problem \eqref{ex:lorenz96} with continuous-time observations, transient time $T_0=2$, reference computed with $M=1.5 \cdot 10^5$, $\Delta t_{branch}= 4\Delta t$.}
	\label{fig: Lorenz96 RMSE vs True sol - co}
\end{figure}

In Figure \ref{fig: Lorenz96 RMSE vs True sol - co} we see how the various filters perform concerning the relative average-in-time \emph{RMSE against the True solution \eqref{eq: avg rel RMSE}} in the case of continuous-time observations. We see that the DLR JMCO is more accurate than all the other DLR filters. Notice that all the proposed DLR strategies are tracking the true state. 

\begin{figure}[!h]
	\centering
	\includegraphics[scale=0.045]{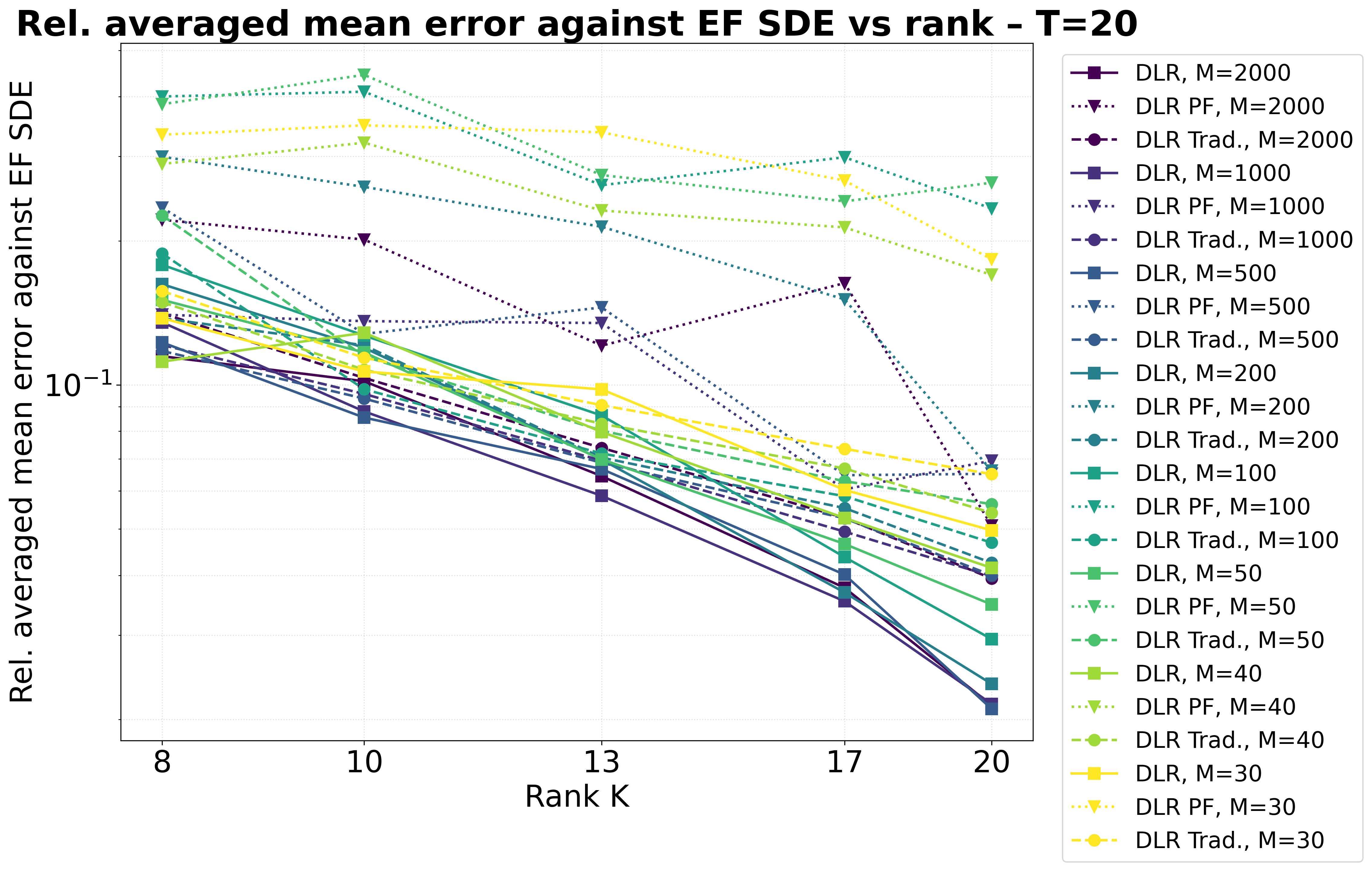}
	\includegraphics[scale=0.045]{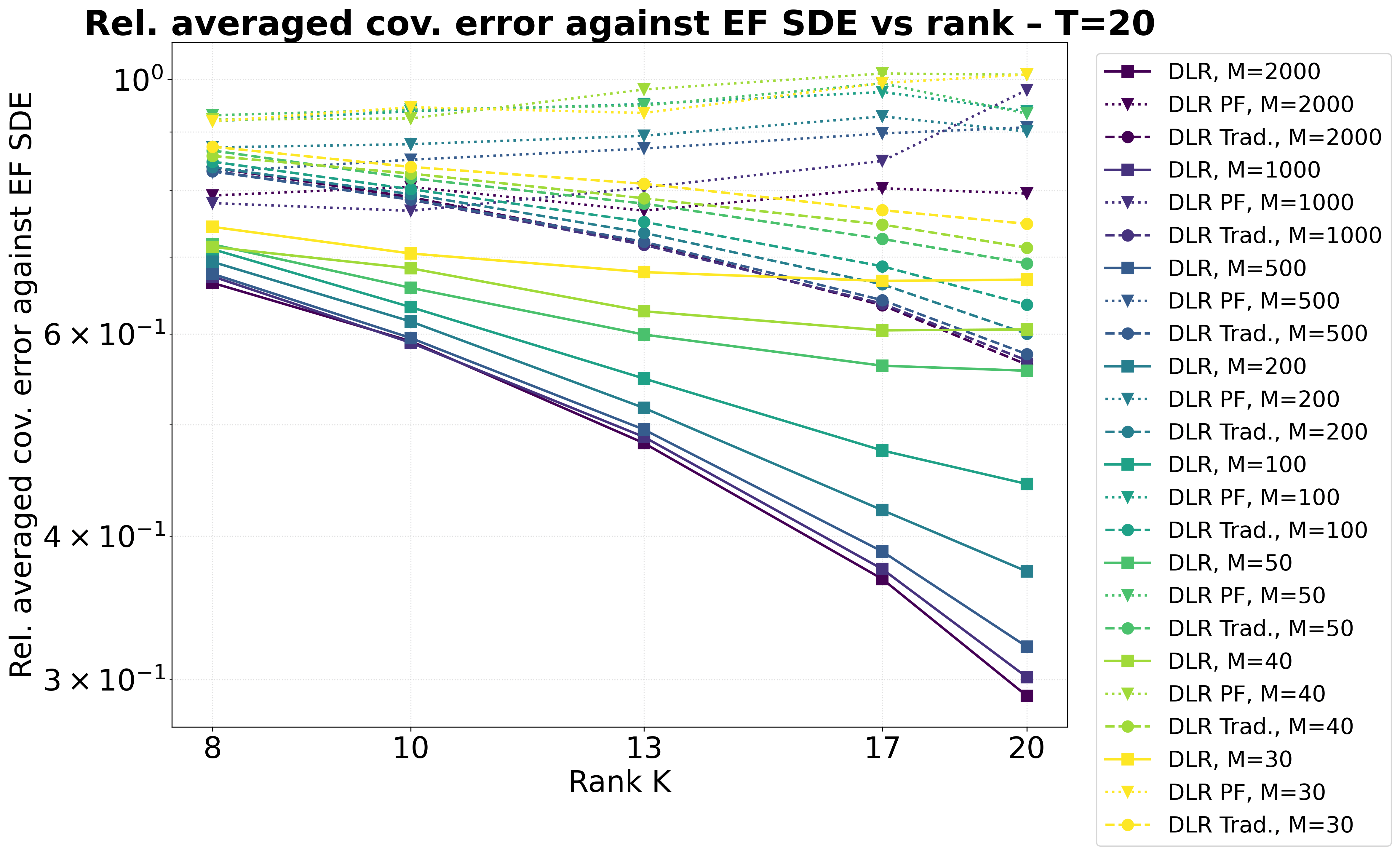} \\
	\includegraphics[scale=0.045]{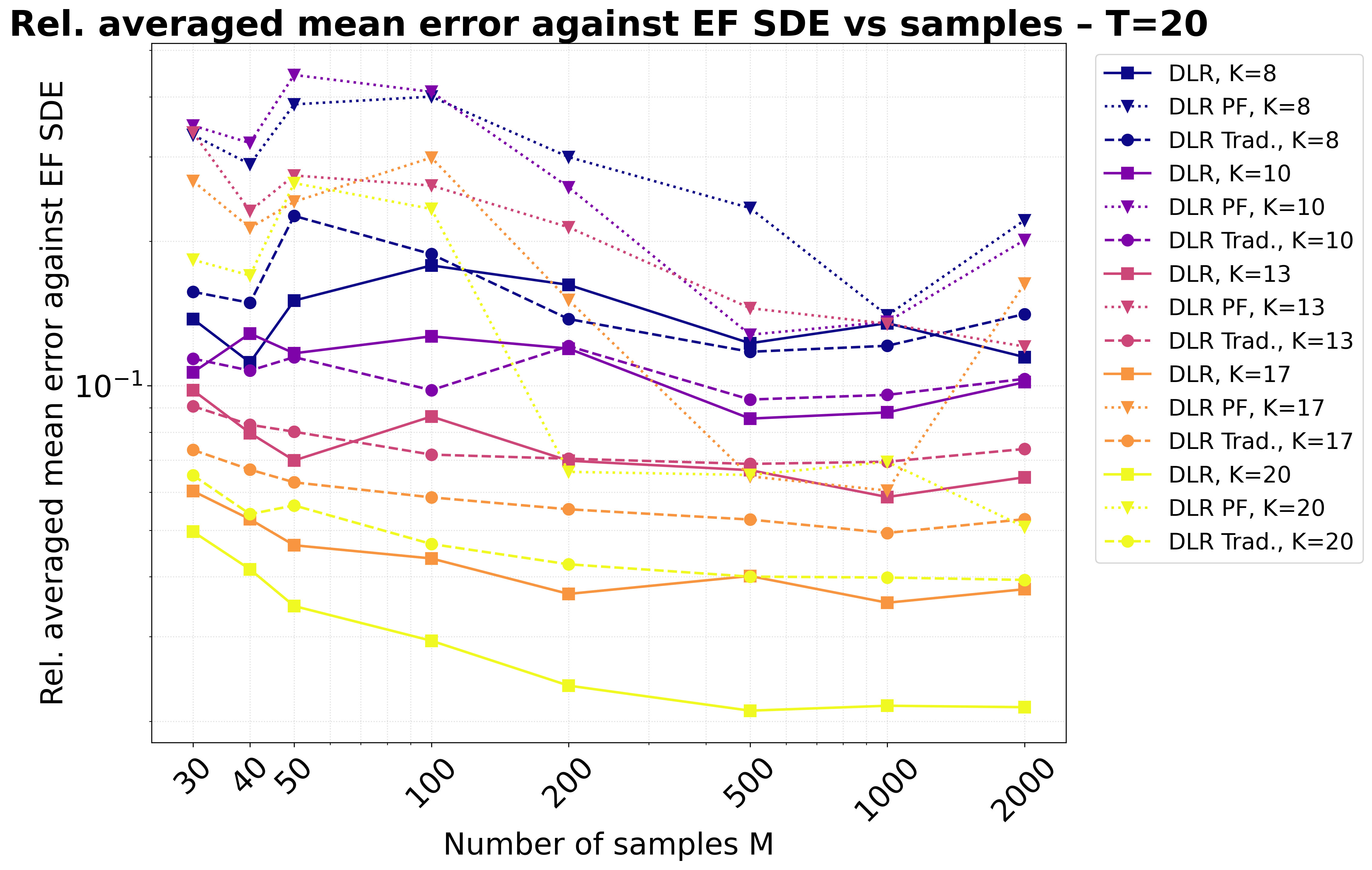}
	\includegraphics[scale=0.045]{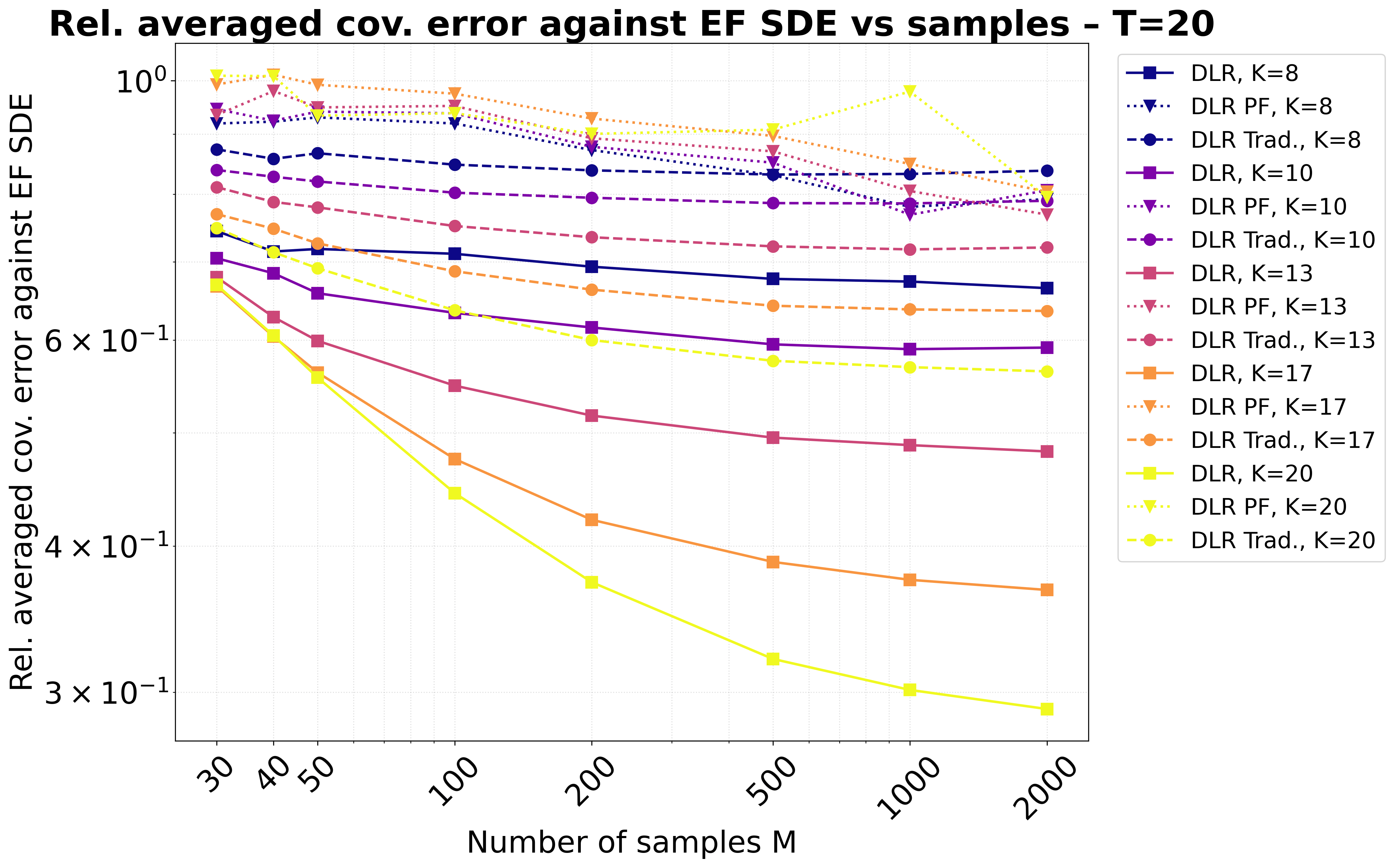} 
	\caption{Rel.\  averaged-over-time \emph{errors against EF SDE}, as a function of the rank $K$ (top row), as a function of the sample size $M$ (bottom row): (left) error on the mean, (right) error on the covariance. Problem \eqref{ex:lorenz96} with continuous-time observations, transient time $T_0=2$, reference computed with $M=1.5 \cdot 10^5$, $\Delta t_{branch}= 4\Delta t$.}
	\label{fig: Lorenz96 EF SDE mean and cov - co}
\end{figure}

\begin{figure}[!h]
	\centering
	\includegraphics[scale=0.045]{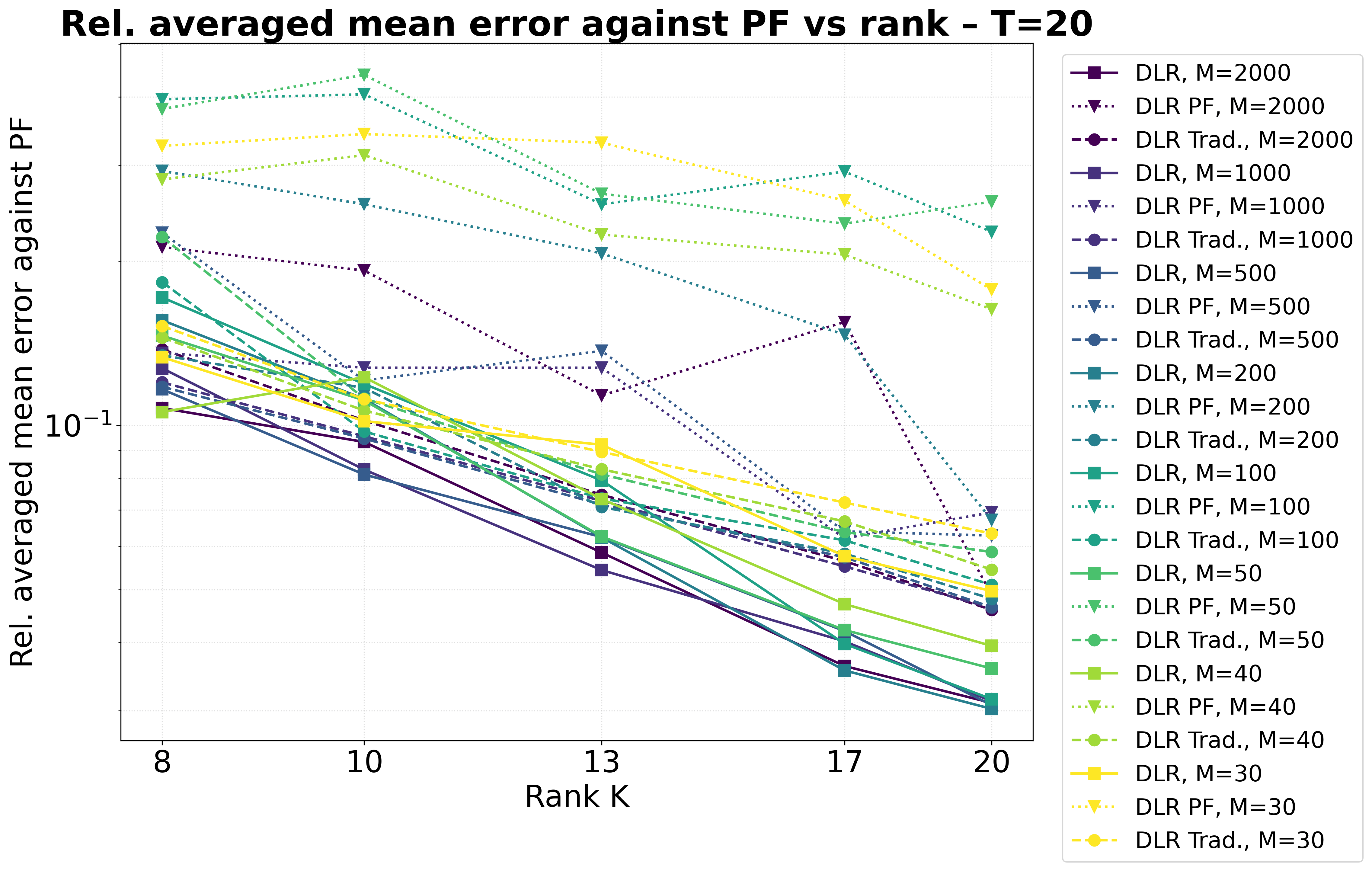}
	\includegraphics[scale=0.045]{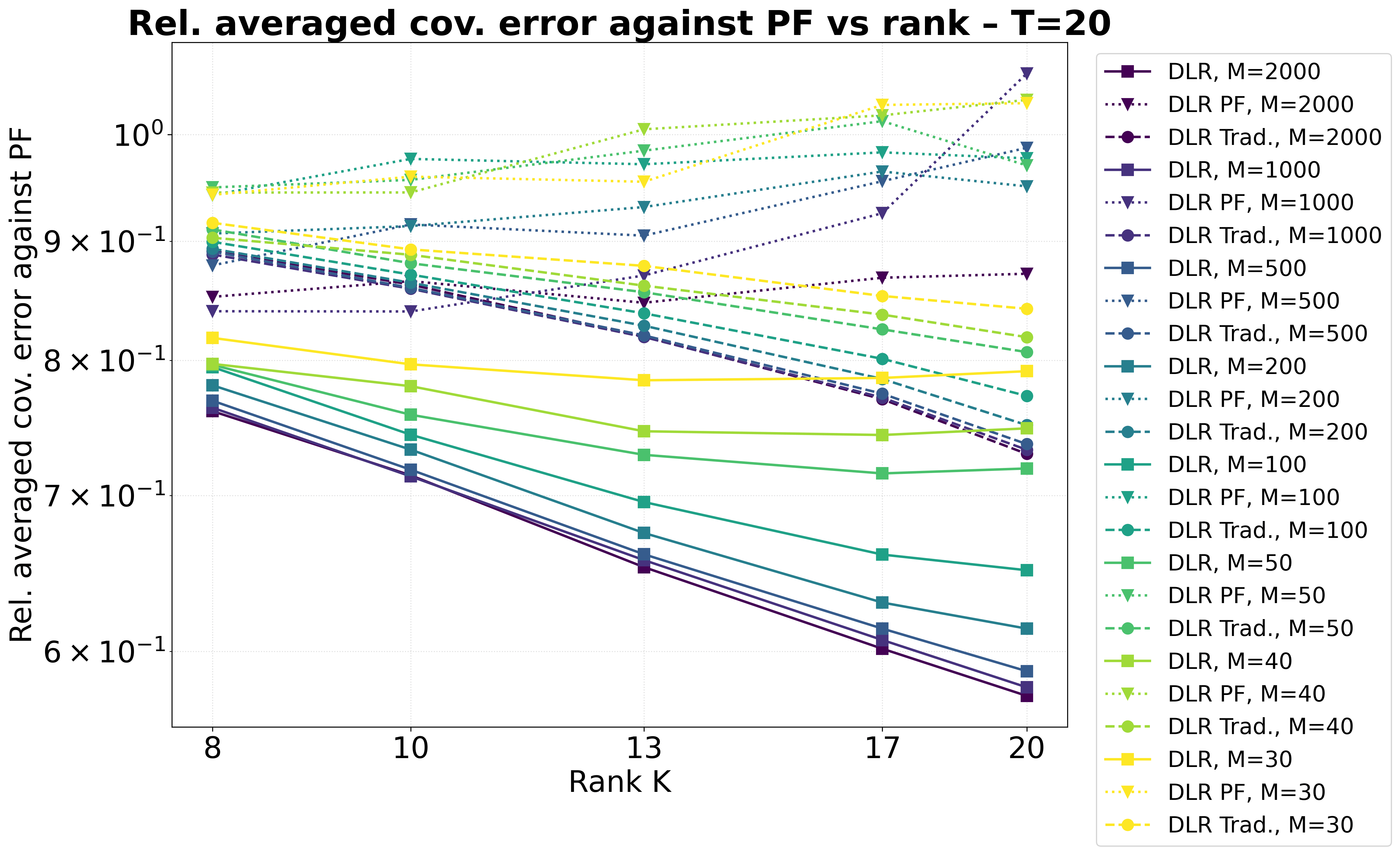} \\
	\includegraphics[scale=0.045]{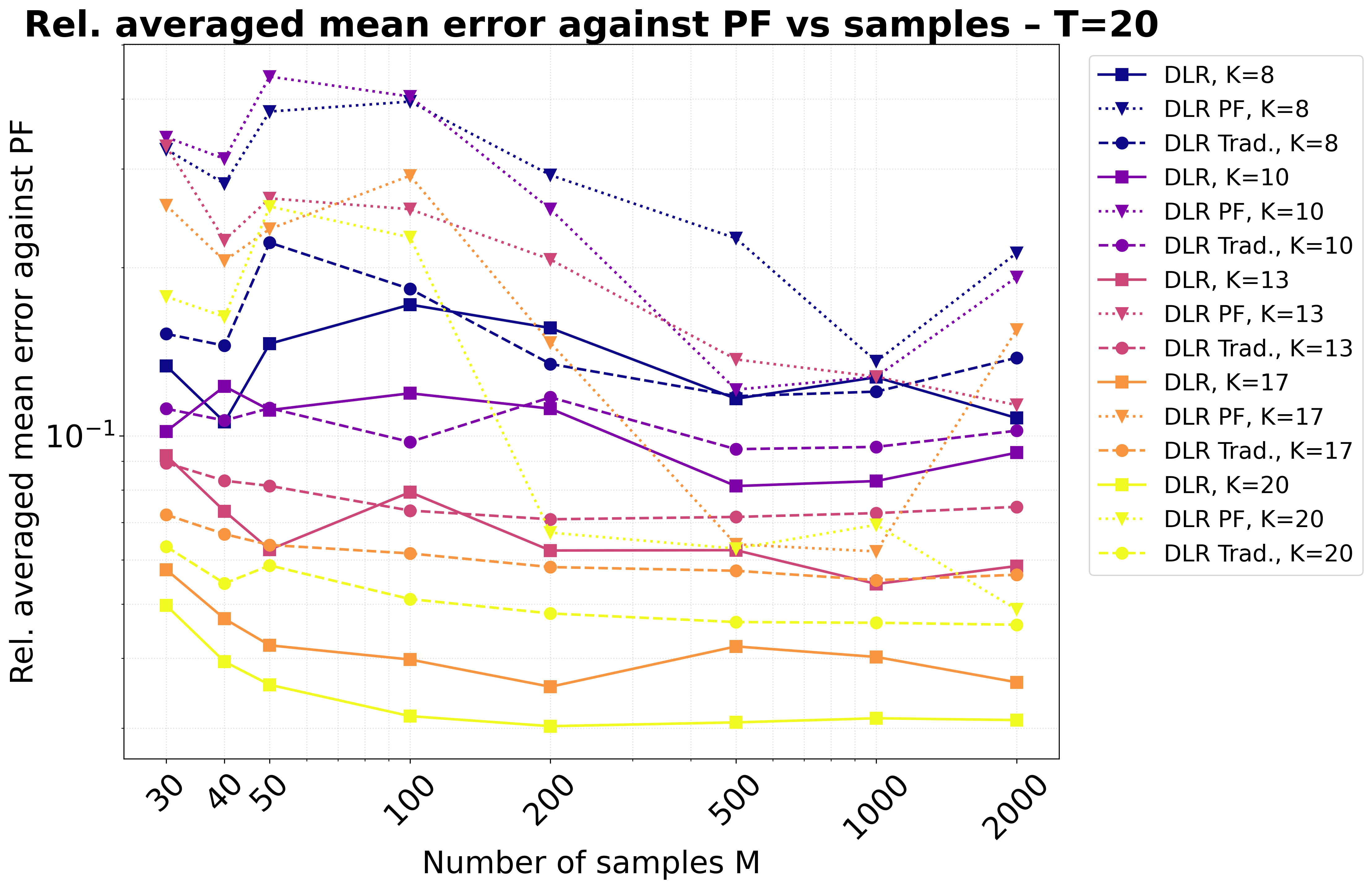}
	\includegraphics[scale=0.045]{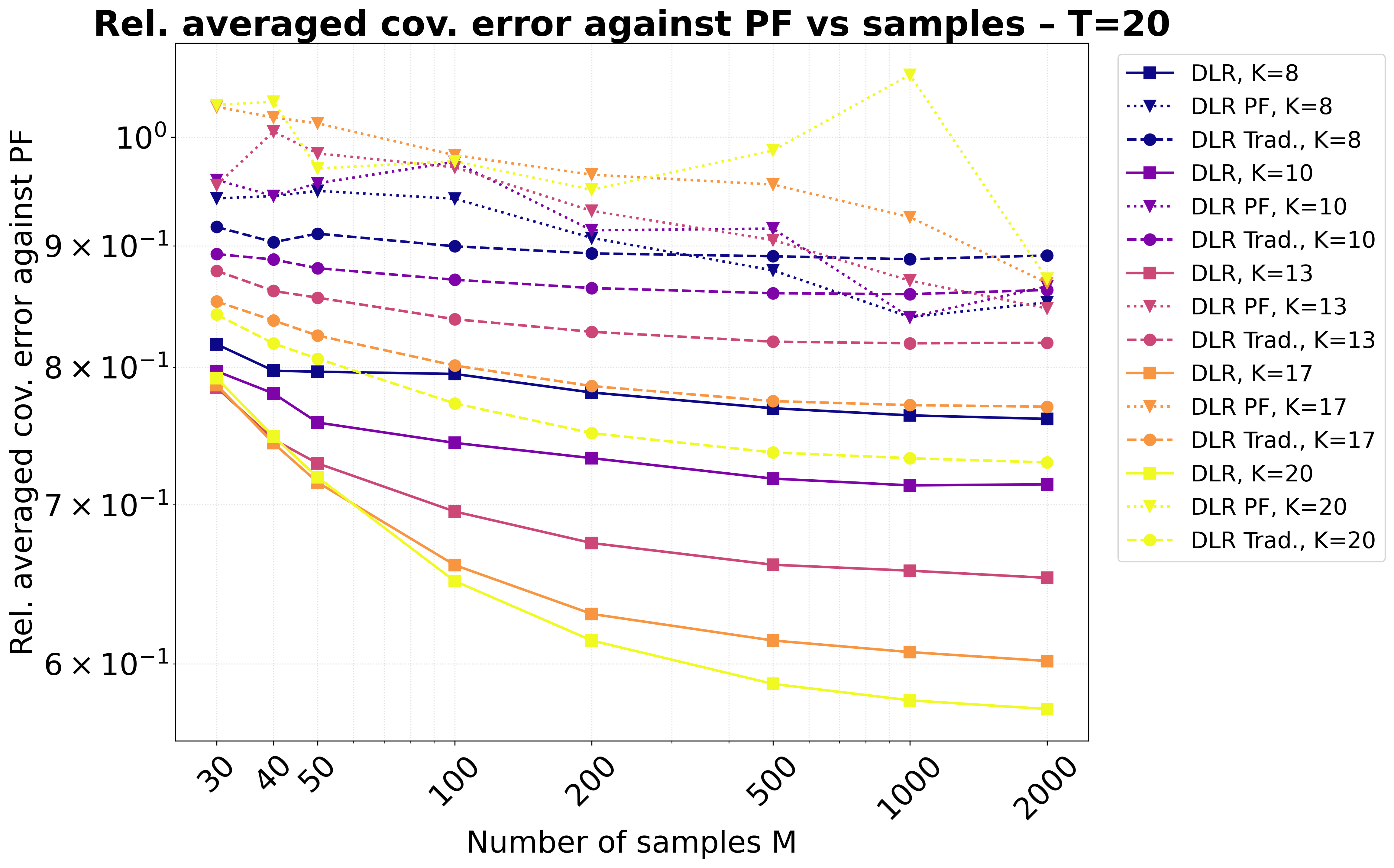} 
	\caption{Rel.\  averaged-over-time \emph{errors against PF}, as a function of the rank $K$ (top row), as a function of the sample size $M$ (bottom row): (left) error on the mean, (right) error on the covariance. Problem \eqref{ex:lorenz96} with continuous-time observations, transient time $T_0=2$, reference computed with $M=1.5 \cdot 10^5$, $\Delta t_{branch}= 4\Delta t$.}
	\label{fig: Lorenz96 PF mean and cov - co}
\end{figure}

From \eqref{fig: Lorenz96 RMSE vs True sol - co}, we see that the PF needs a large number of particles to reach a good accuracy, unlike the EF SDE which performs better with the same number of reference samples ($M=1.5 \cdot 10^5$). We see, as expected, that the DLR JMCO filter also behaves well in approximating the mean and the covariance when compared to the reference filters (see in Figure \ref{fig: Lorenz96 EF SDE mean and cov - co} and \ref{fig: Lorenz96 PF mean and cov - co}). On the other hand, the DLR PF algorithms show bad perfomance with a small number of particles, hence requiring a larger sample size to approximate the same problem.

\begin{figure}[!h]
	\centering
	\includegraphics[scale=0.044]{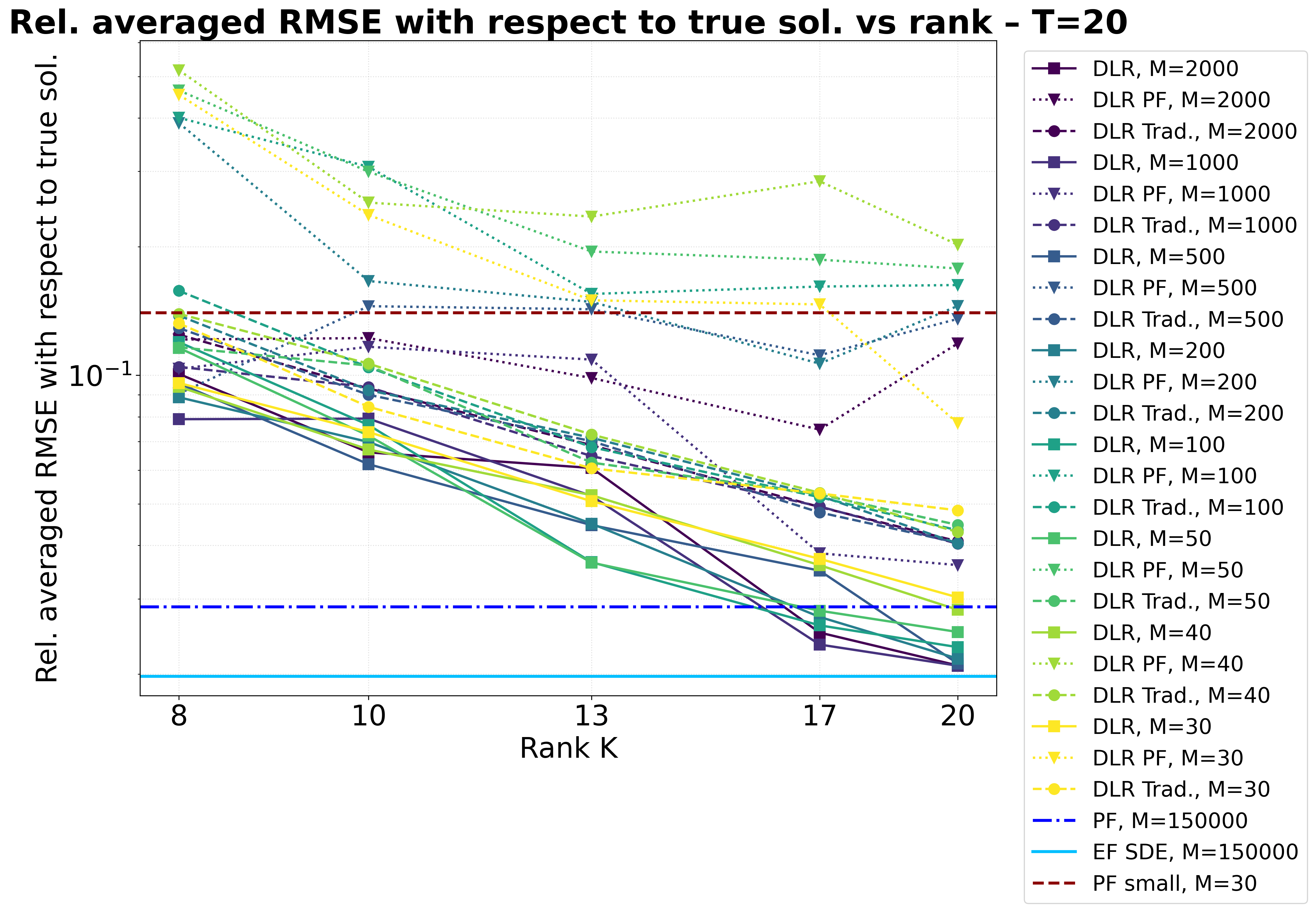}
	\includegraphics[scale=0.044]{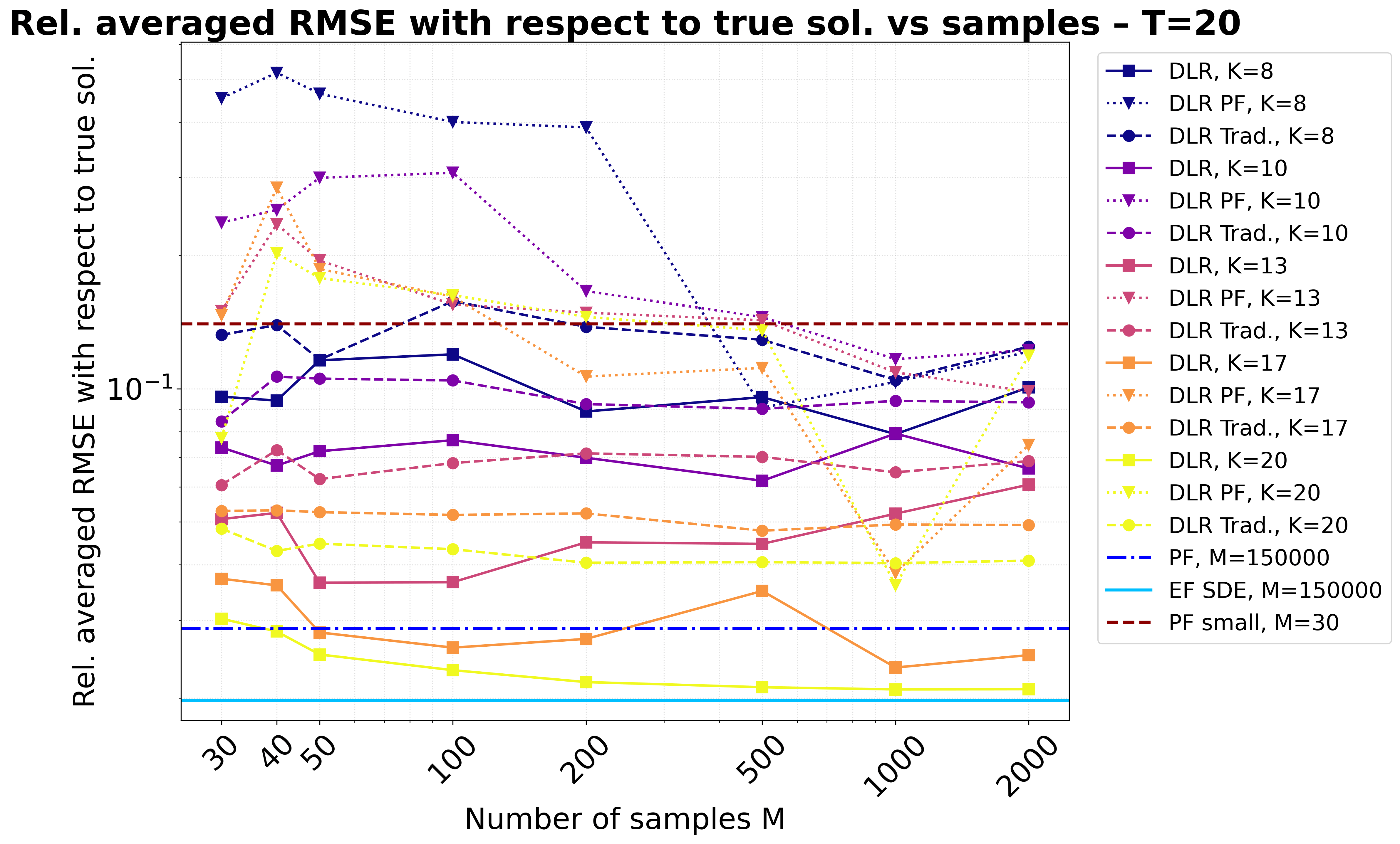}
	\caption{Rel.\  averaged-over-time \emph{RMSE} as a function of the rank $K$ (left) and as a function of the samples sizes $M$ (right). Problem \eqref{ex:lorenz96} with discrete-time observations, transient time $T_0=2$, reference computed with $M=1.5 \cdot 10^5$, $\Delta t_a= 4\Delta t$.}
	\label{fig: Lorenz96 RMSE vs True sol - cd}
\end{figure}

\begin{figure}[!h]
	\centering
	\includegraphics[scale=0.045]{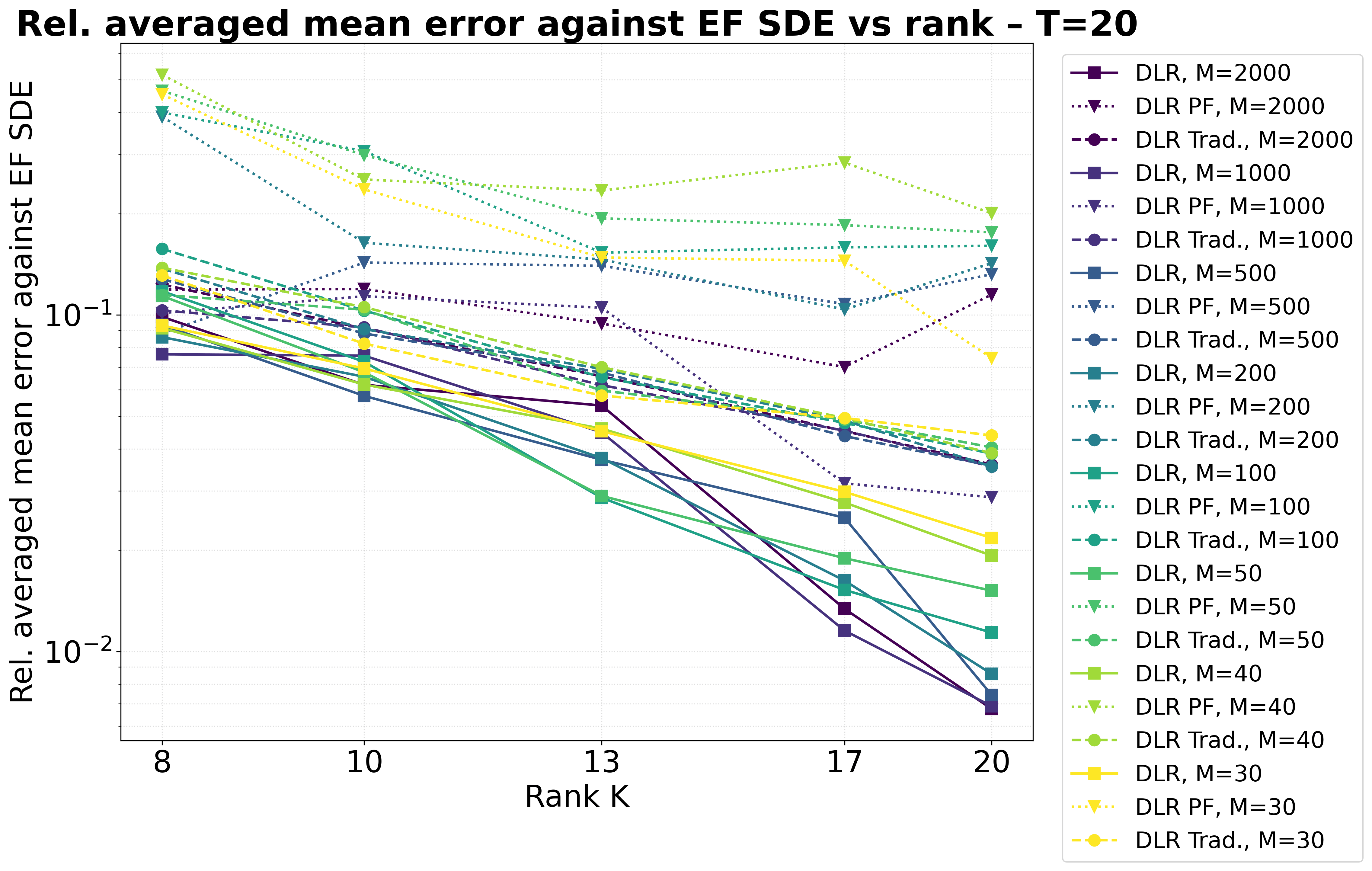}
	\includegraphics[scale=0.045]{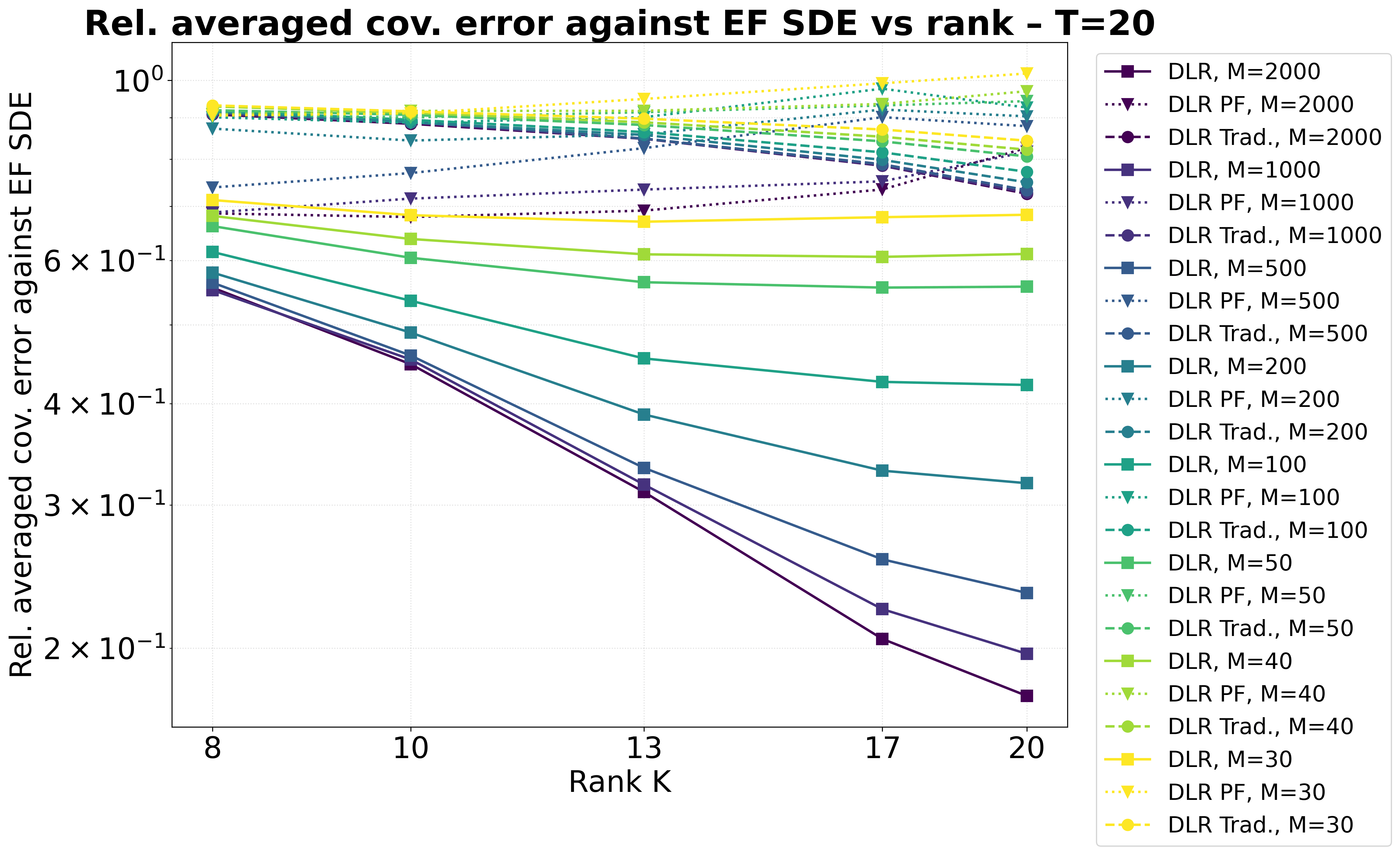} \\
	\includegraphics[scale=0.045]{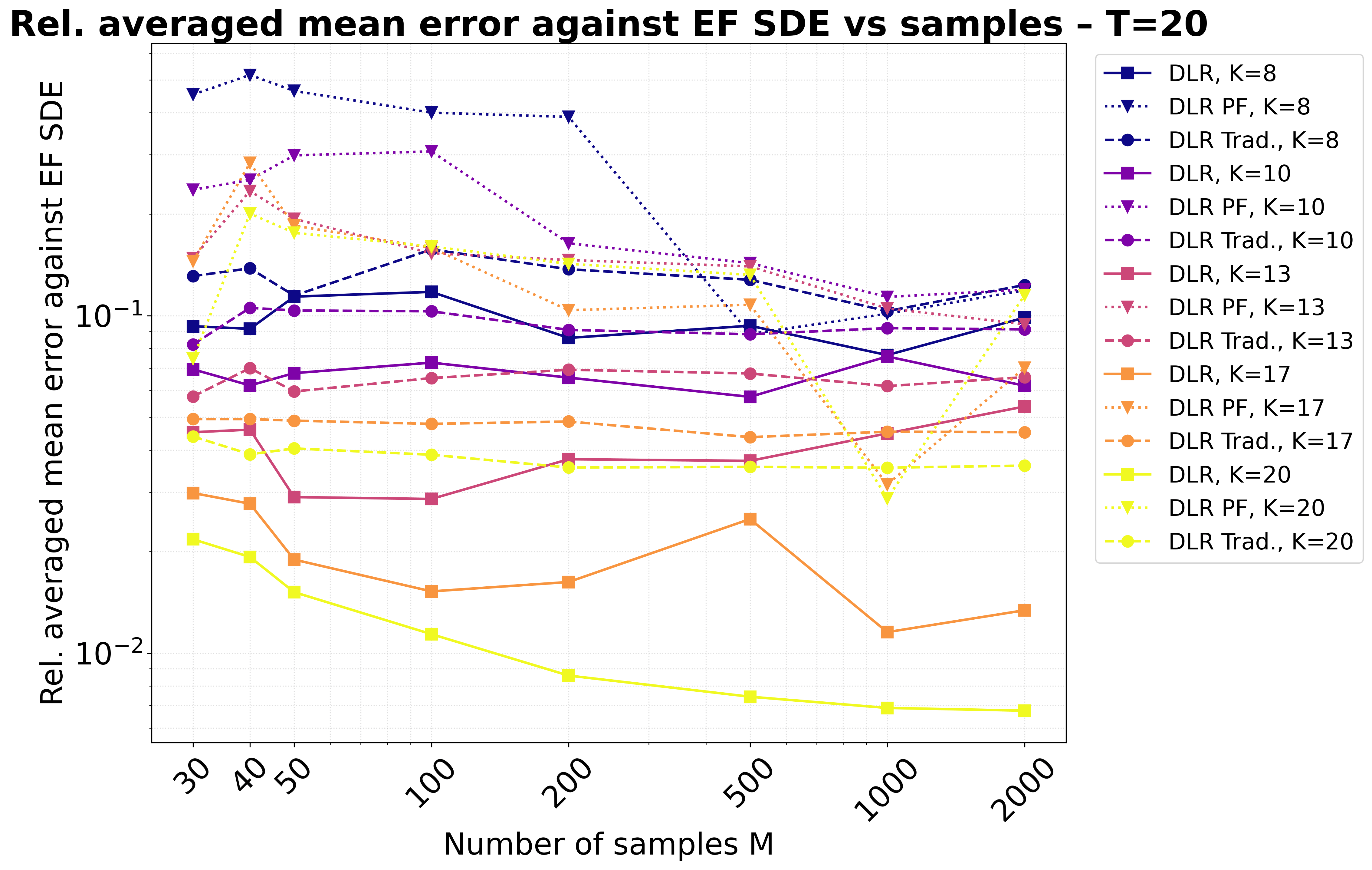}
	\includegraphics[scale=0.045]{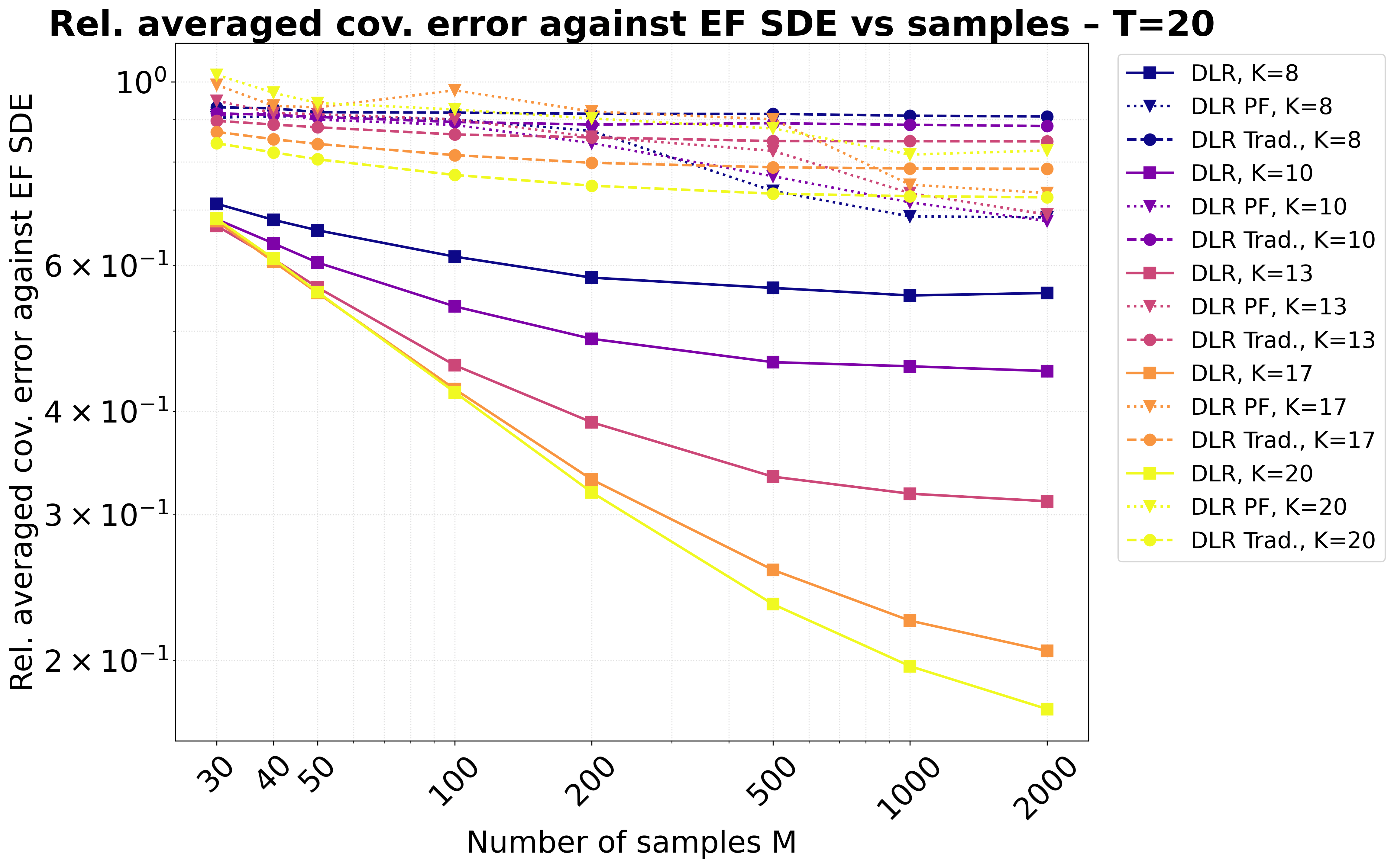} 
	\caption{Rel.\  averaged-over-time \emph{errors against EF SDE}, as a function of the rank $K$ (top row), as a function of the sample size $M$ (bottom row): (left) error on the mean, (right) error on the covariance. Problem \eqref{ex:lorenz96} with discrete-time observations, transient time $T_0=2$, reference computed with $M=1.5 \cdot 10^5$, $\Delta t_a = 10 \Delta t$.}
	\label{fig: Lorenz96 EF SDE mean and cov - cd}
\end{figure}

\begin{figure}[!h]
	\centering
	\includegraphics[scale=0.045]{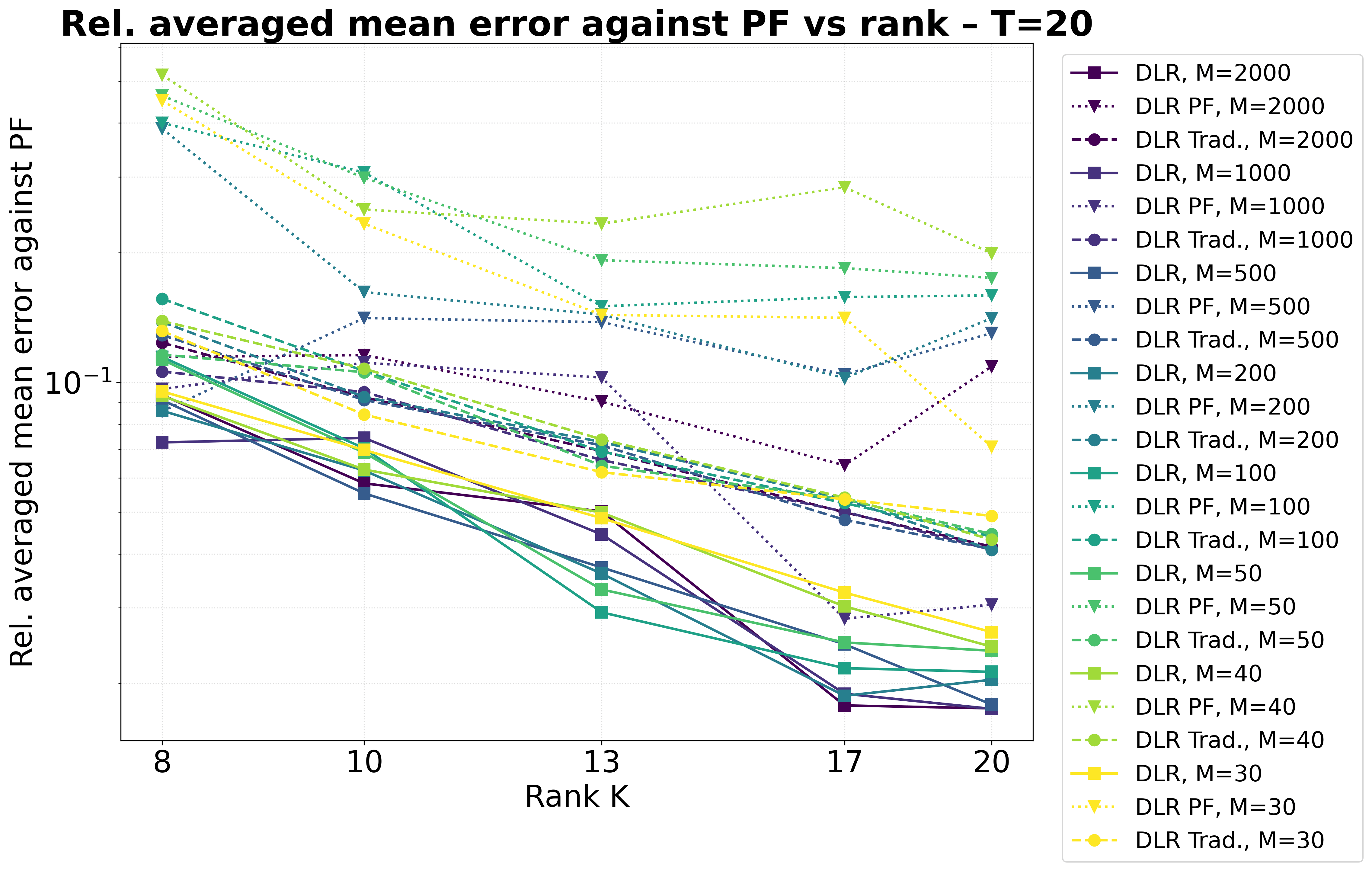}
	\includegraphics[scale=0.045]{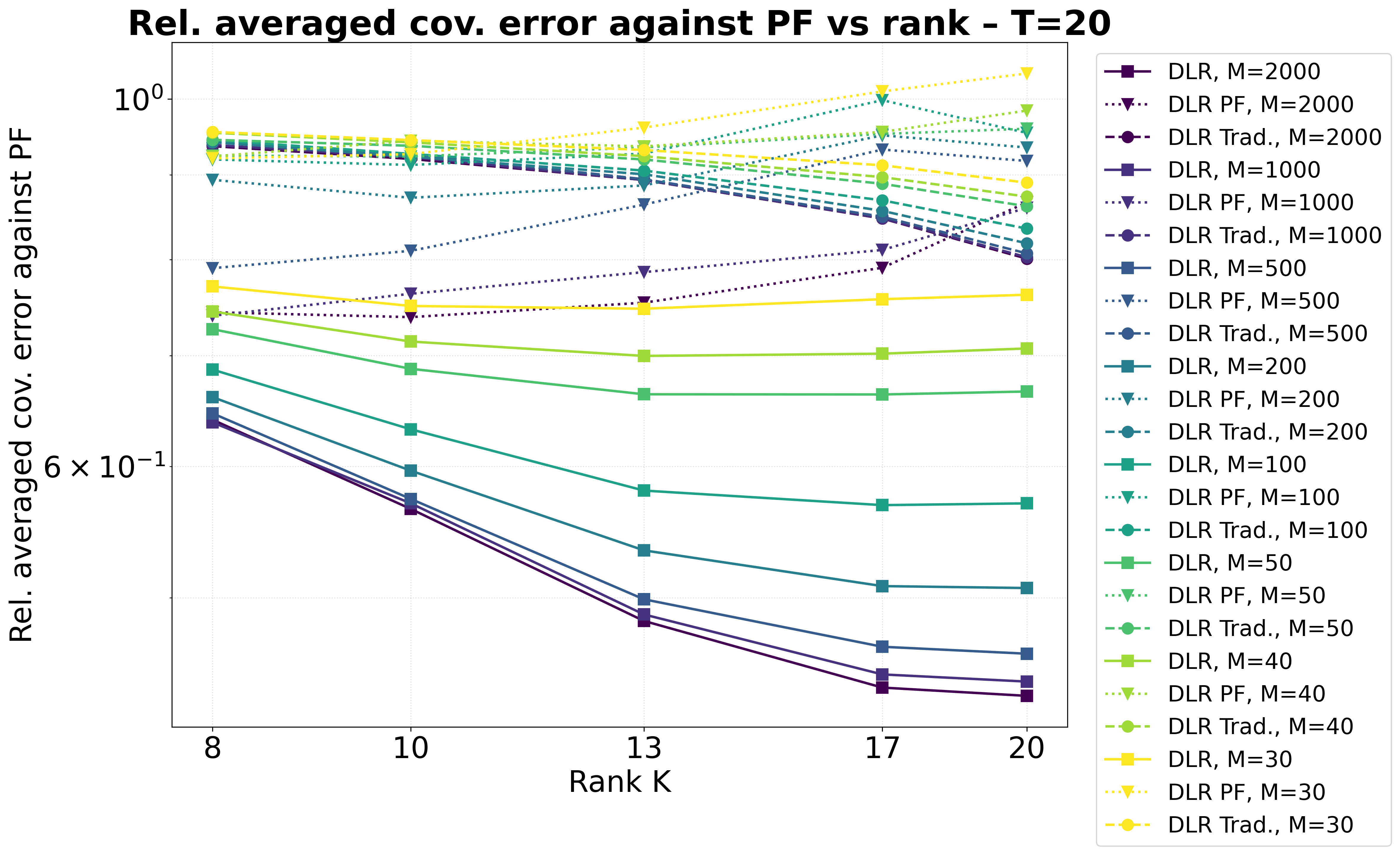} \\
	\includegraphics[scale=0.045]{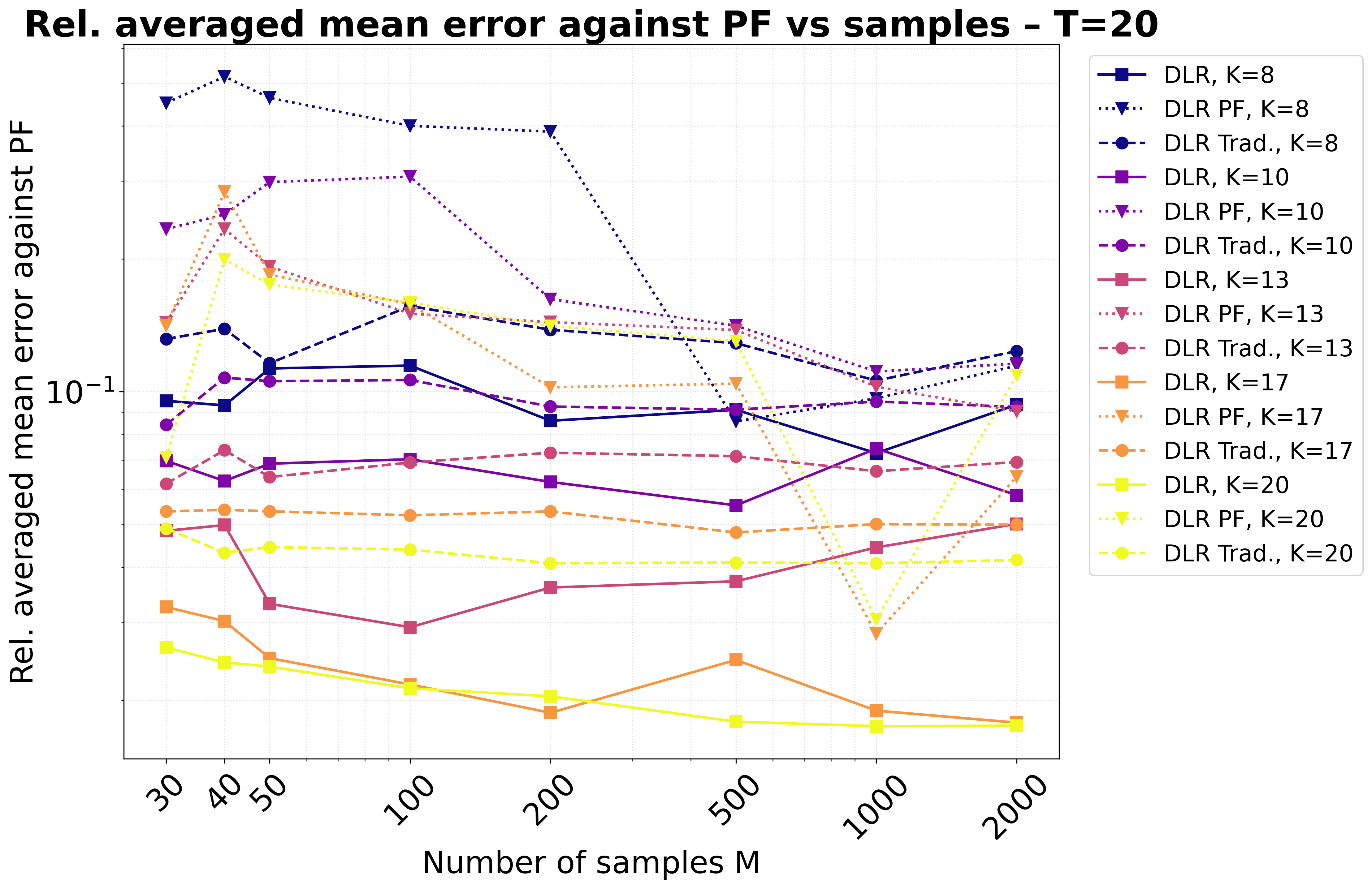}
	\includegraphics[scale=0.045]{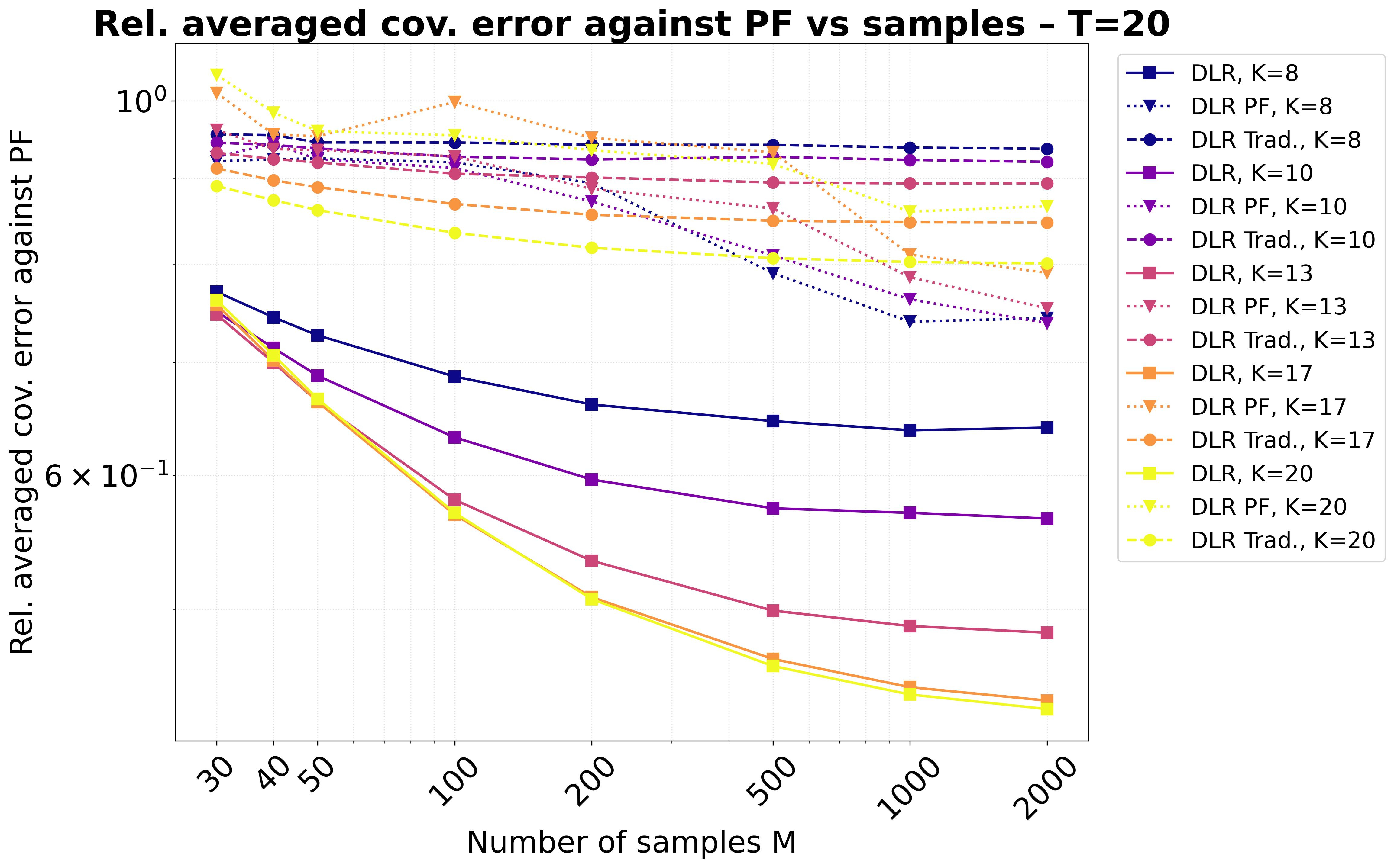} 
	\caption{Rel.\  averaged-over-time \emph{errors against PF}, as a function of the rank $K$ (top row), as a function of the sample size $M$ (bottom row): (left) error on the mean, (right) error on the covariance. Problem \eqref{ex:lorenz96} with discrete-time observations, transient time $T_0=2$, reference computed with $M=1.5 \cdot 10^5$, $\Delta t_a = 10 \Delta t$.}
	\label{fig: Lorenz96 PF mean and cov - cd}
\end{figure}

Similar behavior can be observed for the case of discrete-time observations assimilated every $\Delta t_a = 10 \Delta t$ for the RMSE (Figure \ref{fig: Lorenz96 RMSE vs True sol - cd}) and for the errors on the mean and the covariance against EF SDE reference (Figure \ref{fig: Lorenz96 EF SDE mean and cov - cd}) and the PF reference (Figure \ref{fig: Lorenz96 PF mean and cov - cd}).


\subsection{Nonlinear problem: Two-layer Quasi-Geostrophic System}
The Quasi-Geostrophic model has great relevance in the context of oceanography and atmospheric sciences \cite{majda2012filtering,qi2016low,Vallis2017}, and it is meaningful to study it in the context of data assimilation considering its applications in weather forecasting, where the usual huge dimensionality of the system under study could benefit considerably from using low-rank techniques to reduce the computational cost. The test case considered here describes the geographical flow behavior of two immiscible fluids, a less dense one on the top and a denser one on the bottom. Even though they are immiscible, the two fluids interact in a nonlinear way, and, can produce complex patterns like atmospheric storms or ocean eddies. The two fluids are characterized by the potential vorticity $q$, which is usually considered as the fluid equivalent of angular momentum, and the streamline function $\psi$, whose level curves to which the velocity of the fluid is orthogonal to.

We consider a two-layer quasi-geostrophic (QG) model on a doubly periodic domain $(x,y)\in [0,L_x]\times[0,L_y]$, namely the latitude and the longitude, where the state variables are the potential vorticities $q_1(x,y,t)$ and $q_2(x,y,t)$ associated to the two fluids in the upper and lower layers, respectively. 
The vorticities are related to the streamfunctions $\psi_1$ and $\psi_2$ through the standard two-layer inversion relations
\begin{equation}\label{eq: PV–streamfunction relation}
	\begin{aligned}
		q_1 = \nabla^2 \psi_1 + \alpha(\psi_2 - \psi_1), \qquad
		q_2 = \nabla^2 \psi_2 + \alpha(\psi_1 - \psi_2),
	\end{aligned}
\end{equation}
where $\alpha = \frac{1}{2L_d^2}$, where $L_d=1.0$ is called the Rossby deformation, which is a measure on how strongly the two layers interact, and $\nabla^2$ indicates the (horizontal) Laplacian. The spatial derivatives of the streamfunctions determine the velocity fields in the two layers. Moreover, for the sake of well-posedness, each streamfunction is constrained to have zero spatial mean, i.e.\ $\int \psi_i(x,y,t)\mathrm{d}x \mathrm{d}y = 0,$ for $i=1,2$.

The evolution of the system is governed by the potential vorticity equations
\begin{equation}\label{eq: 2 qg}
	\begin{aligned}
		\frac{\partial q_1}{\partial t}
		+	
		J(\psi_1, q_1) + U \frac{\partial q_1}{\partial x}
		+ (\beta + L_d^2 U) \frac{\partial \psi_1}{\partial x}
		+ \nu \nabla^2 q_1
		+ r_1 q_1	
		&= F_1(x,y,t)
		+\eta_1(x,y,t), \\
		\frac{\partial q_2}{\partial t}
		+
		J(\psi_2, q_2) - U \frac{\partial q_2}{\partial x}
		+ (\beta - L_d^2 U) 
		\frac{\partial \psi_2}{\partial x}
		+ \nu \nabla^2 q_2
		+ r_2 q_2
		&= F_2(x,y,t)
		+ \eta_2(x,y,t),
	\end{aligned}
\end{equation}
where $J(\psi_i,q_i):= \frac{\partial \psi_i}{\partial x}\frac{\partial q_i}{\partial y}-\frac{\partial \psi_i}{\partial y}\frac{\partial q_i}{\partial x}$ denotes the nonlinear Jacobian. In \eqref{eq: 2 qg}, $\nu$ is a hyperviscosity coefficient linked to harmonic dissipation, $U$ is the zonal mean shear, namely the average change in speed of east-west (zonal) winds over distance, and $r_i$ are linear damping coefficients, terms that try to stabilize turbulence. Concerning the possible chaotic behavior of this system, $\beta$ regulates the meridional gradient of the Coriolis parameter, which generates the so-called Rossby waves and large-scale turbulence \cite{Vallis2017}, $F_i$ represents deterministic forcing term, and a colored noise term $\eta_i(x,y,t)$ which is defined as follows.
Each $\eta_i$ is constructed as a low-rank, spatially correlated noise acting through first barotropic  $m_t$ and baroclinic $m_c$ streamfunction modes. Specifically, random perturbations are introduced in the Fourier space using a decomposition into barotropic modes $\psi_1=\psi_2$ and baroclinic modes $\psi_1=-\psi_2$, each driven by independent temporal stochastic processes. These perturbations are then mapped into potential vorticity through the two-layer relation \eqref{eq: PV–streamfunction relation}, ensuring that the stochastic forcing is dynamically consistent with the QG structure. This construction yields a physically informed, low-dimensional noise model that selectively excites large-scale balanced (barotropic) and shear-driven (baroclinic) dynamics \cite{majda2018strategies,Vallis2017}:
\begin{equation}
	\begin{aligned}
		\psi_1^{\mathrm{noise}}(x,y,t) &=
		c_{\mathrm{bt}} \sum_{j=1}^{m_{\mathrm{bt}}} \phi_j^{\mathrm{bt}}(x,y)\,\dot{W}_j^{\mathrm{bt}}(t)
		\;+\;
		c_{\mathrm{bc}} \sum_{j=1}^{m_{\mathrm{bc}}} \phi_j^{\mathrm{bc}}(x,y)\,\dot{W}_j^{\mathrm{bc}}(t),
		\\[6pt]
		\psi_2^{\mathrm{noise}}(x,y,t) &=
		c_{\mathrm{bt}} \sum_{j=1}^{m_{\mathrm{bt}}} \phi_j^{\mathrm{bt}}(x,y)\,\dot{W}_j^{\mathrm{bt}}(t)
		\;-\;
		c_{\mathrm{bc}} \sum_{j=1}^{m_{\mathrm{bc}}} \phi_j^{\mathrm{bc}}(x,y)\,\dot{W}_j^{\mathrm{bc}}(t).
	\end{aligned}
\end{equation}
where $c_{bt}$ and $c_{bc}$ are real amplitudes for the barotropic and the baroclinic modes, respectively. 
The corresponding colored noise is finally obtained through relation \eqref{eq: PV–streamfunction relation}:
\begin{equation}
	\begin{aligned}
		\eta_1(x,y,t) &= \nabla^2 \psi_1^{\mathrm{noise}} +\alpha(\psi_2^{\mathrm{noise}} - \psi_1^{\mathrm{noise}}), \\
		\eta_2(x,y,t) &= \nabla^2 \psi_2^{\mathrm{noise}} +\alpha(\psi_1^{\mathrm{noise}} - \psi_2^{\mathrm{noise}}).
	\end{aligned}
\end{equation}
The difference of amplitude scale between $c_{\mathrm{bt}}$ and $c_{\mathrm{bc}}$, as well as the possible presence of a forcing term, generates turbulence in system \eqref{eq: 2 qg}.

To easily deal with derivatives in relation \eqref{eq: PV–streamfunction relation}, we approximate the spatial dimension in \eqref{eq: 2 qg} via a Fourier spectral method. On the other hand, the integration over time is done using a standard Euler-Maruyama method, denoting $q^{\Delta x}_t$ the state vector containing the spatially discretized components at time $t$. The considered observation process is
\begin{equation}\label{eq: 2 qg - obs proc}
	\mathrm{d} Z_{t} = H q^{\Delta x}_t \mathrm{d}t + R^{\frac{1}{2}} \mathrm{d}B_t.
\end{equation}
where the operator $H$ observes $256$ points in each layer, with equidistant coordinates in both axis, affected by a noise with covariance $R= 10^{-2}\cdot I_{512 \times 512}$. Therefore, the dimension $h$ is equal to $512$.

Our domain is a square of dimension $L_x = L_y = 2 \pi$.
We consider a final time of $T=60$, and we discretized the full-order model with a time step $\Delta t_n=0.02$. 
We choose $m_t=m_c=16$. We then choose the baroclinic modes are using the condition
$7.9 \leq K^2 \leq 13.1$,
where $K^2 = k_x^2 + k_y^2$, with $k_x,k_y$ the Fourier harmonics in $x$ and $y$.
Therefore, the admissible values of \(K^2\) inside this interval are
$K^2 \in \{8,9,10,13\}$, and the corresponding Fourier modes are
\begin{equation}
	\begin{aligned}
		K^2 = 8 &: \qquad (k_x,k_y)=(\pm 2,\pm 2),\\
		K^2 = 9 &: \qquad (k_x,k_y)=(\pm 3,0),(0,\pm 3),\\
		K^2 = 10 &: \qquad (k_x,k_y)=(\pm 3,\pm 1),(\pm 1,\pm 3),\\
		K^2 = 13 &: \qquad (k_x,k_y)=(\pm 3,\pm 2),(\pm 2,\pm 3).
	\end{aligned}
\end{equation}
Thus, the baroclinic forcing acts on a relatively narrow band of
intermediate wavenumbers centred approximately around $K_{\mathrm{BC}} \approx 3$.
The amplitude of the fast barotropic modes is $c_{bt}=0.06$, whereas the one of the slow baroclinic modes are $c_{bc}=3$. This structure of the barotropic-baroclinic modes is chosen to obtain baroclinic instability. We set the rank of our DLRA-filtering algorithms to $k=32$, i.e.\ $k=m_t+m_c$ to well capture this stochastic perturbations. The spectral approximation uses $32 \times 32$ Fourier modes in each layer, hence the dimension fo the discretized system is $d=2048$. The initial condition is also random and follows a mixture of colored Gaussians.

The physical parameters of the model are set to $U=0.6$, $\beta = 0.05$, $\nu = 4e-5$, $r_1=0$, $r_2=0.02$, whereas the forcing term $F$ is deterministic and acts in a slighty asymmetric baroclinic way such that in the first layer we have 
\begin{equation*}
	\begin{aligned}
		F_1(x,y,t)=& -0.002 \left(\sin(\frac{4  \pi  y}{L_y})+0.5\cos(\frac{4  \pi  x}{L_x})\sin(\frac{4  \pi  y}{L_y})\right)
	\end{aligned}
\end{equation*} acting on the second wavenumber, and in the second layer we set $F_2(x,y,t)=-F_1(x,y,t)$.

\begin{figure}[!h]
	\centering 
	\includegraphics[scale=0.052]{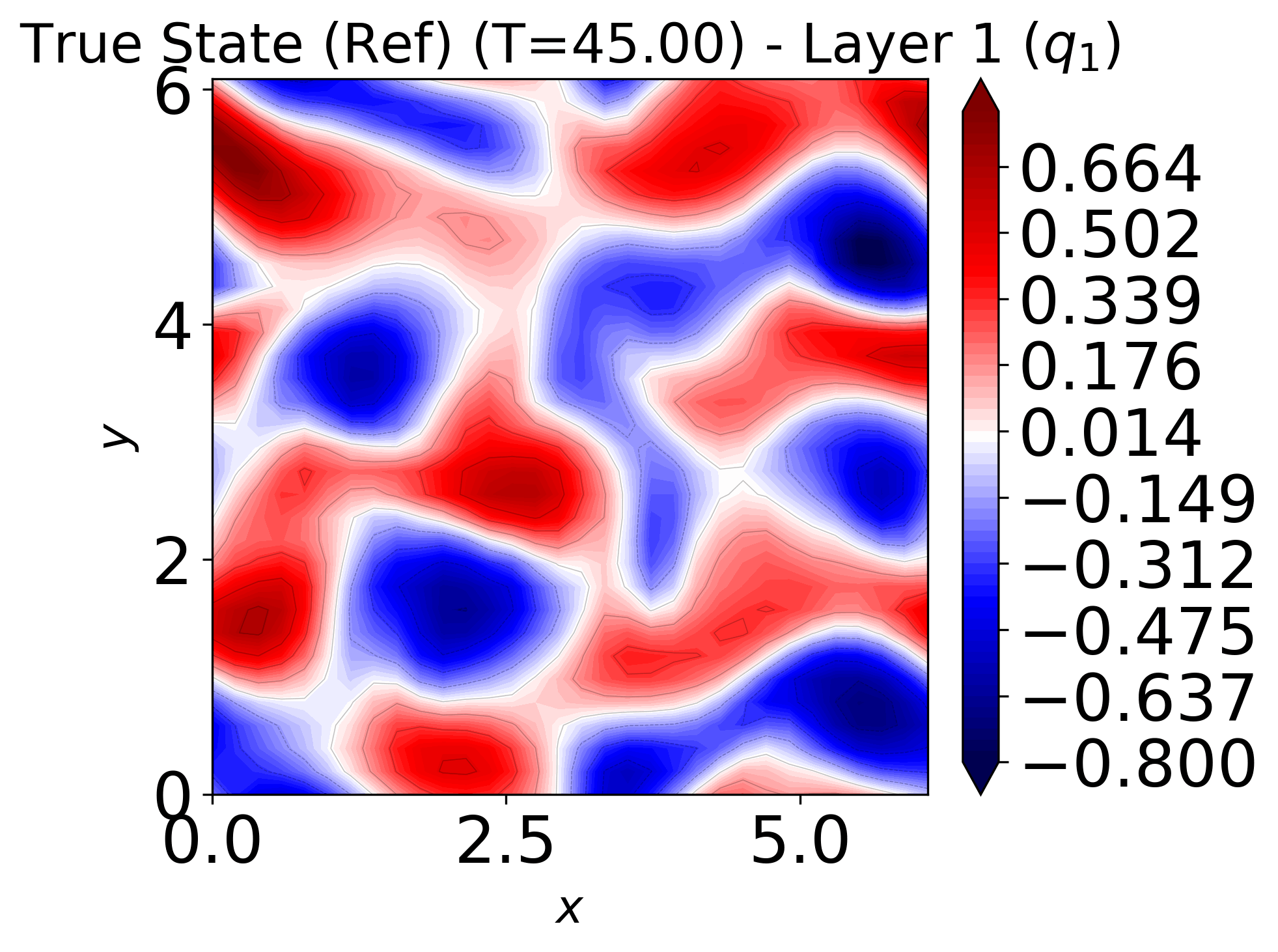} 
	\includegraphics[scale=0.052]{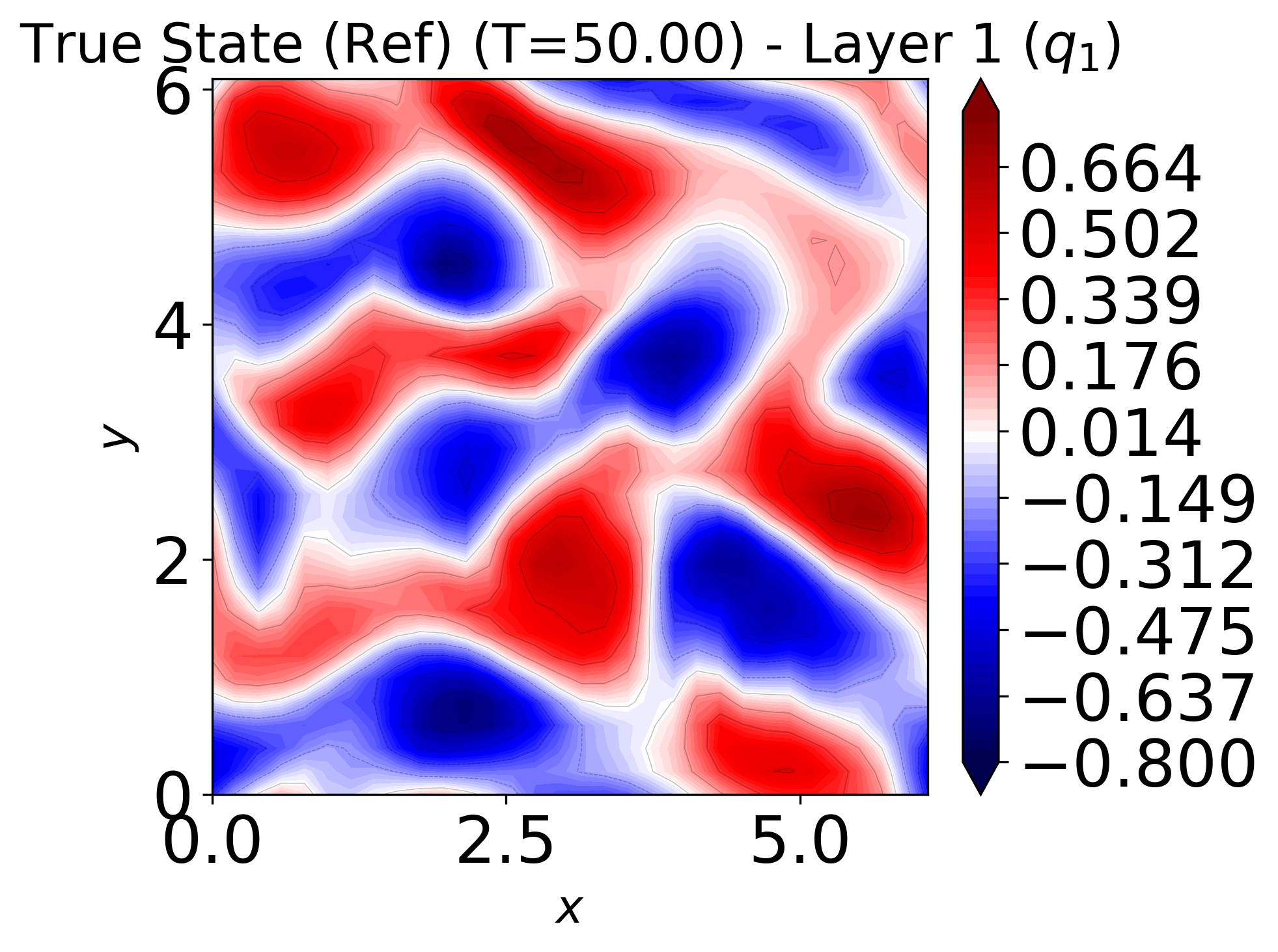} 
	\includegraphics[scale=0.052]{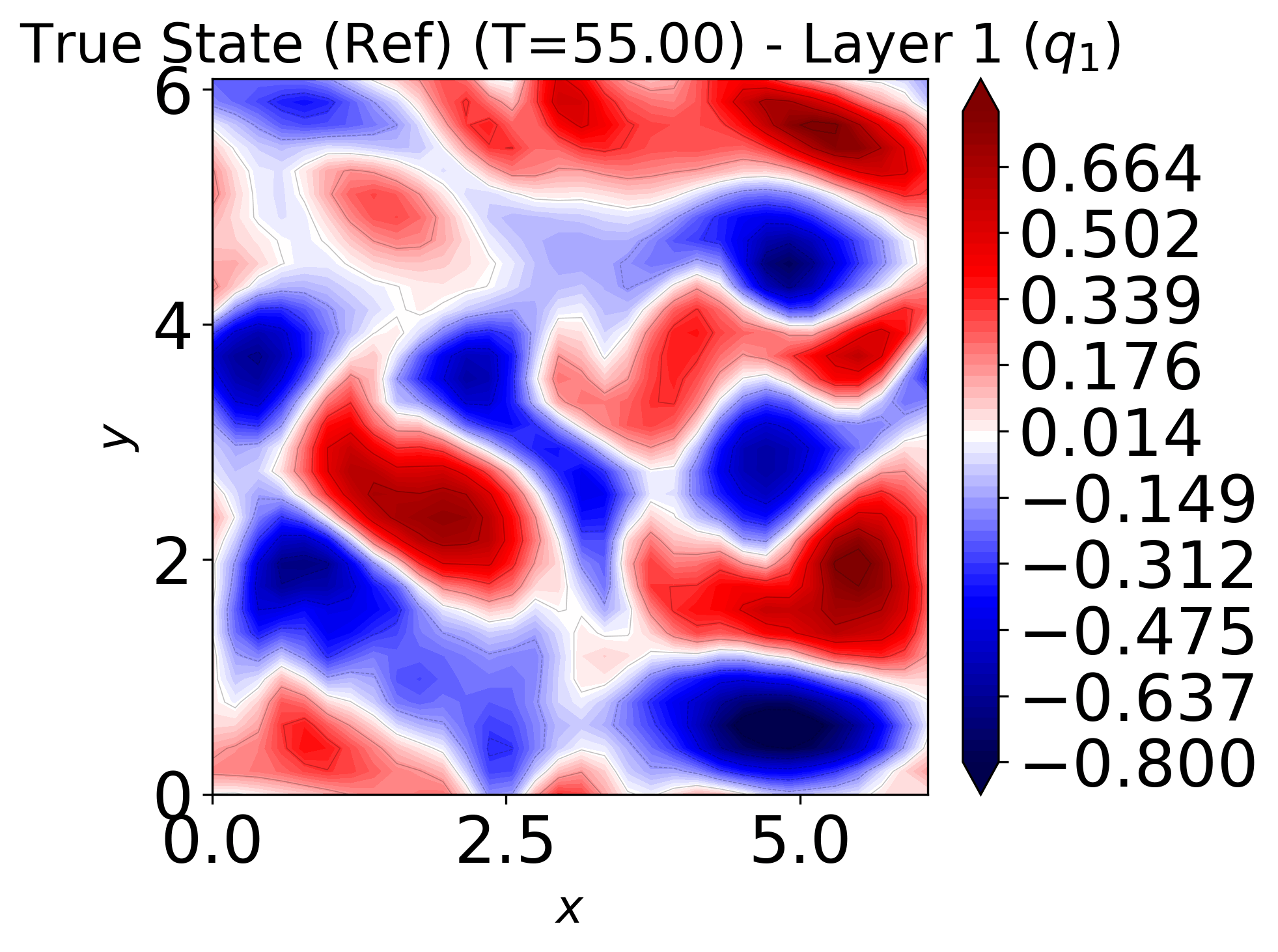} 
	\includegraphics[scale=0.052]{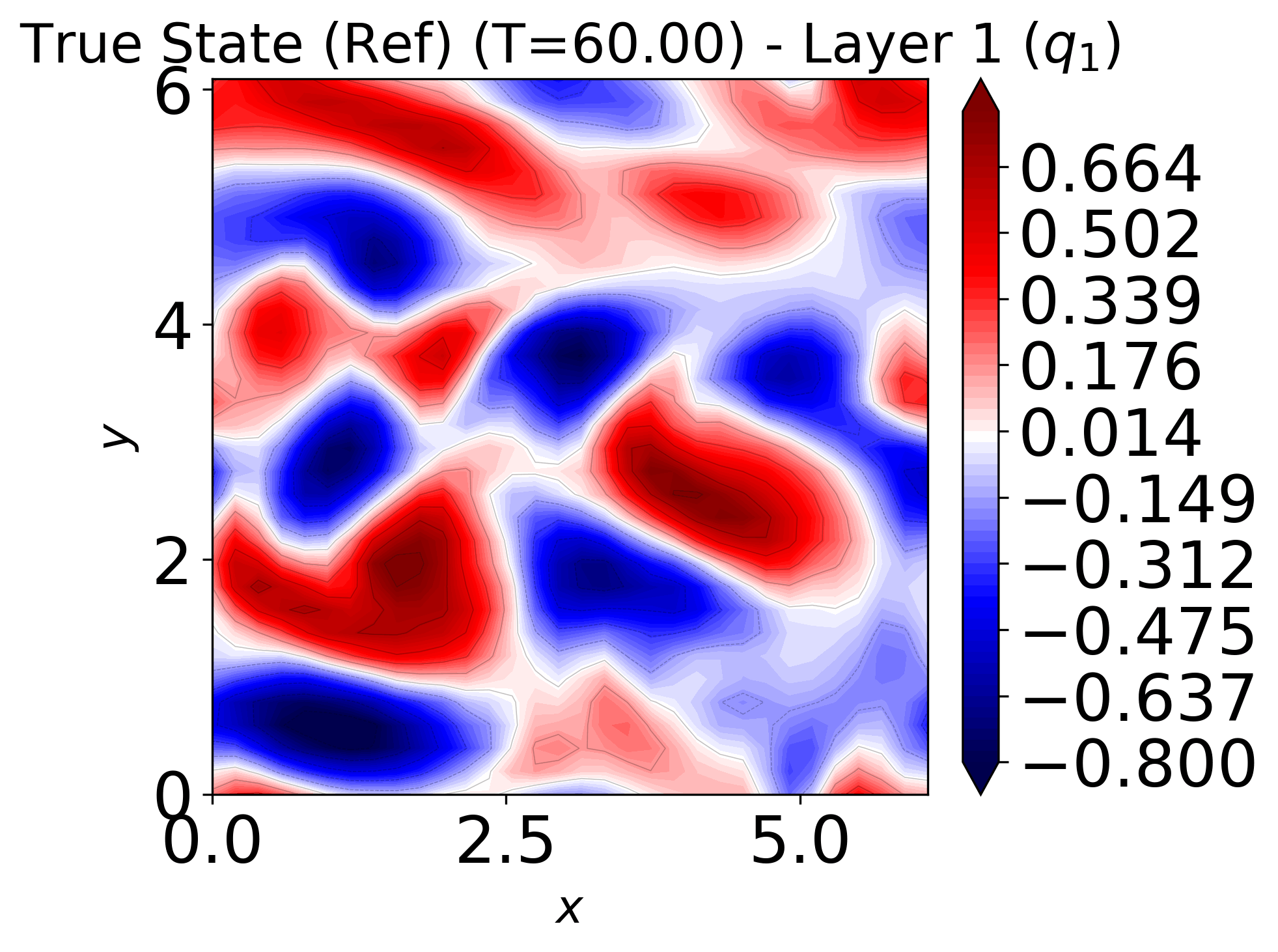}\\
	\includegraphics[scale=0.051]{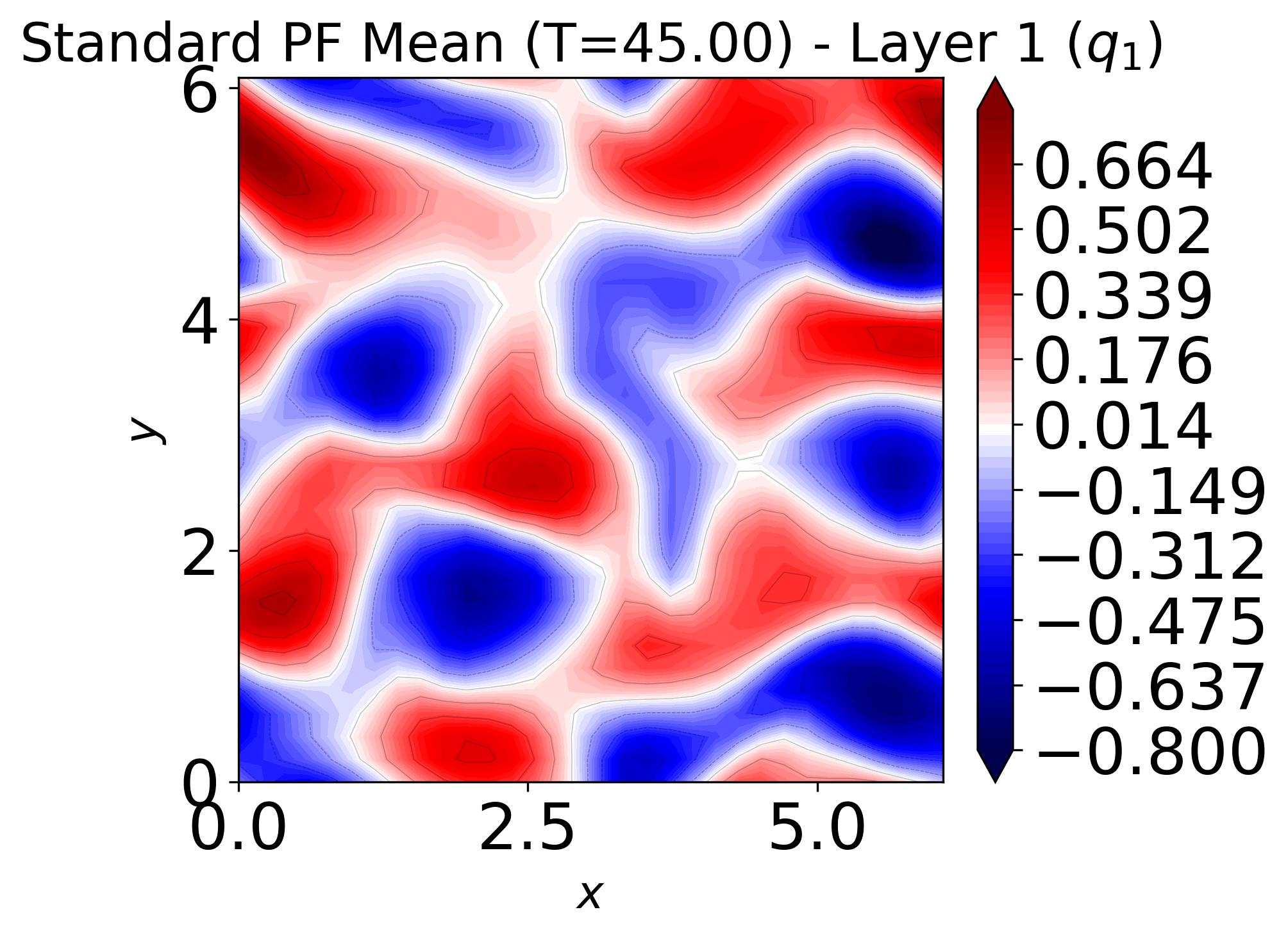} 
	\includegraphics[scale=0.052]{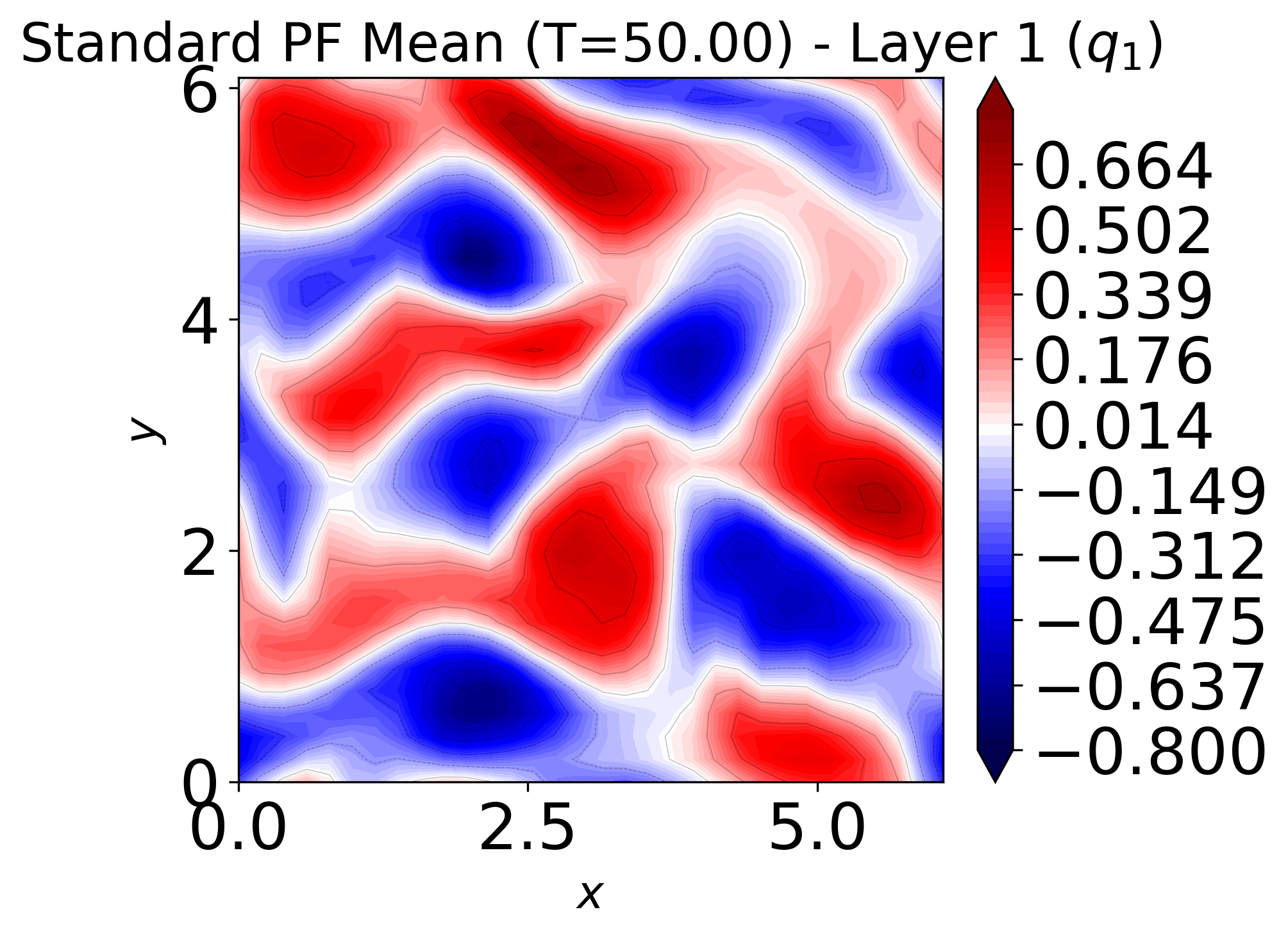} 
	\includegraphics[scale=0.052]{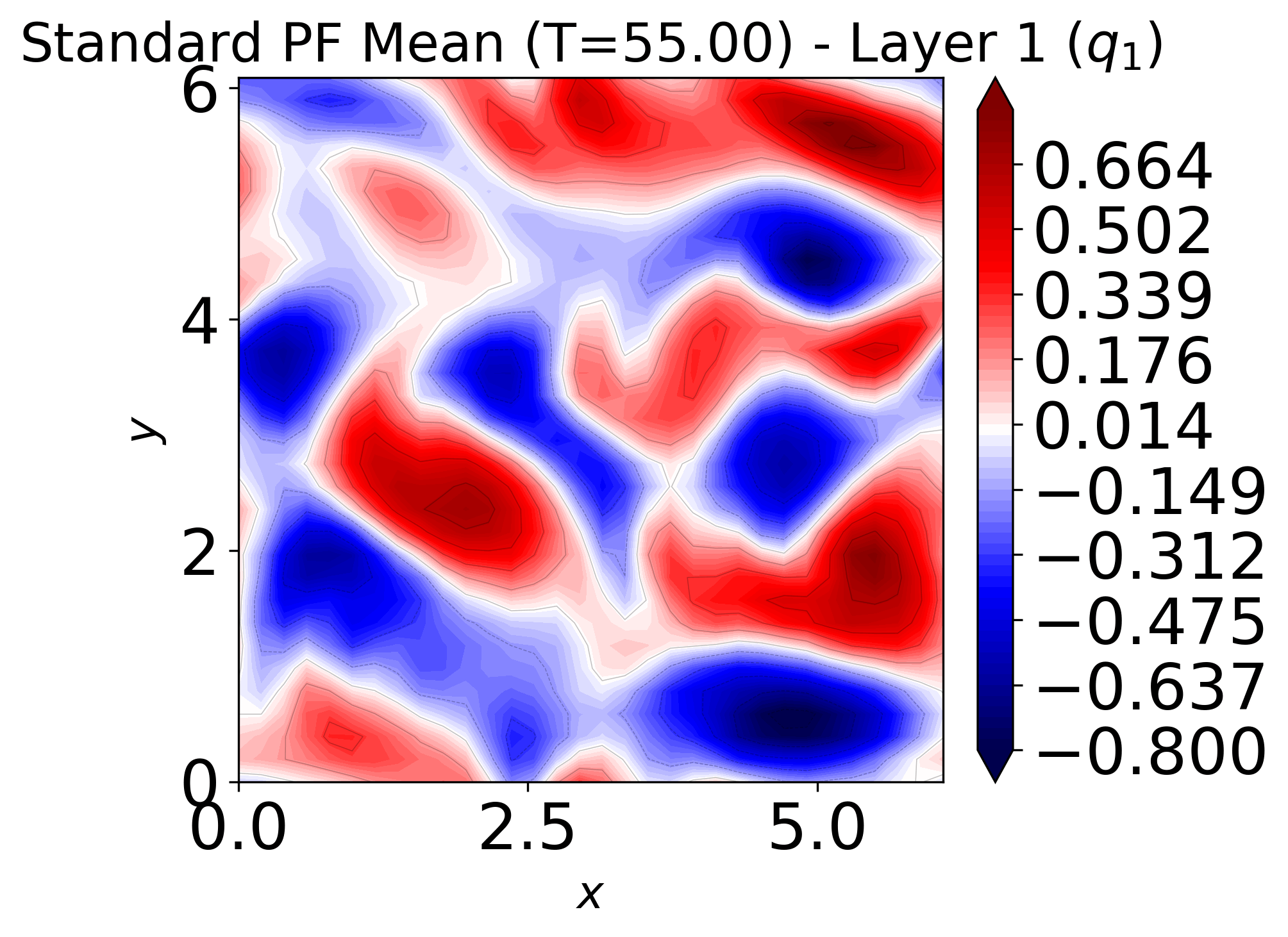} 
	\includegraphics[scale=0.052]{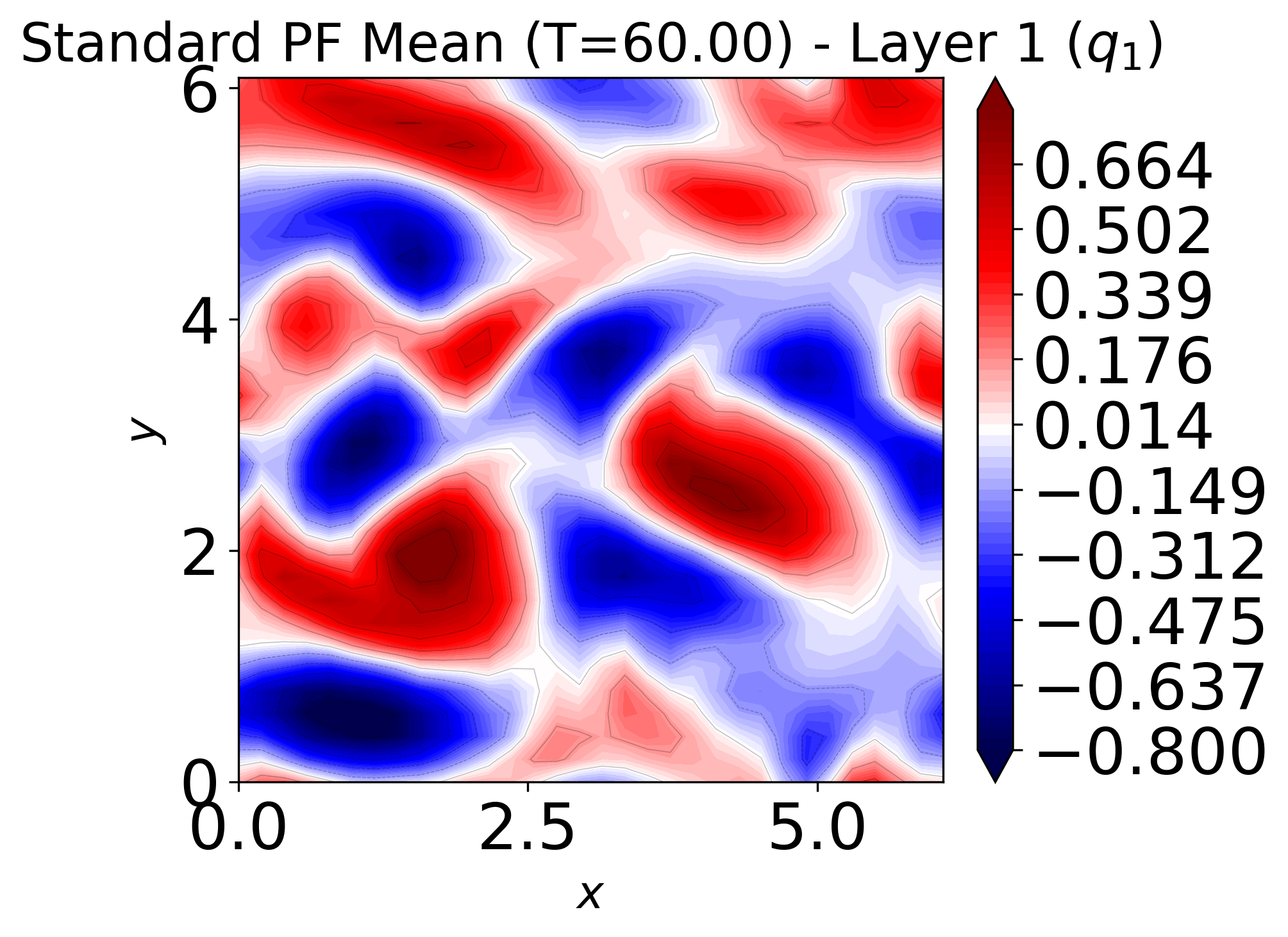}\\
	\includegraphics[scale=0.052]{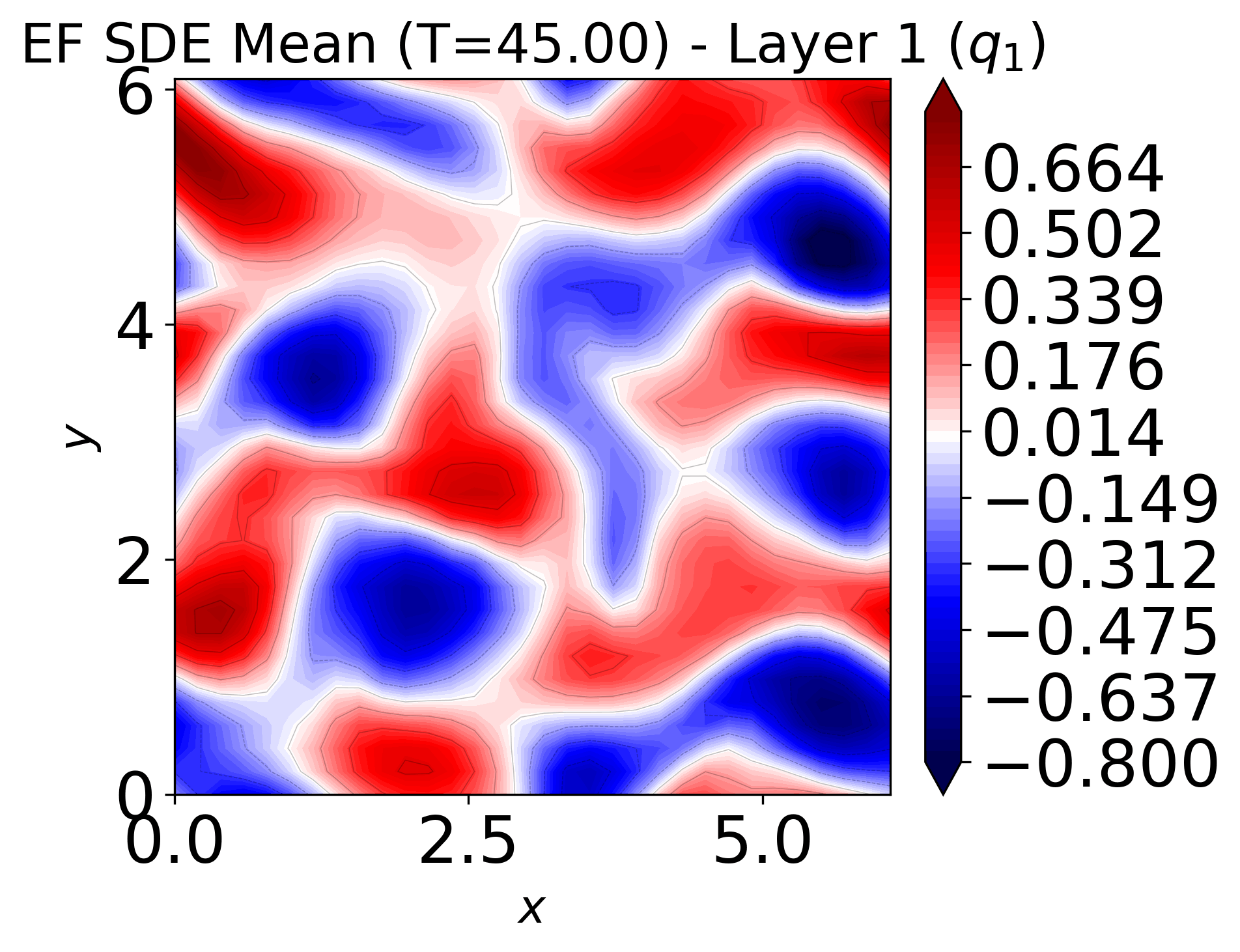}
	\includegraphics[scale=0.052]{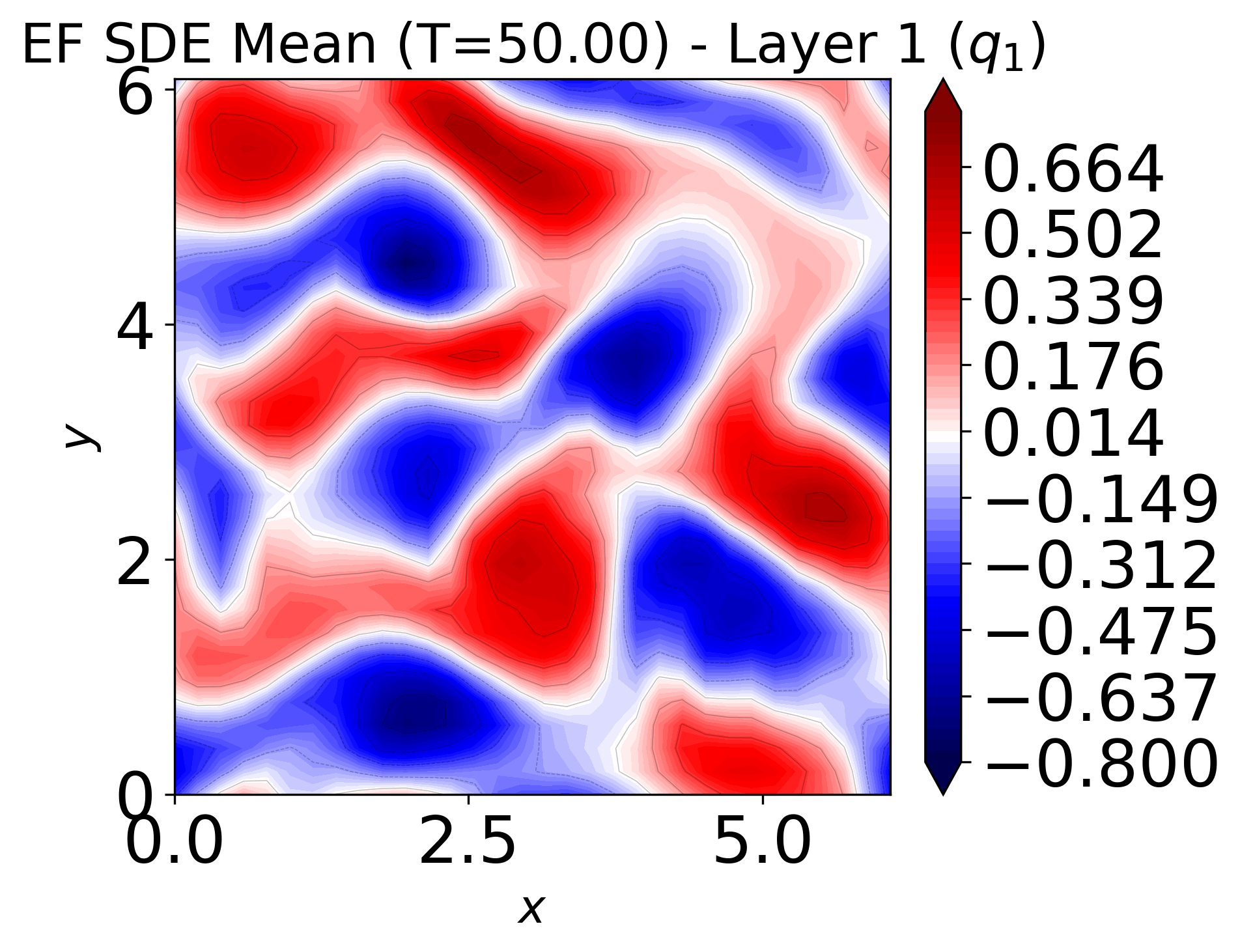} 
	\includegraphics[scale=0.052]{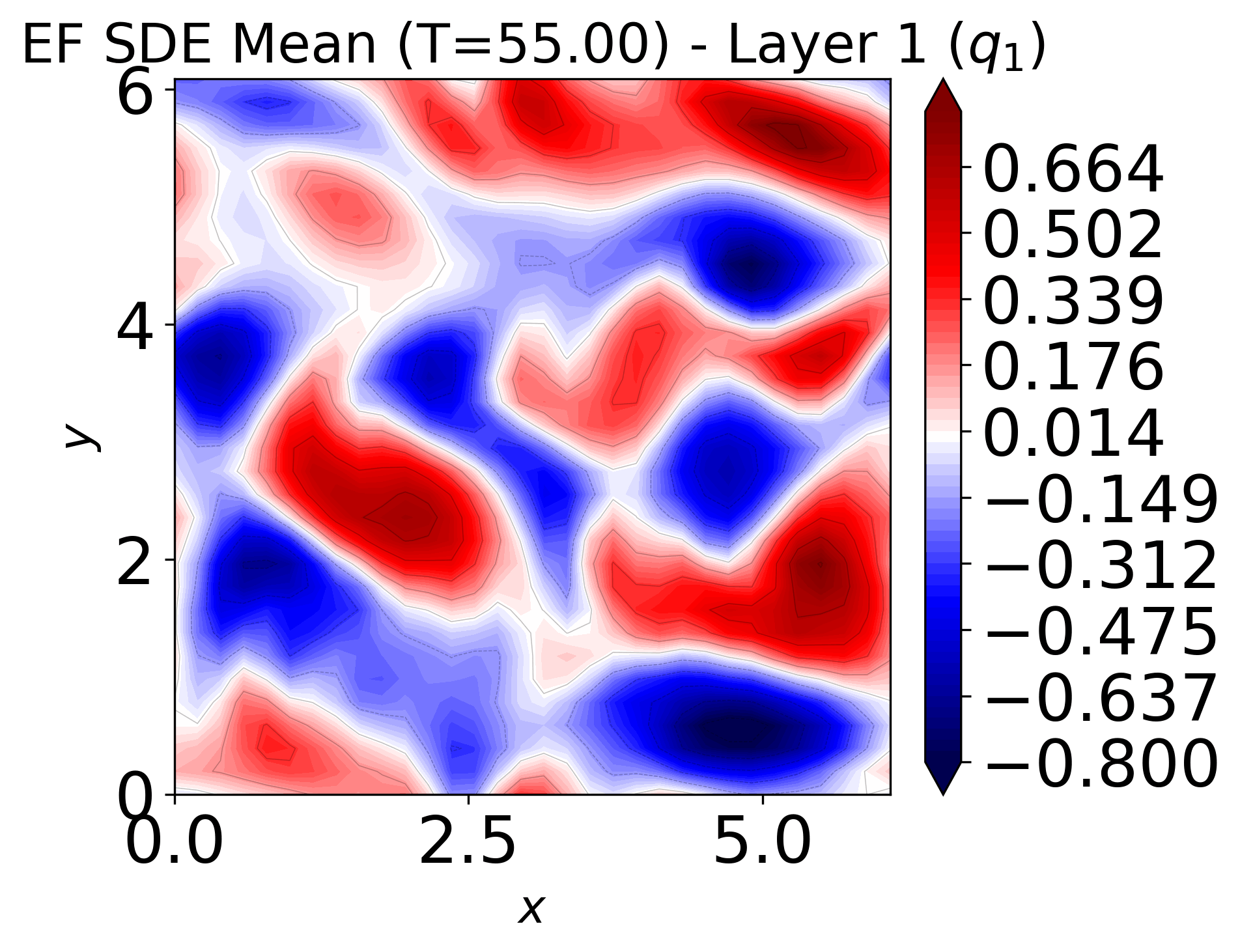} 
	\includegraphics[scale=0.052]{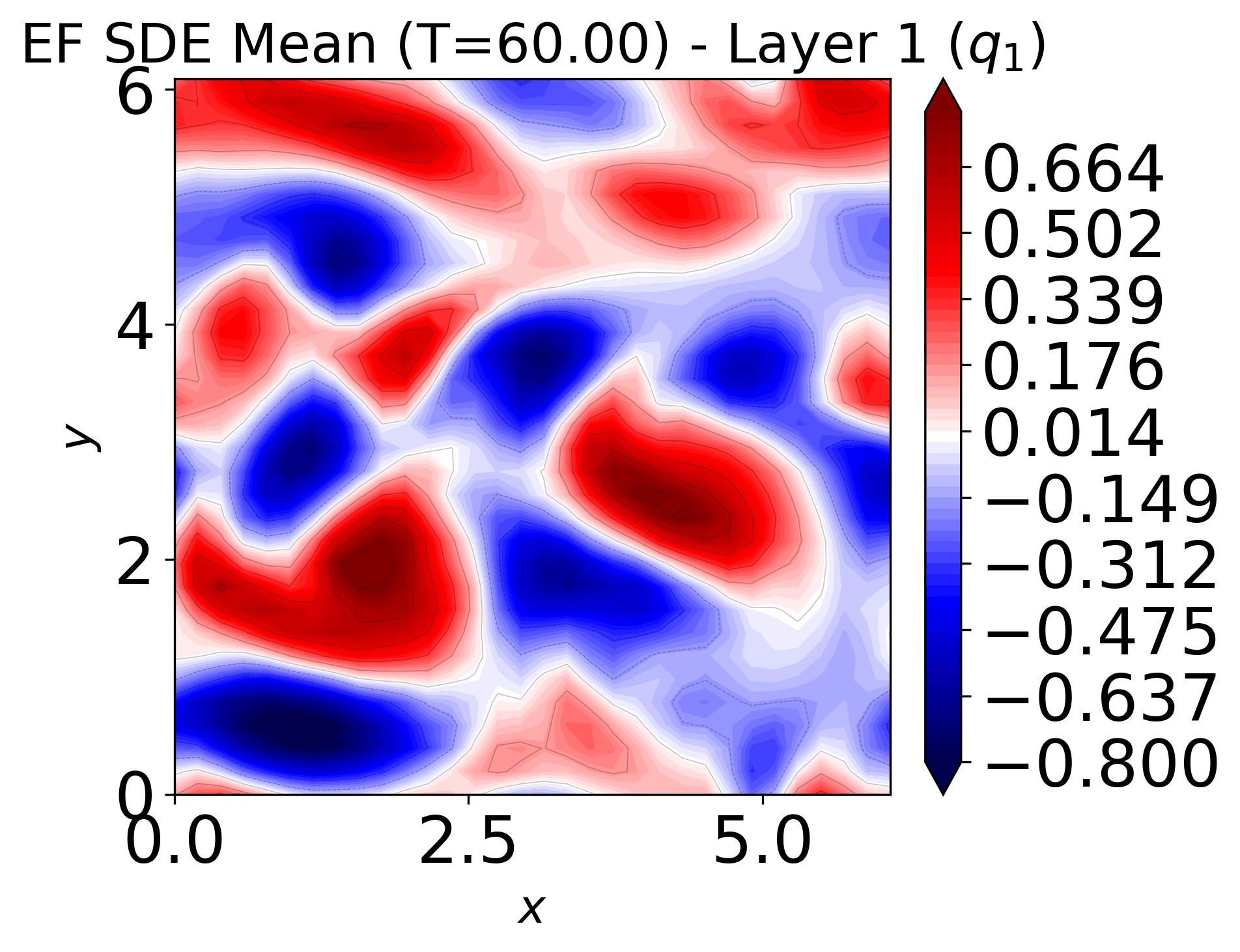}\\
	\includegraphics[scale=0.052]{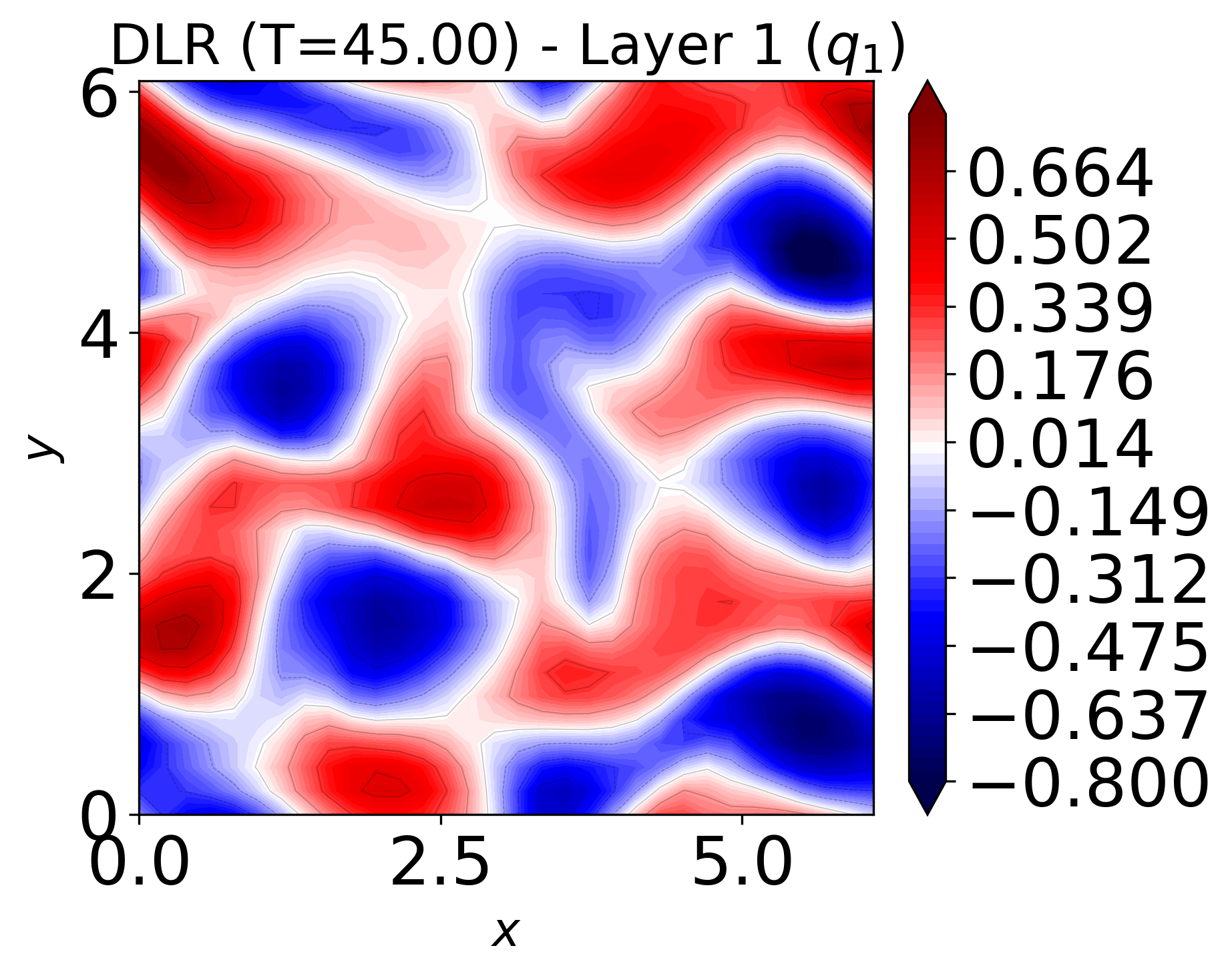} 
	\includegraphics[scale=0.052]{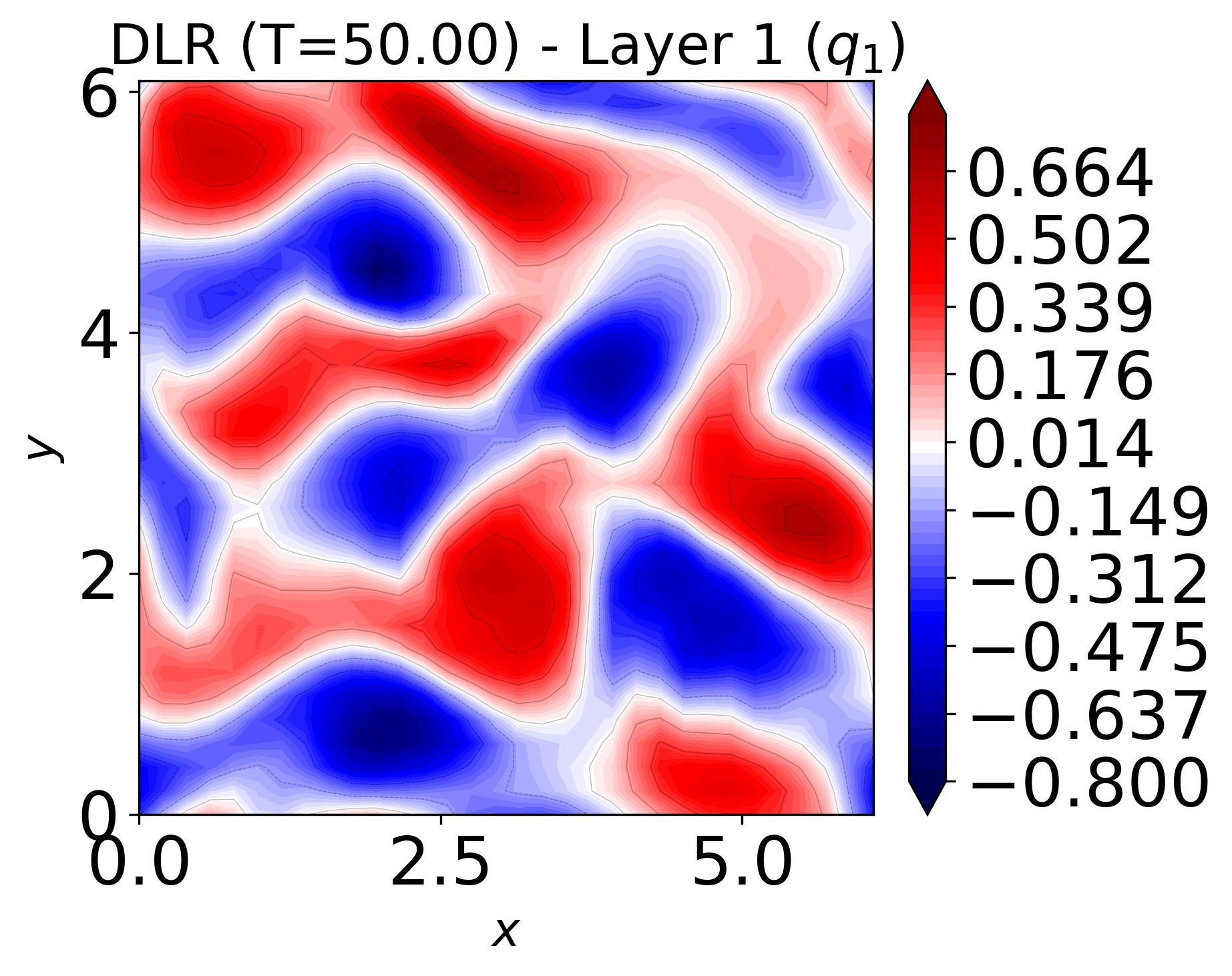} 
	\includegraphics[scale=0.052]{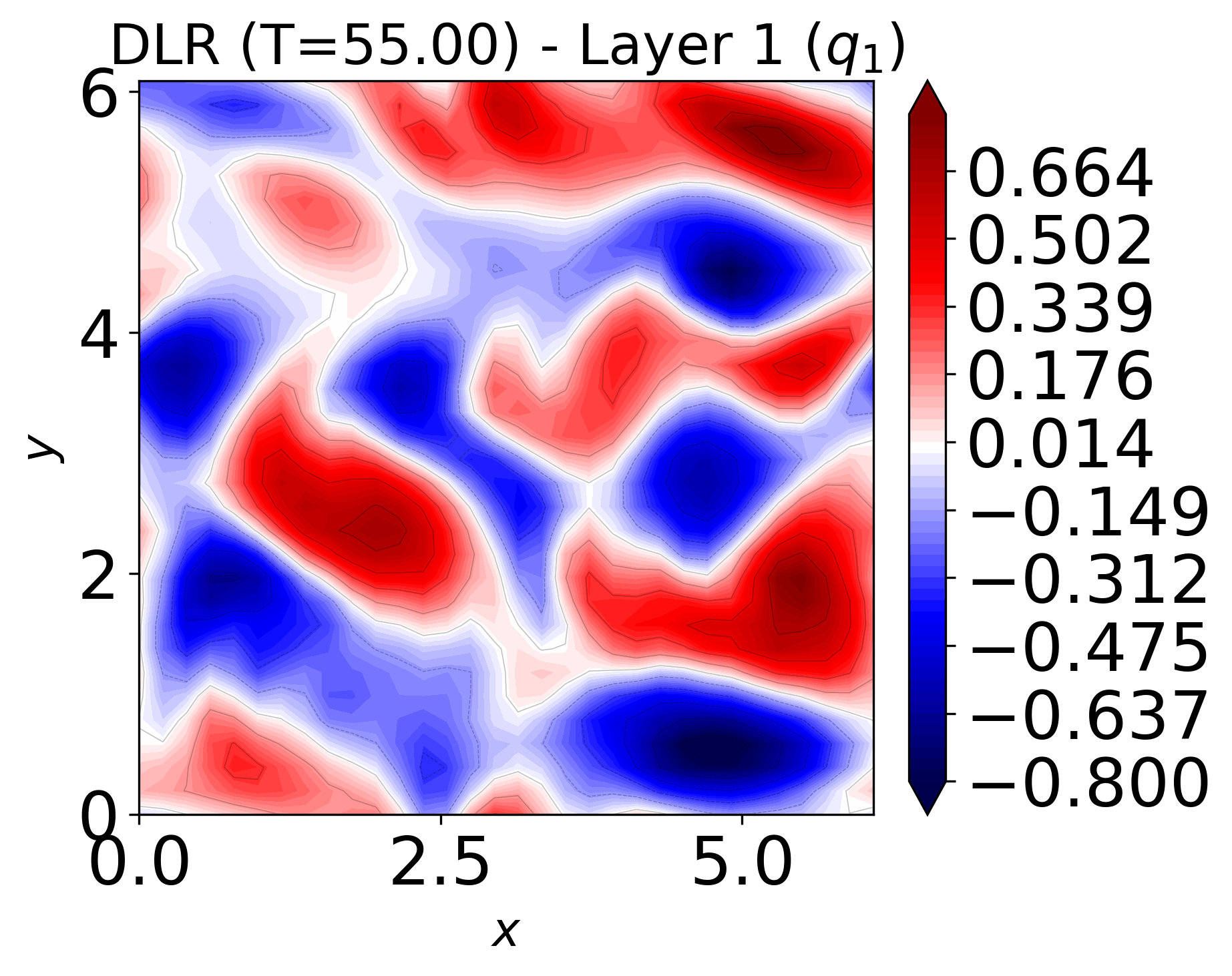} 
	\includegraphics[scale=0.052]{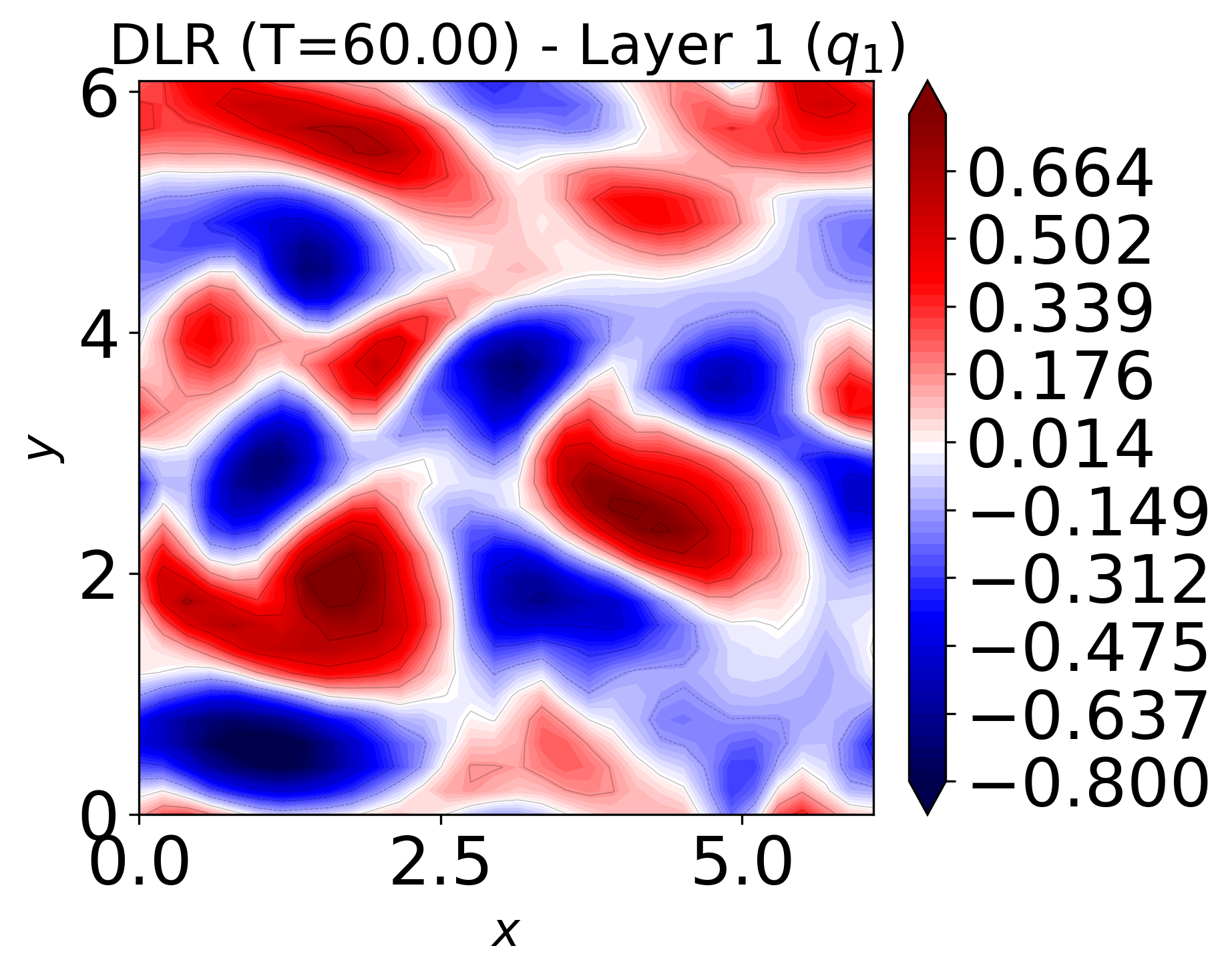}\\
	\caption{ Potential vorticity $q_1(x,y,t)$ in the first layer at different times  ($T=45,50,55,60$). From top to bottom: True state, mean of PF with $M=10^5$ particles, mean of EF SDE with $M=10^5$ particles, mean of DLRA JMCO filter for $M=1000$ and $k=32$. Problem \eqref{eq: 2 qg} with $d=2048$, $h=512$, with observation process \eqref{eq: 2 qg - obs proc}.}
	\label{fig: Oceanography - q1 layer}
\end{figure}

\begin{figure}[!h]
	\centering 
	\includegraphics[scale=0.052]{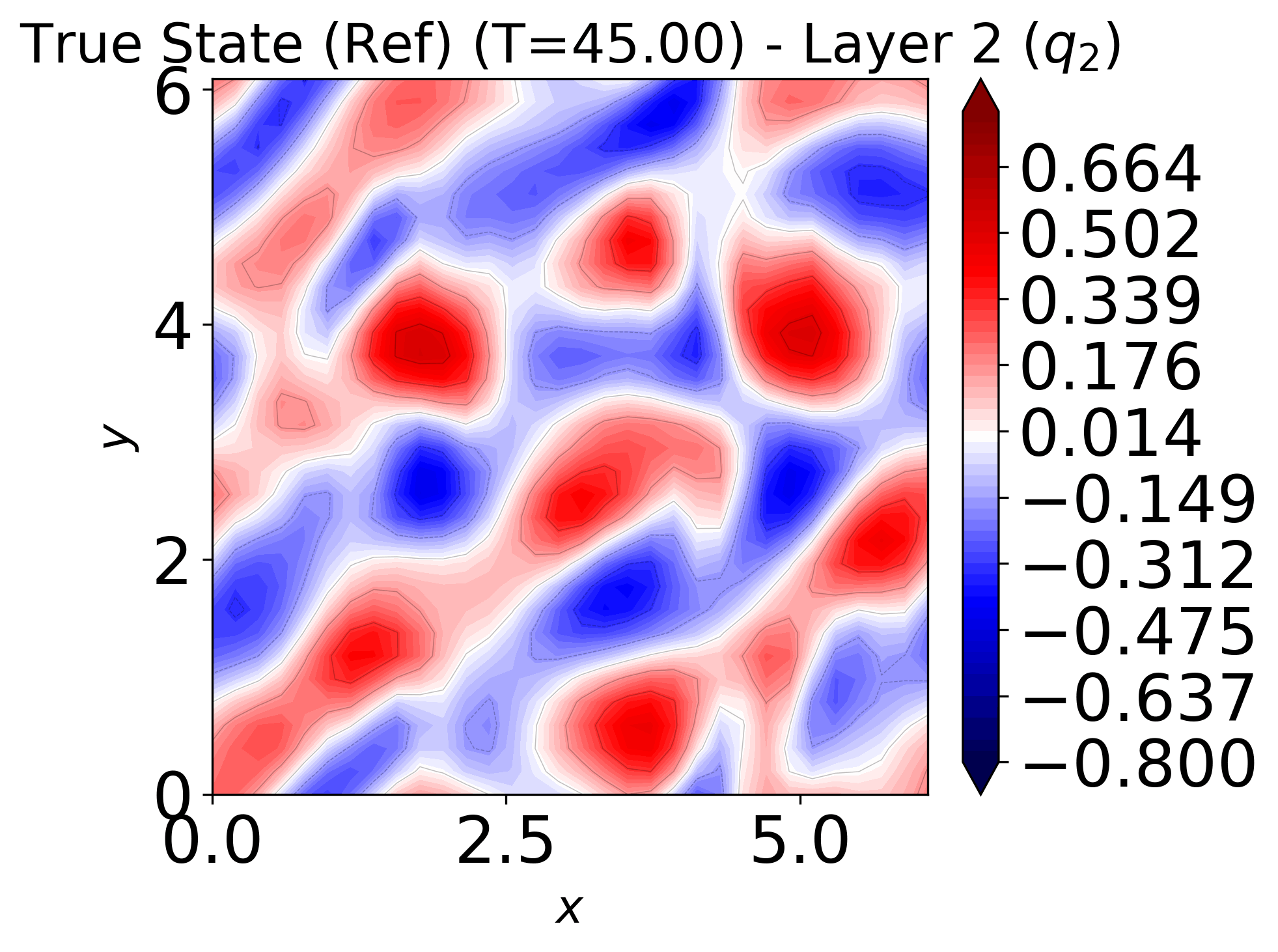} 
	\includegraphics[scale=0.052]{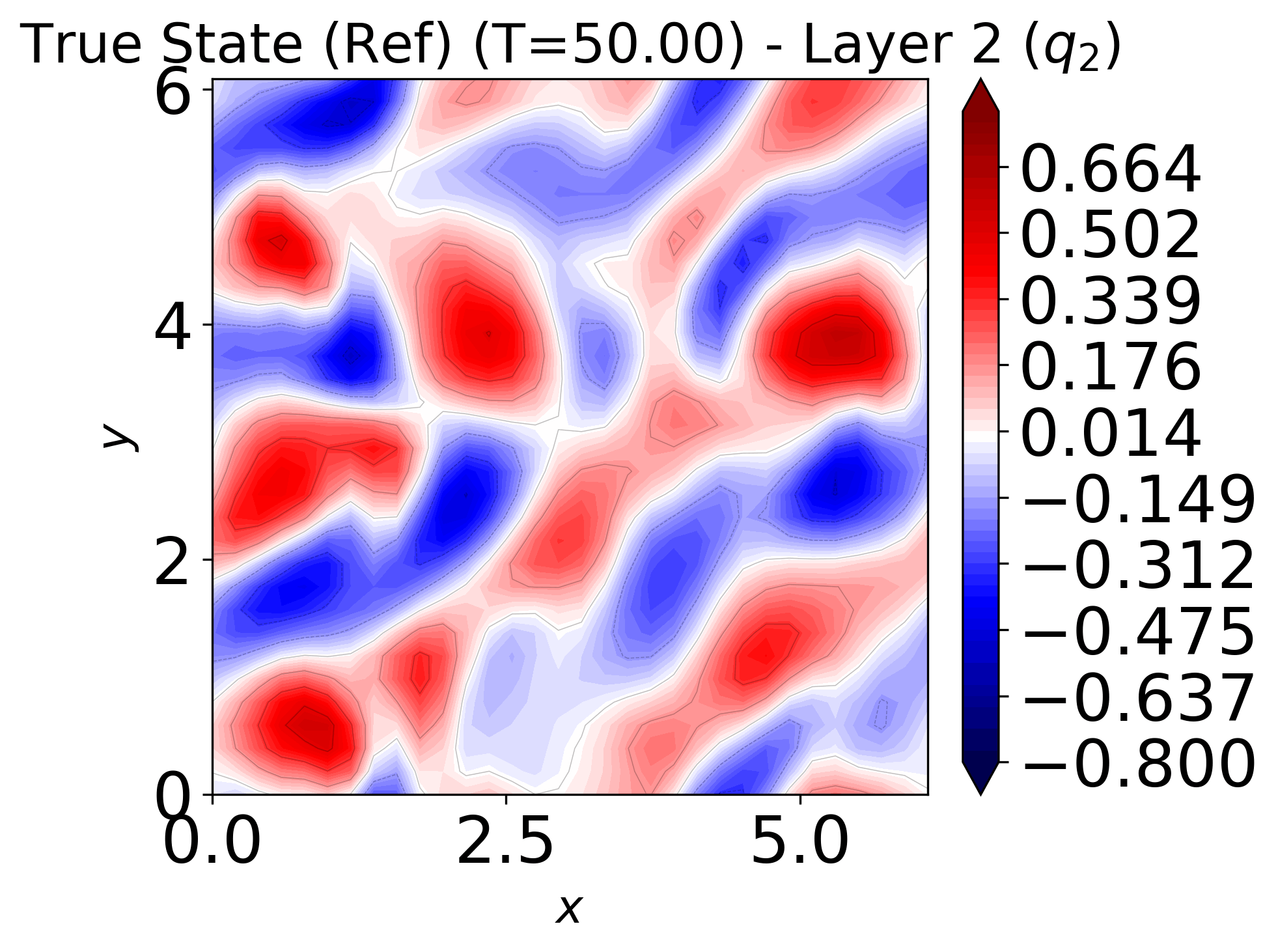} 
	\includegraphics[scale=0.052]{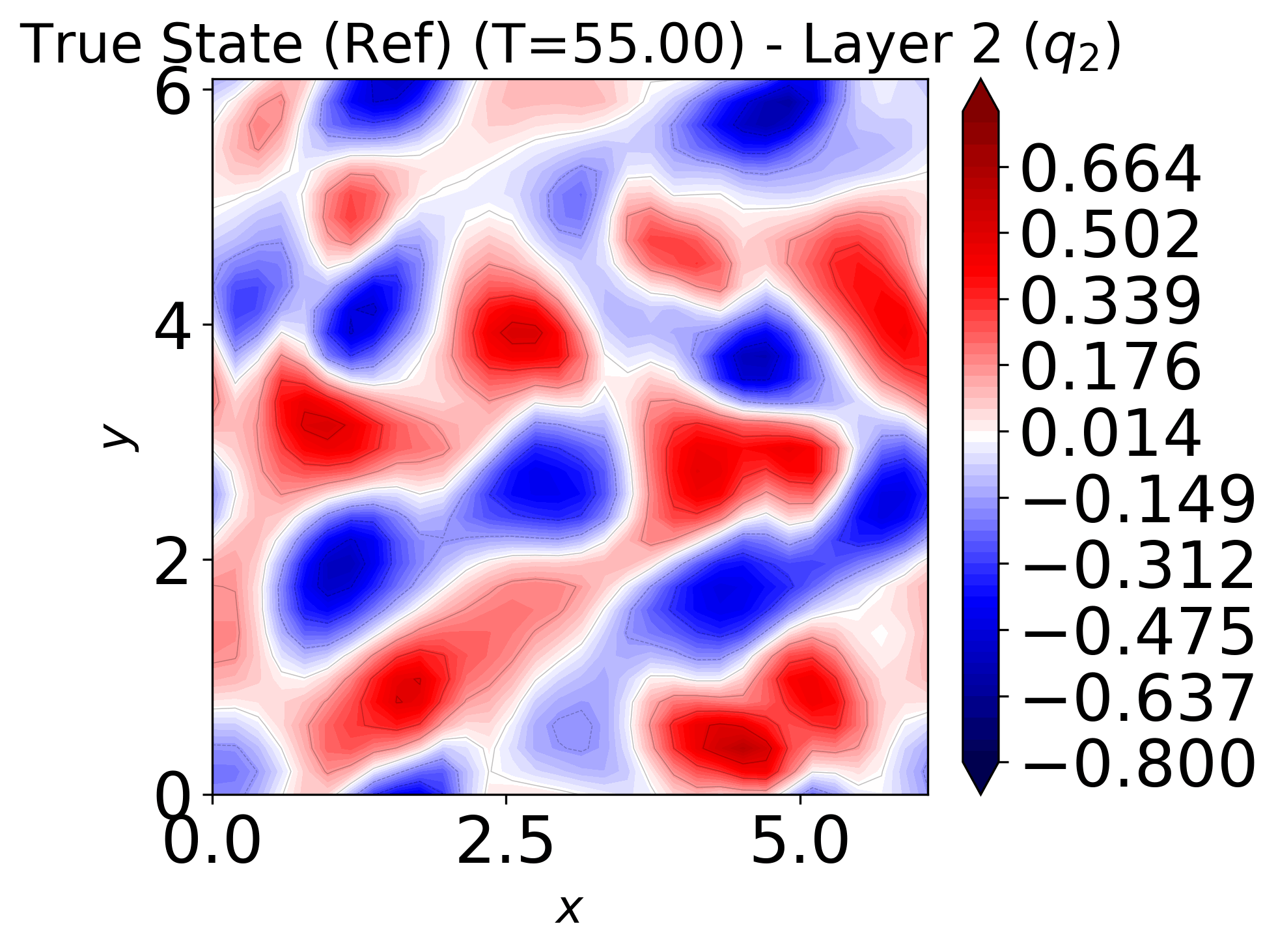} 
	\includegraphics[scale=0.052]{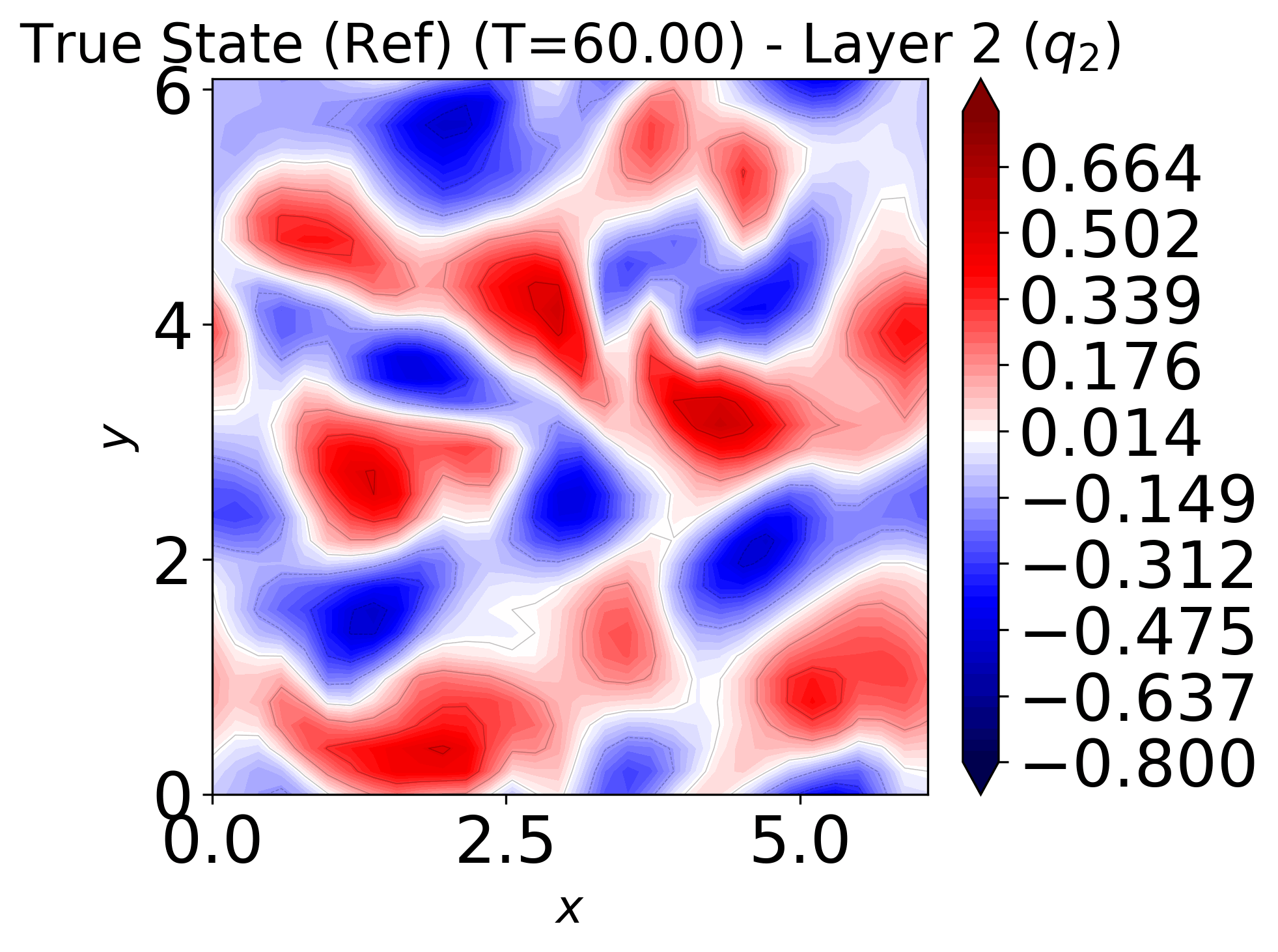}\\
	\includegraphics[scale=0.052]{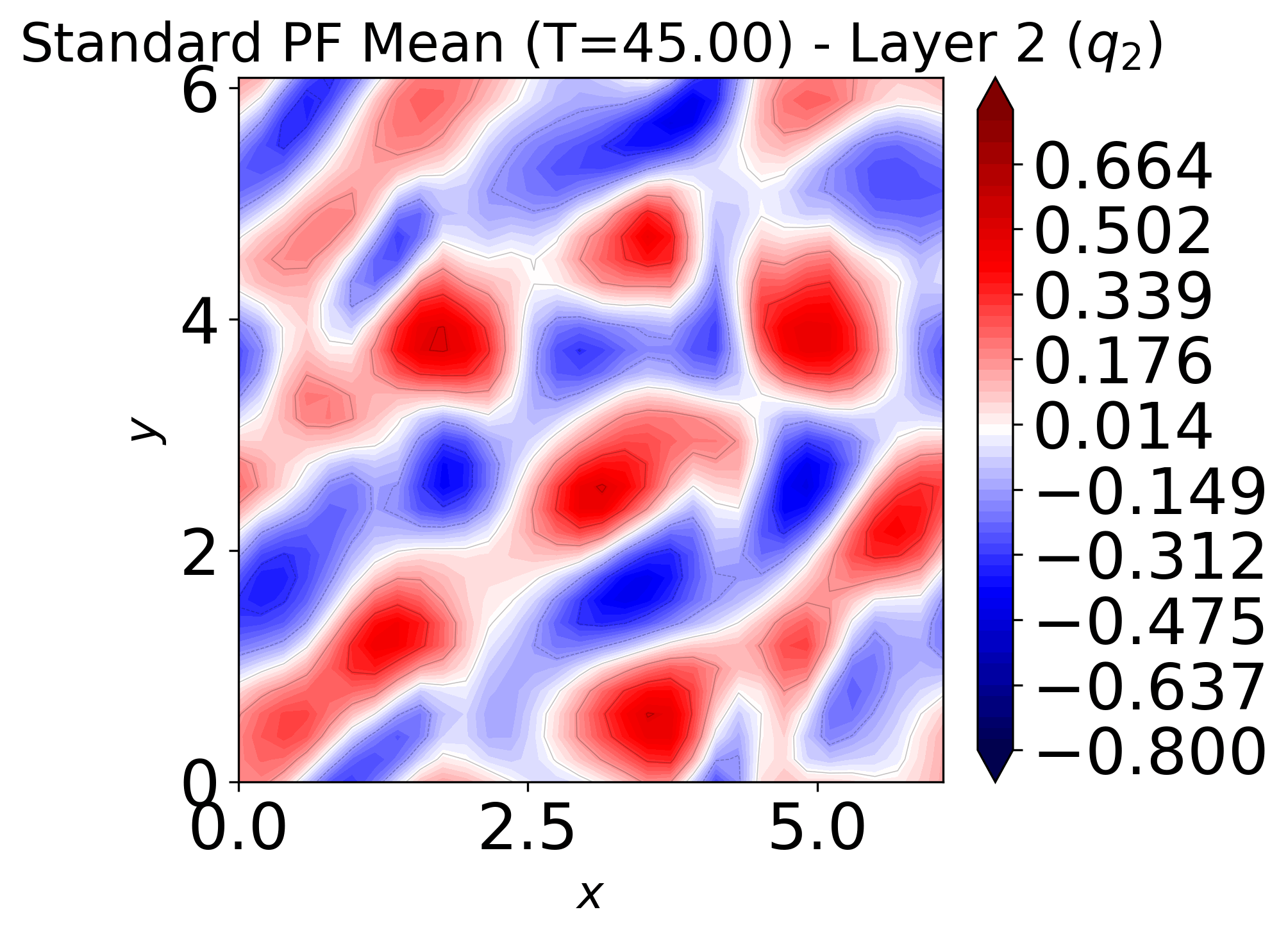} 
	\includegraphics[scale=0.052]{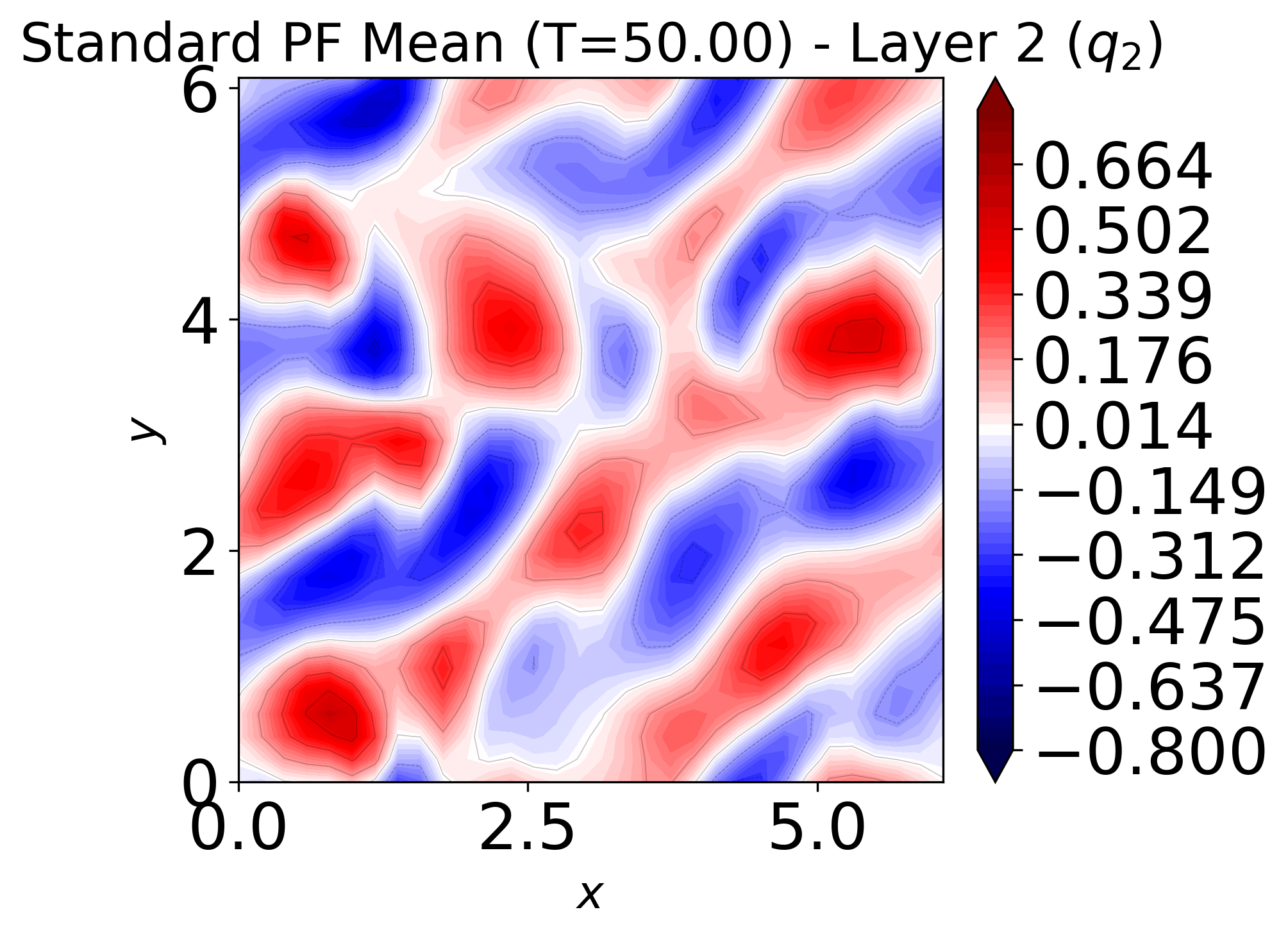} 
	\includegraphics[scale=0.052]{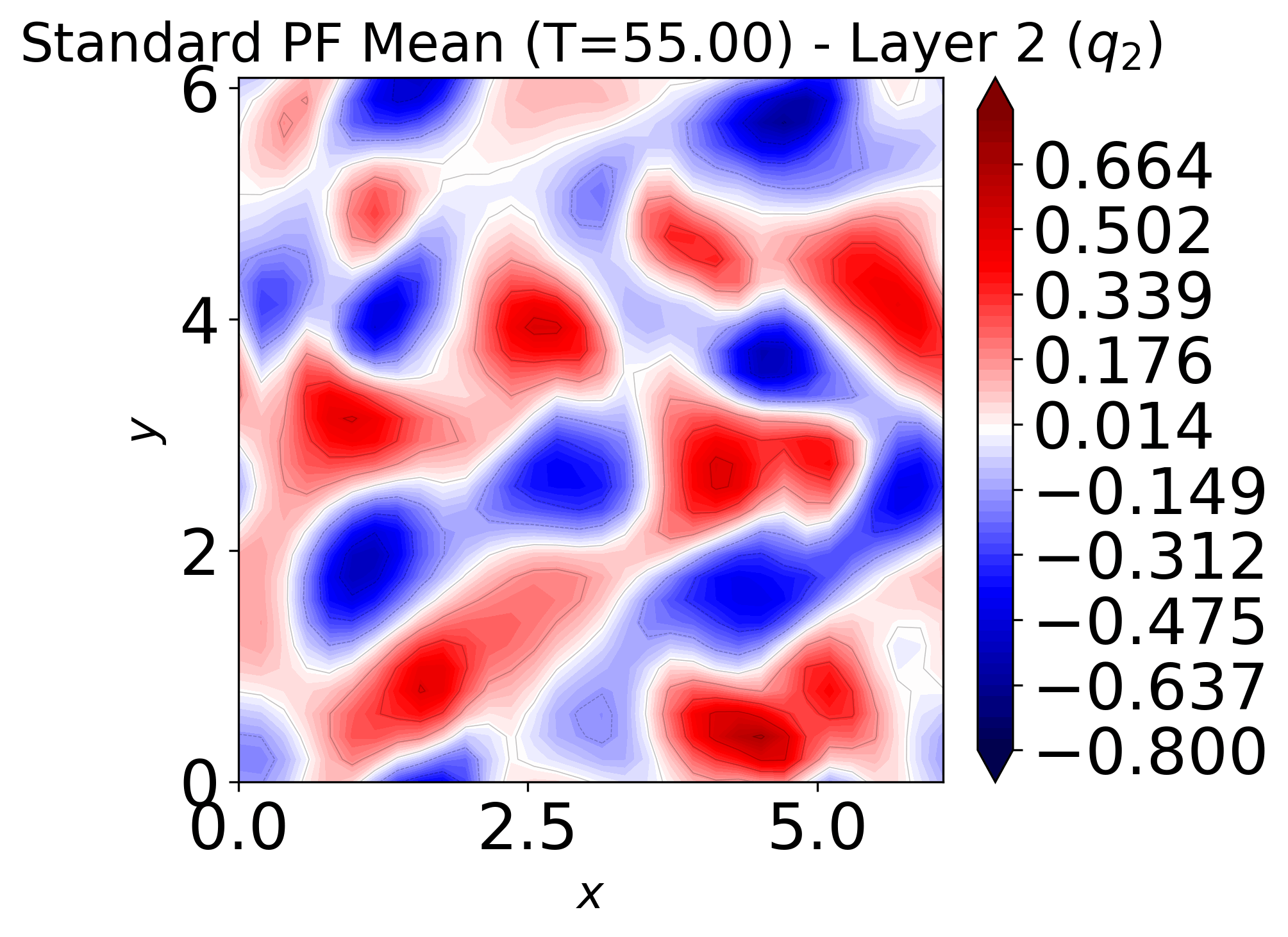} 
	\includegraphics[scale=0.051]{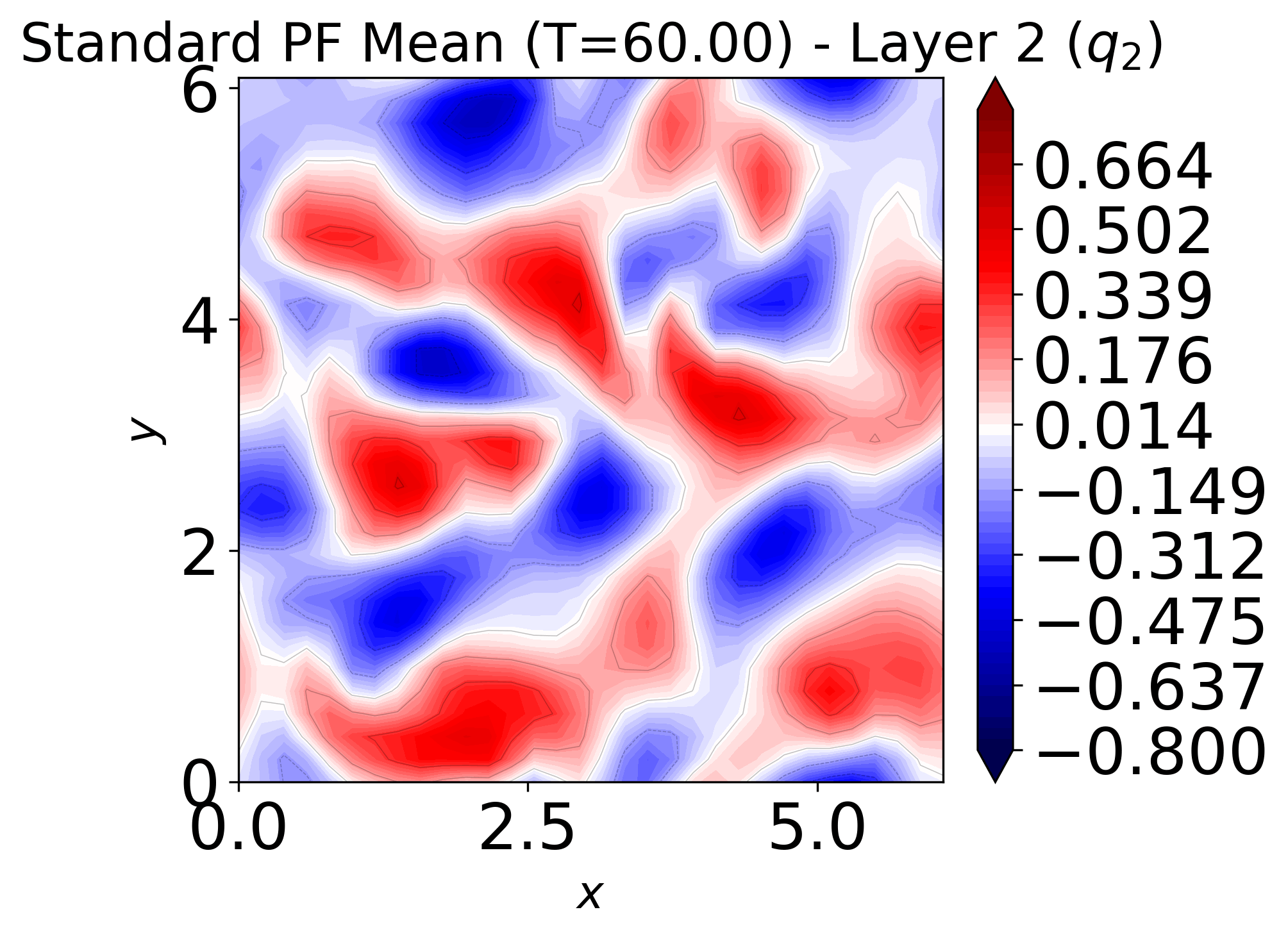}\\
	\includegraphics[scale=0.052]{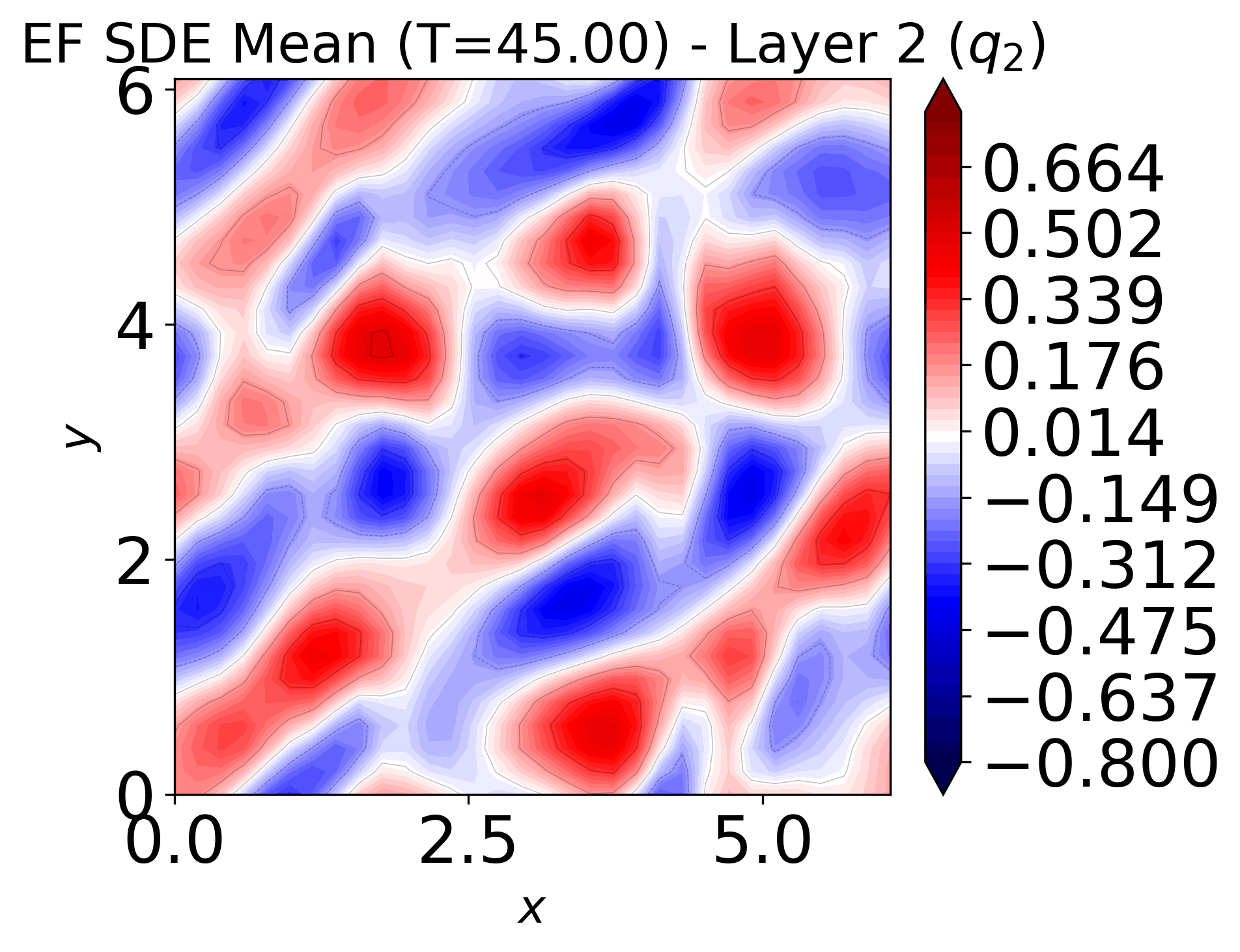}
	\includegraphics[scale=0.052]{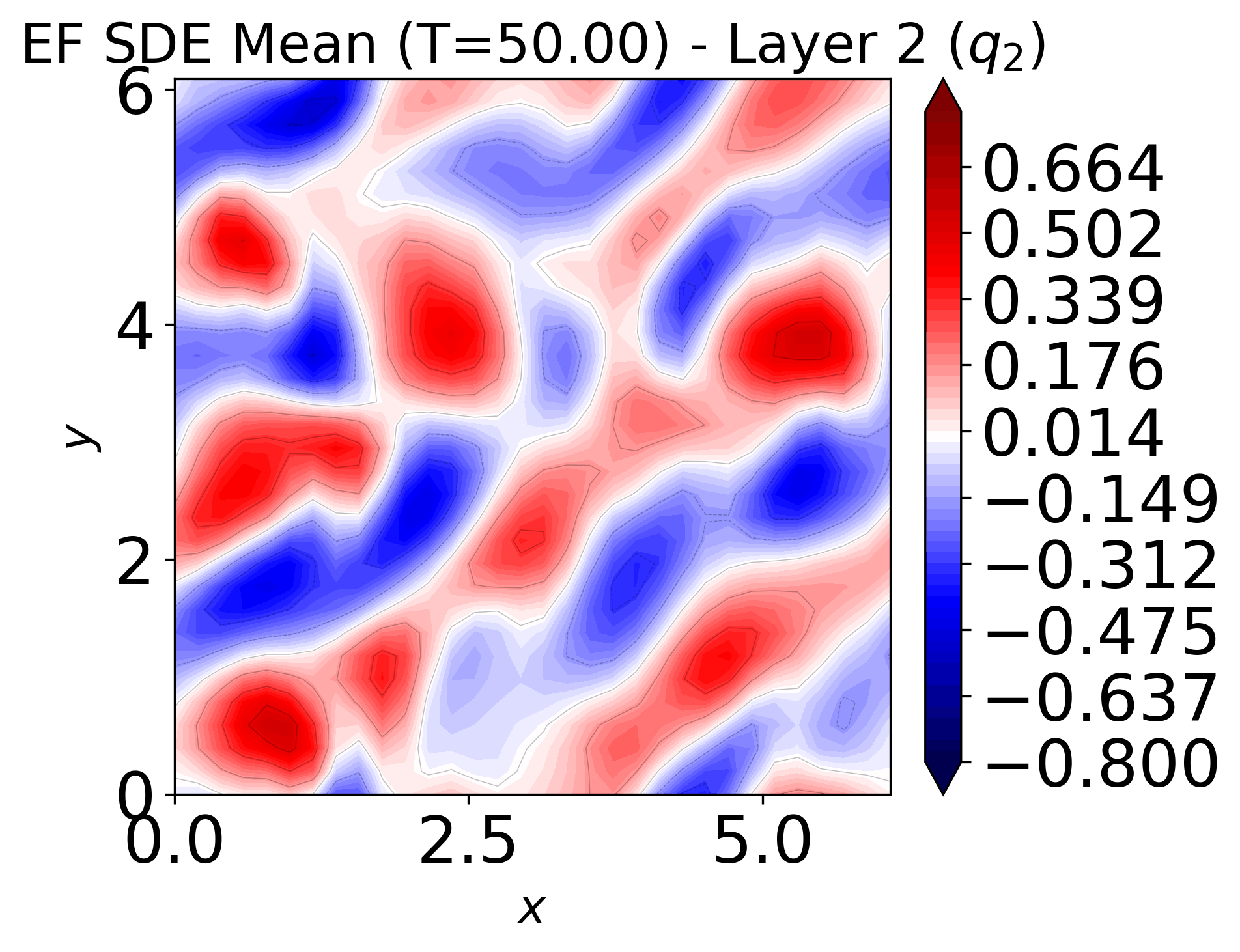} 
	\includegraphics[scale=0.052]{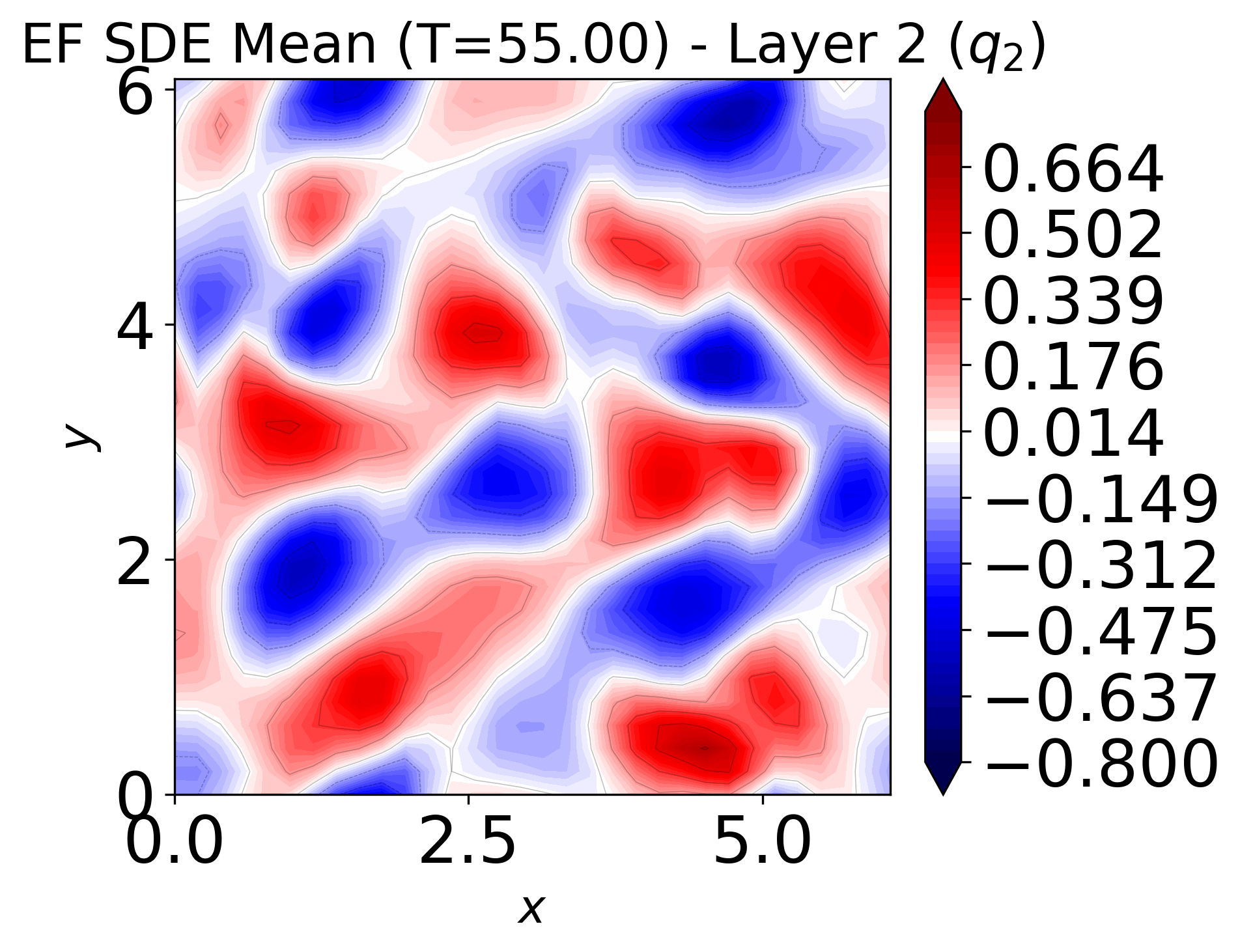} 
	\includegraphics[scale=0.052]{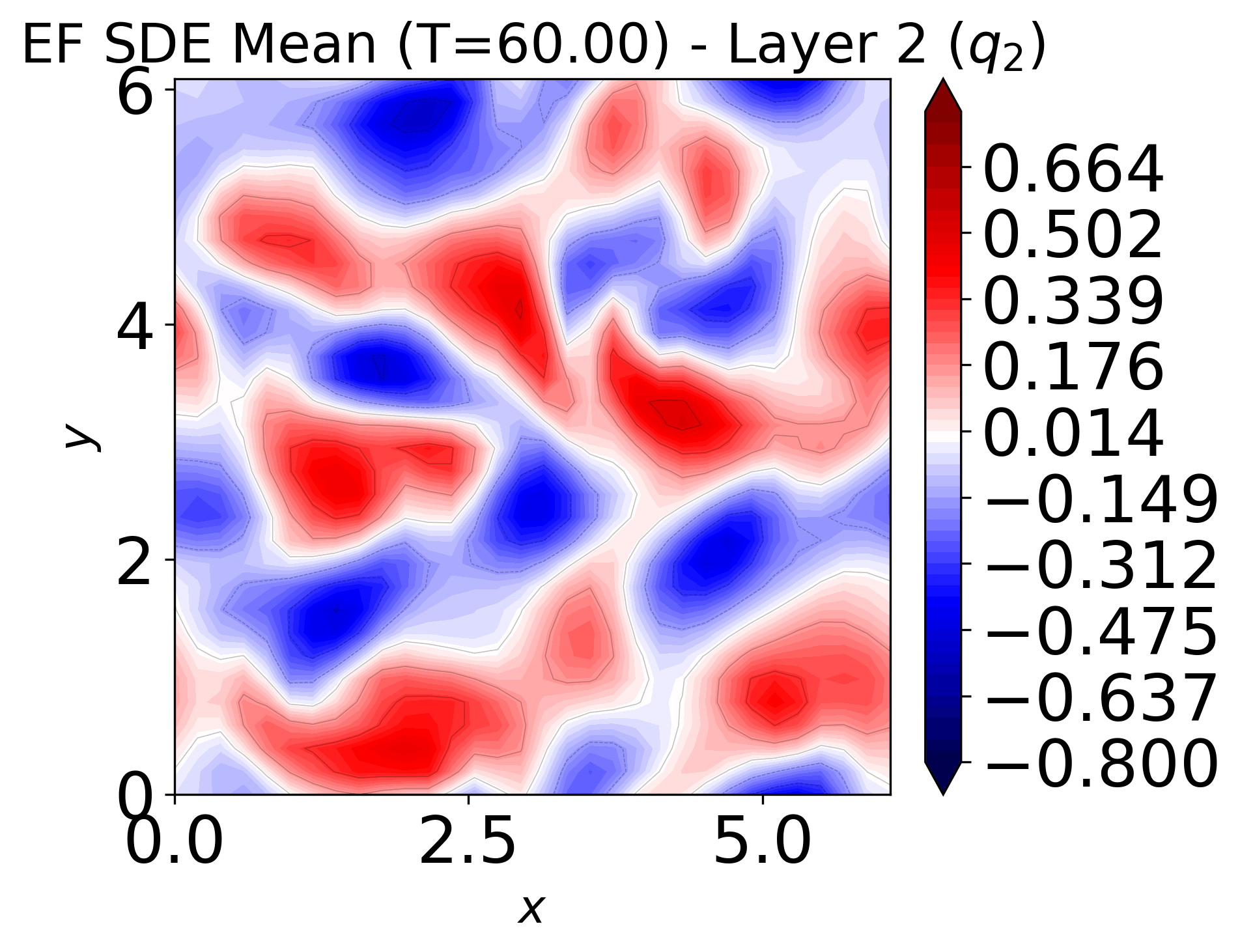}\\
	\includegraphics[scale=0.051]{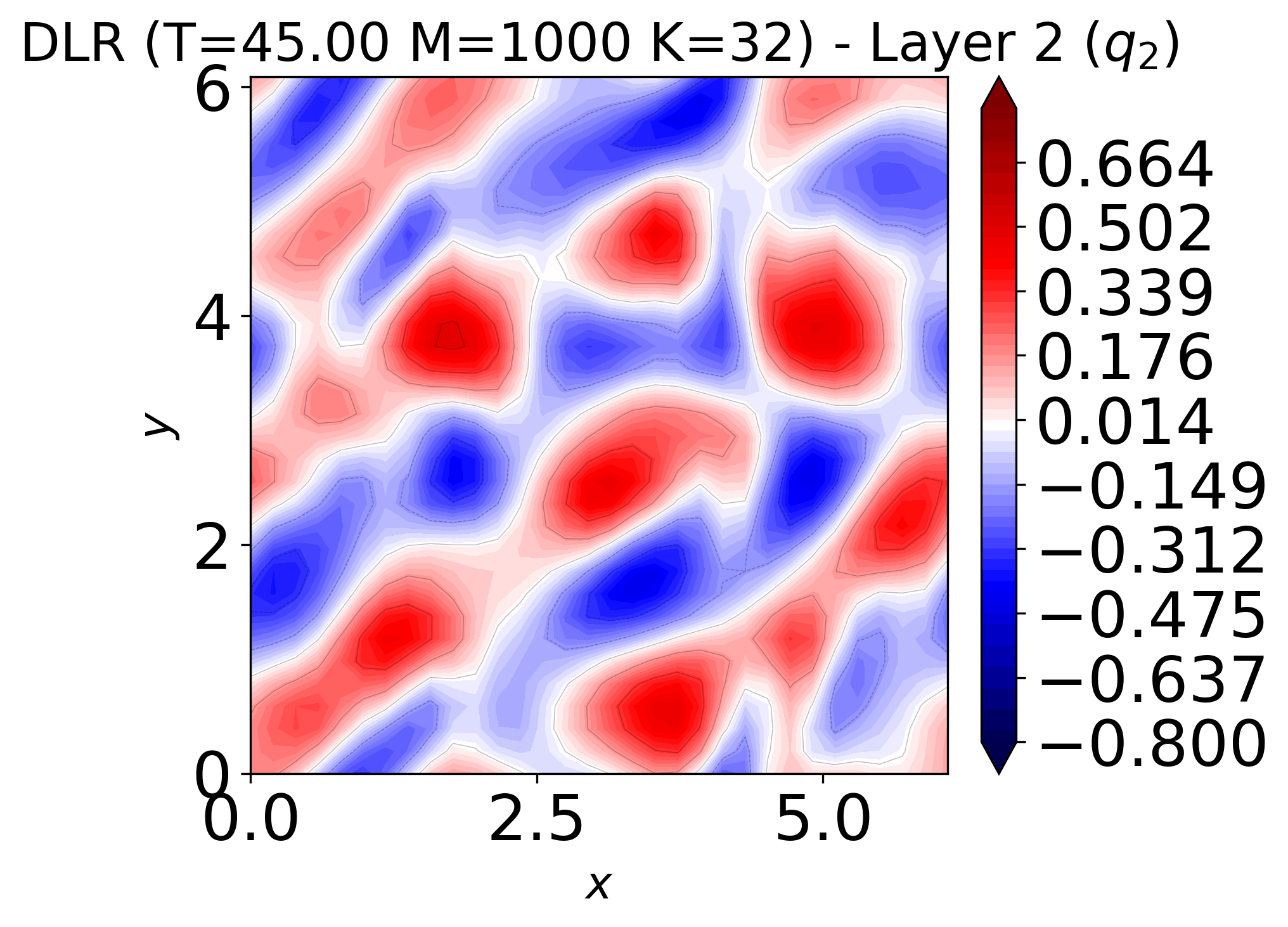} 
	\includegraphics[scale=0.051]{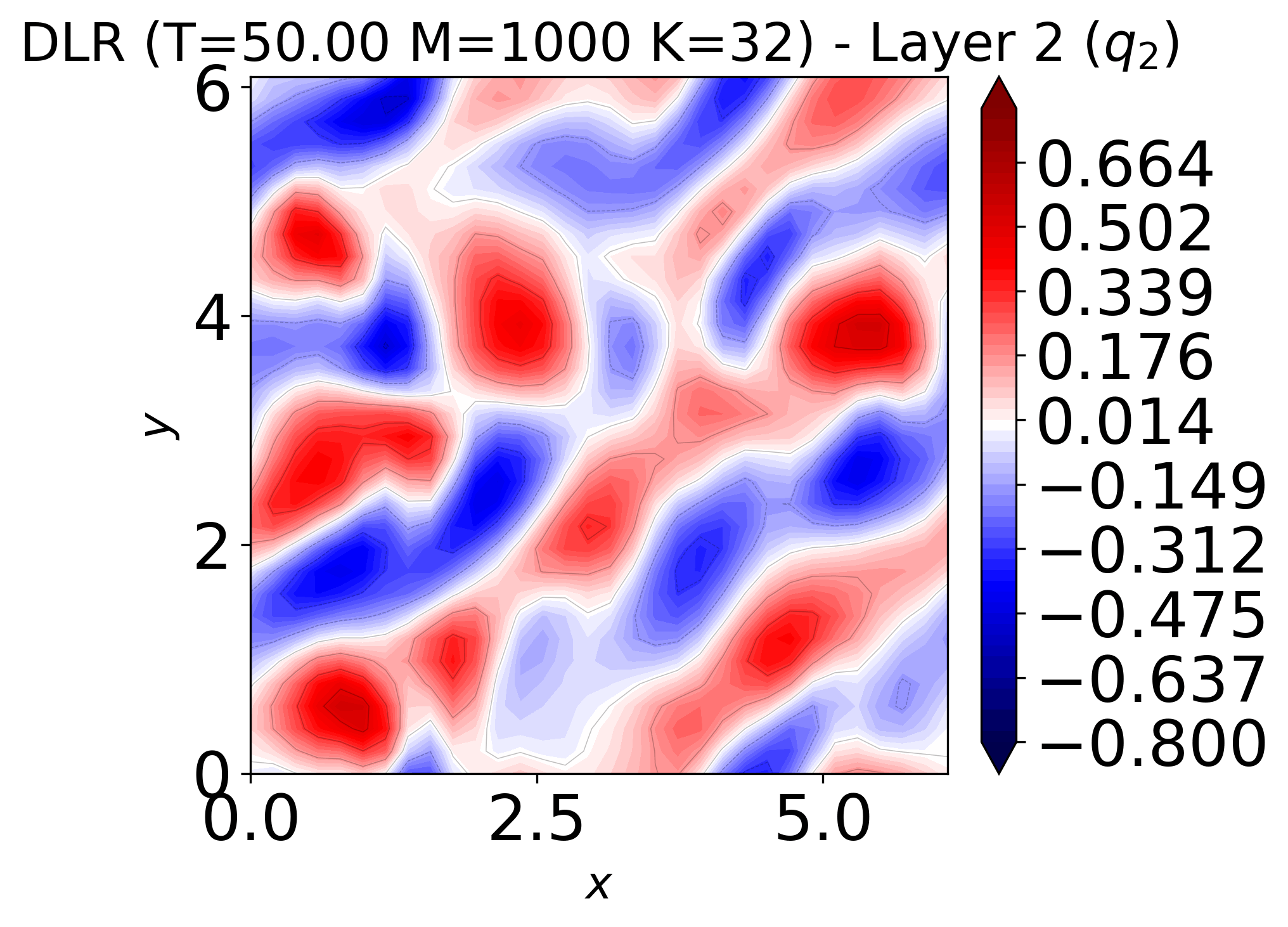} 
	\includegraphics[scale=0.051]{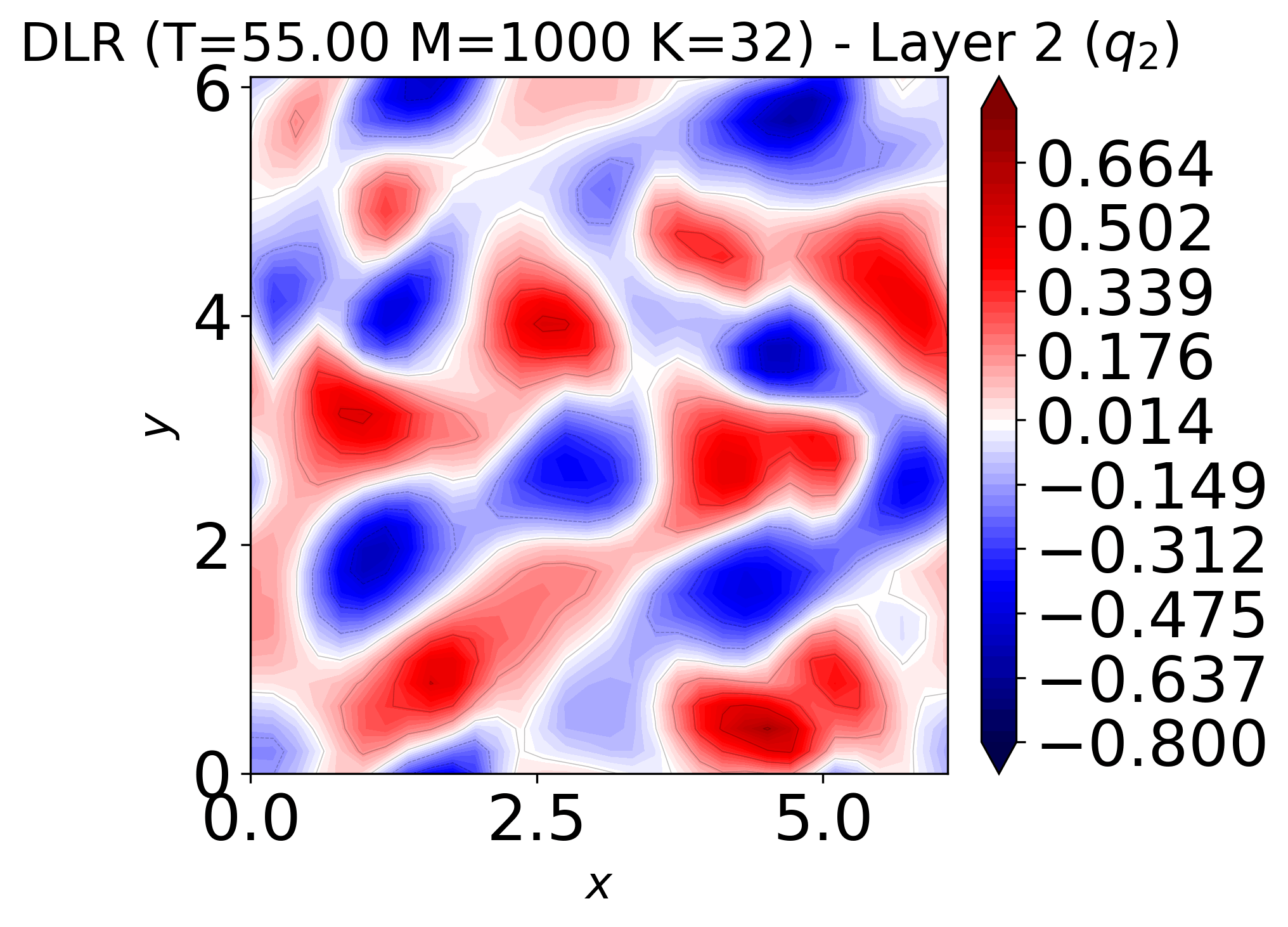} 
	\includegraphics[scale=0.051]{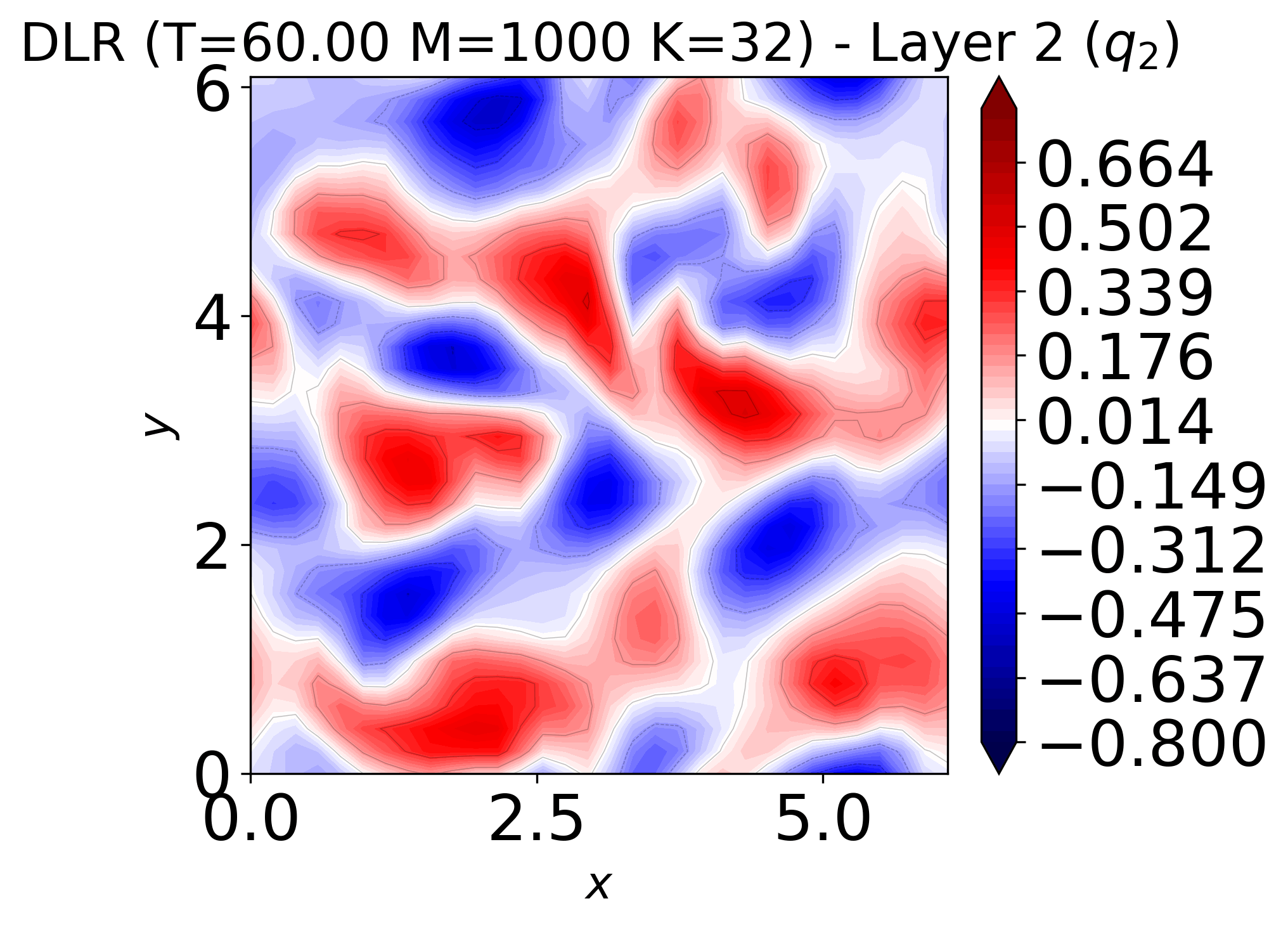}\\
	\caption{ Potential vorticity $q_2(x,y,t)$ in the second layer at different times  ($T=45,50,55,60$). From top to bottom: True state, mean of PF with $M=10^5$ particles, mean of EF SDE with $M=10^5$ particles, mean of DLRA JMCO filter for $M=1000$ and $k=32$. Problem \eqref{eq: 2 qg} with $d=2048$, $h=512$, with observation process \eqref{eq: 2 qg - obs proc}.}
	\label{fig: Oceanography - q2 layer}
\end{figure}

\begin{figure}[!h]
	\centering
	\includegraphics[scale=0.045]{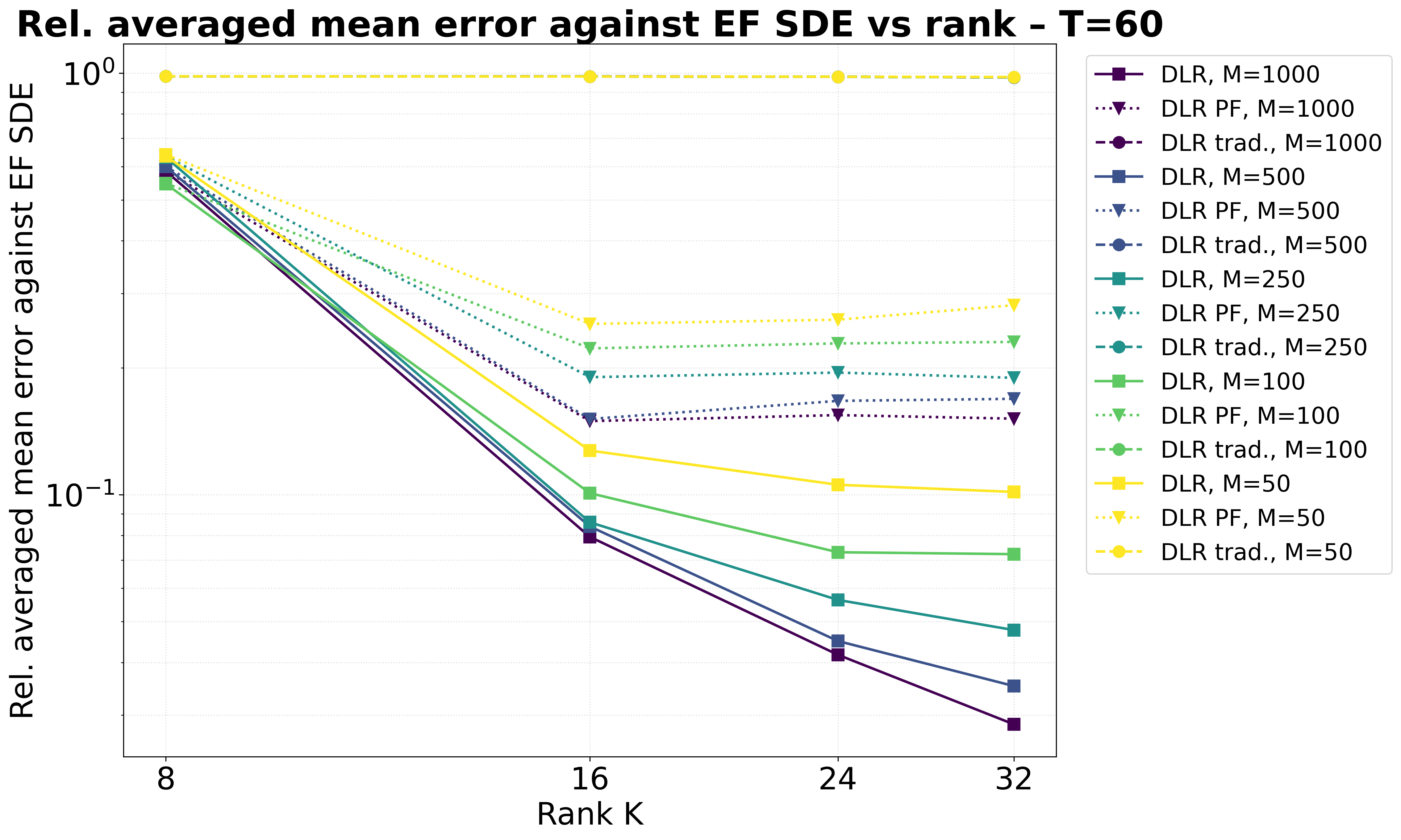}
	\includegraphics[scale=0.045]{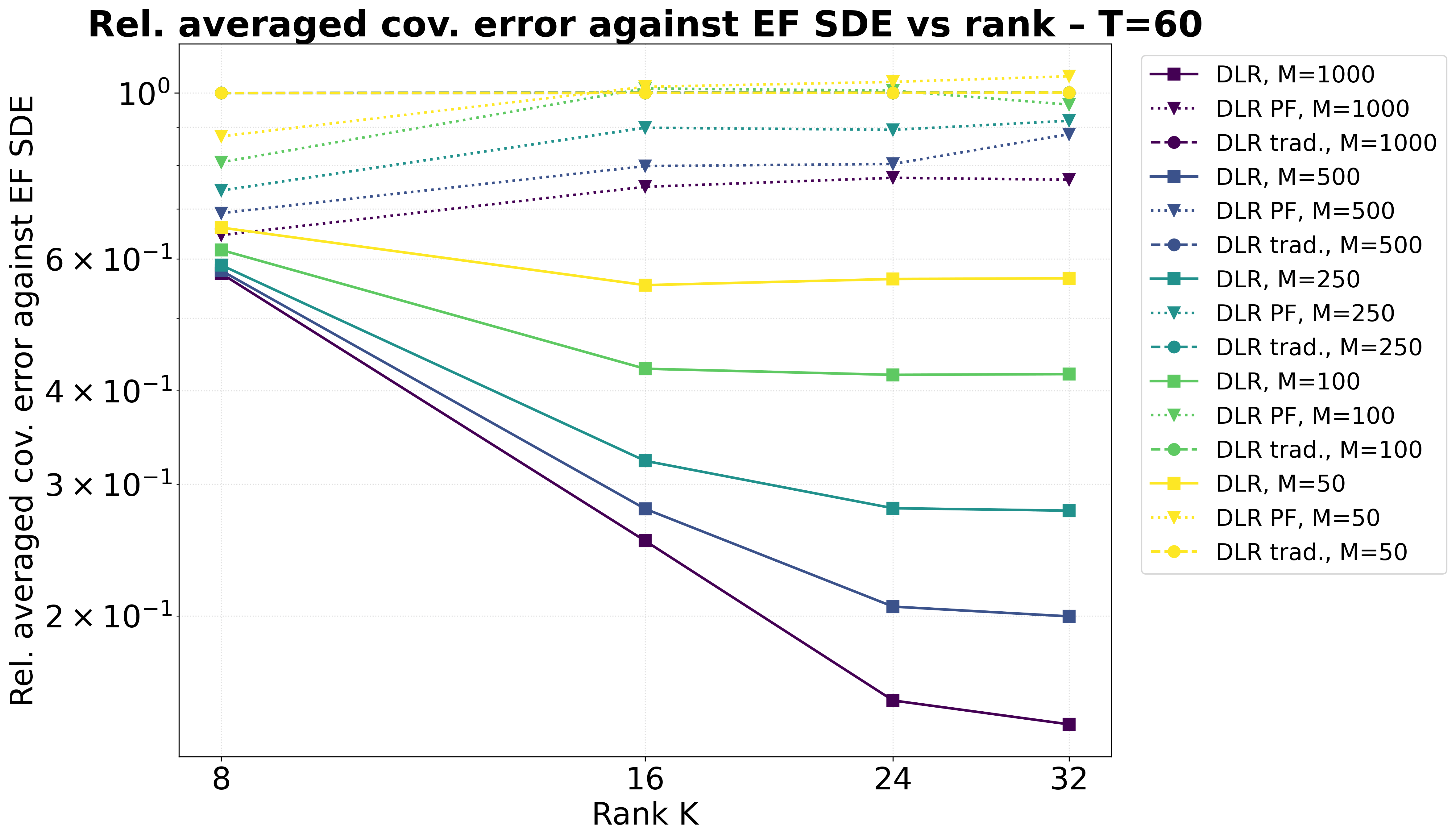} 
	\caption{Rel.\  averaged-over-time \emph{errors against PF} as a function of the rank $K$: (left) error on the mean, (right) error on the covariance. Problem \eqref{eq: 2 qg} with continuous-time observations \eqref{eq: 2 qg - obs proc}, transient time $T_0=2$, reference computed with $M=10^5$, $d=2048$, $h=512$, $\Delta t_{branch}= 4\Delta t$.}
	\label{fig: Oceanography EF SDE rank}
\end{figure}

\begin{figure}[!h]
	\centering
	\includegraphics[scale=0.045]{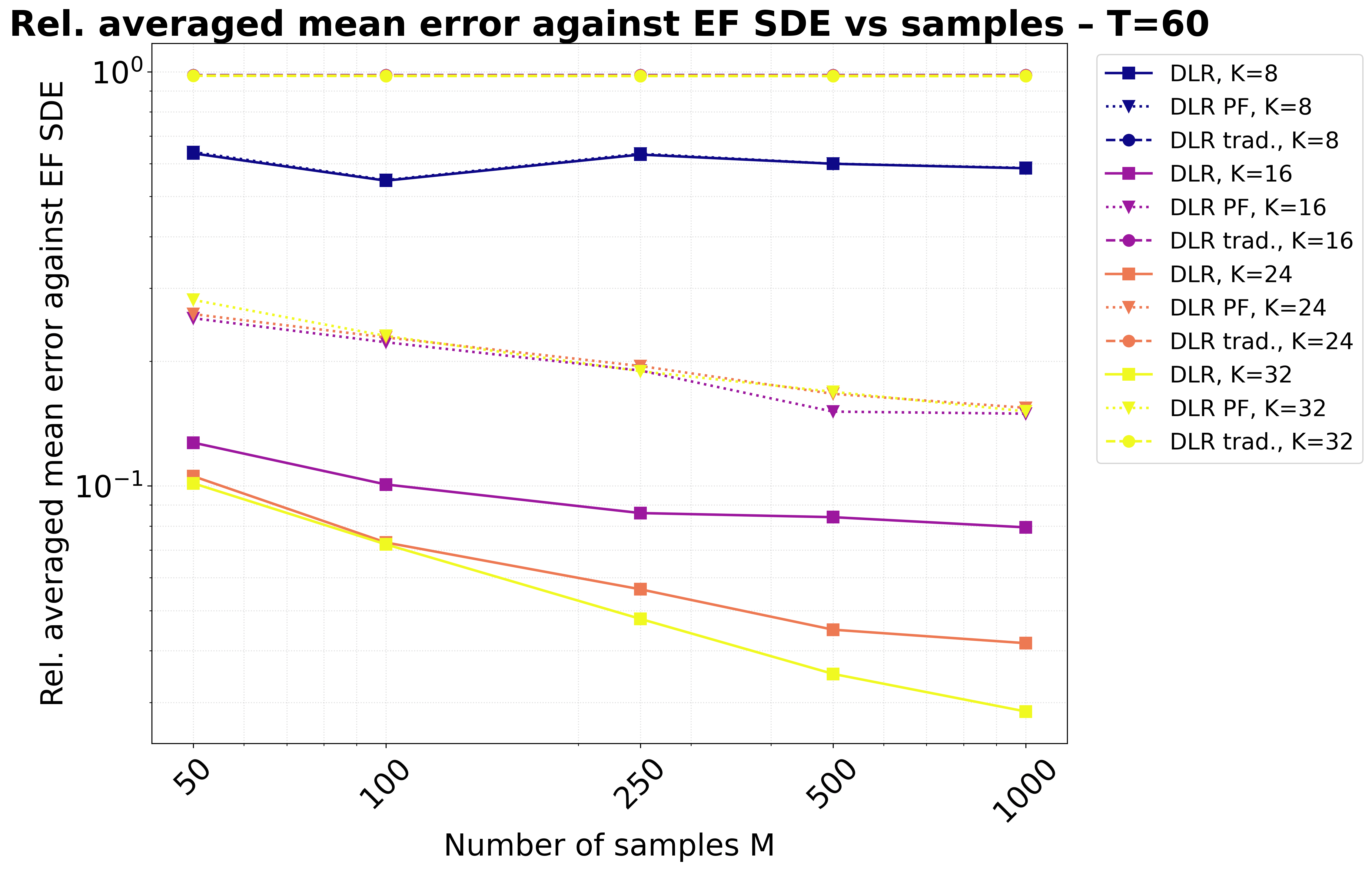}
	\includegraphics[scale=0.045]{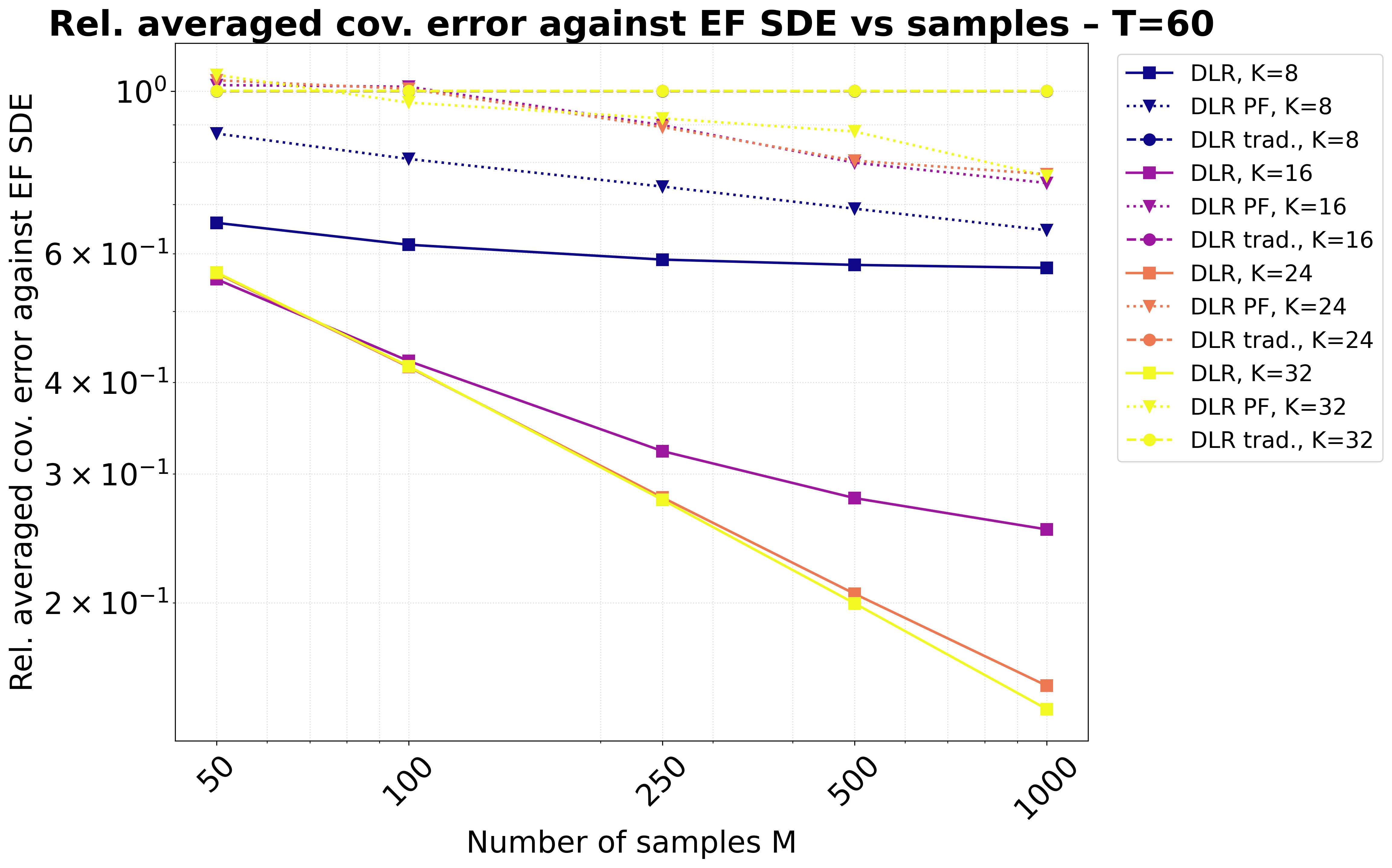} 
	\caption{Rel.\  averaged-over-time \emph{errors against EF SDE} as a function of the sample size $M$: (left) error on the mean, (right) error on the covariance. Problem \eqref{eq: 2 qg} with continuous-time observations \eqref{eq: 2 qg - obs proc}, transient time $T_0=2$, reference computed with $M=10^5$, $d=2048$, $h=512$, $\Delta t_{branch}= 4\Delta t$.}
	\label{fig: Oceanography EF SDE samples}
\end{figure}

We simulate the reference particle filter and the ensemble Kalman filter approximating the full-order system with $M=10^5$ particles and we plot the averaged-over-time errors after a warm-up time $T_0=2$. We also consider a comparison particle filter computed with only $M=60$ particles. The interval of time of the branching procedure for all the particle-filter type algorithms is $\Delta t_{branch}= 4\Delta t$.
Figures \ref{fig: Oceanography EF SDE rank} and  \ref{fig: Oceanography EF SDE samples} show the true potential vorticity in the two layers at four different times, as well as the filters mean for the reference filters PF and EF SDE with $M=10^5$, and DLRA-JMCO with $k=32$ and $M=1000$. There we see that the DLR JMCO filter reaches great performance in approximating the covariance of EF SDE, showing a decreasing error both with respect to the rank and the number of samples. This behavior corroborates the fact that if the true solution lives in a small dimensional manifold, we just need a not so large number of samples to well approximate it.  Moreover, the effect of the correction $P_{U_n}^{\perp} Q U_n C_{Y_n}^{-1}$ in the dynamics of DLR and DLR PF is evident when compared to DLR Trad, whose error stagnates. 

\begin{figure}[!h]
	\centering
	\includegraphics[scale=0.045]{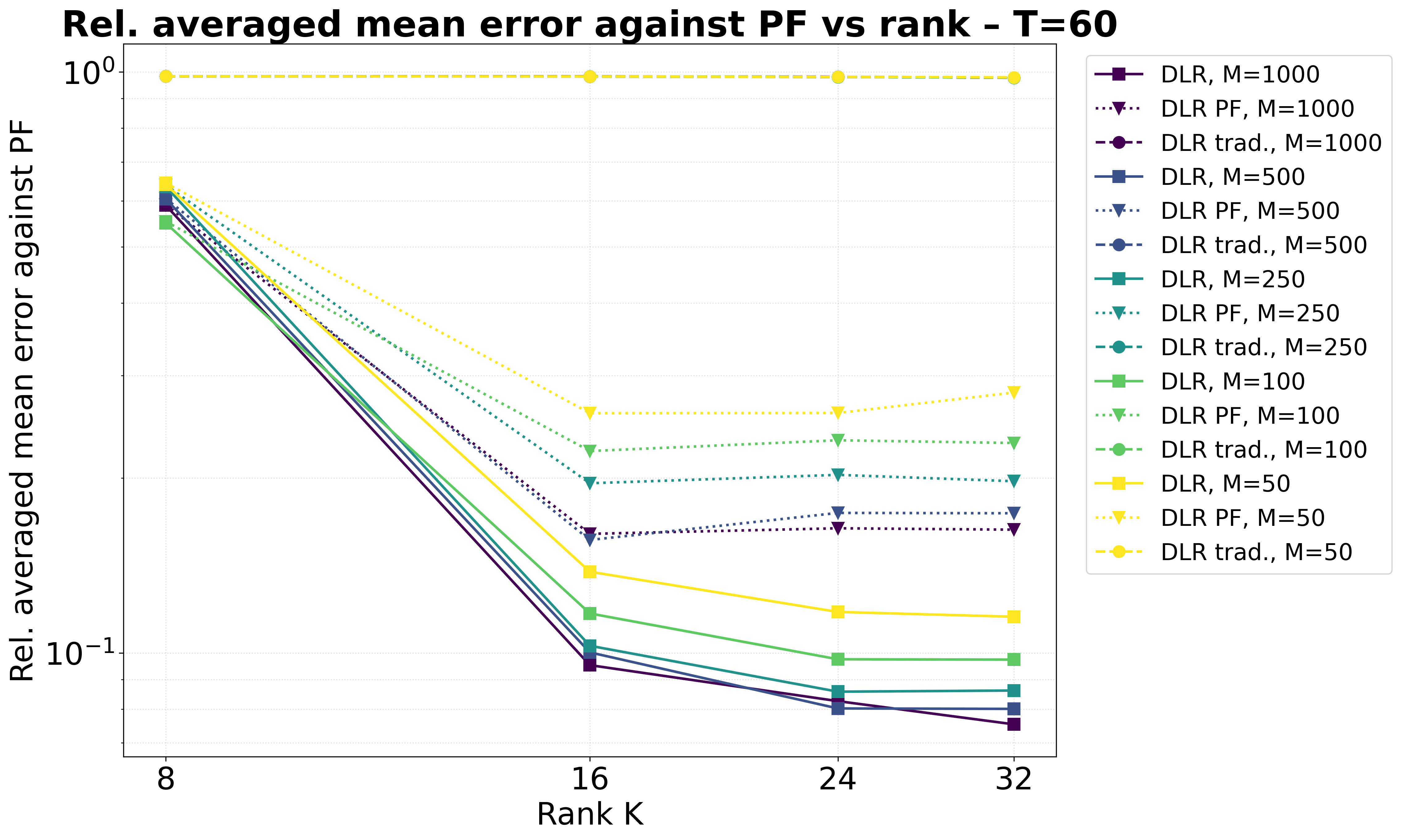}
	\includegraphics[scale=0.045]{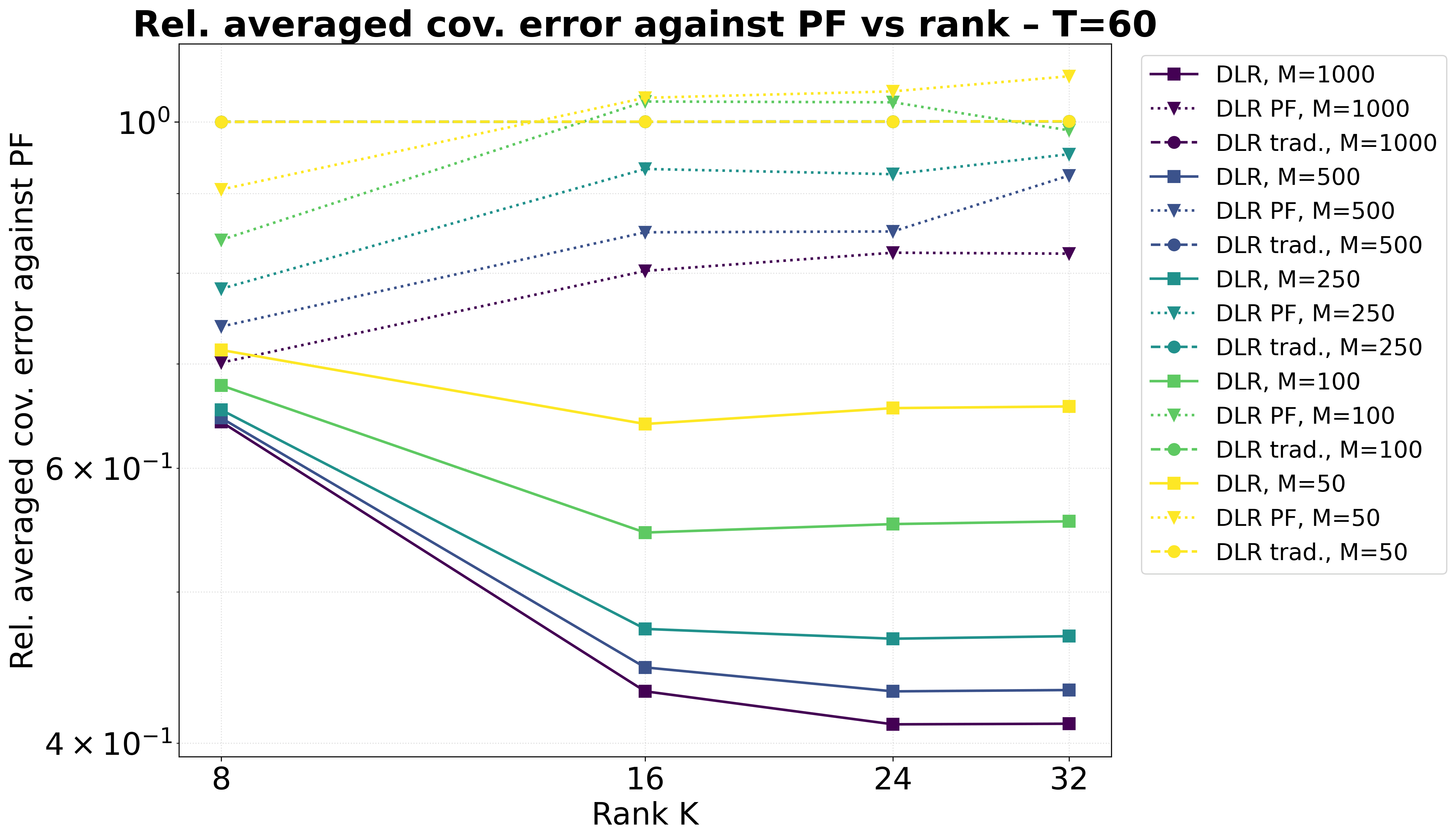} 
	\caption{Rel.\  averaged-over-time \emph{errors against PF} as a function of the rank $K$: (left) error on the mean, (right) error on the covariance. Problem \eqref{eq: 2 qg} with continuous-time observations \eqref{eq: 2 qg - obs proc}, transient time $T_0=2$, reference computed with $M=10^5$, $d=2048$, $h=512$, $\Delta t_{branch}= 4\Delta t$.}
	\label{fig: Oceanography PF rank}
\end{figure}

\begin{figure}[!h]
	\centering
	\includegraphics[scale=0.045]{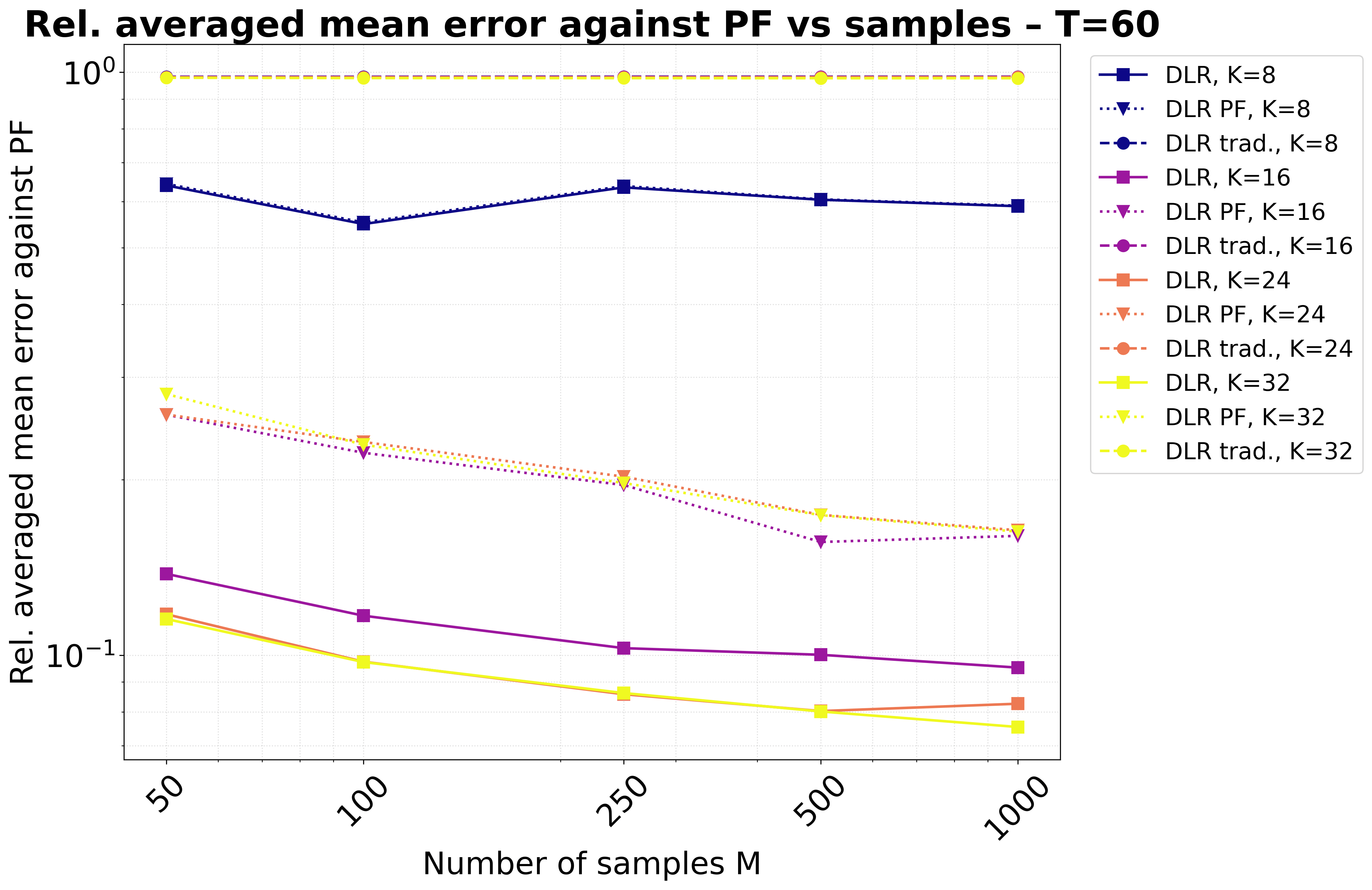}
	\includegraphics[scale=0.045]{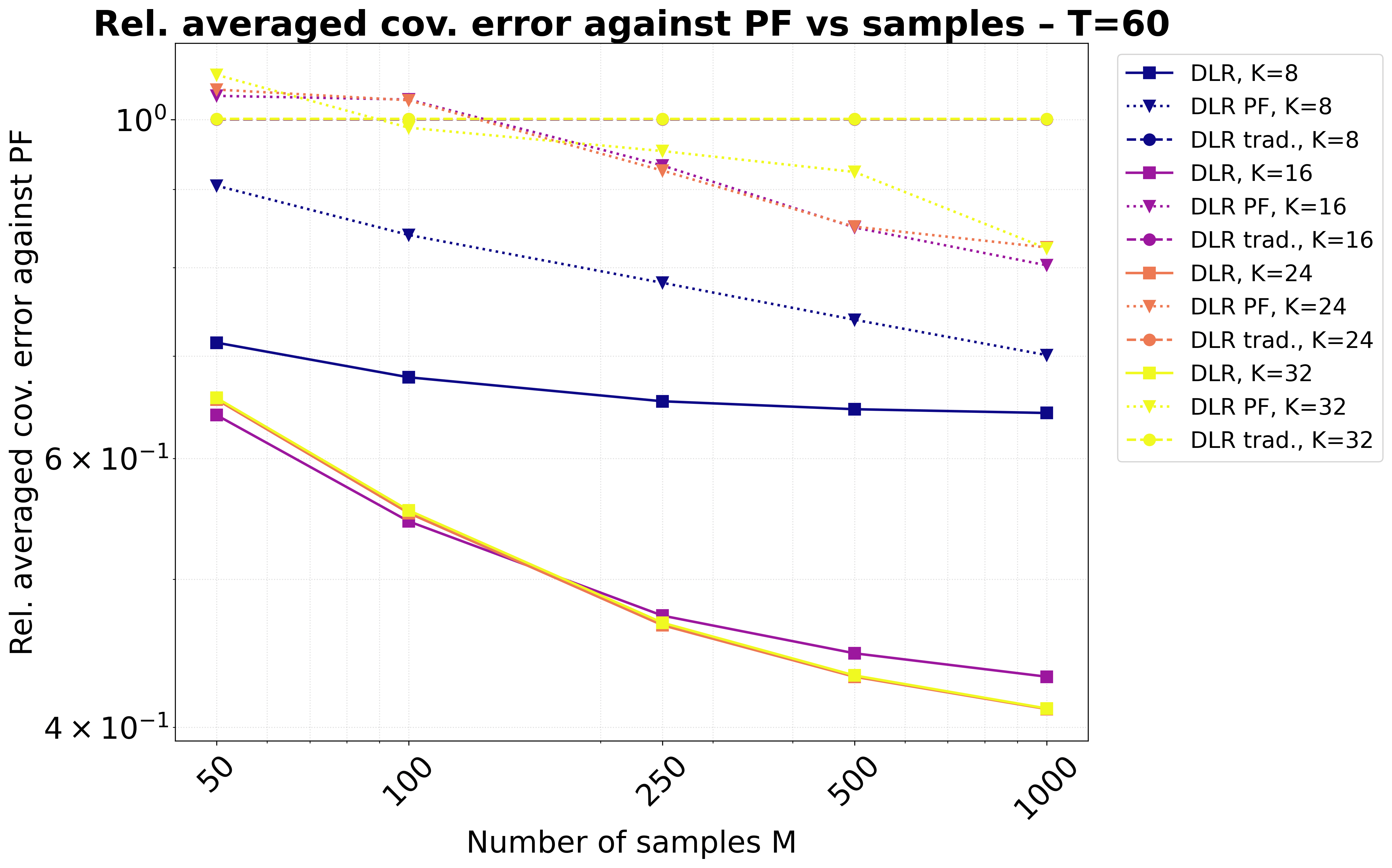} 
	\caption{Rel.\  averaged-over-time \emph{errors against PF}, as a function of the sample size $M$: (left) mean, (right) covariance. Problem \eqref{eq: 2 qg} with continuous-time observations \eqref{eq: 2 qg - obs proc}, transient time $T_0=2$, reference computed with $M=10^5$, $d=2048$, $h=512$, $\Delta t_{branch}= 4\Delta t$.}
	\label{fig: Oceanography PF samples}
\end{figure}

\begin{figure}[!h]
	\centering
	\includegraphics[scale=0.045]{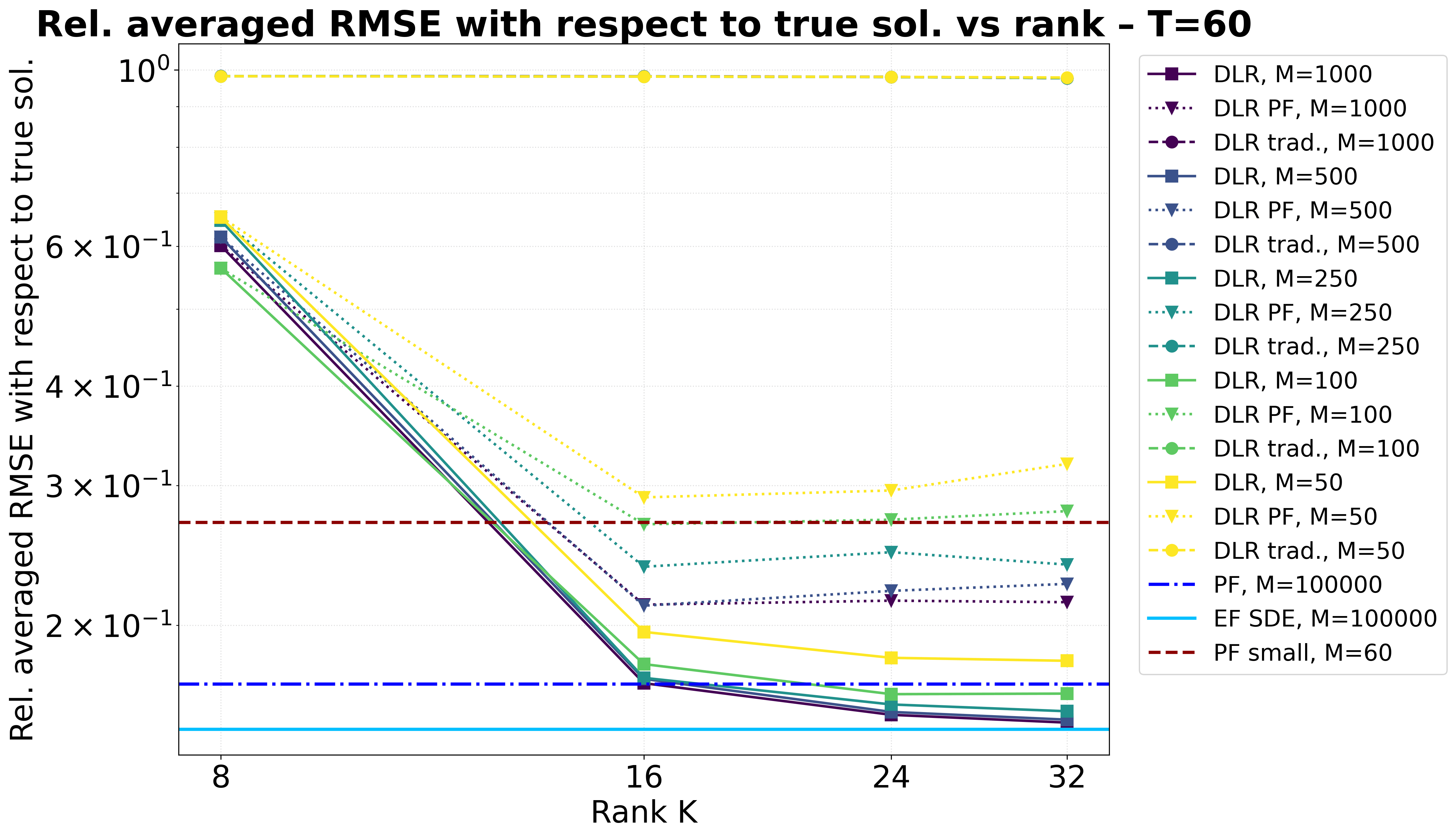}
	\includegraphics[scale=0.043]{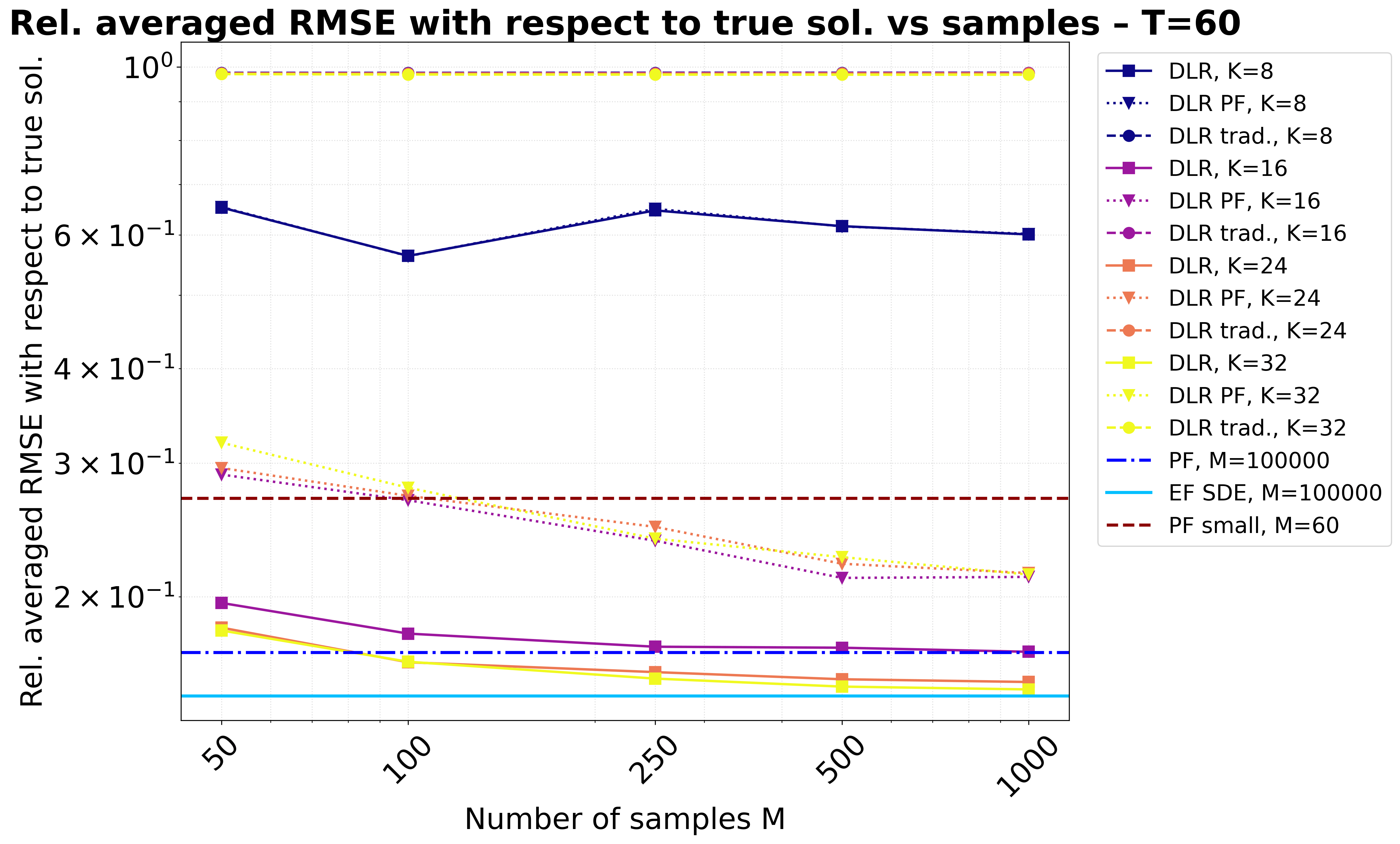}
	\caption{Rel.\  averaged-over-time \emph{RMSE} as a function of the rank $K$ (left) and as a function of the samples sizes $M$ (right). Problem \eqref{eq: 2 qg} with continuous-time observations \eqref{eq: 2 qg - obs proc}, transient time $T_0=2$, reference computed with $M=10^5$, $d=2048$, $h=512$, $\Delta t_{branch}= 4\Delta t$.}
	\label{fig: Oceanography RMSE}
\end{figure}

Similar behaviour can be observed when comparing the DLR filters to the PF, as one can see in Figures  \ref{fig: Oceanography PF rank} and  \ref{fig: Oceanography PF samples}. The results show that the DLR PF requires more samples than its DLR JMCO counterpart to achieve comparable accuracy, further indicating that particle degeneracy is the main limitation of this class of filters. Furthermore, in Figure \ref{fig: Oceanography RMSE} we show the averaged-in-time RMSE. One can notice that all DLRA algorithms perform better than a particle filter with $M=60$ for the same employed number of particles, and they are comparable to the reference PF and EF SDE computed with $M=10^5$.
These results show again the rationale of exploiting DLR strategies in filtering procedures for high-dimensional problems.

\begin{figure}[!h]
	\centering
	\includegraphics[scale=0.19]{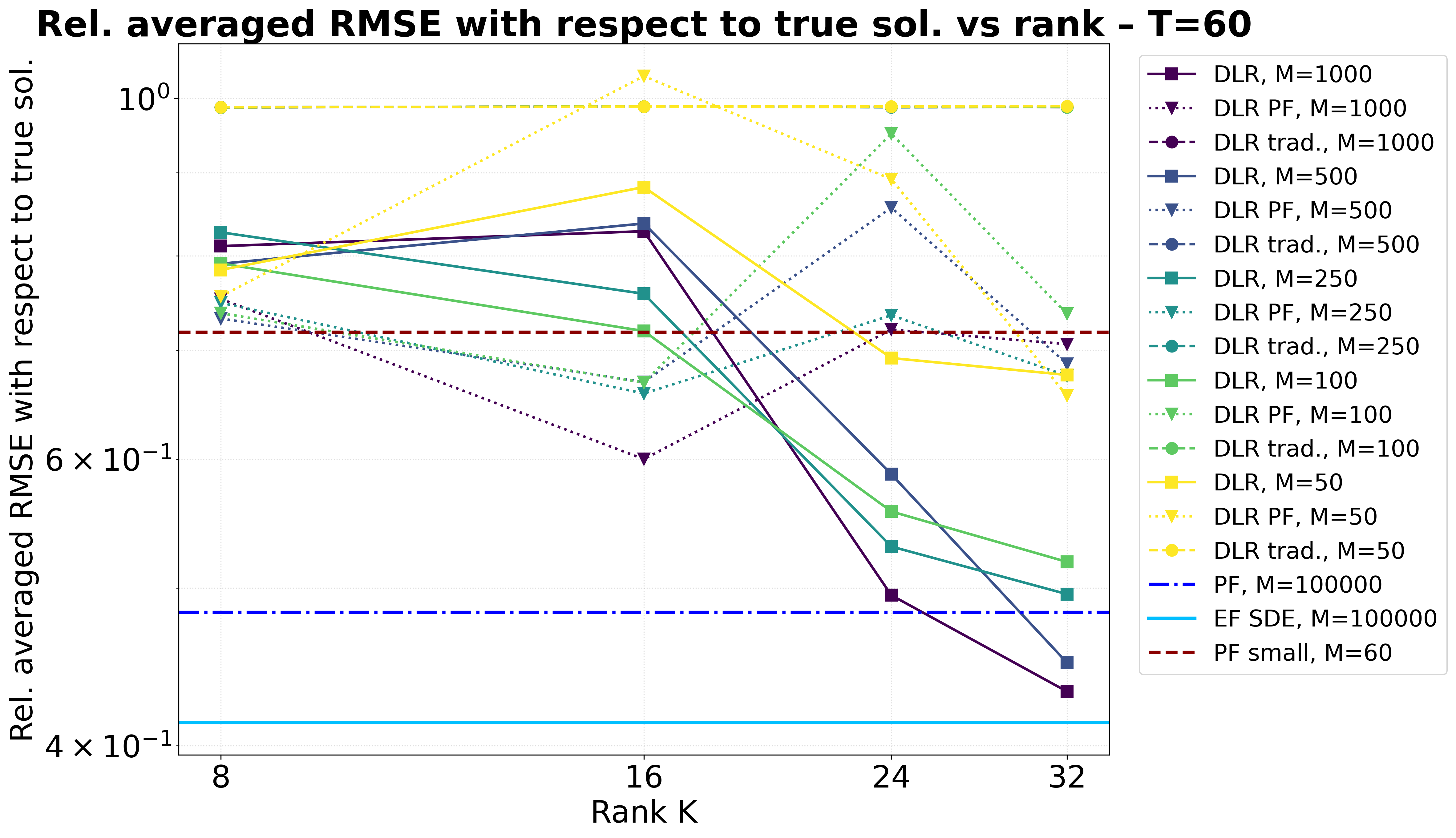}
	\includegraphics[scale=0.18]{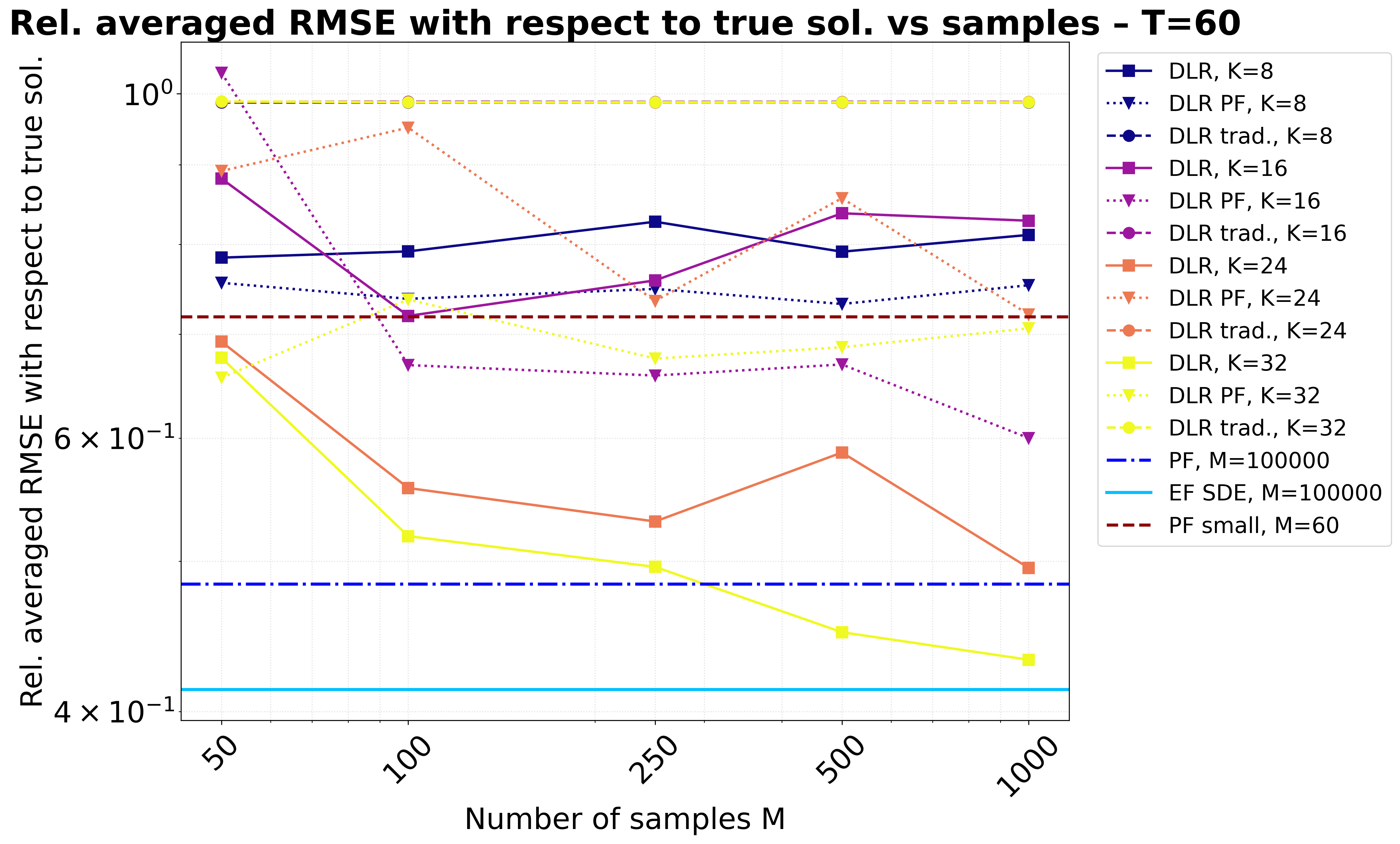}
	\caption{Rel.\  averaged-over-time \emph{RMSE} as a function of the rank $K$ (left) and as a function of the samples sizes $M$ (right). Problem \eqref{eq: 2 qg} with continuous-time observations \eqref{eq: 2 qg - obs proc}, transient time $T_0=2$, reference computed with $M=10^5$, $d=2048$, $h=32$, $\Delta t_{branch}= 4\Delta t$.}
	\label{fig: Oceanography RMSE h 32}
\end{figure}

Finally, we also ran a simulation for the case of $h=32$ observation points, taken with equidistant coordinates in both axis. In Figure \ref{fig: Oceanography RMSE h 32} we see once again that the DLRA JMCO outperforms in terms of RMSE accuracy the other DLR algorithms and approaches the accuracy of the PF and EF SDE for $K=32$ ad $M=500$. This behavior corroborates the fact that DLRA manages to track the main subspace where the solution lives even with a large and uniformly distributed number of observations comparable to the dimension of the noise of the dynamics.

\section{Conclusion and Perspectives}
In this work, we proposed several DLRA approaches for filtering, both for continuous- and discrete-time observations.
These algorithms can be computed completely on-the-fly without the need of an expensive offline phase. Moreover, the time-dependent nature of the DLRA basis allows to cheaply adapt the filtered dynamics to the incoming observations. The benefits of these algorithms were supported by numerical simulations and a more advanced theoretical analysis of these procedures is postponed to future works.

First, we derived a DLRA filter procedure that minimizes jointly the error of mean and covariance between the surrogate and the local full-order approximation. This framework naturally generalizes to ensemble methods to treat general nonlinear and non-Gaussian problems as well as high dimensional linear ones providing a reduced order Kalman-Bucy-type filter.

Then, a complemented DLRA approach is derived targeting to track the directions orthogonal to the main subspace. Unlike other DLRA-type filters, it turns out that this framework provides an equation of the subspace that evolves not only due to the drift of the system, but also due to the observation operator.  This structure is extremely beneficial for data assimilation problems whose observation operator can dramatically change over time or when the subspace individuated by the observation dynamics is different from the one where the state lives.

Finally, we propose a DLRA particle filter algorithm to deal with nonlinear and non-Gaussian systems. This procedure corrects the reduced DLRA dynamics of the prediction by a branching procedure based on the likelihood functional, including the possibility of enriching the predicted subspace via orthogonal components in a cheap way. 

There are several interesting and challenging perspectives of this work, some of which are already under investigation. One is studying the natural extension of the various DLR filters to the smoothing setting, bringing their computational advantages to this more computationally demanding data assimilation procedure.

Concerning the first part of this work, the derivation of the DLR-JMCO filter minimizes some Euclidean-type functionals: in this respect, it might be advantageous to propose filters that minimizes other distances, like the Wasserstein one, exploiting different geometries. Moreover, in order to further deal with nonlinearity in filtering problems, it would be interesting to apply DLRA to the \emph{feedback particle filter} \cite{taghvaei2023survey}, as well as to enhance the deterministic and the stochastic basis via scientific machine learning techniques, like \emph{conditional flow matching}. Last, but not least, in the context of nonlinear problems it seems promising to merge DLRA updates in the prediction step with (triangular) transport maps \cite{parno2018transport} in the analysis procedure, in order to produce a completely online reduced nonlinear filter, resembling the so-called reduced \emph{stochastic map filter} \cite{spantini2022coupling} with the main subspace that is time-evolving.

\section*{Acknowledgements}
This work has also been supported by the Swiss National Science Foundation under the
Project n. 200518 “Dynamical low rank methods for uncertainty quantification and data assimilation” and by the Doc.mobility funding awarded to this project.

\printbibliography

\appendix

\section{Proofs concerning DLRA-JMCO filter}\label{app: proof DLRA JMCO}

\begin{proof}[\bfseries Proof of Proposition \ref{prop: triple JMCO}]
	In order to find the minimizer triple in \eqref{eq: prob DLRA JMCO}, first we compute the quantities
	\begin{equation*}
		|m_{n+1} - m_{n+1}^{\mathrm{DLRA}}|^2 \quad \text{and} \quad
		\| C_{n+1} - C_{n+1}^{\mathrm{DLRA}}  \|_{\mathrm{F}}^2,
	\end{equation*}
	then we set their gradient with respect to the
	unknown triple $(\Delta m_{n}^{\mathrm{DLRA}} ,\Delta U_n,\Delta C_{Y_n})$ equal to $0$
	in a variational formulation. The solution of this variational relation will be our sought DLRA updates. 
	
	\underline{\emph{The equation for $ m_{n+1}^{\mathrm{DLRA}}$.}} Let us first start by minimizing the difference between the means. One has
	\begin{equation}\label{eq: mid m}
		\begin{aligned}
			\big| m_{n+1} - m_{n+1}^{\mathrm{DLRA}}  \big|^{2} 
			= & \Big|
			\mathbb{E}\!\left[A(t_n, X_{n}^{\mathrm{DLRA}} )\right]\Delta t_n + C_{n}^{\mathrm{DLRA}}  H(t_n)^{\top} R^{-1} \\
			& \cdot\left(Z_{n+1}-Z_n- H(t_n)\mathbb{E}[X_{n}^{\mathrm{DLRA}}]\right)´ \Delta t_n 
			- \Delta m_n^{\mathrm{DLRA}}
			\Big|^{2}.
		\end{aligned}
	\end{equation}
	Notice that \eqref{eq: mid m} does not depend on the unknowns $\Delta U_n$ and $\Delta C_{Y_n}$ and we will see later that $\|C_{n+1}- C_{n+1}^{\mathrm{DLRA}}\|_{\mathrm{F}}^2$ does not depend on $\Delta m_n^{\mathrm{DLRA}}$. Therefore, to find $\Delta m_n^{\mathrm{DLRA}}$, one needs to compute only the derivative of \eqref{eq: mid m}, which is
	\begin{equation*}
		\begin{aligned}
			\frac{\partial \big| m_{n+1} - m_{n+1}^{\mathrm{DLRA}}  \big|^{2}}{\partial \Delta m_n^{\mathrm{DLRA}}}
			&=
			\frac{\partial \big| m_{n+1} - m_n^{\mathrm{DLRA}} - \Delta m_n^{\mathrm{DLRA}} \big|^{2}}{\partial \Delta m_n^{\mathrm{DLRA}}}\\
			&= -2\big( m_{n+1} - m_n^{\mathrm{DLRA}} - \Delta m_n^{\mathrm{DLRA}} \big).
		\end{aligned}
	\end{equation*}
Setting the derivative to $0$
	leads to the update
	\begin{equation*}
		\begin{aligned}
			\Delta m_n^{\mathrm{DLRA}}
			&\approx 
			\Big(
			\mathbb{E}[A(t_n, X_{n}^{\mathrm{DLRA}} )]\Delta t_n 
			+ C_{n}^{\mathrm{DLRA}} H(t_n)^{\top} R^{-1}\Big(\Delta Z_n
			- H(t_n)m_n^{\mathrm{DLRA}} \Delta t_n
			\Big).
		\end{aligned}
	\end{equation*}

	\underline{\emph{Minimization of the covariance $C_{n+1}^{\mathrm{DLRA}}$.}}  Now, in order to obtain the update in $\Delta U_n$ and $\Delta C_{Y_n}$, let us minimize the Frobenius norm of the difference between the 
	covariances $C_{n+1}$ and $C_{n+1}^{\mathrm{DLRA}}$.  
	Then, we seek to simultaneously solve the relations
	\begin{equation*}
		\frac{\partial \| C_{n+1} - C_{n+1}^{\mathrm{DLRA}} \|_{\mathrm{F}}^{2} }{\partial \Delta U_n}
		= 0,
		\qquad
		\frac{\partial \| C_{n+1} - C_{n+1}^{\mathrm{DLRA}} \|_{\mathrm{F}}^{2}}{\partial \Delta C_{Y_n}} = 0.
	\end{equation*}
	
	First notice that
	\begin{equation}\label{eq: diff C}
		\| C_{n+1} - C_{n+1}^{\mathrm{DLRA}} \|_{\mathrm{F}}^{2}
		= \operatorname{Tr}( C_{n+1} C_{n+1}^{\!\top} )
		- 2 \operatorname{Tr}( C_{n+1} (C_{n+1}^{\mathrm{DLRA}})^{\!\top} )
		+ \operatorname{Tr}\big( C_{n+1}^{\mathrm{DLRA}} (C_{n+1}^{\mathrm{DLRA}})^{\!\top} \big).
	\end{equation}	
	The first term on the right-hand side of \eqref{eq: diff C} does not depend on $\Delta U_n$ and $\Delta C_{Y_n}$, hence it is not taken into account in the minimization. Moreover, one has	
 \begin{equation*}
		\begin{aligned}
			C_{n+1}^{\mathrm{DLRA}} (C_{n+1}^{\mathrm{DLRA}})^{\!\top}\approx& \left(C_n^{\mathrm{DLRA}}
			+ U_n^{\!\top} C_{Y_n} \Delta U_n 
			+ U_n^{\!\top} \Delta C_{Y_n} U_n
			+ \Delta U_n^{\!\top} C_{Y_n} U_n +  \Delta U_n^{\!\top}  C_{Y_n} \Delta U_n\right) \\
			&\quad  \left(C_n^{\mathrm{DLRA}}
			+ U_n^{\!\top} C_{Y_n} \Delta U_n
			+ U_n^{\!\top} \Delta C_{Y_n} U_n
			+ \Delta U_n^{\!\top} C_{Y_n} U_n +  \Delta U_n^{\!\top}  C_{Y_n} \Delta U_n\right)^{\top} \\
			= &C_{n}^{\mathrm{DLRA}} (C_{n}^{\mathrm{DLRA}})^{\!\top}
			+ U_n^{\top} C_{Y_n} \Delta C_{Y_n} U_n + U_n^{\top} C_{Y_n}^2 \Delta U_n  \\
			&+ U_n^{\top} C_{Y_n} \Delta U_n\Delta U_n^{\top} C_{Y_n} U_n  + U_n^{\top} C_{Y_n} \Delta U_n\Delta U_n^{\top} C_{Y_n} \Delta U_n \\
	& +  U_n^{\top} \Delta C_{Y_n} C_{Y_n} U_n  +  U_n^{\top} (\Delta C_{Y_n})^2 U_n+ U_n^{\top} \Delta C_{Y_n} C_{Y_n} \Delta U_n \\
	&+\Delta U_n^{\top} C_{Y_n}^2 U_n+ \Delta U_n^{\top} C_{Y_n} \Delta C_{Y_n} U_n + \Delta U_n^{\top} C_{Y_n}^2 \Delta U_n \\
	&+ \Delta U_n^{\top} C_{Y_n} \Delta U_n\Delta U_n^{\top} C_{Y_n} U_n +\Delta U_n^{\top} C_{Y_n} \Delta U_n\Delta U_n^{\top} C_{Y_n} \Delta U_n \\
		\end{aligned}
	\end{equation*}
	where in the second line we exploit the discrete gauge condition.
	Then, using properties \eqref{eq: matrix der} one gets
	\begin{equation*}
		\begin{aligned}
			\frac{\partial \operatorname{Tr}\!\left(C_{n+1}^{\mathrm{DLRA}}
				(C_{n+1}^{\mathrm{DLRA}})^{\top}\right)}{\partial \Delta U_n}
			=& 2C_{Y_n}^2 U_n
		+2 C_{Y_n}^2 \Delta U_n
			+ 2 C_{Y_n}^2  \Delta U_n + 2C_{Y_n} \Delta C_{Y_n} U_n\,  \\
			& +2 C_{Y_n} \Delta U_n \Delta U_n^{\top} C_{Y_n} U_n + 2 C_{Y_n} \Delta U_n \Delta U_n^{\top} C_{Y_n} \Delta U_n \\
			& + 2 C_{Y_n} \Delta U_n \Delta U_n^{\top} C_{Y_n} \Delta U_n\\
			=& 2 C_{Y_n}^2 U_n
			+4 C_{Y_n}^2\, \Delta U_n
			+ 2 C_{Y_n} \Delta C_{Y_n} U_n \\
			&+2 C_{Y_n} \Delta U_n \Delta U_n^{\top} C_{Y_n} U_n + 4 C_{Y_n} \Delta U_n \Delta U_n^{\top} C_{Y_n} \Delta U_n,
		\end{aligned}
	\end{equation*}
	where we employed the orthogonality of the rows of $U_n$ and the discrete gauge condition.  
	
	Similarly, for the derivative with respect to $\Delta C_{Y_n}$ we obtain
	\begin{equation*}
		\begin{aligned}
			\frac{\partial \operatorname{Tr}\!\left(C_{n+1}^{\mathrm{DLRA}}
				(C_{n+1}^{\mathrm{DLRA}})^{\top}\right)}{\partial \Delta C_{Y_n}}
			=& C_{Y_n} U_n U_n^{\top}
			+ U_n U_n^{\top} C_{Y_n}  \\
			&+U_n U_n^{\top} \Delta C_{Y_n}
			+ \Delta C_{Y_n} U_n U_n^{\top}\\
			= & 2 C_{Y_n} + 2\, \Delta C_{Y_n},
		\end{aligned}
	\end{equation*}
	where we exploited again the discrete gauge condition.
	
	Concerning the crossed terms between the two covariances from \eqref{eq: diff C}, using again the discrete gauge condition, one obtains
	\begin{equation*}
		\begin{aligned}
		&	\frac{\partial \operatorname{Tr}\!\left( C_{n+1} (C_{n+1}^{\mathrm{DLRA}})^{\top} \right)}{\partial \Delta U_n}\\
			= & \frac{\partial \operatorname{Tr}\!\left(
				C_{n+1}
				\big(
				C_{n}^{\mathrm{DLRA}}
				+ \Delta U_n^{\top} C_{Y_n} U_n
				+ U_n^{\top} C_{Y_n} \Delta U_n
				+ U_n^{\top} \Delta C_{Y_n} U_n
				+ \Delta U_n^{\top} C_{Y_n} \Delta U_n
				\big)^\top
				\right) }{\partial \Delta U_n}\\
			= & \frac{\partial \operatorname{Tr}\!\left(
				C_{n+1}
				(C_{n}^{\mathrm{DLRA}})^{\top}
				+ C_{n+1} U_n^{\top} C_{Y_n} \Delta U_n
				+ C_{n+1}\Delta U_n^{\top} C_{Y_n}  U_n
				+ C_{n+1}U_n^{\top} \Delta C_{Y_n} U_n 
				+ C_{n+1} \Delta U_n^{\top} C_{Y_n} \Delta U_n\right)}{\partial \Delta U_n}\\
			= & C_{Y_n} U_n C_{n+1}
			+ C_{Y_n} U_n C_{n+1} + C_{Y_n} \Delta U_n C_{n+1}
			+ C_{Y_n} \Delta U_n C_{n+1} \\
			= &2\, C_{Y_n} U_n C_{n+1} + 2\, C_{Y_n} \Delta U_n C_{n+1}
			 \end{aligned}
	\end{equation*}
	and, similarly,
	\begin{equation*}
		\begin{aligned}
			\frac{\partial \operatorname{Tr}\!\left(
				C_{n+1} (C_{n+1}^{\mathrm{DLRA}})^{\top}
				\right)}{\partial \Delta C_{Y_n}}
			= U_n C_{n+1} U_n^{\top}.
		\end{aligned}
	\end{equation*}
	
	We are now ready to minimize \eqref{eq: diff C}.  
	With respect to the covariance increment it holds that
	\begin{equation}\label{eq: prel Delta C_y}
		\begin{aligned}
			\frac{\partial \|C_{n+1} - C_{n+1}^{\mathrm{DLRA}} \|_{\mathrm{F}}^2}{\partial \Delta C_{Y_n}}	
			=&
			\frac{\partial
				\operatorname{Tr}\!\left(
				-2 (C_{n+1} (C_{n+1}^{\mathrm{DLRA}})^{\top})
				+ C_{n+1}^{\mathrm{DLRA}} (C_{n+1}^{\mathrm{DLRA}})^{\top}
				\right)}{\partial \Delta  C_{Y_n}} \\
			=& 2C_{Y_n} + 2\Delta C_{Y_n}
			- 2 U_n C_{n+1} U_n^{\top},
		\end{aligned}
	\end{equation}
	and %assuming that $C_{Y_n}$, $\Delta C_{Y_n}$, and $U_n C_{n+1} U_n^{\top}$, which are all $k \times k$ matrices, are always of rank $k$, 
	setting \eqref{eq: prel Delta C_y} to $0$ implies
	\begin{equation}\label{eq: Delta C_y}
		\Delta C_{Y_n} = U_n C_{n+1} U_n^{\top} - C_{Y_n}.
	\end{equation}
	Replacing the expression of $C_{n+1}$ from \eqref{eq: EM cond cov - JMCO} we finally obtain 
		\begin{equation}\label{eq: Delta C_y full}
		\begin{aligned}
			\Delta C_{Y_n}
			=& U_n C_{n+1} U_n^\top - C_{Y_n} \\
			\approx &		
			 \mathbb{E}[Y_n\, \mathring{A}(t_n, X_n^{\mathrm{DLRA}})^\top] \, U_n^\top \Delta t_n
			+ U_n\, \mathbb{E}[\mathring{A}(t_n, X_n^{\mathrm{DLRA}}) Y_n^\top]\, \Delta t_n \\
			&+ U_n Q U_n^\top \Delta t_n
			-C_{Y_n} U_n H(t_n)^\top R^{-1} H(t_n) U_n^{\top} C_{Y_n}\, \Delta t_n.
		\end{aligned}
	\end{equation}	
Notice that $\Delta C_{Y_n}$ does not depend on $\Delta Z_n$.
Notice also that the difference $C_{n+1} - C_{n+1}^{\mathrm{DLRA}}$ does not depend on $\Delta m_n^{\mathrm{DLRA}}$, as it was anticipated.
    
	\underline{\emph{The equation for $U_{n+1}$.}}
	Considering the derivative with respect to the basis increment $\Delta U_n$, one gets
	\begin{equation*}
		\begin{aligned}
			\frac{\partial \|C_{n+1} - C_{n+1}^{\mathrm{DLRA}} \|_{\mathrm{F}}}{\partial \Delta U_{n}}
			=&
			\frac{\partial
				\operatorname{Tr}\!\left(
				-2 (C_{n+1} (C_{n+1}^{\mathrm{DLRA}})^{\top})
				+ C_{n+1}^{\mathrm{DLRA}} (C_{n+1}^{\mathrm{DLRA}})^{\top}
				\right)}{\partial \Delta U_n}. \\
			= & 2\, C_{Y_n}^2 U_n
			+ 4\, C_{Y_n}^2 \Delta U_n
			+ 2\, C_{Y_n} \Delta C_{Y_n} U_n
			- 4\, C_{Y_n} U_n C_{n+1} \\
			&- 4\, C_{Y_n} \Delta U_n C_{n+1}+2 C_{Y_n} \Delta U_n \Delta U_n^{\top} C_{Y_n} U_n + 4 C_{Y_n} \Delta U_n \Delta U_n^{\top} C_{Y_n} \Delta U_n,\\
			= &2\, C_{Y_n}^2 U_n
			+ 4\, C_{Y_n}^2 \Delta U_n
			+ 2\, C_{Y_n} \big( U_n C_{n+1} U_n^{\top} - C_{Y_n} \big) U_n
			- 4\, C_{Y_n} U_n C_{n+1} \\
			&- 4\, C_{Y_n} \Delta U_n C_{n+1}+2 C_{Y_n} \Delta U_n \Delta U_n^{\top} C_{Y_n} U_n + 4 C_{Y_n} \Delta U_n \Delta U_n^{\top} C_{Y_n} \Delta U_n, \\
				= &
			4\, C_{Y_n}^2 \Delta U_n
			+ 2\, C_{Y_n} U_n C_{n+1} P_{U_n} 
			- 4\, C_{Y_n} U_n C_{n+1} - 4\, C_{Y_n} \Delta U_n C_{n+1}\\
			&+2 C_{Y_n} \Delta U_n \Delta U_n^{\top} C_{Y_n} U_n + 4 C_{Y_n} \Delta U_n \Delta U_n^{\top} C_{Y_n} \Delta U_n \\
			\approx & 4\, C_{Y_n}^2 \Delta U_n + 2\, C_{Y_n} U_n C_{n+1} P_{U_n}- 4\, C_{Y_n} U_n C_{n+1} +2 C_{Y_n} \Delta U_n \Delta U_n^{\top} C_{Y_n} U_n
		\end{aligned}
	\end{equation*}
where, in the last line, we have retained only zero-th and first order terms in $\Delta t_n$. Notice, in particular, that $C_{n+1} = U_n^{\top}C_{Y_n}U_n + O(\Delta t_n)$, so that $4C_{Y_n} \Delta U_n C_{n+1}$ is of higher order and has been neglected.
	We want to solve
	\begin{equation*}
		\big\langle \frac{\partial}{\partial \Delta U_n}
		\operatorname{Tr}\!\left( (C_{n+1} - C_{n+1}^{\mathrm{DLRA}})
		(C_{n+1} - C_{n+1}^{\mathrm{DLRA}})^{\top} \right),\,
		V \big\rangle_{\mathrm{F}} = 0,
		\quad
		\forall V \in \mathbb{R}^{k \times d}, \ V^{\top} \in \mathrm{Im}(U_n^{\top})^{\perp}.
	\end{equation*}
	Using the properties of the Frobenius inner product, this corresponds (at first order) to
		\begin{equation*}
		\begin{aligned}
			 \big\langle
			&4\, C_{Y_n} \Delta U_n
			+ 2 U_n C_{n+1} P_{U_n} 
			- 4\, U_n C_{n+1}   \\
			& 
			\qquad + 2 \Delta U_n \Delta U_n^{\top} C_{Y_n} U_n ,\,
			C_{Y_n} V 
			\big\rangle_{\mathrm{F}} =0, \quad
			\forall V^{\top} \in \mathrm{Im}(U_n^{\top})^{\perp},
		\end{aligned}
	\end{equation*}
	and, since, $C_{Y_n}$ is assumed to be invertible by hypothesis (i.e.\ $Y_n$ has rank $k$) , this is equivalent to
	\begin{equation*}
	\begin{aligned}
		\big\langle
		&4\, C_{Y_n} \Delta U_n
		+ 2 U_n C_{n+1} P_{U_n} 
		- 4\, U_n C_{n+1}  \\
		& \qquad
		+ 2\Delta U_n \Delta U_n^{\top} C_{Y_n} U_n ,\,
		V 
		\big\rangle_{\mathrm{F}} =0, \quad
		\forall V^{\top} \in \mathrm{Im}(U_n^{\top})^{\perp},
	\end{aligned}
\end{equation*}
	which implies that
	\begin{equation*}
		\langle 4\, C_{Y_n} \Delta U_n , \,V \rangle_{\mathrm{F}}
		= \left\langle 4 U_n C_{n+1},\, V
		\right\rangle_{\mathrm{F}}, \quad
		\forall V \in \mathbb{R}^{k \times d}, \ V^{\top} \in \mathrm{Im}(U_n^{\top})^{\perp},
	\end{equation*}
	thanks to the fact that $U_n V^{\top} =0$.
	Equivalently, one can write
		\begin{equation*}
		\begin{aligned}
			\langle  \Delta U_n \ ,V \rangle_{\mathrm{F}}
			= &
			\left\langle
			C_{Y_n}^{-1} U_n C_{n+1},\,V
			\right\rangle_{\mathrm{F}}, \quad
			\forall V \in \mathbb{R}^{k \times d}, \ V^{\top} \in \mathrm{Im}(U_n^{\top})^{\perp}.\\
		\end{aligned}
	\end{equation*}	
	Therefore, our deterministic increment is
	\begin{equation*}
		\Delta U_n 
		= C_{Y_n}^{-1} U_n C_{n+1}\,(I_{d \times d} - P_{U_n}).
	\end{equation*}	
	and replacing the expression of $C_{n+1}$ from \eqref{eq: EM cond cov - JMCO} we finally obtain
	\begin{equation*}
		\begin{aligned}
			\Delta U_n
			= &C_{Y_n}^{-1} 
			\, \mathbb{E}\!\left[ Y_n\, \mathring{A}(t_n, X_n^{\mathrm{DLRA}})^\top \right] 
			P_{U_n}^\perp\, \Delta t_n
			\;+\;
			C_{Y_n}^{-1} U_n Q\, P_{U_n}^\perp\, \Delta t_n.
		\end{aligned}
	\end{equation*}
	Notice that the obtained increment $\Delta U_n$ does not depend on $\Delta Z_n$, which implies $f_n=0$ in the representation \eqref{eq: updates}.
	
	\underline{\emph{The equation for $Y_{n+1}$.}} 
	We now seek an increment \(\Delta Y_n\) whose covariance increment matches the previously derived \(\Delta C_{Y_n}\) to first order. Notice that
	\begin{equation*}
		\begin{aligned}
			C_{Y_{n+1}}
			= &\mathbb{E}\!\left[(Y_n + \Delta Y_n)(Y_n + \Delta Y_n)^\top \mid \{Z_{\ell}\}_{\ell = 0, \dots, n}, \Delta Z_n\right]\\
			= &C_{Y_n}
			+ \mathbb{E}[Y_n \Delta Y_n^\top \mid \Delta Z_n]
			+ \mathbb{E}[\Delta Y_n Y_n^\top \mid \Delta Z_n]
			+ \mathbb{E}[\Delta Y_n \Delta Y_n^\top\mid \Delta Z_n ].
		\end{aligned}
	\end{equation*}
	and, hence,
	\begin{equation}\label{eq: Delta C_Y final}
		\begin{aligned}
			\Delta C_{Y_n}
			= &\mathbb{E}[Y_n\, \Delta Y_n^\top \mid \Delta Z_n]
			+ \mathbb{E}[\Delta Y_n\, Y_n^\top \mid \Delta Z_n]
			+ \mathbb{E}[\Delta Y_n \Delta Y_n^\top \mid \Delta Z_n],
		\end{aligned}
	\end{equation}
	One can easily check that the following update of the stochastic modes
	\begin{equation*}
		\begin{aligned}
			\Delta Y_n
			= &U_n\, \mathring{A}(t_n, X_n^{\mathrm{DLRA}})\, \Delta t_n
			+ U_n Q^{\frac{1}{2}} \Delta \widehat{W}_n\\
			&- C_{Y_n} U_n H(t_n)^\top R^{-1} H(t_n) U_n^\top Y_n\, \Delta t_n
			- C_{Y_n} U_n H(t_n)^{\top} R^{-\frac{1}{2}} \Delta \widehat{B}_n,
		\end{aligned}
	\end{equation*}
	allows one to match relations \eqref{eq: Delta C_Y final} and \eqref{eq: Delta C_y full} up to first order.
\end{proof}

\section{Proof concerning enriched DLRA-JMCO filter}\label{app: proof DLRA-JMCO enriched}

\begin{proof}[\bfseries Proof of Lemma \ref{lem: approx error}]
	First, we recall the updates of the DLR Projector Splitting for SDEs, given in \cite{kazashi2025dynamicalpartI}, with the noise correction for the dynamics obtained in the derivation of the DLRA JMCO filter in Section \ref{sec: DLRA for JMCO}. The updates of the bases reads as
	\begin{equation*}
		\begin{aligned}
			\widehat{m}_{n+1} & = m_n + \mathbb{E}\left[A(t_n,X_{n}^{\mathrm{DLRA}})\right] \Delta t_n\\
			\widetilde{Y}_{n+1} &= Y_n + U_n\, \mathring{A}(t_n,X_{n}^{\mathrm{DLRA}}) \Delta t_n +  U_n\, Q^{\frac{1}{2}} \Delta W_n  \\
			C_{\widetilde{Y}_{n+1}} \widetilde{U}_{n+1} &= C_{\widetilde{Y}_{n+1}} U_n + \left( \mathbb{E}\left[\widetilde{Y}_{n+1}\, \mathring{A}(t_n,X_{n}^{\mathrm{DLRA}})^{\top} \right] + U_n Q \right)P_{U_n}^{\perp} \Delta t_n,
		\end{aligned}
	\end{equation*}
	whereas the ambient solution is written as
	\begin{equation}
		\begin{aligned}
			\widehat{X}_{n+1}^{\mathrm{DLRA}}
			= & X_n^{\mathrm{DLRA}} + \mathbb{E}\left[A(t_n,X_{n}^{\mathrm{DLRA}})\right] \Delta t_n + P_{{U}^{\top}_n\widetilde{Y}_{n+1}}[ \mathring{A}(t_n,X_{n}^{\mathrm{DLRA}})] \Delta t_n \\
			&+P_{U_n}^{\perp}C_{\widetilde{Y}_{n+1}}^{-1}\widetilde{Y}_{n+1} Q \Delta t_n + P_{{U}_n}Q^{\frac{1}{2}}\Delta W_n.
		\end{aligned} 
	\end{equation}

	Let us begin from the first inequality. Suppose that $X_{n} = X_{n}^{\mathrm{DLRA}}$. Then one can bound the local error in the
	following way
	\begin{equation*}
		\begin{aligned}
			\mathbb{E}\!\left[\left| \widehat{X}_{n+1} - \widehat{X}^{\mathrm{DLRA}}_{n+1} \right|^2 \right]
			&= \mathbb{E}\Bigg[
			\Bigg| 
			X_{n}^{\mathrm{DLRA}} - X_n^{\mathrm{DLRA}}
			+ 
			\big(
			\mathbb{E}[A(t_n,X_n^{\mathrm{DLRA}})] - \mathbb{E}[A(t_n,X_{n}^{\mathrm{DLRA}})]
			\big)\,\Delta t_n \\
			&\quad
			+  \int_{t_n}^{t_{n+1}}
			\big(
			\mathring{A}(t_n,X_n^{\mathrm{DLRA}})
			- P_{{U}^{\top}_n\widetilde{Y}_{n+1}} \mathring{A}(t_n,X_n)
			- P_{U_n}^\perp Q U_n^{\top}  C_{\widetilde{Y}_{n+1}}^{-1} \widetilde{Y}_{n+1}
			\big)\,\mathrm{d}s \\
			&\quad
			+ \int_{t_n}^{t_{n+1}}
			\big(
			Q^{\frac{1}{2}} - P_{U_n} Q^{\frac{1}{2}}
			\big)\,\mathrm{d}W_s
			\Bigg|^2
			\Bigg]\\
			&= \mathbb{E}\Bigg[
			\Bigg|
			\int_{t_n}^{t_{n+1}}
			\big(
			\mathring{A}(t_n,X_n^{\mathrm{DLRA}})
			- P_{{U}^{\top}_n\widetilde{Y}_{n+1}} \mathring{A}(t_n,X_n^{\mathrm{DLRA}})
			- P_{U_n}^\perp Q U_n^{\top}  C_{\widetilde{Y}_{n+1}}^{-1} \widehat{Y}_{n+1}
			\big)\,\mathrm{d}s \\
			&\quad
			+ \int_{t_n}^{t_{n+1}}
			\big(
			Q^{\frac{1}{2}} - P_{U_n} Q^{\frac{1}{2}}
			\big)\,\mathrm{d}W_s
			\Bigg|^2
			\Bigg]
		\end{aligned}
	\end{equation*}
	Then, via standard inequalities one has that
	\begin{equation*}
		\begin{aligned}
			\mathbb{E}\!\left[\left| \widehat{X}_{n+1} - \widehat{X}^{\mathrm{DLRA}}_{n+1} \right|^2 \right]\leq&
			3 \mathbb{E}\!\left[
			\left|
			\int_{t_n}^{t_{n+1}}
			\mathring{A}(t_n,X_n^{\mathrm{DLRA}})
			- P_{{U}^{\top}_n\widetilde{Y}_{n+1}} \mathring{A}(t_n,X_n^{\mathrm{DLRA}})
			\,\mathrm{d}s
			\right|^2
			\right] \\
			&
			+ 3 \mathbb{E}\!\left[
			\left|
			\int_{t_n}^{t_{n+1}}
			P_{U_n}^\perp Q U_n^{\top}  C_{\widetilde{Y}_{n+1}}^{-1} \widetilde{Y}_{n+1}
			\,\mathrm{d}s
			\right|^2
			\right] 
			+ 3\mathbb{E}\!\left[
			\left|
			\int_{t_n}^{t_{n+1}}
			P_{U_n}^\perp Q^{\frac{1}{2}}
			\,\mathrm{d}W_s
			\right|^2
			\right].
		\end{aligned}
	\end{equation*}
	Using Jensen's inequality, Ito's isometry, and Assumption \ref{ass: eps error} one gets
	\begin{equation*}
		\begin{aligned}
			\mathbb{E}\!\left[\left| \widehat{X}_{n+1} - \widehat{X}^{\mathrm{DLRA}}_{n+1} \right|^2 \right]
			\leq&
			3 \Delta t_n
			\mathbb{E}\!\left[
			\int_{t_n}^{t_{n+1}}
			\left|
			\mathring{A}(t_n,X_n^{\mathrm{DLRA}}) - P_{{U}^{\top}_n\widetilde{Y}_{n+1}} \mathring{A}(t_n,X_n^{\mathrm{DLRA}})
			\right|^2
			\,\mathrm{d}s
			\right] \\
			&
			+ 3 \Delta t_n
			\mathbb{E}\!\left[
			\int_{t_n}^{t_{n+1}}
			\left|
			P_{U_n}^\perp Q U_n^{\top}  C_{\widetilde{Y}_{n+1}}^{-1} \widetilde{Y}_{n+1}
			\right|^2
			\,\mathrm{d}s
			\right] 
			+ 3
			\mathbb{E}\!\left[
			\int_{t_n}^{t_{n+1}}
			\left\|
			P_{U_n}^\perp Q^{\frac{1}{2}}
			\right\|_{\mathrm{F}}^2
			\,\mathrm{d}s
			\right]\\
			\leq&
			3 \Delta t_n
			\mathbb{E}\!\Big[
			\int_{t_n}^{t_{n+1}}
			\Big|
			\mathring{A}(t_n,X_n^{\mathrm{DLRA}}) - A(t_n,{U}^{\top}_n\widetilde{Y}_{n+1})+A(t_n,{U}^{\top}_n\widetilde{Y}_{n+1}) \\
			&- P_{{U}^{\top}_n\widetilde{Y}_{n+1}} \mathring{A}(t_n,{U}^{\top}_n\widetilde{Y}_{n+1})
			+ P_{{U}^{\top}_n\widetilde{Y}_{n+1}} \mathring{A}(t_n,{U}^{\top}_n\widetilde{Y}_{n+1})- P_{{U}^{\top}_n\widetilde{Y}_{n+1}} \mathring{A}(t_n,X_n^{\mathrm{DLRA}})
			\Big|^2
			\,\mathrm{d}s
			\Big] \\
			&
			+ 3 \Delta t_n
			\mathbb{E}\!\left[
			\int_{t_n}^{t_{n+1}}
			\left|
			P_{U_n}^\perp Q U_n^{\top}  C_{\widetilde{Y}_{n+1}}^{-1} \widetilde{Y}_{n+1}
			\right|^2
			\,\mathrm{d}s
			\right] 
			+ 3
			\mathbb{E}\!\left[
			\int_{t_n}^{t_{n+1}}
			\left\|
			P_{U_n}^\perp Q^{\frac{1}{2}}
			\right\|_{\mathrm{F}}^2
			\,\mathrm{d}s
			\right]\\
			&\leq
			3 \Delta t_n^2 \varepsilon^2_n
			+ 3 \varepsilon^2_n  \Delta t_n
			\int_{t_n}^{t_{n+1}}
			\mathbb{E}\!\left[
			\left|
			Q^{\frac{1}{2}} U_n^{\top}  C_{\widetilde{Y}_{n+1}}^{-1} \widetilde{Y}_{n+1}
			\right|^2
			\right]
			\,\mathrm{d}s + 3 \varepsilon^2 \Delta t_n \\
			\leq&
			3 \varepsilon^2_n \Delta t_n
			+ 3 \varepsilon^2_n \Delta t_n^2
			\mathbb{E}\!\left[
			\left| Q^{\frac{1}{2}} U_n^{\top} C_{\widetilde{Y}_{n+1}}^{-\frac{1}{2}} \right|^2
			\right]
			+ 3 \varepsilon^2_n \Delta t_n.
		\end{aligned}
	\end{equation*}
	Now notice that 
\begin{equation*}
	\begin{aligned}
		C_{\widetilde{Y}_{n+1}} = & \mathbb{E}\left[ \left(Y_n +  U_n\mathring{A}(t_n,X_n) \Delta t_n \right)\left(Y_n +  U_n\mathring{A}(t_n,X_n)^{\top} \Delta t_n \right)^{\top}\right] +   U_n Q U_n^{\top} \Delta t_n,
	\end{aligned}
\end{equation*}
which implies that
\begin{equation*}
	\begin{aligned}
		\Delta t_n \left| Q^{\frac{1}{2}} U_n^{\top} C_{\widetilde{Y}_{n+1}}^{-\frac{1}{2}} \right|^2 \leq  &	\Delta t_n \left\| Q^{\frac{1}{2}} U_n^{\top} C_{\widetilde{Y}_{n+1}}^{-\frac{1}{2}} \right\|^2_{\mathrm{F}} \\
		= & \mathrm{Tr}\left(\Delta t_n C_{\widetilde{Y}_{n+1}}^{-\frac{1}{2}} U_n Q^{\frac{1}{2}}  Q^{\frac{1}{2}} U_n^{\top} C_{\widetilde{Y}_{n+1}}^{-\frac{1}{2}}\right)\\
		= & \mathrm{Tr}\left( C_{\widetilde{Y}_{n+1}}^{-1} (U_n Q U_n^{\top} \Delta t_n) \right) \leq k
	\end{aligned}
\end{equation*}
where in the last line we use the cyclic property of the trace and Cauchy-Schwarz inequality, and the thesis follows for the first inequality.
	
To prove the bound on the covariance difference, one can just notice that
		\begin{equation*}
			\begin{aligned}
				\left\| \mathbb{E}\!\left[(\widehat{X}_{n+1}-\widehat{X}^{\mathrm{DLRA}}_{n+1})(\widehat{X}_{n+1}-\widehat{X}^{\mathrm{DLRA}}_{n+1})^{\top}\right]\right\|_{\mathrm{F}}
				&\leq \mathbb{E}\!\left[\left\|(\widehat{X}_{n+1}-\widehat{X}^{\mathrm{DLRA}}_{n+1})(\widehat{X}_{n+1}-\widehat{X}^{\mathrm{DLRA}}_{n+1})^{\top}\right\|_{\mathrm{F}}\right] \\
				&= \mathbb{E}\!\left[\left|\widehat{X}_{n+1}-\widehat{X}^{\mathrm{DLRA}}_{n+1}\right|^2\right] \leq \Lambda_{n+1}.
			\end{aligned}
		\end{equation*}
\end{proof}

	\begin{proof}[\bfseries Proof of Proposition \ref{prop: compl DLRA-JMCO}]
		First notice that by construction, $C_{n+1}$ is not dependent on $m_{n+1}$ and $\Delta m_{n}^c$ is not dependent on $\Delta U_{n}^c$ and $\Delta Y_{n}^c$. Therefore, to find $\Delta U_{n}^c$ and $\Delta Y_{n}^c$ we just need to minimize $C_{n+1}$ with respect to these increments, whereas to obtain the update $\Delta m_{n}^c$ we solve a first order optimality condition for $m_{n+1}$ with respect to $\Delta m_{n}^c$.
		
		We start by deriving equations for $\Delta U_n^c$ and $\Delta C_{Y_n}^c$. First, we approximate \eqref{eq: cov enriched} in a first order of approximation in $O(\Delta t)$ and then we will solve \eqref{eq: DLRA JMCO enriched min}.
		
		Recall that the predicted deterministic basis is
		\begin{equation}\label{eq: pred U incr}
			\begin{aligned}
				\widehat{U}_{n+1}
				= U_n + \Delta \widehat{U}_n
				= U_n
				+ C_{Y_n}^{-1}\mathbb{E}\!\left[Y_n \mathring{A}(t_n, X_n^{\mathrm{DLRA},c})^\top\right]P_{U_n}^{\perp}\Delta t
				+ C_{Y_n}^{-1}U_n QP_{U_n}^{\perp}\Delta t,
			\end{aligned}
		\end{equation}
		and the predicted DLRA covariance is approximated in an order of $O(\Delta t)$ as
		\begin{equation*}
			\begin{aligned}
				\widehat{C}_{n+1}^{\mathrm{DLRA}}
				=&\ \widehat{U}_{n+1}^{\top}C_{\widehat{Y}_{n+1}}\widehat{U}_{n+1} = \widehat{U}_{n+1}^{\top} R_{n+1}C_{\widetilde{Y}_{n+1}} R_{n+1}^{\top}\widehat{U}_{n+1} = \widetilde{U}_{n+1}^{\top}C_{\widetilde{Y}_{n+1}}\widetilde{U}_{n+1} \\
				= &\ \left(
				U_n
				+ C_{Y_n}^{-1}\mathbb{E}\!\left[Y_n \mathring{A}(t_n, X_n^{\mathrm{DLRA},c})^\top\right]P_{U_n}^{\perp}\Delta t
				+ C_{Y_n}^{-1}U_n QP_{U_n}^{\perp}\Delta t
				\right)^\top \\
				&\quad
				\left(
				C_{Y_n}
				+ U_n\mathbb{E}\!\left[\mathring{A}(t_n, X_n^{\mathrm{DLRA},c})Y_n^\top\right]\Delta t
				+ \mathbb{E}\!\left[Y_n \mathring{A}(t_n, X_n^{\mathrm{DLRA},c})^\top\right]U_n\Delta t
				+ U_n^\top Q U_n\Delta t
				\right) \\
				&\quad
				\left(
				U_n
				+ C_{Y_n}^{-1}\mathbb{E}\!\left[Y_n \mathring{A}(t_n, X_n^{\mathrm{DLRA},c})^\top\right]P_{U_n}^{\perp}\Delta t
				+ C_{Y_n}^{-1}U_n QP_{U_n}^{\perp}\Delta t
			\right) + O(\Delta t^2)\\
				\approx&\ C_n^{\mathrm{DLRA}}
				+ P_{U_n}^{\perp}\mathbb{E}\!\left[\mathring{A}(t_n, X_n^{\mathrm{DLRA},c})Y_n^{\top}\right]U_n\Delta t
				+ P_{U_n}^{\perp}QP_{U_n}\Delta t \\
				&\quad
				+ P_{U_n}\mathbb{E}\!\left[A(t_n, X_n^{\mathrm{DLRA},c})Y_n^{\top}\right]U_n\Delta t
				+ U_n^{\top}\mathbb{E}\!\left[Y_n \mathring{A}(t_n, X_n^{\mathrm{DLRA},c})^\top\right]P_{U_n}\Delta t \\
				&\quad
				+ P_{U_n}QP_{U_n}\Delta t 
				+ U_n^{\top} \mathbb{E}\!\left[Y_n\mathring{A}(t_n, X_n^{\mathrm{DLRA},c})^{\top}\right]P_{U_n}^{\perp}\Delta t
				+ P_{U_n}QP_{U_n}^{\perp}\Delta t \\
				=&\ C_n^{\mathrm{DLRA}}
				+ \mathbb{E}\!\left[\mathring{A}(t_n, X_n^{\mathrm{DLRA},c})Y_n^{\top}\right]U_n\Delta t
				+ U_n^{\top}\mathbb{E}\!\left[Y_n\mathring{A}(t_n, X_n^{\mathrm{DLRA},c})^\top\right]\Delta t \\
				&\quad
				+ \left(Q-P_{U_n}^{\perp}QP_{U_n}^{\perp}\right)\Delta t.
			\end{aligned}
		\end{equation*}
		
		On the other hand, 
		the projector onto the predicted deterministic basis can be written as
		\begin{equation*}
			\begin{aligned}
				P_{\widehat{U}_{n+1}}^{\perp}
				=&\ \left(I_{d\times d}-(U_n+\Delta \widehat{U}_n)^\top(U_n+\Delta \widehat{U}_n)\right) \\
				=&\ \left(P_{U_n}^{\perp}-\Delta \widehat{U}_n^\top U_n-U_n^\top\Delta \widehat{U}_n-\Delta \widehat{U}_n^\top\Delta \widehat{U}_n,
				\right)
			\end{aligned}
		\end{equation*}
		where $\Delta \widehat{U}_n$ is the increment of the deterministic basis in the prediction as in \eqref{eq: pred U incr}. 
		\iffalse
		Then,
		one can write
		\begin{equation*}
			\begin{aligned}
				&\left(
				\widehat{C}_{n+1}^{\mathrm{DLRA}}
				+ \widehat{E}_{n+1}
				\right)
				H(t_n)^\top R^{-1} H(t_n)
				\left(
				\widehat{C}_{n+1}^{\mathrm{DLRA}}
				+ \widehat{E}_{n+1}
				\right)\Delta t \\
				\approx& \left(C_n^{\mathrm{DLRA}}+ \widehat{E}_{n+1} \right)H(t_n)^\top R^{-1} H(t_n)\left(C_n^{\mathrm{DLRA}}+ \widehat{E}_{n+1}\right)\Delta t \\
			\end{aligned}
		\end{equation*}
		\fi
		
		As we want to solve the derivative of $C_{n+1}$ with respect to the $\Delta U_n^{c}$, i.e.\
		\begin{equation*}
			\begin{aligned}
				\frac{\partial \left\|C_{n+1}^{\mathrm{DLRA},c}-C_{n+1}\right\|^2_{\mathrm{F}}}{\partial \Delta U_n^c}
				= \frac{\partial \operatorname{Tr}\left(C_{n+1}^{\mathrm{DLRA},c}\left(C_{n+1}^{\mathrm{DLRA},c}\right)^\top\right)}{\partial \Delta U_n^c}
				- \frac{2\partial \operatorname{Tr}\left(C_{n+1}^{\mathrm{DLRA},c}\left(C_{n+1}\right)^\top\right)}{\partial \Delta U_n^c}
			\end{aligned}
		\end{equation*}
		and the derivative with respective $\Delta C_{Y_n}^{c} = \mathbb{E}[\Delta Y_n^{c} (\Delta Y_n^{c})^{\top}]$, we recall that we have the following relations for the derivative with respect to $\Delta U_n^c$
		\begin{equation*}
			\begin{aligned}
				\frac{\partial \operatorname{Tr}\left(\left(C_{n+1}^{\mathrm{DLRA},c}\right)\left(C_{n+1}^{\mathrm{DLRA},c}\right)^\top\right)}{\partial \Delta U_n^c}
				=&\ 2C_{Y_n}^2U_n+4C_{Y_n}^2\Delta U_n^c+2C_{Y_n}\Delta C_{Y_n}^cU_n \\
				&+2C_{Y_n}\Delta U_n^c(\Delta U_n^c)^\top C_{Y_n}U_n+4C_{Y_n}\Delta U_n^c(\Delta U_n^c)^\top C_{Y_n}\Delta U_n^c,
			\end{aligned}
		\end{equation*}
		\begin{equation*}
			\begin{aligned}
				\frac{\partial \operatorname{Tr}\left(\left(C_{n+1}^{\mathrm{DLRA},c}\right)\left(C_{n+1}\right)^\top\right)}{\partial \Delta U_n^c}
				=2C_{Y_n}U_nC_{n+1}+2C_{Y_n}\Delta U_n^cC_{n+1},
			\end{aligned}
		\end{equation*}
		and for $\Delta C_{Y_n}^c$
		\begin{equation*}
			\begin{aligned}
				\frac{\partial \operatorname{Tr}\left(\left(C_{n+1}^{\mathrm{DLRA},c}\right)\left(C_{n+1}^{\mathrm{DLRA},c}\right)^\top\right)}{\partial \Delta C_{Y_n}^c}
				=2C_{Y_n}+2\Delta C_{Y_n}^c,
			\end{aligned}
		\end{equation*}
		and
		\begin{equation*}
			\begin{aligned}
				\frac{\partial \operatorname{Tr}\left(\left(C_{n+1}^{\mathrm{DLRA},c}\right)C_{n+1}^\top\right)}{\partial \Delta C_{Y_n}^c}
				=U_nC_{n+1}U_n^\top.
			\end{aligned}
		\end{equation*}
		
		As computations are very similar to Proposition \eqref{prop: triple JMCO} for the minimization \eqref{eq: DLRA JMCO enriched min} we will streamline them here below. Let us first compute $\Delta C_{Y_n}^c$, whose quantity is given by solving
		\begin{equation*}
			\begin{aligned}
				\frac{\partial \left\|C_{n+1}^{\mathrm{DLRA},c}-C_{n+1}\right\|^2_{\mathrm{F}}}{\partial \Delta C_{Y_n}^{c}}=0.
			\end{aligned}
		\end{equation*}
		This minimization implies that
		\begin{equation*}
			\begin{aligned}
				0=&\ 2C_{Y_n}+2\Delta C_{Y_n}-2U_nC_{n+1}U_n^\top \\
				=&\ 2\Delta C_{Y_n}^c-2\Big(\big(C_{Y_n}+U_n\mathbb{E}\!\left[\mathring{A}(t_n, X_n^{\mathrm{DLRA},c})Y_n^\top\right]\Delta t+\mathbb{E}\!\left[Y_n\mathring{A}(t_n, X_n^{\mathrm{DLRA},c})^\top\right]U_n\Delta t \\
				&+U_nQU_n\Delta t\big)-(C_{Y_n}U_n + U_n \widehat{E}_{n+1} ) H(t_n)^\top R^{-1} H(t_n)( \widehat{E}_{n+1} U_n^{\top} + U_n^\top C_{Y_n})\Delta t\Big) \\
			   =&\ 2\Delta C_{Y_n}^c-2\Big(\big(C_{Y_n}+U_n\mathbb{E}\!\left[\mathring{A}(t_n, X_n^{\mathrm{DLRA},c})Y_n^\top\right]\Delta t+\mathbb{E}\!\left[Y_n\mathring{A}(t_n, X_n^{\mathrm{DLRA},c})^\top\right]U_n\Delta t \\
			   &+U_nQU_n\Delta t\big)-C_{Y_n}U_n H(t_n)^\top R^{-1} H(t_n) U_n^\top C_{Y_n} \Delta t\Big),
			\end{aligned}
		\end{equation*}
		where in the last line we use the fact that 
		\begin{equation*}
			\begin{aligned}
			U_n E_{n+1} = &U_n P_{\widehat{U}_{n+1}}^{\perp} \widehat{E}_{n+1}  \\
			=& U_n \left(P_{U_n}^{\perp}-\Delta \widehat{U}_n^\top U_n-U_n^\top\Delta \widehat{U}_n-\Delta \widehat{U}_n^\top\Delta U_n
			\right) \widehat{E}_{n+1}  = \Delta \widehat{U}_n  \widehat{E}_{n+1} = O(\Delta t).
						\end{aligned}
		\end{equation*}
		This last relations yields the equation for $\Delta Y_n^{c}$.
		
		Finally, we compute $\Delta U_n^c$ from the variational formulation
		\begin{equation*}
			\begin{aligned}
				\left\langle \frac{\partial \left\|C_{n+1}^{\mathrm{DLRA},c}-C_{n+1}\right\|^2_{\mathrm{F}}}{\partial \Delta U_n^c},V\right\rangle=0 \qquad \forall\,V\in\mathbb{R}^{k\times d},\quad V\in\operatorname{Im}(U_n^\perp).
			\end{aligned}
		\end{equation*}
	 This relation implies that
		\begin{equation*}
			\begin{aligned}
				0=&\big\langle 2C_{Y_n}^2U_n+4C_{Y_n}^2\Delta U_n^c+2C_{Y_n}\Delta C_{Y_n}^cU_n+2C_{Y_n}\Delta U_n^c(\Delta U_n^c)^\top C_{Y_n}U_n \\
				&+4C_{Y_n}\Delta U_n^c(\Delta U_n^c)^\top C_{Y_n}\Delta U_n^c-2C_{Y_n}U_nC_{n+1}-2C_{Y_n}\Delta U_n^{c}C_{n+1},V\big\rangle \\
				=&\big\langle 4C_{Y_n}^2\Delta U_n^c+4C_{Y_n}\Delta U_n^c(\Delta U_n^c)^\top C_{Y_n}\Delta U_n^c \\
				&-2C_{Y_n}U_nC_{n+1}P_{U_n}^\perp-2C_{Y_n}\Delta U_n^{c}C_{n+1}P_{U_n}^\perp,V\big\rangle
			\end{aligned}
		\end{equation*}
		which is equivalent in a first order of $O(\Delta t)$ to
		\begin{equation*}
			\begin{aligned}
				C_{Y_n}\Delta U_n^c
				=&\ U_n C_{n+1}P_{U_n}^{\perp}+\Delta U_n^{c} C_{n+1}P_{U_n}^{\perp} \\
				=&\ (U_n+\Delta U_n^{c})\Big(C_n^{\mathrm{DLRA}}+\mathbb{E}\!\left[\mathring{A}(t_n, X_n^{\mathrm{DLRA},c})Y_n^\top\right]U_n\Delta t+U_n^{\top}\mathbb{E}\!\left[Y_n\mathring{A}(t_n, X_n^{\mathrm{DLRA},c})^\top\right]\Delta t \\
				& \qquad \qquad \qquad +\left(Q-P_{U_n}^{\perp}QP_{U_n}^{\perp}\right)\Delta t \\
				&\qquad-\left(C_n^{\mathrm{DLRA}}+ \widehat{E}_{n+1}\right)H(t_n)^\top R^{-1} H(t_n)\left(C_n^{\mathrm{DLRA}}+ \widehat{E}_{n+1}\right)\Delta t\Big)P_{U_n}^{\perp} \\
				\approx& U_n\Bigg(C_n^{\mathrm{DLRA}}+\mathbb{E}\!\left[\mathring{A}(t_n, X_n^{\mathrm{DLRA},c})Y_n^\top\right]U_n\Delta t+U_n^{\top}\mathbb{E}\!\left[Y_n\mathring{A}(t_n, X_n^{\mathrm{DLRA},c})^\top\right]\Delta t \\
				& \qquad \qquad \qquad +\left(Q-P_{U_n}^{\perp}QP_{U_n}^{\perp}\right)\Delta t \\
				&\qquad-\left(C_n^{\mathrm{DLRA}}+\widehat{E}_{n+1}\right)H(t_n)^\top R^{-1} H(t_n)\left(C_n^{\mathrm{DLRA}}+ \widehat{E}_{n+1}\right)\Delta t\Bigg)P_{U_n}^{\perp} + O(\Delta t^2)\\
				=&\ \mathbb{E}\!\left[Y_n\mathring{A}(t_n, X_n^{\mathrm{DLRA},c})^\top\right]P_{U_n}^{\perp}\Delta t+U_nQP_{U_n}^{\perp}\Delta t\\
				&-(C_{Y_n}U_n+ U_n  \widehat{E}_{n+1}) H(t_n)^\top R^{-1} H(t_n)\widehat{E}_{n+1} P_{U_n}^{\perp}\Delta t
			\end{aligned}
		\end{equation*}
		where in the second line we use the fact that $\Delta U_n^{c}U_n^\top=0$ and the ansatz $\Delta U_n^{c}=e_n\Delta t+f_n\Delta Z_n$, which gives the thesis.
		
		Finally, by minimizing $\frac{\partial \left\|m_{n+1}^{\mathrm{DLRA},c}-m_{n+1}\right\|^2_{\mathrm{F}}}{\partial \Delta m_n^{c}}$ one obtains
		\begin{equation}
			\begin{aligned}
			2 \Delta m_{n}^{\mathrm{DLRA},c} =& 2 \widehat{m}_{n+1}^{\mathrm{DLRA}} 
			+ 	\left(\widehat{C}_{n+1}^{\mathrm{DLRA}} + \widehat{E}_{n+1}\right) H^\top R^{-1} \left( \Delta Z_n - H(t_n) \widehat{m}_{n+1}^{\mathrm{DLRA}} \Delta t \right) \\
			=&\widehat{m}_{n+1}^{\mathrm{DLRA}} 
			+ 	\left(C_{n}^{\mathrm{DLRA}} + P_{U_n}^{\perp} \mathbb{E}\left[D_n D_n^{\top}\right] P_{U_n}^{\perp} + O(\Delta t) \right) H^\top R^{-1} \left( \Delta Z_n - H(t_n) \widehat{m}_{n+1}^{\mathrm{DLRA}} \Delta t \right) \\
			\approx & \widehat{m}_{n+1}^{\mathrm{DLRA}} 
			+	\left(C_{n}^{\mathrm{DLRA}} + P_{U_n}^{\perp} \mathbb{E}\left[D_n D_n^{\top}\right] P_{U_n}^{\perp} \right) H^\top R^{-1} \left( \Delta Z_n - H(t_n) \widehat{m}_{n+1}^{\mathrm{DLRA}} \Delta t \right),\\
			\end{aligned}
		\end{equation}
		from which we obtain the thesis.
		\end{proof}
		
\section{Proofs concerning DLRA-PF filter}\label{app: proof DLRA PF}

\paragraph{The Branching Algorithm \cite[Section 9.2.1]{bain2009fundamentals}}

Here we show the branching process employed at the end of the analysis step of the particle filter only after certain intervals of time $\Delta t_{\mathrm{branch}}$, both for the full-order version and the DLRA one. 
By $r_i$ we denote the number of offspring assigned to the $i$-th particle during the branching step, i.e.\ the integer number of copies of the $i$-th realization of the state produced by the branching. Furthermore $\lambda^{i}$ denote the normalized weight produced by the analysis step through the likelihood, associate to the $i$-th realization of the state.
The number of particles in
the system remains constant at $M$; that is, we require that for each $i$,

\begin{center}
	$\displaystyle
	\sum_{i=1}^{M} r_{i} = M.
	$
\end{center}

Let $u_{i}^{M}$, $i = 1,\ldots,M-1$ be $M-1$ mutually independent
random variables, uniformly distributed on $[0,1]$, which are independent of
all other random variables in the system. Then the \emph{branching algorithm} is summarized below

\begin{algorithmic} \label{alg: branching}
	\State $g := M \qquad h := M$
	
	\For{$i := 1$ to $M-1$}
	
	\If{$\{M\lambda^{i}\}
		+ \{g-M\lambda^{i}\} < 1$}
	
	\If{$u_{i}^{M}
		< 1-\left(
		\{M\lambda_{i}^{M}\}/\{g\}
		\right)$}
	
	\State $r_{i}
	:= [M\lambda^{i}]$
	
	\Else
	
	\State $r_{i}
	:= [M\lambda^{i}]
	+ (h-[g])$
	
	\EndIf
	
	\Else
	
	\If{$u_{i}^{M}
		< 1-
		\left(
		1-\{M\lambda^{i}\}
		\right)
		/
		\left(
		1-\{g\}
		\right)$}
	
	\State $r_{i}
	:= [M\lambda^{i}] + 1$
	
	\Else
	
	\State $r_{i}
	:= [M\lambda^{i}]
	+ (h-[g])$
	
	\EndIf
	
	\EndIf
	
	\State $g := g-M\lambda^{i}$
	\State $h := h-r_{i}$
	
	\EndFor
	
	\State $r_M := h$,
\end{algorithmic}
where $[ \cdot ]$ denotes the integer part of a number.

\begin{proof}[\bfseries Proof of Lemma \ref{lem: invar ran}]
	According to \cite[Section 9.2.1]{bain2009fundamentals}, the branching procedure applied
	to a matrix $\widehat{\mathbb{X}} = [\widehat{X}^1 \mid \dots \mid \widehat{X}^M] \in \mathbb{R}^{d \times M}$, where $d$ is the physical dimension and $M$ is
	the number of particles, has as output a new matrix $\mathbb{X} = [X^1 \mid \dots \mid X^M] \in \mathbb{R}^{d \times M}$ whose columns
	are extracted, and possibly repeated, from the ones of $\widehat{\mathbb{X}}$, i.e.\
	$X^i \in \{\widehat{X}^i\}_{i=1,..,M}$, therefore
	\begin{equation*}
		\mathbb{X}= \widehat{\mathbb{X}} B,
	\end{equation*}
	where $B \in \mathbb{R}^{M \times M}$, with $B = [ \underbrace{e_1 | e_1 | e_1}_{r_1} | ...\underbrace{ | e_i | }_{r_i}... \underbrace{| e_M | e_M}_{r_M} ] ,$
	$e_i$ is the i-th canonical basis, $r_i$ are random variables with multinomial distribution with values $r_i \in \mathbb{N} \cup \{0\}$,
	and $\sum_{i=1}^{M} r_i =  M.$
	Therefore, the analysis step of the DLRA-PF at time $n+1$ can be
	obtained as
	\begin{equation*}
		\mathbb{X}_{n+1} = \widehat{\mathbb{X}}_{n+1} B,
	\end{equation*}
	where $\widehat{\mathbb{X}}_{n+1}  = \left(\widehat{\mathbb{X}}_{n+1}^{1}, \dots, \widehat{\mathbb{X}}_{n+1}^{M} \right)$ are the $M$ predicted particles, $\widehat{\mathbb{X}}_{n+1} =  \widehat{m}_{n+1}^{\mathrm{DLRA}} \mathbbm{1}^\top + \widehat{U}_{n+1}^{\top} \widehat{\mathbb{Y}}_{n+1}$ and $\widehat{\mathbb{Y}}_{n+1}= [\widehat{Y}_{n+1}^1 \mid \dots \mid \widehat{Y}_{n+1}^M] \in \mathbb{R}^{k \times M}$. Hence $X_{n+1}$ still admits a three-term decomposition with $m_{n+1}= \mathbb{E}[X_{n+1}]= \widehat{m}_{n+1}^{\mathrm{DLRA}}+ \widehat{U}_{n+1}^{\top} \widehat{\mathbb{Y}}_{n+1} B \mathbbm{1} \frac{1}{M}$ and $\mathbb{Y}_{n+1} = \widehat{\mathbb{Y}}_{n+1}B -  \widehat{\mathbb{Y}}_{n+1}B \frac{1}{M}$. We see then that $\mathrm{rank}(U_{n+1}^{\top}) =\mathrm{Ran}(\mathbb{X}_{n+1}- \mathbb{E}[\mathbb{X}_{n+1}]) \subset \mathrm{Ran}( U_{n+1}^{\top})$.
\end{proof}

\end{document}